\documentclass[12pt]{article}
\usepackage{amsmath,mathtools} 
\usepackage{amsthm}
\usepackage{dsfont}
\usepackage{thmtools}

\usepackage{booktabs} 
\usepackage{xspace}
\usepackage{graphicx}
\usepackage{xspace}
\usepackage{url}
\usepackage{enumerate}

\usepackage{caption} 
\usepackage{subcaption}

\usepackage[algo2e,ruled,vlined,linesnumbered]{algorithm2e} 
\SetAlgoSkip{}

\usepackage{balance} 

\allowdisplaybreaks[3] 

\usepackage{tikz}
\usetikzlibrary{shapes,arrows}
\usetikzlibrary{decorations}
\usetikzlibrary{plotmarks}
\usetikzlibrary{calc}
\usetikzlibrary{mindmap}
\usetikzlibrary{shadows}
\usetikzlibrary{backgrounds}
\usetikzlibrary{patterns}
\usetikzlibrary{shapes.symbols}

\usepackage{pgfplots}
\usepackage{pgfplotstable}
\usepgfplotslibrary{statistics}

\newtheorem{theorem}    		 {Theorem}
\newtheorem{lemma}      [theorem]{Lemma}
\newtheorem{corollary}  [theorem]{Corollary}

\newtheorem{fact}		[theorem]{Fact}
\newtheorem{claim}    [theorem]{Claim}
\theoremstyle{definition}
\newtheorem{definition} [theorem]{Definition}
\newtheorem{remark}     [theorem]{Remark}
\newtheorem{example}    [theorem]{Example}

\newtheorem{convention}    [theorem]{Convention}

\usepackage[english]{babel}
\usepackage{amssymb}
\usepackage{amsmath}
\usepackage{graphicx}
\usepackage{bbm}
\usepackage[right]{eurosym}
\usepackage[labelfont={bf,sf},font={small},%
labelsep=space]{caption}
\usepackage{chngcntr}
\usepackage{makeidx}
\usepackage{xcolor}
\makeindex

\counterwithin{figure}{section}

\def\vs{\vspace}

\def\exp{\mathrm{exp}}

\def\arg{\mathrm{arg}}

\usepackage{imakeidx}
\makeindex[columns=2, options=-s mystyle]

\usepackage{filecontents}

\begin{document}

\begin{center}
	{\Large\bf Global Complexification of Real Analytic Restricted Log-Exp-Analytic Functions}
	
	\vs{0.3cm}
	Andre Opris
\end{center}

{\small {\bf Abstract.}
We show that a real analytic restricted log-exp-analytic function has a holomorphic extension which is again restricted log-exp-analytic. We also establish a parametric version of this result.}

\section*{Introduction}
This paper contributes to analysis in the framework of o-minimal structures. O-minimality is a concept from mathematical logic with connections and applications to geometry, analysis, number theory and other areas. Sets and functions definable in an o-minimal structure (i.e. ‘belonging to’) exhibit tame geometric, combinatorial and topological behavior. This includes for example the cell decomposition theorem, the fact that every definable set is homeomorphic to a finite union of hypercubes and that every definable set has finitely many connected components which are again definable. We refer to the book of Van den Dries~\cite{11} for the general properties of o-minimal structures.

The major line of research was based on discovering o-minimal structures on the reals. For example the structure on the real field $\mathbb{R}$ is o-minimal. The definable sets are exactly the semialgebraic sets since they are closed under projections (see Tarski~\cite{36} which was also popularized by Seidenberg~\cite{32}). But there are much bigger o-minimal structures on the real field. Together with a theorem of Gabrielov Van den Dries could prove in~\cite{10} for example the remarkable result that $\mathbb{R}_{\textnormal{an}}$, the structure generated by all restricted analytic functions, is o-minimal. The definable sets in $\mathbb{R}_{\textnormal{an}}$ are precisely the globally subanalytic ones. It is also shown in~\cite{10} that $\mathbb{R}_{\textnormal{an}}$ is polynomially bounded. So the global exponential function is not definable in there. Van den Dries et al. even showed in~\cite{15} that we stay in the o-minimal context if we add the global exponential function (see also van den Dries et al.~\cite{13}). Then we obtain the o-minimal structure $\mathbb{R}_{\textnormal{an,exp}}$. This is one of the most important o-minimal structures, because a huge class of elementary functions like polynomial functions, the global arctangent function, the global logarithm, the global exponential function and hyperbolic functions are definable in there. They have crucial application in diophantine geometry (see for example Wilkie~\cite{37} where definable functions with special diophantine properties and growth rates are investigated or Pila~\cite{31} for the treatment of the Andr\'{e}-Oort conjecture for $\mathbb{C}^n$ from the viewpoint of o-minimality).

Kobi Peterzil and Sergei Stepanovich Starchenko designed the development of complex analysis within the o-minimal framework also in non-standard setting (see Peterzil et al.~\cite{30}). One obtains analogies of classical results as well as strong differences, due to the o-minimality assumption. For example the maximum principle, the open mapping theorem, the theorem of Liouville, the identity theorem and Riemann's removeable theorem hold in this setting. However, the key feature is that a definable holomorphic function in an o-minimal expansion of $\mathbb{R}$ does not have an essential singularity. So the entire definable holomorphic functions are exactly the polynomials. Further results of Peterzil and Starchenko concern analytic geometry. For example Peterzil et al.~\cite{29} gave some stronger version of Chow's theorem which states that a definable analytic subset of $\mathbb{C}^n$ is already an algebraic variety. All this theory can be formulated for an arbitrary real closed field $R$ with algebraic closure $K=R[\sqrt{-1}]$. But to describe such features in this non-standard setting it is necessary to use ''topological analysis'' instead of power series and integration.
So it would be desirable to understand definability in the complex setting outgoing from definability on the reals. A very interesting question in this context is the following.\\

{\bf Question}

{\it
Are real analytic functions definable in an o-minimal structure $\mathcal{M}$ on the reals reducts of definable holomorphic functions, i.e. does $\mathcal{M}$ have global complexification?}\\

To investigate this question one has to find a definable holomorphic extension for a given definable real analytic function (see Krantz et al.~\cite{23} for elementary properties of real analytic functions). Of course finding a holomorphic extension is not difficult simply by Taylor series expansion, but the crucial point is that this extension needs not to be definable. Tobias Kaiser gave a positive answer for this question for the o-minimal structure on the real field $\mathbb{R}$ (see Kaiser~\cite{18}) and $\mathbb{R}_{\textnormal{an}}$ (see Kaiser~\cite{17}). He could even establish parametric versions of this results. However, for the structure $\mathbb{R}_{\exp}$ (which is also o-minimal by a theorem of Khovanskii~\cite{22} together with the result from Wilkie~\cite{38} that the theory of $\mathbb{R}_{\exp}$ is model complete) it is easy to see that the above question has a negative answer (see Example 1.3 in~\cite{17}). The global exponential function $f:\mathbb{R} \to \mathbb{R}, x \mapsto e^x$ is $\mathbb{R}_{\text{exp}}$-definable, real analytic, and every holomorphic extension of $f$ is of the form $F:V \to \mathbb{C}, x+iy \mapsto e^x(\cos(y)+i\sin(y))$ where $\mathbb{R} \subset V \subset \mathbb{C}$ is open in $\mathbb{C}$. However, by a result of Bianconi in~\cite{1}), no restriction of the global sine function on a non-empty open set is $\mathbb{R}_{\exp}$-definable and therefore, $F$ is not $\mathbb{R}_{\exp}$-definable. 
For $\mathbb{R}_{\textnormal{an,exp}}$ issues regarding global complexifixation are not solved at present. But there is the deep model theoretical fact that every definable function in $\mathbb{R}_{\textnormal{an,exp}}$ is piecewise given by $\mathcal{L}_{\textnormal{an}}(\exp,\log)$-terms ($\mathcal{L}_{\textnormal{an}}(\exp,\log)$ is the language of ordered rings augmented by all restricted analytic functions, the global exponential and the global logarithm, see~\cite{13}). Kaiser used this fact to show in \cite{17} that univariate real analytic functions which are definable in $\mathbb{R}_{\textnormal{an,exp}}$ extend in a definably way to a holomorphic function. Also a quantitative version of this result was given by Kaiser et al.~\cite{21} (see also Wilkie~\cite{37}). But especially for this big structure it would be desirable to have an analoguous result for multiple variables as in the globally subanalytic case, because there are many fields of application. 

But there are strong differences between the structures $\mathbb{R}_{\textnormal{an}}$ and $\mathbb{R}_{\textnormal{an,exp}}$ from the point of analysis: There are $\mathbb{R}_{\textnormal{an,exp}}$-definable functions which are $C^\infty$, but not real analytic (for example flat functions) or violate Tamm's Theorem (see van den Dries et al.~\cite{14}, Kaiser et al~\cite{19} or Opris~\cite{28} for versions of Tamm's Theorem) in contrast to globally subanalytic functions. So there is the natural question if there is a big non-trivial class of $\mathbb{R}_{\textnormal{an,exp}}$-definable functions which share some analytic properties with globally subanalytic ones, but behaves completely differently from the point of o-minimality. This applies for example for restricted log-exp-analytic functions which have been defined by Opris (see~\cite{28}). They are given piecewise by compositions of globally subanalytic functions, the global logarithm and exponentials of locally bounded functions. For example these functions satisfy a parametric version of Tamm's theorem (see also~\cite{28}), as well as globally subanalytic ones. Since the global exponential function $\mathbb{R} \to \mathbb{R}, x \mapsto e^x,$ is restricted log-exp-analytic (because it is the exponential of the identity and the latter is locally bounded), this class of functions generates the whole structure $\mathbb{R}_{\textnormal{an,exp}}$ in the sense of definability, but is not polynomially bounded in contrast to globally subanalytic functions.\\

{\bf Example.} The function
$$g: \textnormal{}]0,1[^2 \to \mathbb{R}, (t,x) \mapsto \arctan(\log(e^{ \log^2(1/x)/t}+\log(e^{e^{1/t}}+2))),$$ 
is restricted log-exp-analytic since $g$ is definable and the real functions $(t,x) \mapsto \log^2(1/x)/t$, $(t,x) \mapsto 1/t$ and $(t,x) \mapsto e^{1/t}$ with domain $]0,1[^2$ are locally bounded.\\

We will not answer the question above for all real analytic $\mathbb{R}_{\textnormal{an,exp}}$-definable functions, but for real analytic restricted log-exp-analytically ones which is the main goal for this paper. Hence, we obtain a further analytic property  that restricted log-exp-analytic functions and globally subanalytic ones have in common.\\

{\bf Theorem A.} {\it Let $n \in \mathbb{N}$ and $U \subset \mathbb{R}^n$ be an $\mathbb{R}_{\text{an,exp}}$-definable open set. Let $f:U \to \mathbb{R}$ be a real analytic restricted log-exp-analytic function. Then there is a definable open $V \subset \mathbb{C}^n$ with $U \subset V$ and a holomorphic restricted log-exp-analytic $F:V \to \mathbb{C}$ such that $F|_U=f$ (i.e. the real and imaginary part of $F$ considered as real functions are restricted log-exp-analytic).}\\

The idea for the proof of Theorem A is the same procedure as in the globally subanalytic case: Outgoing from a preparation theorem for restricted log-exp-analytic functions in~\cite{28} (see also Theorem C in Opris~\cite{27} for a more general preparation theorem for $\mathbb{R}_{\textnormal{an,exp}}$-definable functions or Lion et al.~\cite{25} for original versions) we construct a suitable holomorphic extension of a parametric family of real analytic restricted log-exp-analytic one variable functions and then we do a technical induction on the number of variables. This gives the following version of Theorem A for parameters.\\

{\bf Theorem B.} {\it Let $n \in \mathbb{N}_0$, $m \in \mathbb{N}$ and $X \subset \mathbb{R}^n \times \mathbb{R}^m$ be $\mathbb{R}_{\text{an,exp}}$-definable such that $X_t:=\{x \in \mathbb{R}^m \mid (t,x) \in X\}$ is open for all $t \in \mathbb{R}^n$. Let $f:X \to \mathbb{R}, (t,x) \mapsto f(t,x)$, be restricted log-exp-analytic in $x$ such that $f_t:X_t \to \mathbb{R}, x \mapsto f(t,x),$ is real analytic for every $t \in \mathbb{R}^n$. Then there is a definable $Z \subset \mathbb{R}^n \times \mathbb{C}^m$ with $X \subset Z$ such that $Z_t:=\{z \in \mathbb{C}^m \mid (t,z) \in Z\}$ is open for all $t \in \mathbb{R}^n$ and a function $F:Z \to \mathbb{C}, (t,z) \mapsto F(t,z),$ which is restricted log-exp-analytic in $z$ such that $F_t:Z_t \to \mathbb{C}, z \mapsto F(t,z),$ is holomorphic for every $t \in \mathbb{R}^n$ and $F|_X=f$ (i.e. the real and imaginary part of $F$ considered as real functions are restricted log-exp-analytic in $(x,y)$ where $z=x+iy$).}\\

Here a restricted log-exp-analytic function $f:X \to \mathbb{R}, (t,x) \mapsto f(t,x),$ in $x$ where $X \subset \mathbb{R}^n \times \mathbb{R}^m$ is definable such that $X_t$ is open for $t \in \mathbb{R}^n$ is piecewise given by compositions of globally subanalytic functions, the global logarithm and exponentials of locally bounded functions $g:X \to \mathbb{R}, (t,x) \mapsto g(t,x),$ in $x$ (where $g$ is locally bounded in $x$ if $g_t$ is locally bounded for every $t \in \mathbb{R}^n$).\\
However, obtaining a parametric version for unary functions requires a rather sophisticated setup. One needs definability results of special parameterized integrals (similar as in Kaiser~\cite{17} for the globally subanalytic case), several results about the asymptotical behavior of restricted log-exp-analytically prepared functions to establish restricted log-exp-analyticity of the constructed holomorphic extension and methods from complex analysis.\\

This paper is organized as follows. After a short introduction in o-minimality (see Section 1) we introduce global complexification in Section 2 and give examples. In Section 3 we pick up the definition of log-analytic functions, the exponential number and restricted log-exp-analytic functions from~\cite{28}. Section 4 is about a preparation theorem for definable functions from~\cite{27} and for restricted log-exp-analytic functions from~\cite{28}. 
In Section~5 we construct a holomorphic extension of a parametric family of restricted log-exp-analytically prepared one variable functions and investigate their asymptotical behavior under taking definable curves in Section~6. Then we establish a further preparation theorem (Section~7) to obtain the desired result on integration (Section~8). Finally, Section~9 contains the proof of our main results. \\

{\bf Notation:} The empty sum is by definition $0$ and the empty product is by definition $1$. By $\mathbb{N}=\{1,2, \ldots \}$\index{$\mathbb{N}$} we denote the set of natural numbers and by $\mathbb{N}_0=\{0,1,2, \ldots \}$\index{$\mathbb{N}_0$} the set of nonnegative integers. 
For $a \in \mathbb{R}$ we set $\mathbb{R}_{>a}:=\{x \in \mathbb{R} \mid x > a\}$\index{$\mathbb{R}_{>a}$}. Given $x \in \mathbb{R}$ let $\lceil{x}\rceil$ \index{$\lceil{x}\rceil$} be the smallest integer which is not smaller than $x$ and let $\textnormal{sign}(x) \in \{\pm 1\}$ be its \index{$\textnormal{sign}(x)$}sign if $x \neq 0$. For $a \in \mathbb{R}$ we have $a<\infty$, $-\infty < a$ and we set $a+\infty:=\infty$ and $a-\infty:=-\infty$. We set $\textnormal{sup}(\emptyset)=-\infty$ and $\textnormal{inf}(\emptyset)=\infty$. For $m,n \in \mathbb{N}$ we denote by $M(m \times n,\mathbb{Q})$\index{$M(m \times n,\mathbb{Q})$} the set of $(m \times n)$-matrices with rational entries. For $P \in M(m \times n, \mathbb{Q})$ we denote by $^tP \in M(n \times m,\mathbb{Q})$ its transpose\index{$^tP$}. Given a subset $A$ of $\mathbb{R}^n$ we denote by $\overline{A}$ its closure and by $\partial{A}$ its boundary. Given $x=(x_1, \ldots ,x_n) \in \mathbb{R}^n$ and $r>0$, we define the box
$$Q^n(x,r):=\{(y_1, \ldots ,y_n) \in \mathbb{R}^n \mid \vert{y_j-x_j}\vert<r \textnormal{ for all }j \in \{1, \ldots ,n\}\}.$$  
\index{$Q^n(x,r)$}For $m \in \mathbb{N}$ and $X \subset \mathbb{R}^m$ we say for two functions $f,g:X \to \mathbb{R}$ that $f<g$ ($f \leq g$) if $f(x)<g(x)$ ($f(x) \leq g(x)$) for every $x \in X$. For $x \in \mathbb{R}^m$ let\index{$\text{dist}(x,X)$}  
$$\text{dist}(x,X):=\inf\{\vert{x-w}\vert \mid w \in X\}.$$
A \textbf{curve}\index{curve} is a continuous mapping $\gamma:\text{}]0,1[\text{} \to X$, a \textbf{subcurve}\index{subcurve} of a curve $\gamma: \textbf{}]0,1[\textnormal{} \to X$ is a curve $\hat{\gamma}: \textnormal{}]0,1[\textnormal{} \to X$ of the following form: There is $\delta > 1$ such that $\hat{\gamma}(y) = \gamma(y/\delta)$ for $y \in \textbf{}]0,1[$.\\

For $X \subset \mathbb{R}^n \times \mathbb{R}^m$ and $t \in \mathbb{R}^n$ we set for the \textbf{fiber} $X_t:=\{x \in \mathbb{R}^m \mid (t,x) \in X\}$ \index{$X_t$}\index{fiber}and for a function $f:X \to \mathbb{C}, (t,x) \mapsto f(t,x),$ let $f_t:X_t \to \mathbb{C}, x \mapsto f(t,x)$.\index{$f_t$}\\

For this paragraph let $t$ ranges over $\mathbb{R}^n$ and $x$ over $\mathbb{R}$ and let $\pi:\mathbb{R}^n \times \mathbb{R} \to \mathbb{R}^n, (t,x) \mapsto t$. Then for $X \subset \mathbb{R}^n \times \mathbb{R}$ let $X_{\neq 0}:=\{(t,x) \in X \mid x \neq 0\}$ \index{$X_{\neq 0}$} and for a function $h:\pi(X) \to \mathbb{R}$ and a set $C \subset \mathbb{R}^n \times \mathbb{R}$ we say that $h<C$ ($h>C$) if for every $(t,x) \in X$: $h(t,x)<a$ ($h(t,x)>a$) for all $a \in C_t$\index{$h>C$}\index{$h<C$}. Further we say for $X_1,X_2 \subset \mathbb{R}^n \times \mathbb{R}$ with $\pi(X_1)=\pi(X_2)$ that $X_1<X_2$ \index{$X_1<X_2$} if for every $t \in \pi(X_1)$ we have that $x_1<x_2$ for every $x_1 \in (X_1)_t$ and $x_2 \in (X_2)_t$.\\

We set $\mathbb{C}^-:=\mathbb{C} \setminus \mathbb{R}_{\leq 0}$ \index{$\mathbb{C}^-$}, $\mathbb{C}^-_a:=\mathbb{C} \setminus \mathbb{R}_{\leq a}$ \index{$\mathbb{C}^-_a$} and $\mathbb{C}^+_a:=\mathbb{C} \setminus \mathbb{R}_{\geq a}$ \index{$\mathbb{C}^+_a$} for $a \in \mathbb{R}$. We denote by $\arg:\mathbb{C} \setminus \{0\} \to ]-\pi,\pi[$ the standard argument function \index{$\textnormal{arg}$} and define $\log:\mathbb{C}^- \to \mathbb{C}, z \mapsto \log(z):= \log(\vert{z}\vert)+i\textnormal{arg}(z)$ (i.e. $z$ is mapped on the principal value of the complex logarithm), and $z^q:=e^{q \log(z)}$ for $z \in \mathbb{C}^-$ and $q \in \mathbb{Q} \setminus \mathbb{Z}$\index{$z^q$}\index{$\log(z)$}. By $\log_k$\index{$\log_k$} we denote the $k$-times iterated of the complex logarithm 
and by $\exp_k$\index{$\exp_k$} the $k$-times iterated of the complex exponential 
where $\log_0=\exp_0=\textnormal{id}$. 
Let $z$ range over $\mathbb{C}$. 
For $r \in \mathbb{R}_{>0}$ and $c \in \mathbb{C}$ we define the open ball \index{$B(c,r)$}
$$B(c,r):=\{z \in \mathbb{C} \mid \vert{z-c}\vert<r\}.$$
Given $z=(z_1, \ldots ,z_n) \in \mathbb{C}^n$ and $r>0$, we define the polydisc \index{$D^n(z,r)$}
$$D^n(z,r):=\{(w_1, \ldots ,w_n) \in \mathbb{C}^n \mid \vert{w_j-z_j}\vert<r \textnormal{ for all }j \in \{1, \ldots ,n\}\}.$$

\index{$\textnormal{dist}$}
Let $m \in \mathbb{N}$ and $Z \subset \mathbb{C}^m$. For a set $E$ of functions on $Z$ with values in $\mathbb{C}^-$ we set $\log(E):=\{\log(g) \mid g \in E\}$\index{$\log(E)$}. For $Y \subset \mathbb{C}^m$ with $Y \subset Z$ and a set $E$ of functions on $Z$ we set $E|_Y:=\{G|_Y \mid G \in E\}$\index{$E|_Y$}.\\


{\bf Terminologies from model theory:} By $\mathcal{L}:=(\leq,+,\cdot,-,0,1)$ we denote the language of ordered rings, by $\mathcal{L}_{\textnormal{an}}$ the language of ordered rings expanded by symbols for all restricted analytic functions (compare with Definition~\ref{7} for the notion of a restricted analytic function), by $\mathcal{L}_{\textnormal{an}}(\exp,\log)$ the language $\mathcal{L}_{\textnormal{an}}$ expanded by symbols for the global exponential function and the global logarithm and by $\mathcal{L}_{\textnormal{an}}(^{-1},(\sqrt[n]{ \ldots })_{n=2,3, \ldots }, \log)$ the language $\mathcal{L}_{\textnormal{an}}$ augmented by a function symbol $^{-1}$ for taking reciprocals with respect to ''$\cdot$'', by function symbols $\sqrt[n]{ \ldots }$ for taking the $n$'th root for $n \in \mathbb{N}$ with $n \geq 2$, and by a symbol for the global logarithm. \index{$\mathcal{L}_{\textnormal{an}}$}\index{$\mathcal{L}$}\index{$\mathcal{L}_{\textnormal{an}}(\exp,\log)$}\index{$\mathcal{L}_{\textnormal{an}}(^{-1},(\sqrt[n]{ \ldots })_{n=2,3, \ldots }, \log)$}\\

An $\mathcal{L}'$-term is inductively defined as follows for a formal language $\mathcal{L}'$. \index{term!$\mathcal{L}_{\textnormal{an}}$}\index{term!$\mathcal{L}_{\textnormal{an}}(^{-1},(\sqrt[n]{ \ldots })_{n=2,3, \ldots }, \log)$}\index{term!$\mathcal{L}_{\textnormal{an}}(\exp,\log)$}
\begin{itemize}
	\item [(i)] $0$ and $1$ are $\mathcal{L}'$-terms.
	\item [(ii)] Let $x_1,x_2, \ldots $ be the infinite list of variables in $\mathcal{L}'$. Then $x_1,x_2, \ldots $ are $\mathcal{L}'$-terms.
	\item [(iii)] If $f$ is an $m$-ary function symbol in $\mathcal{L}'$ for $m \in \mathbb{N}$ and $t_1, \ldots ,t_m$ are $\mathcal{L}'$-terms then $f(t_1, \ldots ,t_m)$ is an $\mathcal{L}'$-term.
\end{itemize}

{\bf Conventions:} Definable means always $\mathbb{R}_{\textnormal{an,exp}}$-definable if not otherwise mentioned. We identify $\mathbb{C}$ with $\mathbb{R}^2$ via $x+iy \mapsto (x,y)$. So $Z \subset \mathbb{C}^m$ is definable if $Z$ is definable considered as a subset of $\mathbb{R}^{2m}$ for $m \in \mathbb{N}$. \index{definable!in $\mathbb{C}^m$} We say that a function $F:Y \to \mathbb{C}$ where $Y \subset \mathbb{R}^n \times \mathbb{C}$ (such that $Y_t$ is open for $t \in \mathbb{R}^n$) is continuous (holomorphic) in $z$ if $F_t$ is continuous (holomophic) for every $t \in \mathbb{R}^n$.\index{continuous in $z$}\index{holomorphic in $z$} For $Z \subset \mathbb{R}^m \times \mathbb{C}^n$ we say that a function $F:Z \mapsto \mathbb{C}$ is bounded at \index{bounded at $z$}$(u_0,z_0) \in \mathbb{R}^m \times \mathbb{C}^n$ if there is an open neighborhood $U$ of $(u_0,z_0)$ in $\mathbb{R}^m \times \mathbb{C}^n$ such that either $U \cap Z = \emptyset$ or $f|_{U \cap Z}$ is bounded. \\
Given $a,c \in \mathbb{C}$, $s \in \mathbb{R}_{>0}$, $D \subset \mathbb{C}$ with $\partial B(a,s) \subset D \cup \{c\}$ and a function $F:D \to \mathbb{C}$, we set
$$\int_{\partial B(a,s)} F(\xi) d \xi :=\int_{\partial B(a,s)} \tilde {F}(\xi) d\xi$$
if the integral to the right is well-defined where $\tilde{F}:D \cup \{c\} \to \mathbb{C}$ is any function with $\tilde{F}|_D=F$.
\section{O-Minimal Structures on the Real Field}

\begin{definition}
\label{1}
A subset $A$ of $\mathbb{R}^n$, $n \geq 1$, is called \textbf{semialgebraic}\index{semialgebraic} if there are $k,l \in \mathbb{N}_0$ and real polynomials $f_i,g_{i,1}, \ldots ,g_{i,k} \in \mathbb{R}[X_1, \ldots ,X_n]$ for $1 \leq i \leq l$ such that 
$$A = \bigcup_{i=1}^l \{x \in \mathbb{R}^n \mid f_i(x)=0, g_{i,1}(x)>0, \ldots ,g_{i,k}(x)>0\}.$$
A map is called semialgebraic if its graph is semialgebraic.
\end{definition}

\begin{definition}
\label{2}
A subset $A$ of $\mathbb{R}^n$, $n \geq 1$, is called \textbf{semianalytic}\index{semianalytic} if for each $a \in \mathbb{R}^n$ there are open neighborhoods $U,V$ of $a$ with $\overline{U} \subset V$, $k,l \in \mathbb{N}_0$ and real analytic functions $f_i,g_{i,1}, \ldots ,g_{i,k}$ on $V$ for $1 \leq i \leq l$, such that
$$A \cap U = \bigcup_{i=1}^l\{x \in U \mid f_i(x)=0, g_{i,1}(x)>0, \ldots ,g_{i,k}(x)>0\}.$$
A map is called semianalytic if its graph is semianalytic.
\end{definition}

\begin{definition}
\label{3}
A subset $B$ of $\mathbb{R}^n$, $n \geq 1$, is called \textbf{subanalytic}\index{subanalytic} if for each $a \in \mathbb{R}^n$ there is an open neighborhood $U$ of $a$, some $p \geq n$ and some bounded semianalytic set $A \subset \mathbb{R}^p$ such that $B \cap U=\pi_n(A)$ where $$\pi_n:\mathbb{R}^p \to \mathbb{R}^n, (x_1, \ldots ,x_p) \mapsto (x_1, \ldots ,x_n),$$
is the projection on the first $n$ coordinates.
A map is called subanalytic if its graph is subanalytic.
\end{definition}
	
\begin{remark}
\label{4}
A semialgebraic set is semianalytic. A semianalytic set is subanalytic.
\end{remark}

See Bierstone et al.~\cite{2} and Shiota~\cite{33} for geometrical descriptions of semianalytic resp. subanalytic sets and functions.

\begin{definition}
\label{5}
A subset $B$ of $\mathbb{R}^n$, $n \geq 1$, is called \textbf{globally subanalytic} \index{globally subanalytic!set} if it is subanalytic after applying the semialgebraic homeomorphism 
$$\mathbb{R}^n \to \textnormal{}]-1,1[^n, x_i \mapsto \frac{x_i}{\sqrt{1+x_i^2}},$$
for $i \in \{1, \ldots ,n\}$. A map is called globally subanalytic \index{globally subanalytic!map} if its graph is globally subanalytic.
\end{definition}

\begin{example}
\label{6}
A semialgebraic set is globally subanalytic. The restriction of the global sine function on a compact interval is globally subanalytic.	
\end{example}

\begin{definition}
\label{7}
A function $f:\mathbb{R}^n \to \mathbb{R}$ is called \textbf{restricted analytic}\index{restricted analytic} if there is a real convergent power series $p$ in $n$ variables which converges absolutely on an open neighborhood of $[-1,1]^n$ such that
$$f(x)=\left\{\begin{array}{ll} p(x), & x \in [-1,1]^n, \\
0, & \textnormal{ else.}\end{array}\right.$$
\end{definition}

\begin{definition}
\label{8}
A \textbf{structure}\index{structure} on the real field is axiomatically defined as follows. For $n \in \mathbb{N}$ let $M_n$ be a set of subsets of $\mathbb{R}^n$ and let $\mathcal{M}:=(M_n)_{n \in \mathbb{N}}$. Then $\mathcal{M}$ is a structure on the real field if the following holds for all $m,n,p \in \mathbb{N}$.
\begin{itemize}
	\item [(S1)] If $A,B \in \mathcal{M}^n$ then $A \cup B$, $A \cap B$ and $\mathbb{R}^n \setminus A \in M_n$. (So $M_n$ is a Boolean algebra of subsets of $M_n$.)
	\item [(S2)] If $A \in M_n$ and $B \in M_m$ then $A \times B \in M_{n+m}$.
	\item [(S3)] If $A \in M_p$ and $p \geq n$ then $\pi_n(A) \in M_n$ where $\pi_n:\mathbb{R}^p \to \mathbb{R}^n, (x_1, \ldots ,x_p) \mapsto (x_1, \ldots ,x_n)$, denotes the projection on the first $n$ coordinates.
	\item [(S4)] $M_n$ contains the semialgebraic subsets of $\mathbb{R}^n$. 
\end{itemize}
The structure $\mathcal{M}=(M_n)_{n \in \mathbb{N}}$ on the real field is called \textbf{o-minimal}\index{o-minimal structure} if additionally the sets in $M_1$ are exactly the finite unions of intervals and points.
\end{definition}

\begin{definition}
\label{9}
Let $\mathcal{M}=(M_n)_{n \in \mathbb{N}}$ be a structure. Let $n \in \mathbb{N}$.
\begin{itemize}
	\item [(a)]
	A subset $A$ of $\mathbb{R}^n$ is called \textbf{definable in $\mathcal{M}$}\index{definable!set} if $A \in M_n$.
	\item [(b)]
	Let $B \subset \mathbb{R}^n$ and $m \in \mathbb{N}$. A function $f:B \to \mathbb{R}^m$ is \textbf{definable in $\mathcal{M}$}\index{definable!function} if its graph $\{(x,f(x)) \mid x \in B\}$ is definable in $\mathcal{M}$. 
\end{itemize}
\end{definition}

\begin{remark}
\label{10}
The global sine or cosine function is not definable in an o-minimal expansion of $\mathbb{R}$.
\end{remark}

\begin{example}
\label{11}
\begin{itemize}
	\item [(1)] The smallest o-minimal structure on the real field is given by the semialgebraic sets (see~\cite{36} and~\cite{32}). We call this the structure on the real field \index{structure on the real field} and it is denoted by $\mathbb{R}$. \index{$\mathbb{R}$}
	\item [(2)] $\mathbb{R}_{\textnormal{exp}}$\index{$\mathbb{R}_{\textnormal{exp}}$}, the structure generated on the real field by the global exponential function $\exp:\mathbb{R} \to \mathbb{R}_{>0}$ (i.e. the smallest structure containing the semialgebraic sets and the graph of the exponential function), is o-minimal. (With a theorem of Khovanskii~\cite{22} together with the result of Wilkie~\cite{38} that the theory of $\mathbb{R}_{\textnormal{exp}}$ is model complete we obtain o-minimality.)
	\item [(3)] $\mathbb{R}_{\textnormal{an}}$\index{$\mathbb{R}_{\textnormal{an}}$}, the structure generated on the real field by the restricted analytic functions, is o-minimal (see~\cite{10} together with a theorem of Gabrielov~\cite{16} that the complement of a globally subanalytic function is again globally subanalytic; see also \cite{25}, Section 1 for an analytic proof). The sets definable in $\mathbb{R}_{\textnormal{an}}$ are precisely the globally subanalytic ones (see~\cite{11}).
	\item [(4)] $\mathbb{R}_{\textnormal{an,exp}}$\index{$\mathbb{R}_{\textnormal{an,exp}}$}, the structure generated by $\mathbb{R}_{\textnormal{an}}$ and $\mathbb{R}_{\textnormal{exp}}$, is o-minimal (see for example~\cite{13} for a model-theoretic proof and~\cite{25}, Section 2 for an analytic proof).
\end{itemize}
\end{example}

Especially the o-minimal structure $\mathbb{R}_{\textnormal{an,exp}}$ is very important, because a huge class of elementary functions are definable in there. With the following fact we understand how $\mathbb{R}_{\textnormal{an,exp}}$-definable functions look like.

\begin{fact}[Van den Dries et al.~\cite{13} Corollary 4.7]
	\label{12}
	Let $n \in \mathbb{N}$ and let $f:\mathbb{R}^n \to \mathbb{R}$ be definable in $\mathbb{R}_{\textnormal{an,exp}}$. Then there are $s \in \mathbb{N}$ and $\mathcal{L}_{\textnormal{an}}(\exp,\log)$-terms $t_1, \ldots ,t_s$ such that for every $x \in \mathbb{R}^n$ there is $j \in \{1, \ldots ,s\}$ with $f(x)=t_j(x)$. So $f$ is \textbf{piecewise given by $\mathcal{L}_{\textnormal{an}}(\exp,\log)$-terms}.\index{piecewise}
\end{fact}

{\textbf{This means:} Definable functions are piecewise compositions of globally subanalytic functions, the global logarithm and the global exponential. So we call them also \textbf{log-exp-analytic} in some context. \index{log-exp-analytic}\\

Subsequently, we exhibit the tame geometric behavior of o-minimal structures. An essential role in this context plays the cell decomposition theorem which states that we can decompose every definable set in finitely many disjoint subsets of a special form called cells. This concept is very helpful for technical proofs: The idea is to do calculation on every single cell at first and finally obtain the result by considering all of them.

For the rest of the section ''definable'' means always ''definable in $\mathcal{M}$'' if the underlying o-minimal structure $\mathcal{M}$ on the real field is clear from the context.

\begin{definition}
\label{13}
Let $\mathcal{M}$ be an o-minimal structure on the real field. By induction on $n$ we define a \textbf{definable cell $C \subset \mathbb{R}^n$} \index{definable!cell} as follows.

$n=1$: $C$ is either a singleton or an open interval.

$n \to n+1$: $C$ has one of the following form.
\begin{itemize}
	\item [(i)] $C=\textnormal{graph}(f)$,
	\item [(ii)] $C=\{(t,x) \in B \times \mathbb{R} \mid f(t)<x<g(t)\}:= \textnormal{}]f,g[_B$,
	\item [(iii)] $C=\{(t,x) \in B \times \mathbb{R} \mid f(t)<x<\infty\}:= \textnormal{}]f,\infty[_B$,
	\item [(iv)] $C=\{(t,x) \in B \times \mathbb{R} \mid -\infty<x<f(x)\}:= \textnormal{}]-\infty,f[_B$,
	\item [(v)] $C=\{(t,x) \in B \times \mathbb{R} \mid -\infty<x<\infty\}:= \textnormal{}]-\infty,\infty[_B$,
\end{itemize}
where $B \subset \mathbb{R}^n$ is a definable cell called the \textbf{base}\index{base} of $C$ in $\mathbb{R}^n$ and $f,g:B \to \mathbb{R}$ are continuous definable functions with additionally $f<g$ in ii).
\end{definition}

\begin{figure}[h]
	\begin{center}\includegraphics[width=9cm,height=6cm,keepaspectratio]{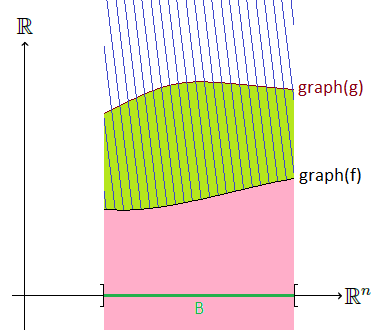}\end{center}
	\caption{Different cell types in $\mathbb{R}^{n+1}$ on the base $B$}
\end{figure}

\begin{definition}
\label{14}
Let $\mathcal{M}$ be an o-minimal structure on the real field. Let $m \in \mathbb{N}$ and $A \subset \mathbb{R}^m$ be non-empty and definable. A finite partition $\mathcal{C}$ of $A$ into definable cells is called a \textbf{definable cell decomposition $\mathcal{C}$ of $A$}\index{definable!cell decomposition}.
\end{definition}

\begin{fact}[Van den Dries~\cite{11}, Chapter 3]
\label{15}
Let $\mathcal{M}$ be an o-minimal structure on the real field. Let $n \in \mathbb{N}$ and $A \subset \mathbb{R}^n$ be definable. The following holds.
	\begin{itemize}
		\item [(1)] There is a definable cell decomposition $\mathcal{C}$ of $A$.
		\item [(2)] Let $m \in \mathbb{N}$. If $f:A \to \mathbb{R}^m$ is a definable function then the definable cell composition $\mathcal{C}$ of $A$ can be chosen in a way that $f|_C$ is continuous for every $C \in \mathcal{C}$.	
	\end{itemize}
\end{fact}

Since a definable cell is definably connected (i.e. cannot be written as the disjoint union of two definable open subsets, see~\cite{11}) an immediate consequence from the cell decomposition theorem is that every definable set is the union of finitely many connected components and these are again definable. A connected component $X_{\text{comp}}$ of a non-empty definable set $X \subset \mathbb{R}^m$ is a connected subset of $X$ which is maximal with respect to inclusion (i.e. if $X_{\text{comp}} \subset Y \subset X$ and $Y$ is connected then $X_{\text{comp}}=Y$)\index{definably connected}. Note also that the notion "definably connected" is equivalent to "definably path connected" in our context since we work with o-minimal structures on the reals which all contain the addition (see for example Chapter~6,~(3.2) in~\cite{11} for a formal proof). A set $X \subset \mathbb{R}^m$ is definably path connected\index{definably path connected} if for every $x_1,x_2 \in X$ there is a definable continuous $f:[0,1] \to X$ with $f(0)=x_1$ and $f(1)=x_2$. We refer to the book of Van den Dries~\cite{11} or Coste~\cite{9} for more on the general properties of o-minimal structures. \\

To close this section we pick up results on definability of parameterized integrals over the reals which will be important for our purposes in Section~8 below. The class of constructible functions discussed by Cluckers et al.~\cite{7} plays a major role in this context.

\begin{definition}
\label{16}
Let $n \in \mathbb{N}$. Let $A \subset \mathbb{R}^n$ be a globally subanalytic set. A function $f:A \to \mathbb{R}$ is called \textbf{constructible}\index{constructible} if there are $k \in \mathbb{N}$, $l_i \in \mathbb{N}$ for $i \in \{1, \ldots ,k\}$ and globally subanalytic functions $g_i:A \to \mathbb{R}$ and $h_{ij}:A \to \mathbb{R}_{>0}$ for $i \in \{1, \ldots ,k\}$ and $j \in \{1, \ldots ,l_i\}$ such that 
$$f(x)=\sum_{i=1}^k g_i(x) \prod_{j=1}^{l_i} \log(h_{ij}(x))$$
for every $x \in A$. 
\end{definition}

In other words, a constructible function is the sum of products of globally subanalytic functions and their logarithm. 

\begin{remark}
\label{17}
The following properties hold.
\begin{itemize}
\item[(1)] A globally subanalytic function is constructible. 
\item[(2)] A constructible function is log-analytic (compare with Lion et al.~\cite{25} or Definition 3.1 below for the notion of a log-analytic function).
\end{itemize}
\end{remark}

Let $n \in \mathbb{N}_0$ and $m \in \mathbb{N}$. Let $u$ range over $\mathbb{R}^n$ and $v$ over $\mathbb{R}^m$.

\begin{definition}
\label{18}
Let $f:\mathbb{R}^n \times \mathbb{R}^m \to \mathbb{R}$ be globally subanalytic. We set
$$\textnormal{Fin}(f):=\Bigl\{u \in \mathbb{R}^n \mid \int_{\mathbb{R}^m} \vert{f(u,v)}\vert dv < \infty \Bigl\}.$$
\index{$\textnormal{Fin}(f)$}
\end{definition}

Comte, Lion and Rolin showed that parameterized integrals of globally subanalytic functions are constructible.

\begin{fact}[Comte et al.~\cite{8} Theorem 1]
\label{19}
Let $f:\mathbb{R}^n \times \mathbb{R}^m \to \mathbb{R}$ be globally subanalytic. The following holds.
	\begin{itemize}
		\item [(1)] The set $\textnormal{Fin}(f)$ is globally subanalytic.
		\item [(2)] The function
		$$h:\textnormal{Fin}(f) \to \mathbb{R}, u \mapsto \int_{\mathbb{R}^m} f(u,v)dv,$$
		is constructible. 
	\end{itemize}
\end{fact}

Cluckers and Miller generalized this result on constructible functions.

\begin{fact}[Cluckers et al.~\cite{6} Theorem 2.5]
\label{20}
The class of constructible functions is stable under parametric integration.
\end{fact}

This means that the class of constructible functions is the smallest subclass of all $\mathbb{R}_{\textnormal{an,exp}}$-definable functions which contains all the globally subanalytic ones and is stable under parametric integration. Further and more deep theory on integration of constructible functions where loci of integrability and Lebesque classes of constructible functions are studied can be found in Cluckers et al.~\cite{5} and Cluckers et al.~\cite{6}.
However, it is in general not possible to extend this result within $\mathbb{R}_{\textnormal{an,exp}}$ beyond the constructible setting as the following fact indicates.

\begin{fact}[Van den Dries et al.~\cite{12} Theorem 5.11]
\label{21}
The function $\mathbb{R} \to \mathbb{R}, u \mapsto \int_0^u e^{-v^2} dv,$ is not definable in $\mathbb{R}_{\textnormal{an,exp}}$.
\end{fact}

An important consequence of this fact is that even parameterized integrals of log-analytic functions are in general not definable in $\mathbb{R}_{\textnormal{an,exp}}$ which we define in Definition~\ref{3.1} below.

\begin{example}
\label{22}
Consider
$$f:\mathbb{R} \times \mathbb{R} \to \mathbb{R}_{\geq 0}, (u,v) \mapsto \left\{\begin{array}{lll} \frac{1}{2v^2 \sqrt{\log(v)}},&& 1<v<u, \\
&\textnormal{if}&\\
0,&& \textnormal{else.} \end{array}\right.$$
Then $\mathbb{R} \to \mathbb{R}, u \mapsto \int_{\mathbb{R}} f(u,v) dv,$ is not definable in $\mathbb{R}_{\textnormal{an,exp}}$, i.e. $f$ is not constructible.
\end{example}

\begin{proof}
We obtain for $1 < u$
\begin{eqnarray*}
	\int_0^{\sqrt{\log(u)}} e^{-v^2} dv &=& \int_0^{\log(u)}\frac{e^{-v}}{2 \sqrt{v}} dv\\
	&=& \int_1^u \frac{1}{2v^2 \sqrt{\log(v)}} dv.\\
\end{eqnarray*}
For $u \in \mathbb{R}$ let 
$$g(u):=\int_{\mathbb{R}} f(u,v) dv.$$
We have $g(u)=0$ for $u \leq 1$ and
$$g(u) = \int_0^{\sqrt{\log(u)}} e^{-v^2} dv \leq \sqrt{\pi}$$
for $u > 1$. So we see that $\textnormal{Fin}(f)=\mathbb{R}$. By Fact~\ref{21} we obtain that $\mathbb{R} \to \mathbb{R}, u \mapsto \int_0^u e^{-v^2} dv,$ is not definable in $\mathbb{R}_{\textnormal{an,exp}}$. Therefore $\mathbb{R} \to \mathbb{R}, u \mapsto \int_{\mathbb{R}} f(u,v) dv,$ is also not definable in $\mathbb{R}_{\textnormal{an,exp}}$.
\end{proof}

Note that the antiderivative of a continuous function $f:\mathbb{R} \to \mathbb{R}$ can be considered as a special parameterized integral:
$$\int_0^x f(t)dt = \int_\mathbb{R} \tilde{f}(x,t)dt$$
where $\tilde{f}(x,t)=f(t) \mathbbm{1}_{[0,x]}(t)$. For an o-minimal structure $\mathcal{M}$ on the real field consider $\mathcal{P}(\mathcal{M})$, the Pfaffian closure of $\mathcal{M}$, which is an o-minimal expansion of $\mathcal{M}$. It goes beyond the scope of this paper to give the exact definition. For more details we refer to Speissegger~\cite{35}.

\begin{fact}[Speissegger~\cite{35} Corollary on p.1]
\label{23}
Let $n=1$ and $m=1$. Let $\mathcal{M}$ be an o-minimal structure on the real field and let $f:\mathbb{R} \to \mathbb{R}$ be definable in $\mathcal{M}$. Then the function $F:\mathbb{R} \to \mathbb{R}, x \mapsto \int_0^x f(t)dt$, is definable in $\mathcal{P}(\mathcal{M})$.
\end{fact}

This shows that the parameterized integrals from Fact~\ref{21} and Example~\ref{22} are indeed definable in $\mathcal{P}(\mathbb{R}_{\textnormal{an,exp}})$. 
But questions about this matter regarding definability of such integrals in more than one variable are quite open. (See for example Kaiser et al.~\cite{20} for definability results and open questions regarding parameterized exponential integrals given by the Brownian motion on globally subanalytic sets.) 

\section{Global Complexification}

This section is based on Kaiser~\cite{17} where $\mathcal{M}$ denotes a fixed o-minimal structure on the real field. In this section the expression "definable" means always "definable in $\mathcal{M}$" if not otherwise mentioned.

\begin{definition}
\label{2.1}
We say that $\mathcal{M}$ has \textbf{complexification}\index{complexification} if the following holds. Let $l \in \mathbb{N}$. Let $X \subset \mathbb{R}^l$ be open and let $f:X \to \mathbb{R}$ be a definable real analytic function. Then for every $x \in X$ there is a definable open neighborhood $U \subset \mathbb{R}^l$ of $x$, an open set $V \subset \mathbb{C}^l$ with $U \subset V$, and a definable holomorphic function $F:V \to \mathbb{C}$ such that $F|_U=f|_U$.
\end{definition}

So $\mathcal{M}$ has complexification if every definable real analytic function has locally a definable holomorphic extension.

\begin{example}
\label{2.2}
The following holds.
\begin{itemize}
	\item[(1)] Let $\mathcal{M}$ be an o-minimal expansion of $\mathbb{R}_{\textnormal{an}}$ (e.g. $\mathcal{M}=\mathbb{R}_{\textnormal{an,exp}}$). Then $\mathcal{M}$ has complexification.
	\item[(2)] The o-minimal structure $\mathbb{R}_{\textnormal{exp}}$ does not have complexification.
\end{itemize}
\end{example}

\begin{proof}
(1): This follows simply by the Taylor series expansion.\\
(2): We denote by $f$ the real exponential function which is definable in $\mathbb{R}_{\textnormal{exp}}$. Consider
$$\exp:\mathbb{C} \to \mathbb{C}, x+iy \mapsto f(x)(\cos(y)+i\sin(y)),$$
which is holomorphic. Let $x \in \mathbb{R}$ and let $V$ be an open ball in $\mathbb{C}$ around $x$. By the identity theorem we see that $\exp|_V$ is the unique holomorphic extension of $f|_{V \cap \mathbb{R}}$ on $V$. But Bianconi showed in~\cite{1} that no restriction of the global sine function on a non-empty open interval is definable in $\mathbb{R}_{\textnormal{exp}}$.
\end{proof}

We mention that such questions about complexification are only poorely answered. A field of open problems in this context are for example o-minimal structures which are generated by so-called convergent Weiherstrass systems. We refer to Miller~\cite{26} for the details.

\begin{definition}
\label{def:complexification-non-open}
Let $l \in \mathbb{N}$, $X \subset \mathbb{R}^l$ and let $f:X \to \mathbb{R}$ be a definable function\index{global complexification}. Then a definable function $F:V \to \mathbb{C}$ where $V \subset \mathbb{C}^l$ is open is called \textbf{global complexification} of $f$ if $X \subset V$, $F$ is holomorphic,  and $F|_X=f$.
\end{definition}

\begin{definition}
\label{2.3}
Let $l \in \mathbb{N}$. We say that $\mathcal{M}$ has \textbf{$l$-ary global complexification}\index{global complexification} if the following holds. Let $X \subset \mathbb{R}^l$ be open and let $f:X \to \mathbb{R}$ be a definable real analytic function. Then f has a global complexification.
\end{definition}

So $\mathcal{M}$ has $l$-ary global complexification for $l \in \mathbb{N}$ if every definable real analytic function $F:X \to \mathbb{R}$ where $X \subset \mathbb{R}^l$ is open has a definable holomorphic extension.

\begin{remark}
\label{2.4}
Let $l \in \mathbb{N}$. If $\mathcal{M}$ has $(l+1)$-ary global complexification then it has $l$-ary global complexification.
\end{remark}

By Fact~\ref{12} every definable function in $\mathbb{R}_{\textnormal{an,exp}}$ is piecewise given by $\mathcal{L}_{\textnormal{an}}(\exp,\log)$-terms. A consequence is the following.

\begin{fact}[Kaiser~\cite{17} Theorem C]
\label{2.5}
The o-minimal structure $\mathbb{R}_{\textnormal{an,exp}}$ has unary global complexification.
\end{fact}

\begin{definition}
\label{2.6}
We say that $\mathcal{M}$ has \textbf{global complexification}\index{global complexification} if $\mathcal{M}$ has $l$-ary global complexification for every $l \in \mathbb{N}$.
\end{definition}

So $\mathcal{M}$ has global complexification if every definable real analytic function has a definable holomorphic extension.

\begin{remark}
\label{2.7}
If $\mathcal{M}$ has global complexification then it has complexification. Consequently, $\mathbb{R}_{\textnormal{exp}}$ does not have global complexification.
\end{remark}

It is not known to us whether there is an o-minimal structure which has complexification but no global complexification. The problem is that there is no general concept how definable functions in o-minimal structures look like. But for the case that $\mathcal{M}$ is the structure on the real field all definable real analytic functions are so-called Nash functions. For this class of functions there are strong results like the implicit function theorem or the Artin-Mazur description (see Bochnak et al.~\cite{3}, Chapter 8) which are sufficient to show the following (see also Shiota~\cite{34}, Chapter I.6.7.).

\begin{fact}[Kaiser~\cite{18} Theorem B]
\label{2.8}
The structure on the real field $\mathbb{R}$ has global complexification.
\end{fact}

For the o-minimal structure $\mathbb{R}_{\textnormal{an}}$ the situation is more complicated since much more functions are definable in there than in $\mathbb{R}$. Fortunately there are deep geometrical results for globally subanalytic functions. These are the preparation theorems of Lion and Rolin (see \cite{25}) which give a nice representation of a globally subanalytic function in one variable and are precise enough to establish global complexification: The idea is to compute a suitable holomorphic extension of a parametric family of real analytic globally subanalytic one-variable functions at first and then do a non-trivial induction on the number of variables. Therefore in~\cite{17} Tobias Kaiser set up the notion of parametric global complexification.

For $n \in \mathbb{N}_0$ and $l \in \mathbb{N}$ let $t$ range over $\mathbb{R}^n$ and $x$ over $\mathbb{R}^l$. Let $z$ range over $\mathbb{C}^l$. Here $x$ and $z$ are serving as tuples of independent
variables of families of functions parameterized by $t$.

\begin{definition}
	\label{def:global-complex-parametric}
	Let $l \in \mathbb{N}$, $n \in \mathbb{N}_0$, $X \subset \mathbb{R}^n \times \mathbb{R}^l$ and let $f:X \to \mathbb{R}, (t,x) \mapsto f(t,x),$ be definable. Then we call a definable function $F:Z \to \mathbb{C}, (t,z) \mapsto F(t,z),$ where $Z \subset \mathbb{R}^n \times \mathbb{C}^l$ is such that $Z_t=\{z \in \mathbb{C}^l \mid (t,z) \in Z\}$ is open for every $t \in \mathbb{R}^n$, an \textbf{$l$-ary parametric global complexification} of $f$ if $X \subset Z$, $F_t$ is holomorphic for all $t \in \mathbb{R}^n$, and $F|_X=f$.\index{global complexification!parametric}
\end{definition}

\begin{definition}
\label{2.9}
Let $l \in \mathbb{N}$. We say that $\mathcal{M}$ has \textbf{$l$-ary parametric global complexification}\index{global complexification!parametric} if the following holds. Let $n \in \mathbb{N}_0$, $X \subset \mathbb{R}^n \times \mathbb{R}^l$ be definable such that $X_t$ is open for all $t \in \mathbb{R}^n$. Let $f:X \to \mathbb{R}, (t,x) \mapsto f(t,x),$ be definable such that $f_t:X_t \to \mathbb{R}$ is real analytic for all $t \in \mathbb{R}^n$. Then there is an $l$-ary parametric global complexification $F:Z \to \mathbb{C}, (t,z) \mapsto F(t,z),$ of $f$.
\end{definition}

\begin{remark}
\label{2.10}
Let $l \in \mathbb{N}$.
\begin{itemize}
	\item[(1)]
	If $\mathcal{M}$ has $(l+1)$-ary parametric global complexification then it has $l$-ary parametric global complexification.
	
	\item[(2)]
	If $\mathcal{M}$ has $l$-ary parametric global complexification then it has $l$-ary global complexification.
\end{itemize}
\end{remark}

\begin{definition}
\label{2.11}
We say that $\mathcal{M}$ has \textbf{parametric global complexification}\index{global complexification!parametric} if $\mathcal{M}$ has $l$-ary parametric global complexification for every $l \in \mathbb{N}$.
\end{definition}

\begin{remark}
\label{2.12}
If there is an open non-empty subset of $\mathbb{R}^l$ that is contained in every $X_t$, there may be no $Z \subset \mathbb{R}^n \times \mathbb{C}^l$ such that there is an open and non-empty subset of $\mathbb{C}^l$ that is contained in every $Z_t$ as shown by the example 
$$\mathbb{R}_{>0} \times \mathbb{R} \to \mathbb{R}, (t,x) \mapsto 1/(t^2 + x^2).$$
\end{remark}

\begin{remark}
\label{2.13}
If $\mathcal{M}$ has parametric global complexification then it has global complexification.
\end{remark}

When establishing results about $l$-ary parametric global complexification by an induction on the independent variables there is the problem to ''lift'' the open property of $Z_t \subset \mathbb{C}^l$ from $l=1$ to higher $l$ since "fiberwise open" does not imply open. So we set up the concept of high parametric global complexification (see also~\cite{17} for the ideas).
For $m \in \mathbb{N}$ let $u$ range over $\mathbb{R}^m$. In the following $(u,x)$ and $(u,z)$ are serving as tuples of independent variables of families of functions parameterized by $t$. 

\begin{definition}
\label{2.14}
We define the following.
\begin{itemize}
	\item [(a)] Let $l \in \mathbb{N}$. We say that $\mathcal{M}$ has \textbf{$l$-ary high parametric global complexification}\index{global complexification!high parametric} if the following holds. Let $n,m \in \mathbb{N}_0$ and $X \subset \mathbb{R}^n \times \mathbb{R}^m \times \mathbb{R}^l$ be definable such that $X_t=\{(u,x) \in \mathbb{R}^m \times \mathbb{R}^l \mid (t,u,x) \in X\}$ is open in $\mathbb{R}^m \times \mathbb{R}^l$ for every $t \in \mathbb{R}^n$. Let $f:X \to \mathbb{R}$ be definable such that $f_t:X_t \to \mathbb{R}, (u,x) \mapsto f(t,u,x),$ is real analytic for every $t \in \mathbb{R}^n$. Then there is an $l$-ary parametric global complexification $F:Z \to \mathbb{C}, (t,u,z) \mapsto F(t,u,z),$ of $f$ such that $Z_t=\{(u,z) \in \mathbb{R}^m \times \mathbb{C}^l \mid (t,u,z) \in Z\}$ is open in $\mathbb{R}^m \times \mathbb{C}^l$ 
	for $t \in \mathbb{R}^n$. We call $F$ a \textbf{high parametric global complexification of $f$ with respect to $(u,x)$}.
	
	\item [(b)] We say that $\mathcal{M}$ has \textbf{high parametric global complexification}\index{global complexification!high parametric} if $\mathcal{M}$ has $l$-ary high parametric global complexification for every $l \in \mathbb{N}$.    
\end{itemize}
\end{definition}

\begin{remark}
\label{2.15}
If $\mathcal{M}$ has high parametric global complexification then it has parametric global complexification and therefore global complexification.
\end{remark}

\begin{fact}[Kaiser~\cite{17} Theorem 2.10]
\label{2.16}
The o-minimal structure $\mathbb{R}_\textnormal{an}$ has high parametric global complexification. Therefore $\mathbb{R}_\textnormal{an}$ has parametric global complexification and global complexification.
\end{fact}

Finally we give a small example to illustrate the difference between parametric global complexification from Definition~\ref{2.9} and high parametric global complexification from Definition~\ref{2.14}. Let $\mathcal{M}=\mathbb{R}_{\textnormal{an,exp}}$ for the rest of this section.

\begin{example}
\label{2.17}
Let $X:=\mathbb{R} \times \mathbb{R}$ and consider the definable function
$$f:X \to \mathbb{R}, \textnormal{ } (u,x)\mapsto 
\left\{\begin{array}{lll} \exp(\frac{x}{u}),&& u > 0, \\
&\mbox{if}&\\
0,&& u \leq 0. \end{array}\right.$$
Let
$$Z:=\{(u,z) \in \mathbb{R} \times \mathbb{C} \mid \vert{\textnormal{Im}(z)}\vert <\pi u \textnormal{ if } u>0\}.$$
Then the function
$$G:Z \to \mathbb{C}, (u,z) \mapsto 
\left\{\begin{array}{lll} \exp(\frac{z}{u}),&& u > 0, \\
&\mbox{if}&\\
0,&& u \leq 0, \end{array}\right.$$
is a unary parametric global complexification of $f$.
\end{example}

\begin{proof}
	Note that $Z_u$ is open and $G_u$ is holomorphic for every $u \in \mathbb{R}$. Further we have for $u \in \mathbb{R}_{>0}$ and $x,y \in \mathbb{R}$ with $x+iy \in Z_u$ 
	\[
	G(u,x+iy) = \exp((x+iy)/u) = \exp(x/u)(\cos(y/u)+i\sin(y/u)).
	\]
	which implies definability of $G$ since $-\pi<y/u<\pi$. 
\end{proof}

\begin{lemma}
\label{2.18}
Let $X$ and $f$ be as in Example~\ref{2.17}. Let $Z \subset \mathbb{R} \times \mathbb{C}$ be definable and let $F:Z \to \mathbb{C}$ be a unary parametric global complexification of $f$. Then there is no $r>0$ such that $]0,r[ \textnormal{} \times (]0,r[ + i]0,r[) \subset Z$, i.e. $F$ is no unary high parametric global complexification of $f$ with respect to $(u,x)$. 
\end{lemma}

\begin{proof}
Note that $Z_u$ is open for all $u \in \mathbb{R}$. Assume the contrary. Fix $r>0$ such that 
$$Q:= \textnormal{} ]0,r[ \times \textnormal{}(]0,r[ \textnormal{}+ i \textnormal{} ]0,r[) \subset Z.$$ 
By the identity theorem we have for $(u,z) \in Z$ that $F(u,z)=\exp(z/u)$ if $u>0$ and $F(u,z)=0$ otherwise. Note that $Q$ is definable, but $F|_Q$ is not definable, because for every $(u,x+iy) \in Q$ we have
$$F(u,x+iy)=\exp(x/u)(\cos(y/u) + i\sin(y/u))$$
and the global sine function (or cosine function) is not definable at infinity, a contradiction.
\end{proof}

Nevertheless it is absolutely possible that $\mathbb{R}_{\textnormal{an,exp}}$ has global complexification since the underlying function must be real analytic on its domain and $f$ in Example~\ref{2.17} is not real analytic at zero.

\section{$\mathbb{R}_{\textnormal{an,exp}}$-Definable Functions}

\subsection{Log-Analytic Functions and the Exponential \\ Number}

We fix $m \in \mathbb{N}$ and a definable $X \subset \mathbb{R}^m$.

\begin{definition}[Lion et al.~\cite{25}]
\label{3.1}
Let $f:X \to \mathbb{R}$ be a function.\index{log-analytic}
\begin{itemize}
	\item [(a)] Let $r \in \mathbb{N}_0$. By induction on $r$ we define that $f$ is \textbf{log-analytic of order at most} $r$.

	\textbf{Base case}: The function $f$ is log-analytic of order at most $0$ if there is a partition $\mathcal{C}$ of $X$ into finitely many definable sets such that for $C \in \mathcal{C}$ there is a globally subanalytic function $F:\mathbb{R}^m \to \mathbb{R}$ such that $f|_C = F|_C$.
	
	\textbf{Inductive step}: The function $f$ is log-analytic of order at most $r$ if the following holds: There is a partition $\mathcal{C}$ of $X$ into finitely many definable sets such that for $C \in \mathcal{C}$ there are $k,l \in \mathbb{N}_{0}$, a globally subanalytic function $F:\mathbb{R}^{k+l} \to \mathbb{R}$, and log-analytic functions $g_1, \ldots ,g_k:C \to \mathbb{R}, h_1, \ldots ,h_l:C \to \mathbb{R}_{>0}$ of order at most $r-1$ such that
	$$f|_C=F(g_1, \ldots ,g_k,\log(h_1), \ldots ,\log(h_l)).$$
	
	\item[(b)] Let $r \in \mathbb{N}_0$. We call $f$ \textbf{log-analytic of order}\index{log-analytic} $r$ if $f$ is log-analytic of order at most $r$ but not of order at most $r-1$.
	
	\item[(c)] We call $f$ \textbf{log-analytic}\index{log-analytic} if $f$ is log-analytic of order $r$ for some $r \in \mathbb{N}_0$.
\end{itemize}
\end{definition}

\begin{definition}
\label{3.2}
Let $f:X \to \mathbb{R}$ be a function. Let $E$ be a set of positive definable functions on $X$.
\begin{itemize}
	\item [(a)] By induction on $e \in \mathbb{N}_0$ we define that $f$ has \textbf{exponential number at most $e$ with respect to $E$}\index{exponential!number}.
	
	{\bf Base Case}: The function $f$ has exponential number at most $0$ with respect to $E$ if $f$ is log-analytic.
	
	{\bf Inductive Step}: The function $f$ has exponential number at most $e$ with respect to $E$ if the following holds: There are $k,l \in \mathbb{N}_0$, functions $g_1, \ldots ,g_k:X \to \mathbb{R}$ and $h_1, \ldots ,h_l:X \to \mathbb{R}$ with exponential number at most $e-1$ with respect to $E$ function $F:\mathbb{R}^{k+l} \to \mathbb{R}$ such that
	$$f=F(g_1, \ldots ,g_k,\exp(h_1), \ldots ,\exp(h_l))$$
	and $\exp(h_1), \ldots ,\exp(h_l) \in E$.
	
	\item [(b)] Let $e \in \mathbb{N}_0$. We say that $f$ has \textbf{exponential number $e$ with respect to $E$}\index{exponential!number} if $f$ has exponential number at most $e$ with respect to $E$ but not at most $e-1$ with respect to $E$.
	
	\item [(c)] We say that $f$ \textbf{can be constructed from $E$}\index{constructed from} if there is $e \in \mathbb{N}_0$ such that $f$ has exponential number $e$ with respect to $E$. 
\end{itemize}
\end{definition}

\begin{example}
\label{3.3}
The function $f:\mathbb{R} \to \mathbb{R}, x \mapsto \log(x^2+\vert{x}\vert+1),$ is log-analytic (and even constuctible). The function $g:\mathbb{R} \to \mathbb{R}, x \mapsto \log(\exp(x+\exp(x))+1),$ has exponential number at most $2$ with respect to $E:=\{g_1,g_2\}$ where
$$g_1:\mathbb{R} \to \mathbb{R}_{>0}, x \mapsto \exp(x), \textnormal{ } g_2: \mathbb{R} \to \mathbb{R}_{>0}, x \mapsto \exp(x+\exp(x)).$$
\end{example}
	
\begin{remark}
\label{3.4}
Let $e \in \mathbb{N}_0$. Let $E$ be a set of positive definable functions on $X$.
\begin{itemize}
	\item[(1)] Let $f:X \to \mathbb{R}$ be a function with exponential number at most $e$ with respect to $E$. Then $\exp(f)$ has exponential number at most $e+1$ with respect to $E\cup \{\exp(f)\}$.
	\item[(2)] Let $s \in \mathbb{N}_0$. Let $f_1, \ldots ,f_s:X \to \mathbb{R}$ be functions with exponential number at most $e$ with respect to $E$ and let $F:\mathbb{R}^s \to \mathbb{R}$ be log-analytic. Then $F(f_1, \ldots ,f_s)$ has exponential number at most $e$ with respect to $E$. 
\end{itemize}
\end{remark}

See also~\cite{27}, Section 1 for Definition~\ref{3.1}, Remark~\ref{3.4} and also for some more facts and examples in this context.

\begin{remark}
\label{3.5}
Let $X_1,X_2 \subset \mathbb{R}^m$ be definable and disjoint. Let $X = X_1 \cup X_2$. For $j \in \{1,2\}$ let $E_j$ be a set of positive definable functions on $X_j$ and $f_j:X_j \to \mathbb{R}$ be a function. Let $e \in \mathbb{N}_0$ be such that $f_j$ has exponential number at most $e$ with respect to $E_j$ for $j \in \{1,2\}$. Let 
$$E:=\{g \mid g:X \to \mathbb{R} \textnormal{ is function with } g|_{X_j} \in E_j \textnormal{ for } j \in \{1,2\}\}.$$
Then 
$$f:X \to \mathbb{R}, x \mapsto \left\{\begin{array}{ll} f_1(x) , & x \in X_1,  \\
f_2(x), & x \in X_2, \end{array}\right.$$
has exponential number at most $e$ with respect to $E$.
\end{remark}

We need the concept of the exponential number also in the complex setting.\\

For the rest of Section~3.1 we fix a definable $Z \subset \mathbb{C}^m$. Let $z:=(z_1, \ldots ,z_m)$ range over $\mathbb{C}^m$.  Let $x:=(x_1, \ldots ,x_m)$ and $y:=(y_1, \ldots ,y_m)$ range over $\mathbb{R}^m$ such that $z=x+iy$. We fix a function $F:Z \to \mathbb{C}$ with $F(z)=u(x,y)+iv(x,y)$ where $u$ is the real and $v$ the imaginary part of $F$ considered as real functions.

\begin{definition}
\label{3.6}
We call $F$ log-analytic if $u$ and $v$ are log-analytic\index{log-analytic}.  
\end{definition}

\begin{lemma}
\label{3.7}
The complex logarithm function 
$$\log: \mathbb{C}^- \to \mathbb{R} \times \text{}]-\pi,\pi[, w \mapsto \log(w)=\log(\vert{w}\vert) + i \arg(w),$$ is log-analytic.
\end{lemma}

\begin{proof}
Note that for every $w \in \mathbb{C}^-$
$$\log(w)=\log(\vert{w}\vert) + i \arg(w).$$
We see that 
$$\mathbb{R}^2 \setminus \{0\} \to \mathbb{R}, (u,v) \mapsto \log(\sqrt{u^2+v^2}),$$
is log-analytic and that
$$\mathbb{R}^2 \setminus \{(a,0) \mid a \leq 0\} \to \mathbb{R}, (u,v) \mapsto \arg(u+iv),$$
is log-analytic (since it is globally subanalytic).
\end{proof}

Note that for a function $h:\mathbb{C} \to \mathbb{C}^-$ the composition 
$$\log(h):\mathbb{C} \to \mathbb{R} \times \text{} ]-\pi,\pi[, z \mapsto \log(\vert{h(z)}\vert) + i \text{arg}(h(z)),$$
is well-defined. Hence we are able to define the following.

\begin{definition}
\label{3.8}
Let $E$ be a set of definable functions on $Z$ with values in $\mathbb{C}^-$. We say that $F$ can be constructed from $E$ \index{constructed from} if $u$ and $v$ can be constructed from
$$E^{\text{Re}}:=\{\exp(\textnormal{Re}(G)) \mid G \in \log(E)\}$$
where $\log(E):=\{\log(H) \mid H \in E\}$.
\end{definition}

\begin{example}
	\label{ex:constructed-complex}
	 Let $Z:=\{z \in \mathbb{C} \mid \vert{\text{Im(z)}\vert}<\pi \text{ and } \vert{\text{Im}(z+\exp(z))}\vert<\pi\}$. The function $G:Z \to \mathbb{C}, z \mapsto \log(\exp(z+\exp(z))+1),$ can be constructed from $E:=\{G_1,G_2\}$ where
	\[
	G_1:Z \to \mathbb{C}^-, z \mapsto \exp(z), \textnormal{ } G_2: Z \to \mathbb{C}^-, z \mapsto \exp(z+\exp(z)).
	\]
\end{example}

\begin{proof}
	Note that $\exp(x+iy)=\exp(x)(\cos(y)+i\sin(y))$ for $x,y \in \mathbb{R}$. 
	Let $\hat{G}:Z \to \mathbb{C}, z \mapsto \exp(z+\exp(z))+1$. Then $\hat{G}=g_1+ig_2$ for 
	$$g_1:Z \mapsto \mathbb{R}, x+iy \mapsto \exp(x+\exp(x)\cos(y)) \cos(y+\exp(x) \sin(y))+1,$$
	and 
	$$g_2:Z \mapsto \mathbb{R}, x+iy \mapsto \exp(x+\exp(x)\cos(y)) \sin(y+\exp(x) \sin(y)).$$
	Since $E^{\text{Re}} = \{H_1,H_2\}$ for 
	$$H_1:Z \mapsto \mathbb{R}, z \mapsto \exp(\text{Re}(z)) = \exp(x),$$ 
	and 
	$$H_2:Z \mapsto \mathbb{R}, z \mapsto \exp(\text{Re}(z+\exp(z))) = \exp(x+\exp(x)\cos(y)),$$
	where $z=x+iy$, we see with Definition~\ref{3.2}(c) that $g_1$ and $g_2$ can be constructed from $E^{\text{Re}}$. Since 
	$$G = \log(g_1+ig_2) = \log \Bigl(\sqrt{g_1^2+g_2^2}\Bigl) + i\arg(g_1+ig_2)$$
	we see that $\text{Re}(G)$ and $\text{Im}(G)$ can also be constructed from $E^{\text{Re}}$ which shows that $G$ can be constructed from $E$.
\end{proof}

Definition~\ref{3.8} is a generalization of Definition~\ref{3.2}(c) in the following sense.

\begin{lemma}
\label{3.9}
Let $E$ be a set of positive real valued definable functions on $Z$. Assume that $F$ takes only real values. Then Definition~\ref{3.8} coincides with Definition~\ref{3.2}(c) if one considers $F$ and every $G \in E$ as a real function.
\end{lemma}

\begin{proof}
We have $F=u$, $v=0$ and $E^{\textnormal{Re}}=\{\exp(\textnormal{Re}(G)) \mid G \in \log(E)\} = E$ since $\textnormal{Re}(G)=G$ for every $G \in \log(E)$.
\end{proof}

We can formulate Remark~\ref{3.4} and Remark~\ref{3.5} also for the complex setting and give short proofs.

\begin{lemma}
\label{3.10}
Let $e \in \mathbb{N}_0$. Let $E$ be a set of definable functions on $Z$ with values in $\mathbb{C}^-$.
\begin{itemize}
	\item[(1)] Assume that $\vert{\textnormal{Im}(F)}\vert<\pi$ and that $F:Z \to \mathbb{C}$ can be constructed from $E$. Then $\exp(F)$ can be constructed from $D:=E \cup \{\exp(F)\}$.
	\item[(2)] Let $F_1, \ldots ,F_l:Z \to \mathbb{C}$ be functions which can be constructed from $E$ and let $G:\mathbb{C}^l \to \mathbb{C}$ be log-analytic. Then $G(F_1, \ldots ,F_l)$ can be constructed from $E$. 
\end{itemize}
\end{lemma}

\begin{proof}
(1): Note that $\exp(F)$ is definable, $D$ is a set of functions which takes only values in $\mathbb{C}^-$ and $D^{\text{Re}}=E^{\textnormal{Re}} \cup \{\exp(u)\}$. We have
$$\exp(F(x+iy))=\exp(u(x,y))(\cos(v(x,y))+i\sin(v(x,y)))$$
for $x+iy \in Z$. By Remark~{3.4}(1) $\exp(u)$ can be constructed from $E^{\textnormal{Re}} \cup \{\exp(u)\}$. 
Let $T \in \{\cos,\sin\}$. Then 
$$T^*:\mathbb{R} \to \mathbb{R}, x \mapsto \left\{\begin{array}{ll} T(x) , & \vert{x}\vert < \pi,  \\
0, & \textnormal{else}, \end{array}\right.
$$
is globally subanalytic. Therefore by Remark~\ref{3.4}(2) 
$$Z \mapsto \mathbb{R}, x+iy \mapsto \exp(u(x,y))T^*(v(x,y)),$$ 
can be constructed from $D^{\textnormal{Re}}$ since $F$ can be constructed from $E$ and therefore $v$ from $E^{\textnormal{Re}}$ by Definition~\ref{3.8}.

(2): Let $u_j$ be the real part and $v_j$ be the imaginary part of $F_j$ considered as real functions, i.e. 
$$F_j(x+iy)=u_j(x,y)+iv_j(x,y)$$
for $j \in \{1, \ldots ,l\}$ and $x+iy \in Z$. Let $\tilde{u}$ be the real part and $\tilde{v}$ be the imaginary part of $G$ considered as real functions. Then 
$$G(F_1, \ldots ,F_l)=\tilde{u}(u_1, \ldots ,u_l,v_1, \ldots ,v_l) + i\tilde{v}(u_1, \ldots ,u_l,v_1, \ldots ,v_l).$$ 
We are done with Remark~\ref{3.4}(2). 
\end{proof}

\begin{lemma}
	\label{lem:construction-complex}
	Let $Z_1,Z_2 \subset \mathbb{C}^m$ be definable and disjoint. Let $Z = Z_1 \cup Z_2$. For $j \in \{1,2\}$ let $E_j$ be a set of definable functions on $Z_j$ with values in $\mathbb{C}^-$ and let $F_j:Z_j \to \mathbb{C}$ be a function. Let 
	$$E:=\{G \mid G:Z \to \mathbb{C} \textnormal{ is function with } G|_{Z_j} \in E_j \textnormal{ for } j \in \{1,2\}\}.$$
	Then 
	$$F:Z \to \mathbb{C}, z \mapsto \left\{\begin{array}{ll} F_1(z) , & z \in Z_1,  \\
		F_2(z), & z \in Z_2, \end{array}\right.$$
	can be constructed from $E$.
\end{lemma}

\begin{proof}
	Let $u_j:=\text{Re}(F_j)$ and $v_j:=\text{Im}(F_j)$ for $j \in \{1,2\}$. By Definition~\ref{3.8} the functions $u_j$ and $v_j$ can be constructed from $(E_j)^{\text{Re}}=\{\exp(\text{Re}(G)) \mid G \in \log(E_j)\}$ and hence, by Remark~\ref{3.5}, 
	$$u:Z \to \mathbb{R}, z \mapsto \left\{\begin{array}{ll} u_1(z) , & z \in Z_1,  \\
		u_2(z), & z \in Z_2, \end{array}\right.$$
	can be constructed from 
	$$\hat{E}:=\{e^{\text{Re}(H)} \mid H:Z \to \mathbb{C} \textnormal{ is function with } H|_{Z_j} \in \log(E_j) \textnormal{ for } j \in \{1,2\}\}.$$
	Note that $\hat{E}=E^{\text{Re}}$ which proves the lemma since $v$ is handled analoguously. 
\end{proof}

\subsection{Restricted Log-Exp-Analytic Functions}

We establish a result on global complexification also in the parametric setting. Hence, we set up the concept of parametric families of restricted log-exp-analytic functions.

For this section let $n \in \mathbb{N}_0$ and $m \in \mathbb{N}$. Let $t$ range over $\mathbb{R}^n$ and $x$ over $\mathbb{R}^m$. We fix definable sets $C,X \subset \mathbb{R}^n \times \mathbb{R}^m$ with $C \subset X$. Suppose that $X_t$ is open for every $t \in \mathbb{R}^n$.  Here $x=(x_1, \ldots, x_m)$ is serving as the tuple of independent variables of families of functions parameterized by $t=(t_1, \ldots, t_n)$.

\begin{definition}
\label{3.11}
Let $f:C \to \mathbb{R}, (t,x) \mapsto f(t,x),$ be a function.
\begin{itemize}
	\item [(a)] We call $f$ \textbf{locally bounded in $x$ with reference set $X$}\index{locally bounded!with reference set} if the following holds. For $t \in \mathbb{R}^n$ and $w \in X_t$ there is an open neighborhood $U$ of $w$ in $X_t$ such that $U \cap C_t = \emptyset$ or $f_t|_{U \cap C_t}$ is bounded. 
	\item [(b)] Suppose that $C_t$ is open for every $t \in \mathbb{R}^m$. We call $f$ \textbf{locally bounded in $x$}\index{locally bounded} if $f$ is locally bounded in $x$ with reference set $C$.	
\end{itemize}
\end{definition}

\begin{remark}
\label{3.12}
The set of locally bounded functions in $x$ with reference set $X$ on $C$ is an $\mathbb{R}$-algebra with respect to pointwise addition and multiplication. 
\end{remark}

\begin{remark}
\label{3.13}
Let $C_1,C_2 \subset X$ be disjoint definable sets such that $C = C_1 \cup C_2$. Let $g_j:C_j \to \mathbb{R}$ be locally bounded in $x$ with reference set $X$ for $j \in \{1,2\}$. Then 
$$g:C \to \mathbb{R}, (t,x) \mapsto \left\{\begin{array}{ll} g_1(t,x) , & (t,x) \in C_1,  \\
g_2(t,x), & (t,x) \in C_2, \end{array}\right.$$
is locally bounded in $x$ with reference set $X$.
\end{remark}

\begin{definition}
\label{3.14}
Let $f:C \to \mathbb{R}, (t,x) \mapsto f(t,x),$ be a function.
\begin{itemize}
	\item [(a)]
	Let $e \in \mathbb{N}_0$. We say that $f$ is \textbf{restricted log-exp-analytic in $x$ of order (at most) $e$ with reference set $X$}\index{restricted log-exp-anal.!with reference set} if $f$ has exponential number (at most) $e$ with respect to a set $E$ of positive definable functions on $C$ such that every $g \in \log(E)$ is locally bounded in $x$ with reference set $X$.
	
	\item [(b)]
	We say that $f$ is \textbf{restricted log-exp-analytic in $x$ with reference set $X$}\index{restricted log-exp-anal.!with reference set} if $f$ can be constructed from a set $E$ of positive definable functions on $C$ such that every $g \in \log(E)$ is locally bounded in $x$ with reference set $X$.
\end{itemize}
\end{definition}

\begin{remark}
\label{3.15}
Let $l \in \mathbb{N}$. For $j \in \{1, \ldots ,l\}$ let $f_j:C \to \mathbb{R}$ be a function which is restricted log-exp-analytic in $x$ with reference set $X$. Let $F:\mathbb{R}^l \to \mathbb{R}$ be log-analytic. Then 
$$C \to \mathbb{R}, (t,x) \mapsto F(f_1(t,x), \ldots ,f_l(t,x)),$$
is restricted log-exp-analytic in $x$ with reference set $X$.
\end{remark}

\begin{definition}
\label{3.16}
A function $f:X \to \mathbb{R}, (t,x) \mapsto f(t,x),$ is called \textbf{restricted log-exp-analytic in $x$}\index{restricted log-exp-anal.} if $f$ is restricted log-exp-analytic in $x$ with reference set $X$.
\end{definition}

One of the main results is formulated in the non-parametric setting in which Definition~\ref{3.11} and~\ref{3.14} simplify as follows. 

\begin{definition}
\label{3.17}
Let $Y \subset \mathbb{R}^m$ be open. 
\begin{itemize}
\item[(a)] A function $g:Y \to \mathbb{R}, x \to g(x),$ is called \textbf{locally bounded}\index{locally bounded} if $g$ is locally bounded in $x$.
\item[(b)] A function $f:Y \to \mathbb{R}, x \mapsto f(x),$ is called \textbf{restricted log-exp-analytic}\index{restricted log-exp-anal.} if $f$ is restricted log-exp-analytic in $x$.
\end{itemize}
\end{definition}

(See also~\cite{28}, Section 2 for Definition~\ref{3.11} to Definition~\ref{3.17}.)

\begin{example}
\label{3.18}
The function $g:\mathbb{R} \to \mathbb{R}, x \mapsto \log(\exp(x+\exp(x))+1),$ is restricted log-exp-analytic (of order at most $2$). The function
$$h: \mathbb{R} \to \mathbb{R}, x \mapsto \left\{\begin{array}{ll} e^{-1/x}, & x > 0, \\
0, & \textnormal{ else, } \end{array}\right.$$
is not restricted log-exp-analytic, but the function $h|_{\mathbb{R}_{>0}}$ is. 
\end{example}

\begin{proof}
It is easy to see with Definition~\ref{3.2}(c) that the function $g$ can be constructed from $E:=\{h_1,h_2\}$ where $h_1:\mathbb{R} \to \mathbb{R}, x \mapsto \exp(x)$ and $h_2:\mathbb{R} \to \mathbb{R}, x \mapsto \exp(x+\exp(x))$. Clearly $h_1$ and $h_2$ are locally bounded. 

For a restricted log-exp-analytic function $f:\mathbb{R} \to \mathbb{R}$ the following holds for all sufficiently small positive $y$: Either $f(y)=0$ or there is $a \in \mathbb{R} \setminus \{0\}$, a non-negative integer $r \in \mathbb{N}_0$ and $q_0,\ldots,q_r \in \mathbb{Q}$ such that 
$\lim_{y \searrow 0} f(y)/h(y) = a$ where $h(y):=y^{q_0} \cdot (-\log(y))^{q_1} \cdot \ldots \cdot \log_{r-1}(-\log(y))^{q_r}$ (see for example Definition~1.13 and Proposition~3.16 in \cite{28}). Hence, $h$ cannot be restricted log-exp-analytic since $\lim_{y \searrow 0} h(y) / e^{-1/y} = 1$.
\end{proof}

So we see that the class of restricted log-exp-analytic functions is a proper subclass of the class of all definable functions. Since the global logarithm function is log-analytic, but not globally subanalytic and the global exponential function is restricted log-exp-analytic, but not log-analytic we obtain the following chain.

\begin{center}
$\textbf{globally subanalytic} \subsetneq \textbf{log-analytic} \subsetneq \textbf{restricted log-exp-analytic}$\\
$\subsetneq \textbf{log-exp-analytic} = \textbf{definable}.$
\end{center}

Finally we describe briefly the class of restricted log-exp-analytic functions in the complex setting.

For the rest of Section~3.2 let $z:=(z_1, \ldots ,z_m)$ range over $\mathbb{C}^m$. Let $x:=(x_1, \ldots ,x_m)$ and $y:=(y_1, \ldots ,y_m)$ range over $\mathbb{R}^m$ such that $z=x+iy$. 

\begin{definition}
\label{3.19}
Let $D \subset \mathbb{C}^m$ be definable and let $F:D \to \mathbb{C}$ be a function with $F(z)=u(x,y)+iv(x,y)$ where $u$ is the real and $v$ is the imaginary part of $F$ considered as real functions.
\begin{itemize}
	\item [(a)]  Let $Z \subset \mathbb{C}^m$ be open and definable with $D \subset Z$. We say that $F$ is \textbf{restricted log-exp-analytic with reference set $Z$}\index{restricted log-exp-anal.!with reference set} if there is a set $E$ of functions on $D$ with values in $\mathbb{C}^-$ such that $F$ can be constructed from $E$ and $\log(E^{\text{Re}})$ is a set of locally bounded functions on $D$ with reference set $Z$ where $D$ and $Z$ are considered as subsets of $\mathbb{R}^{2m}$.
	\item [(b)] Assume that $D$ is open in $\mathbb{C}^m$. We say that $F$ is \textbf{restricted log-exp-analytic}\index{restricted log-exp-anal.} if $F$ is restricted log-exp-analytic with reference set $D$. 
\end{itemize}
\end{definition}

\begin{example}
	Let $Z:=\{z \in \mathbb{C} \mid \vert{\text{Im(z)}}\vert<\pi \text{ and } \vert{\text{Im}(z+\exp(z))}\vert<\pi\}$. The function $G:Z \to \mathbb{C}, z \mapsto \log(\exp(z+\exp(z))+1),$ is restricted log-exp-analytic.
\end{example}

\begin{proof}
	Note that $Z$ is open. As in Example~\ref{ex:constructed-complex} we see that $G$ can be constructed from $E:=\{G_1,G_2\}$ where
	\[
	G_1:Z \to \mathbb{C}, z \mapsto \exp(z), \textnormal{ } G_2: Z \to \mathbb{C}, z \mapsto \exp(z+\exp(z)).
	\]
	Then with the proof of Example~\ref{ex:constructed-complex} it is easy to see that $E^{\text{Re}}$ is a set of locally bounded functions on $Z$ considered as subset of $\mathbb{R}^2$.
\end{proof}

We close this section by looking at the parametric setting: Let $l \in \mathbb{N}_0$, let $w$ range over $\mathbb{R}^l$, let $D \subset \mathbb{R}^n \times \mathbb{R}^l \times \mathbb{C}^m$ be a definable set and let $F:D \to \mathbb{C}$ be a function with $F(t,w,z)=u(t,w,x,y)+iv(t,w,x,y)$ where $u$ is the real part and $v$ is the imaginary part of $F$ considered as real functions. Here $(w,z)$ and $(w,x,y)$ are serving as tuples of independent variables of families of functions parameterized by $t$.

\begin{definition}
\label{3.20}
We define the following.
\begin{itemize}
	\item [(a)]
	Let $Z \subset \mathbb{R}^n \times \mathbb{R}^l \times \mathbb{C}^m$ be definable with $D \subset Z$ such that $Z_t$ is open in  $\mathbb{R}^l \times \mathbb{C}^m$ for every $t \in \mathbb{R}^n$. We say that $f$ is \textbf{restricted log-exp-analytic in $(w,z)$ with reference set $Z$} \index{restricted log-exp-anal.!with reference set} if there is a set $E$ of functions on $D$ with values in $\mathbb{C}^-$ such that $F$ can be constructed from $E$ and $\log(E^{\text{Re}})$ is a set of locally bounded functions in $(w,x,y)$ on $D$ with reference set $Z$ where $D$ and $Z$ are considered as subsets of $\mathbb{R}^n \times \mathbb{R}^l \times \mathbb{R}^{2m}$.
	
	
	\item [(b)]
	Suppose that $D_t$ is open for every $t \in \mathbb{R}^n$. We say that $f$ is \textbf{restricted log-exp-analytic in $(w,z)$} if $f$ is restricted log-exp-analytic in $(w,z)$ with reference set $D$. \index{restricted log-exp-anal.} 
\end{itemize}
\end{definition}

\begin{remark}
	\label{rem:rest-log-exp-real-imaginary}
	Let $Z \subset \mathbb{R}^n \times \mathbb{R}^l \times \mathbb{C}^m$ be definable with $D \subset Z$ such that $Z_t$ is open in $\mathbb{R}^l \times \mathbb{C}^m$ for every $t \in \mathbb{R}^n$. Let $F:D \to \mathbb{C}$ be restricted log-exp-analytic in $(w,z)$ with reference $Z$. Then $u$ and $v$ are restricted log-exp-analytic in $(w,x,y)$ with reference set $Z$ considered as subset of $\mathbb{R}^n \times \mathbb{R}^l \times \mathbb{R}^{2m}$.
\end{remark}


\section{A Preparation Theorem for Restricted \\ Log-Exp-Analytic Functions}

In this section we recall a preparation theorem for $\mathbb{R}_{\textnormal{an,exp}}$-definable functions from~\cite{27}. For its formulation technical tools like logarithmic scales, $C$-heirs and $C$-nice functions are needed (compare also with~\cite{27}). Then we formulate a preparation theorem for restricted log-exp-analytic functions from~\cite{28}. 

For this whole section let $n \in \mathbb{N}_0$, $t$ range over $\mathbb{R}^n$ and $x$ over $\mathbb{R}$. 
We fix a definable set $C \subset \mathbb{R}^n \times \mathbb{R}$. Let $\pi:\mathbb{R}^n \times \mathbb{R} \to \mathbb{R}^n, (t,x) \mapsto t$ be the projection on the first $n$ coordinates. 

\begin{definition}[Van den Dries et al.~\cite{100}]
\label{4.1}
Let $r \in \mathbb{N}_0$. A tuple \index{$y_l$} $\mathcal{Y}:=(y_0, \ldots ,y_r)$ of functions on $C$ is called an \textbf{$r$-logarithmic scale}\index{logarithmic scale} on $C$ with \textbf{center} $\Theta=(\Theta_0, \ldots ,\Theta_r)$ if the following holds:
\begin{itemize}
	\item[(a)] $y_j>0$ or $y_j<0$ for every $j \in \{0, \ldots ,r\}$.
	\item[(b)] $\Theta_j$ is a definable function on $\pi(C)$ for every $j \in \{0, \ldots ,r\}$.
	\item[(c)] For $(t,x) \in C$ we have $y_0(t,x)=x-\Theta_0(t)$ and $y_j(t,x)=\log(\vert{y_{j-1}(t,x)}\vert) - \Theta_j(t)$  for every $j \in \{1, \ldots ,r\}$.
	\item[(d)] We have $\Theta_0=0$ or there is $\varepsilon_0 \in \textnormal{}]0,1[$ such that $0<\vert{y_0(t,x)}\vert < \varepsilon_0\vert{x}\vert$ for all $(t,x) \in C$. For $j \in \{1, \ldots ,r\}$ the following holds: We have $\Theta_j=0$ or there is $\varepsilon_j \in \textnormal{}]0,1[$ such that $0<\vert{y_j(t,x)}\vert<\varepsilon_j\vert{\log(\vert{y_{j-1}(t,x)}\vert)}\vert$ for all $(t,x) \in C$.
\end{itemize}
\end{definition}

\begin{definition}
\label{4.2}
Let $r \in \mathbb{N}_0$, $\mathcal{Y}:=(y_0, \ldots ,y_r)$ be a logarithmic scale on $C$.
\begin{itemize}
	\item[(a)] We set $$\textrm{sign}(\mathcal{Y}):=\sigma:=(\sigma_0, \ldots ,\sigma_r):=(\textrm{sign}(y_0), \ldots, \textrm{sign}(y_r)) \in \{-1,1\}^{r+1}$$
	which defines the \textbf{sign} of $\mathcal{Y}$. \index{sign of $\mathcal{Y}$}
	\item[(b)] For $q:=(q_0, \ldots ,q_r) \in \mathbb{Q}^{r+1}$ let \index{$\vert{\mathcal{Y}}\vert^{\otimes q}$}
	$$\vert{\mathcal{Y}}\vert^{\otimes q}:C \to \mathbb{R}, (t,x) \mapsto \prod_{j=0}^r \vert{y_j(t,x)}\vert^{q_j} = \prod_{j=0}^r (\sigma_jy_j(t,x))^{q_j}.$$ 
\end{itemize}
\end{definition}

\begin{example}
\label{4.3}
Consider
$$C:=\{(t,x) \in \mathbb{R} \times \mathbb{R} \mid t \in \textnormal{} ]0,1[ \textnormal{}, \textnormal{}\tfrac{1}{1+t}+e^{-t-1/t} < x < \tfrac{1}{1+t} + e^{-t/2-1/t}\}.$$
Let $\Theta_0:\pi(C) \to \mathbb{R}, t \mapsto \tfrac{1}{1+t}$, $\Theta_1:\pi(C) \to \mathbb{R}, t \mapsto -\tfrac{1}{t}$ and $\Theta_2=0$. Let $y_0:C \to \mathbb{R}, (t,x) \mapsto x-\Theta_0(t)$, and inductively for $j \in \{1,2\}$ let $y_j:C \to \mathbb{R}, (t,x) \mapsto \log(\vert{y_{j-1}(t,x)}\vert)-\Theta_j(t)$. Then $\mathcal{Y}:=(y_0,y_1,y_2)$ is a $2$-logarithmic scale with center $(\Theta_0,\Theta_1,\Theta_2)$ and $\textnormal{sign}(\mathcal{Y})=(1,-1,-1) \in \{-1,1\}^3$.
\end{example}
	
\begin{proof}
Since $\Theta_0<C$ we see that $y_0>0$. So we have 
$$y_1(t,x)=\log(y_0(t,x))-\Theta_1(t)$$
for $(t,x) \in C$. We have $\Theta_0+e^{\Theta_1}>C$ and therefore $y_1<0$ since $y_1$ is strictly monotone increasing in $x$. This gives 
$$y_2(t,x)=\log(-y_1(t,x))-\Theta_2(t) = \log(-y_1(t,x))$$
for $(t,x) \in C$. We have $\Theta_0+e^{\Theta_1-1}<C$ and therefore $y_2<0$ and $\sigma_2=-1$ since $y_2$ is strictly monotone decreasing in $x$. An easy calculation shows for $(t,x) \in C$ that
$$x-\Theta_0(t)<e^{-1/t} < 2/(e+et) < \varepsilon_0x$$
for $\varepsilon_0:=2/e$ and that 
$$\vert{\log(x-\Theta_0(t))-\Theta_1(t)}\vert < t < \varepsilon_1(t/2 + 1/t) < \varepsilon_1\vert{\log(x-\Theta_0(t))}\vert$$
for $\varepsilon_1:=2/3$. So one sees that $(y_0,y_1,y_2)$ is a $2$-logarithmic scale with center $(\Theta_0,\Theta_1,\Theta_2)$. We have $\textnormal{sign}(\mathcal{Y})=(1,-1,-1) \in \{-1,1\}^3$.
\end{proof}

\begin{definition}[Lion et al.~\cite{25}, Section 0.4]
\label{4.6}
Let $r \in \mathbb{N}_0$. Let $f:C \to \mathbb{R}, (t,x) \mapsto f(t,x),$ be a function. We say that $g$ is \textbf{$r$-log-analytically prepared in $x$ with center $\Theta$}\index{$r$-log-analytically prepared} if
$$g(t,x)=a(t) \vert{\mathcal{Y}}\vert^{\otimes q}(t,x)u(t,x)$$
for all $(t,x) \in C$ where $a$ is a definable function on $\pi(C)$ which vanishes identically or has no zero, $\mathcal{Y}=(y_0, \ldots ,y_r)$ is an $r$-logarithmic scale with center $\Theta$ on $C$, $q \in \mathbb{Q}^{r+1}$ and the following holds for $u$ which we call \textbf{unit}\index{unit}. There is $s \in \mathbb{N}$ such that $u=v \circ \phi$ where $v$ is a power series which converges on an open neighborhood of $[-1,1]^s$ with $v([-1,1]^s) \subset \mathbb{R}_{>0}$ and $\phi:=(\phi_1, \ldots ,\phi_s):C \to [-1,1]^s$ is a function of the form 
$$\phi_j(t,x):=b_j(t)\vert{\mathcal{Y}}\vert^{\otimes p_j}(t,x)$$
for every $j \in \{1, \ldots ,s\}$ and $(t,x) \in C$ where $b_j:\pi(C) \to \mathbb{R}$ is definable without a zero and $p_j:=(p_{j0}, \ldots ,p_{jr}) \in \mathbb{Q}^{r+1}$. 
We call $a$ \textbf{coefficient}\index{coefficient} and $b:=(b_1, \ldots ,b_s)$ a tuple of \textbf{base functions}\index{base functions} for $f$. An \textbf{LA-preparing tuple}\index{LA-preparing tuple} for $f$ is then
$$(r,\mathcal{Y},a,q,s,v,b,P)$$
where
$$P:=\left(\begin{array}{cccc}
p_{10}&\cdot&\cdot&p_{1r}\\
\cdot&& &\cdot\\
\cdot&& &\cdot\\
p_{s0}&\cdot&\cdot&p_{sr}\\
\end{array}\right)\in M\big(s\times (r+1),\mathbb{Q}).$$
\end{definition}

Compare also Van den Dries et al.~\cite{100} for Definition~\ref{4.6} without a precise statement how the unit $u$ looks like.\\

Globally subananlytically prepared functions (see~\cite{25}) can be considered as special log-analytically prepared functions.

\begin{definition}[Lion et al.~\cite{25}, Section 0.3]
	\label{def:globally-subanalytical-preparation}
	We call an $r$-log-analytically prepared function $f:C \to \mathbb{R}$ with LA-preparing tuble $(r,\mathcal{Y},a,q,s,v,b,P)$ \textbf{globally subanalytically prepared in $x$}\index{globally subanalytically prepared} if $r=0$, and if $a$, $b$ and the center $\theta:\pi(C) \to \mathbb{R}$ of $\mathcal{Y}$ are globally subanalytic.
\end{definition}

\begin{fact}[Lion et al.~\cite{25}, Section 0.3]
	\label{fact:globally-subanalytical-preparation}
	Let $m \in \mathbb{N}$ and let $X \subset \mathbb{R}^n \times \mathbb{R}$. 
	Let $f:X \to \mathbb{R}$ be globally subanalytic. Then there is a globally subanalytic cell decomposition $\mathcal{D}$ of $X_{\neq 0}$ such that $f_1|_D, \ldots , f_m|_D$ are globally subanalytically prepared in $x$ with common center $\theta$.
\end{fact}

Now we proceed on formulating a preparation theorem for log-analytic functions established in~\cite{27}.

\begin{definition}
	\label{4.4}
	We call $g:\pi(C) \to \mathbb{R}$ \textbf{$C$-heir}\index{heir} if there is $l \in \mathbb{N}_0$, an $l$-logarithmic scale $\hat{\mathcal{Y}}$ with center $(\hat{\Theta}_0, \ldots ,\hat{\Theta}_l)$ on $C$, and $j \in \{1, \ldots ,l\}$ such that $g=\exp(\hat{\Theta}_j)$.
\end{definition}

\begin{definition}
	\label{4.5}
	We call $g:\pi(C) \to \mathbb{R}$ \textbf{$C$-nice}\index{nice} if there is a set $E$ of $C$-heirs such that $g$ can be constructed from $E$.
\end{definition}

Note that the class of log-analytic functions on $\pi(C)$ is a proper subclass of the class of $C$-nice functions. (See~\cite{27}, Section 2.3 for further properties and examples for $C$-heirs and $C$-nice functions.)

\begin{fact}[Opris~\cite{27} Theorem A]
\label{4.7}
Let $m \in \mathbb{N}$, $r \in \mathbb{N}_0$. Let $X \subset \mathbb{R}^n \times \mathbb{R}$ be definable. Let $f_1, \ldots .,f_m:X \to \mathbb{R}$ be log-analytic functions of order at most $r$. Then there is a definable cell decomposition $\mathcal{D}$ of $X_{\neq 0}$ such that $f_1|_D, \ldots ,f_m|_D$ are $r$-log-analytically prepared in $x$ with $D$-nice coefficient, $D$-nice base functions and common $D$-nice center for $D \in \mathcal{D}$.
\end{fact}

This gives roughly that the function $f(t,-)$ behaves cellwise as iterated logarithms independently of $t$ where the order of iteration is bounded in terms of $f$ (compare also with Van den Dries et al.~\cite{100} and Lion et al.~\cite{24}). However, our version is more precise than those from~\cite{100} and~\cite{24} since the coefficient, center and base functions of the preparation is $C$-nice (i.e. compositions of log-analytic functions and $C$-heirs) and not only definable. For example if $C$ is simple (i.e. for every $t \in \pi(C)$ there is $d_t>0$ such that $C_t=]0,d_t[$, see Definition 2.15 in Kaiser et al.~\cite{19}) we have that $a,\Theta_0, \ldots ,\Theta_r$ and $b_1, \ldots ,b_s$ are log-analytic (compare with Remark~3.10 in~\cite{28}) which cannot be obtained by using only the preparation theorems from~\cite{100} or~\cite{24}. However, $C$-heirs are in general not log-analytic as the following example shows. 

\begin{example}[Opris~\cite{27} Example 2.37]
\label{4.8}
Consider the definable cell
$$C:=\{(t,x) \in \mathbb{R}^2 \mid 0<t<1, \textnormal{} \tfrac{1}{1+t}+e^{-2/t+2e^{-1/t}}<x<\tfrac{1}{1+t}+e^{-1/t}\}.$$
There it is shown that $g:\pi(C) \to \mathbb{R}, t \mapsto e^{-1/t},$ is a $C$-heir. So the function 
$$f:C \to \mathbb{R}, (t,x) \mapsto e^{-1/t}x,$$
is $0$-log-analytically prepared in $x$ with center $0$ with a $C$-nice coefficient which is not log-analytic. 
\end{example}

This example also shows that the logarithms of $C$-heirs occuring in the compositions of the coefficient, center and base functions in a log-analytical preparation may differ to the components of the center of the logarithmic scale $\mathcal{Y}$ occuring in the preparation itself (see also Opris~\cite{27} for a further example).\\ 

Now we proceed with log-exp-analytic functions.

\begin{definition}[Opris~\cite{27} Definition 3.10]
\label{4.9}
Let $r \in \mathbb{N}_0$ and $e \in \mathbb{N}_0 \cup \{-1\}$. Let $f:C \to \mathbb{R}$ be a function. Let $E$ be a set of positive definable functions on $C$. By induction on $e$ we define that $f$ is \textbf{log-exp-analytically $(e,r)$-prepared}\index{log-exp-analytically prepared!$(e,r)$ with respect to} in $x$ with center $\Theta$ with respect to $E$. To this preparation we assign a \textbf{preparing tuple} for $f$.\index{preparing tuple}\\

$e=-1$: We say that $f$ is log-exp-analytically $(-1,r)$-prepared in $x$ with center $\Theta$ with respect to $E$ if $f$ is the zero function. A preparing tuple for $f$ is $(0)$.\\

$e-1 \to e$: We say that $f$ is log-exp-analytically $(e,r)$-prepared in $x$ with center $\Theta$ with respect to $E$ if 
$$f(t,x)=a(t)\vert{\mathcal{Y}}\vert^{\otimes q}(t,x)\exp(c(t,x)) u(t,x)$$
for every $(t,x) \in C$ where $a:\pi(C) \to \mathbb{R}$ is $C$-nice which vanishes identically or has no zero, $\mathcal{Y}:=(y_0, \ldots ,y_r)$ is an $r$-logarithmic scale with $C$-nice center $\Theta$, $q:=(q_0, \ldots ,q_r) \in \mathbb{Q}^{r+1}$, $\exp(c) \in E$ where $c$ is log-exp-analytically $(e-1,r)$-prepared in $x$ with center $\Theta$ with respect to $E$ and the following holds for $u$ which we call \textbf{unit}\index{unit}. There is $s \in \mathbb{N}$ such that $u=v \circ \phi$ where $\phi:=(\phi_1, \ldots ,\phi_s):C \to [-1,1]^s$ with 
$$\phi_j(t,x)=b_j(t)\vert{\mathcal{Y}}\vert^{\otimes p_j}(t,x)\exp(d_j(t,x))$$
for every $j \in \{1, \ldots ,s\}$ where $p_j:=(p_{j0}, \ldots ,p_{jr}) \in \mathbb{Q}^{r+1}$, $b_j:\pi(C) \to \mathbb{R}$ is $C$-nice without a zero, $\exp(d_j) \in E$, $d_j:C \to \mathbb{R}$ is log-exp-analytically $(e-1,r)$-prepared in $x$ with center $\Theta$ with respect to $E$ and $v$ is a power series which converges absolutely on an open neighborhood of $[-1,1]^s$ with $v([-1,1]^s) \subset \mathbb{R}_{>0}$. We call $a$ \textbf{coefficient} and $b:=(b_1, \ldots ,b_s)$ a tuple of \textbf{base functions}\index{coefficient}\index{base functions}. A preparing tuple for $f$ is then
$$\mathcal{J}:=(r,\mathcal{Y},a,\exp(c),q,s,v,b,\exp(d),P)$$
where $\exp(d):=(\exp(d_1), \ldots ,\exp(d_s))$ and 
$$P:=\left(\begin{array}{cccc}
	p_{10}&\cdot&\cdot&p_{1r}\\
	\cdot&& &\cdot\\
	\cdot&& &\cdot\\
	p_{s0}&\cdot&\cdot&p_{sr}\\
\end{array}\right)\in M\big(s\times (r+1),\mathbb{Q}).$$
\index{preparing tuple}\index{log-exp-analytically prepared!$(e,r)$ with respect to}
\end{definition}
Theorem B above is formulated in the parametric setting. So we recall a preparation theorem for restricted log-exp-analytic functions in the variables $(w_1, \ldots ,w_l,u_1, \ldots ,u_m,x)$, where $(u_1, \ldots ,u_m,x)$ is serving as the tuple of independent variables of families of functions parameterized by $w:=(w_1, \ldots ,w_l)$. The preparation theorem is formulated with respect to the variable $x$. \\

For the rest of Section 4 let $m,l \in \mathbb{N}_0$ be with $n=l+m$. Let $w:=(w_1, \ldots ,w_l)$ range over $\mathbb{R}^l$ and $u:=(u_1, \ldots ,u_m)$ over $\mathbb{R}^m$. Let $X \subset \mathbb{R}^l \times \mathbb{R}^m \times \mathbb{R}$ be a definable set with $C \subset X$. Assume that $X_w$ is open for every $w \in \mathbb{R}^l$.

\begin{definition}
\label{4.10}
Let $e \in \mathbb{N}_0 \cup \{-1\}$ and $r \in \mathbb{N}_0$. We call $f:C \to \mathbb{R}, (w,u,x) \mapsto f(w,u,x),$ \textbf{$(m+1,X)$-restricted log-exp-analytically $(e,r)$-prepared in $x$}\index{log-exp-analytically prepared!$(m+1,X)$-restricted $(e,r)$} if $f$ is log-exp-analytically $(e,r)$-prepared in $x$ with respect to a set $E$ of positive definable functions such that every $g \in \log(E)$ is locally bounded in $(u,x)$ with reference set $X$.
\end{definition}

\begin{example}
\label{4.11}
Let $C=\textnormal{}]0,1[^2$ and $X=\mathbb{R}_{>0} \times \mathbb{R}_{>0}$. The function 
$$C \to \mathbb{R}, (u,x) \mapsto u \cdot e^{x \cdot \log(x)},$$
is $(2,X)$-restricted log-exp-analytically $(1,1)$-prepared with center $0$. The function 
$$C \to \mathbb{R}, (u,x) \mapsto u \cdot e^{\frac{x \cdot \log(x)}{1-u}},$$
is $(1,X)$-restricted log-exp-analytically $(1,1)$-prepared with center $0$.
\end{example}

\begin{fact}[Opris~\cite{28} Proposition 2.12]
\label{4.12}
Let $e \in \mathbb{N}_0 \cup \{-1\}$. Let $f:X \to \mathbb{R}$ be a restricted log-exp-analytic function in $(u,x)$ of order at most $e$. Then there is $r \in \mathbb{N}_0$ and a definable cell decomposition $\mathcal{C}$ of $X_{\neq 0}$ such that $f|_C$ is $(m+1,X)$-restricted log-exp-analytically $(e,r)$-prepared in $x$ for every $C \in \mathcal{C}$.
\end{fact}

\section{Holomorphic Extensions of Log-Exp-Analy-tically Prepared Functions}

Fact~\ref{4.12} indicates that a restricted log-exp-analytic function can be prepared cellwise. So the first step in the construction of a unary parametric global complexification of such a function is to find a suitable holomorphic extension of a parametric family of restricted log-exp-analytically prepared one-variable functions on a single cell $C$, as discussed in this section. As indicated with Definitions~\ref{4.9},~\ref{4.10} and Fact~\ref{4.12} logarithmic scales play a crucial role in this context. 

The procedure of this section is as follows. In the first subsection we construct a suitable unary parametric global complexification of a logarithmic scale, and in the second one we investigate log-analytically prepared functions. Finally, we extend these results on log-exp-analytically prepared functions which are more general than restricted log-exp-analytically prepared ones.\\

For the whole Section let $n \in \mathbb{N}_0$, let $t$ range over $\mathbb{R}^n$ and $x$ over $\mathbb{R}$. We fix a non-empty definable cell $C \subset \mathbb{R}^n \times \mathbb{R}_{\neq 0}$ and an $r$-logarithmic scale $\mathcal{Y}:=(y_0, \ldots ,y_r)$ on $C$ with center $\Theta:=(\Theta_0, \ldots ,\Theta_r)$. We fix $\textrm{sign}(\mathcal{Y}):=\sigma:=(\sigma_0, \ldots ,\sigma_r) \in \{-1,1\}^{r+1}$ (i.e. $y_0=x-\Theta_0$ and inductively $y_j=\log(\sigma_{j-1}y_{j-1})-\Theta_j$ for $j \in \{1, \ldots , r\}$).
Let $z$ range over $\mathbb{C}$. Let $\pi:\mathbb{R}^n \times \mathbb{C} \to \mathbb{R}^n, (t,z) \mapsto t,$ be the projection on the first $n$ real coordinates (which includes also the projection $\mathbb{R}^n \times \mathbb{R} \to \mathbb{R}^n, (t,x) \mapsto t$)\index{$\pi(A)$}\index{projection on the first $n$ real coordinates}.

We say that a function $f:C \to \mathbb{R}, (t,x) \mapsto f(t,x),$ is strictly monotone increasing in $x$ if $f_t$ is strictly monotone increasing for all $t \in \pi(C)$. Similarly, $f$ is strictly monotone decreasing in $x$ if $f_t$ is strictly monotone decreasing for all $t \in \pi(C)$ and $f$ is strictly monotone in $x$ if $f_t$ is strictly monotone for all $t \in \pi(C)$.

\subsection{Logarithmic Scales}

The goal for this subsection is to compute a unary parametric global complexification of a logarithmic scale $\mathcal{Y}$. Particulary, we construct $C \subset H \subset \mathbb{R}^n \times \mathbb{C}$ such that $H_t$ is open for all $t \in \pi(C)$ and a definable $\mathcal{Z}:H \to \mathbb{C}$ such that $\mathcal{Z}|_C=\mathcal{Y}$ and $\mathcal{Z}_t$ is holomorphic for every $t \in \pi(C)$. 
The first step in the construction is to naturally extend $y_0, \ldots ,y_r:C \to \mathbb{R}$ to functions $g_0, \ldots ,g_r:D \to \mathbb{R}$ beyond the cell $C$ where $\pi(C)=\pi(D)$.\\

We set
$D_0:=\pi(C) \times \mathbb{R}$ and 
$$g_0:D_0 \to \mathbb{R}, (t,x) \mapsto x-\Theta_0(t).$$
For $l \in \{1, \ldots ,r\}$ we set inductively
$$D_l:=\{(t,x) \in D_{l-1} \mid \sigma_{l-1}g_{l-1}(t,x)>0\}$$
and
$$g_l:D_l \to \mathbb{R}, (t,x) \mapsto \log(\sigma_{l-1}g_{l-1}(t,x))- \Theta_l(t).$$
Set 
$$D:=\{(t,x) \in D_r \mid \sigma_rg_r(t,x)>0\}.$$

\begin{remark}
	\label{5.1}
	We have $C \subset D \subset D_r \subset  \ldots  \subset D_1$ and for $l \in \{0, \ldots ,r\}$ we have that $g_l|_C=y_l$ (i.e. $g_l$ is an extension of $y_l$), and $(g_l)_t:(D_l)_t \to \mathbb{R}$ is injective and real analytic for every $t \in \pi(C)$. 
\end{remark}

To be able to holomorphically extend the iterated logarithms occuring in $\mathcal{Y}$ one has to determine the zero set $Z(g_l)$ of $g_l$ for $l \in \{0, \ldots ,r\}$. 
In Lemma~\ref{5.4} and Corollary~\ref{5.5} we show that $Z(g_l)$ is exactly the graph of the function $\mu_l$ which we introduce in the next definition. The functions $\mu_{j,l}$ are needed for technical reasons in this context for $0 \leq j < l$. 

\begin{definition}
\label{5.2}
Let $l \in \{0, \ldots ,r\}$. Let $\mu_{0,l}:\pi(C) \to \mathbb{R}, t \mapsto \Theta_l(t),$ and inductively for $j \in \{1, \ldots ,l\}$ we set \index{$\mu_{j,l}$}
$$\mu_{j,l}:\pi(C) \to \mathbb{R}, t \mapsto \Theta_{l-j}(t)+\sigma_{l-j}e^{\mu_{j-1,l}(t)}.$$
Set $\mu_l:=\mu_{l,l}$\index{$\mu_\kappa$,$\mu_l$}.
\end{definition}

Since the pointwise addition is globally subanalyic and therefore log-analytic, we obtain with Lemma~\ref{3.4}(1),(2) the following.

\begin{remark}
\label{5.3}
Let $l \in \{0, \ldots ,r\}$ and let $\mathcal{E}$ be a set of positive definable functions on $\pi(C)$ such that $\Theta$ can be constructed from $\mathcal{E}$. Then $\mu_{j,l}$ for $j \in \{0, \ldots ,l\}$ can be constructed from $\mathcal{E} \cup \{e^{\mu_{k,l}} \mid k \in \{0, \ldots , j-1\}\}$.
\end{remark}


\begin{lemma}
\label{5.4}
Let $t \in \pi(C)$. Let $l \in \{0, \ldots ,r\}$ and $j \in \{0, \ldots ,l\}$. We have $(t,\mu_l(t)) \in D_j$ and
$$g_j(t,\mu_l(t)) = \sigma_je^{\mu_{l-j-1,l}(t)}$$
if $j < l$.
\end{lemma}

\begin{proof}
That $(t,\mu_l(t)) \in D_0$ is clear. So suppose that $l>0$. We show the statement by induction on $j \in \{0, \ldots ,l\}$.

$j=0$: With Definition~\ref{5.2} we obtain
$$g_0(t,\mu_l(t)) = \mu_l(t) - \Theta_0(t) = \sigma_0e^{\mu_{l-1,l}(t)}.$$

$j-1 \to j$:
By the inductive hypothesis $(t,\mu_l(t)) \in D_{j-1}$ and 
$$g_{j-1}(t,\mu_l(t))=\sigma_{j-1}e^{\mu_{l-j,l}(t)}.$$
So we see that $(t,\mu_l(t)) \in D_j$. If $j<l$ we obtain
\begin{align*}
	 g_j(t,\mu_l(t)) &= \log(\sigma_{j-1}g_{j-1}(t,\mu_l(t)))- \Theta_j(t) \\
	&=\mu_{l-j,l}(t) - \Theta_j(t) = \sigma_je^{\mu_{l-j-1,l}(t)}. \qedhere
\end{align*}
\end{proof}

\begin{corollary}
\label{5.5}	
For $l \in \{0, \ldots ,r\}$ and $t \in \pi(C)$ we have $g_l(t,\mu_l(t)) = 0$ and
$$Z(g_l):=\{(t,x) \in D_l \mid g_l(t,x)=0\}=\{(t,x) \in D_l \mid x=\mu_l(t)\}.$$
\end{corollary}

\begin{proof}
For $l=0$ this is clear. So assume $l>0$. With Lemma~\ref{5.4} we obtain  $(t,\mu_l(t)) \in D_l$,
$$g_{l-1}(t,\mu_l(t))=\sigma_{l-1}e^{\Theta_l(t)},$$
and therefore
$$g_l(t,\mu_l(t)) = \log(\sigma_{l-1}g_{l-1}(t,\mu_l(t)))-\Theta_l(t) = 0$$
for every $t \in \pi(C)$. Due to Remark~\ref{5.1} we see that $(g_l)_t$ is injective for every $t \in \pi(C)$ and we are done with the proof.
\end{proof}

Since $g_l|_C=y_l$, we see with Definition~\ref{4.1}(a) that either $g_l|_C>0$ or $g_l|_C<0$. Hence, the graph of $\mu_l$ is disjoint to $C$. Since $(g_l)_t$ is also injective for $t \in \pi(C)$, we obtain with Lemma~\ref{5.4} that the graphs of $\mu_l$ and $\mu_j$ for $j \neq l$ are disjoint.\\
In the next lemma we will strengthen these results and even
show that $\text{sign}(\mathcal{Y})$ uniquely determines whether the graph of $\mu_l$ is above or below both $C$ and the graph of $\mu_j$ for $j \neq l$.

\begin{lemma}
\label{5.6}
Let $l \in \{0, \ldots ,r\}$. The following properties hold.
	\begin{itemize}
		\item[(1)] The function $g_l$ is strictly monotone in $x$. It is strictly monotone increasing in $x$ if and only if $\prod_{j=0}^{l-1}\sigma_j=1$.
		
		\item[(2)]
		Either $\mu_l<C$ or $\mu_l>C$. We have $\mu_l<C$ if and only if $\prod_{j=0}^l\sigma_j=1$.	
		
		\item[(3)] Let $k \in \{l+1, \ldots ,r\}$.
		Either $\mu_l<\mu_k$ or $\mu_k<\mu_l$. we have $\mu_l<\mu_k$ if and only if $\prod_{j=0}^l \sigma_j=1$.
\end{itemize}
\end{lemma}

\begin{proof}
(1): We do an induction on $l \in \{0, \ldots ,r\}$. For $l=0$ there is nothing to show.

$l-1 \to l$: We assume $\sigma_{l-1}=1$. Then $g_l$ is strictly monotone increasing in $x$ if and only if $g_{l-1}$ is strictly monotone increasing in $x$. By the inductive hypothesis we have that $g_{l-1}$ is strictly monotone increasing in $x$ if and only if $\prod_{j=0}^{l-2}\sigma_j=1$. This gives the result. The case ''$\sigma_{l-1}=-1$'' is done similarly.\\

(2): Because $g_l|_C=y_l$ we obtain with Definition~\ref{4.1} $g_l|_C>0$ or $g_l|_C<0$. We assume $g_l|_C>0$, i.e. $\sigma_l=1$. By (1) we have that $g_l$ is either strictly monotone increasing in $x$ or strictly monotone decreasing in $x$. If the former holds then $\mu_l<C$ by Corollary~\ref{5.5}. Additionally, $\prod_{j=0}^{l-1} \sigma_j=1$ by (1) and therefore $\prod_{j=0}^l \sigma_j=1$. If the latter holds then $\mu_l>C$ by Corollary~\ref{5.5}. Additionally, $\prod_{j=0}^{l-1} \sigma_j=-1$ by (1) and therefore $\prod_{j=0}^l \sigma_j=-1$. The case ''$g_l|_C<0$'' is done similarly.\\

(3): With Lemma~\ref{5.4} we obtain 
$$g_l(t,\mu_k(t))=\sigma_le^{\mu_{k-l-1,k}(t)}$$
for every $t \in \pi(C)$. We assume that $g_l$ is strictly monotone increasing in $x$. By (1) we get $\prod_{j=0}^{l-1}\sigma_j=1$. If $\sigma_l=1$ then $g_l(t,\mu_k(t))=e^{\mu_{k-l-1,k}(t)}>0$ for every $t \in \pi(C)$ and consequently $\mu_l<\mu_k$ by Corollary~\ref{5.5}. If $\sigma_l=-1$ then $g_l(t,\mu_k(t))=-e^{\mu_{k-l-1,k}(t)}<0$ for every $t \in \pi(C)$ and consequently $\mu_k<\mu_l$ by Corollary~\ref{5.5}. The case ''$g_l$ is strictly monotone decreasing'' is done similarly.
\end{proof}

An immediate consequence of Lemma~\ref{5.6} is the following.

\begin{corollary}
\label{only-r}
If $\textnormal{sign}(\mathcal{Y})=(1, \ldots ,1)$ then $\mu_0<\mu_1 < \ldots < \mu_r < C$. If $\textnormal{sign}(\mathcal{Y})=(-1,1,\ldots,1)$ then $\mu_0>\mu_1> \ldots > \mu_r>C$. If $\textnormal{sign}(\mathcal{Y}) \neq (1, \ldots ,1)$ and $\textnormal{sign}(\mathcal{Y}) \neq (-1,1 \ldots ,1)$ there is $j,l \in \{0, \ldots, r\}$ with $j \neq l$ such that $\mu_j < C < \mu_l$.
\end{corollary}

\begin{proof}
	In the first case we have $\prod_{j=0}^l \sigma_j=1$ for all $l \in \{0, \ldots, r\}$  and in the second  $\prod_{j=0}^l \sigma_j=-1$ for all $l \in \{0, \ldots, r\}$. Hence, the desired properties follow from Lemma~\ref{5.6}(2),(3).  In the third case choose a maximum $l \in \{1, \ldots , r\}$ with $\sigma_l=-1$. Then $\prod_{j=0}^{l-1} \sigma_j \neq \prod_{j=0}^{l} \sigma_j$ and we obtain the result with Lemma~\ref{5.6}(2). 
\end{proof}


When $\text{sign}(\mathcal{Y}) \neq (1, \ldots, 1)$ and $\text{sign}(\mathcal{Y}) \neq (-1, \ldots, 1)$, we will prove the existence of two unique zeros such that the graph of one zero is closest to $C$ from above, and the graph of the other one is closest to $C$ from below.

\begin{definition}
\label{5.7}
We define the \textbf{change index}\index{change index} $k^{\textnormal{ch}}$ \index{$k^{\textnormal{ch}}$} for $\mathcal{Y}$ as follows: If there is $l \in \{0, \ldots ,r\}$ with $\sigma_l=-1$ set
$$k^{\textnormal{ch}}:=\max\{l \in \{0, \ldots ,r\} \mid \sigma_l=-1\}-1.$$
Otherwise set $k^{\textnormal{ch}}:=-2$.
\end{definition}

\begin{remark}
	\label{rem:change:globally-subanalytical}
	If $r=0$ then $k^{\text{ch}}<0$.
\end{remark}

\begin{remark}
\label{5.8}
The following holds.
\begin{itemize}
	\item[(1)] If $k^{\textnormal{ch}}=-2$ then $\sigma_0= \ldots =\sigma_r=1$, i.e. $\textnormal{sign}(\mathcal{Y})=(1,1, \ldots ,1)$.
	\item[(2)] If $k^{\textnormal{ch}}=-1$ then $\sigma_0=-1$ and $\sigma_1= \ldots =\sigma_r=1$, i.e. $\textnormal{sign}(\mathcal{Y})=(-1,1, \ldots ,1)$.
\end{itemize}
\end{remark}

An immediate consequence of Lemma~\ref{5.6} and Definition~\ref{5.7} is that $\mu_r$ and $\mu_{k^{\textnormal{ch}}}$ are indeed the two unique zeros with the desired properties.

\begin{corollary}
	\label{k-and-r}
	Suppose that $k:=k^{\textnormal{ch}} \geq 0$.
	\begin{itemize}
		\item[(1)] Let $\prod_{j=0}^r \sigma_j=1$.  Then $\mu_r<C<\mu_k$ and for any zero $\mu_l$ with $l \neq r$ and $\mu_l<C$ we have that $\mu_l < \mu_r$. For any zero $\mu_l$ with $l \neq k$ and $\mu_l>C$ we have that $\mu_l > \mu_k$.
		\item[(2)] Let $\prod_{j=0}^r \sigma_j=-1$. Then $\mu_k<C<\mu_r$ and for any zero $\mu_l$ with $l \neq r$ and $\mu_l>C$ we have that $\mu_l > \mu_r$. For any zero $\mu_l$ with $l \neq k$ and $\mu_l<C$ we have that $\mu_l < \mu_k$.
	\end{itemize}
\end{corollary}

\begin{proof}
	We show only (1). The second property is proven analoguously. If $\mu_l<C$ for $l<r$ then $\prod_{j=0}^l \sigma_j=1$ by Lemma~\ref{5.6}(1) and hence $\mu_l<\mu_r$ by Lemma~\ref{5.6}(3). If $\mu_l>C$ then $\prod_{j=0}^l \sigma_j=-1$ and $l \leq k$ by Definition~\ref{5.7} since $\prod_{j=0}^m \sigma_j=1$ for $m > k$. If also $l \neq k$ then $l<k$ and hence, $\mu_k<\mu_l$ by Lemma~\ref{5.6}(3).
\end{proof}

Now we can fully describe $D$ (defined at the beginning of this subsection) using $\mu_r$ and $\mu_{k^{\text{ch}}}$ where the latter is required only if $k^{\text{ch}} \geq 0$. 

\begin{lemma}
\label{5.9}
Let $k:=k^{\textnormal{ch}}$. The following holds.

\begin{itemize}
	\item [(1)] If $\prod_{j=0}^r\sigma_j=1$ and $k<0$ then
	$$D=\{(t,x) \in \pi(C) \times \mathbb{R} \mid x > \mu_r(t)\}.$$
	\item [(2)]
	If $\prod_{j=0}^r\sigma_j=-1$ and $k<0$ then
	$$D=\{(t,x) \in \pi(C) \times \mathbb{R} \mid x < \mu_r(t)\}.$$
	\item [(3)]
	If $\prod_{j=0}^r\sigma_j=1$ and $k \geq 0$ then
	$$D=\{(t,x) \in \pi(C) \times \mathbb{R} \mid \mu_r(t) < x < \mu_k(t) \}.$$
	\item [(4)]
	If $\prod_{j=0}^r\sigma_j=-1$ and $k \geq 0$ then
	$$D=\{(t,x) \in \pi(C) \times \mathbb{R} \mid \mu_k(t) < x < \mu_r(t)\}.$$
\end{itemize}
\end{lemma}

\begin{proof}
We show property (3). The rest is proven similarly. Note that $\prod_{j=0}^k\sigma_j=-1$, $\mu_r<C<\mu_k$ and $\sigma_\kappa g_\kappa(t,x)=\sigma_\kappa y_\kappa(t,x)>0$ for $(t,x) \in C$ and $\kappa \in \{k,r\}$.

Let $(t,x) \in \pi(C) \times \mathbb{R}$ be with $\mu_r(t)<x<\mu_k(t)$. If $\prod_{j=0}^l \sigma_j=1$ then $\mu_l<\mu_r$ if $l \neq r$ by Lemma~\ref{5.6}(3). If $\prod_{j=0}^l \sigma_j=-1$ then $\mu_k<\mu_l$ if $l \neq k$ (since then $l < k$). So we obtain with the intermediate value theorem, Definition~\ref{4.1} and Corollary~\ref{5.5} that $\sigma_jg_j(t,x) > 0$ for every $j \in \{0, \ldots ,r\}$ since $C \subset D$ and $(g_j)_t$ is continuous . Therefore $(t,x) \in D$. 

Let $(t,x) \in D$. Then we have $\sigma_\kappa g_\kappa(t,x)>0$ for $\kappa \in \{k,r\}$. Again with Lemma~\ref{5.6}(2), Corollary~\ref{5.5}, the continuity of $g_\kappa$ in $x$ and the monotony property of $g_\kappa$ for $\kappa \in \{k,r\}$ we obtain that $\mu_r(t)<x<\mu_k(t)$.
\end{proof}

\begin{figure}[h]
	\begin{center}\includegraphics[width=9cm,height=6cm,keepaspectratio]{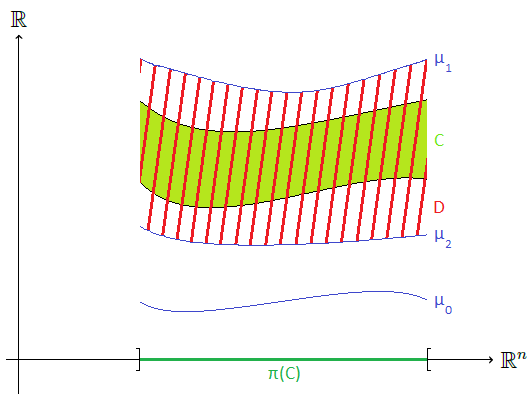}\end{center}
	\caption{The set $D$ and the functions $\mu_j$  for a $2$-logarithmic scale $\mathcal{Y}$ on $C$ with $\textnormal{sign}(\mathcal{Y})=(1,-1,-1)$ and continuous center. Here we have $r=2$ and $k^{\textnormal{ch}}=1$.}
\end{figure}

\begin{example}
	\label{ex:C-simple}
	If the cell $C$ is \emph{simple} (i.e. $C_t$ of the form $]0,d_t[$ for every $t \in \pi(C)$, see also Kaiser et al.~\cite{19} or Opris~\cite{28}) then there is only one unique logarithmic scale $\mathcal{Y}$ on $C$ with center $0$ and $\text{sign}(\mathcal{Y}) =(1,-1,1, \ldots ,1)$ (see for example Proposition 2.17 in~\cite{19}). Hence, we have $k^{\text{ch}}=0$. In this situation $\mu_0=\Theta_0=0$ and $\mu_l=1/\exp_{l-1}(1)$ for every $l \in \{1, \ldots , r\}$.
\end{example}

\begin{example}
\label{5.10}
Consider 
$$C:=\{(t,x) \in \mathbb{R} \times \mathbb{R} \mid t \in \textnormal{} ]0,1[ \textnormal{}, \textnormal{}\tfrac{1}{1+t}+e^{-t-1/t} < x < \tfrac{1}{1+t} + e^{-t/2-1/t}\}.$$
Let $\Theta_0:\pi(C) \to \mathbb{R}, t \mapsto \tfrac{1}{1+t}$, $\Theta_1:\pi(C) \to \mathbb{R}, t \mapsto -\tfrac{1}{t}$ and $\Theta_2=0$. Let $y_0:C \to \mathbb{R}, (t,x) \mapsto x-\Theta_0(t)$, and inductively for $j \in \{1,2\}$ let $y_j:C \to \mathbb{R}, (t,x) \mapsto \log(\vert{y_{j-1}(t,x)}\vert)-\Theta_j(t)$. Then $\mathcal{Y}:=(y_0,y_1,y_2)$ is a $2$-logarithmic scale with center $(\Theta_0,\Theta_1,\Theta_2)$ and $\textnormal{sign}(\mathcal{Y}):=(\sigma_0,\sigma_1,\sigma_2)=(1,-1,-1) \in \{-1,1\}^3$. Therefore $k^{\textnormal{ch}}=1$, $\mu_0=\Theta_0$, $\mu_1=\Theta_0+e^{\Theta_1}$, $\mu_2=\Theta_0+e^{\Theta_1-1}$ and $\mu_{0,1}=\Theta_1$, $\mu_{0,2}=0$ and $\mu_{1,2}=\Theta_1-1$.
\end{example}

\begin{proof}
Compare with the proof of Example~\ref{4.3}.
\end{proof}

Now we are able to compute the desired unary parametric global complexification of $\mathcal{Y}$. So we merge into the complex setting.

\begin{definition}
\label{5.11}
We define the following.
\begin{itemize}
	\item [(a)] Consider $H_0:=\pi(C) \times \mathbb{C}$ and 
	$$z_0:H_0 \to \mathbb{C}, (t,z) \mapsto z-\Theta_0(t).$$
	Inductively for $l \in \{1, \ldots ,r\}$ we set \index{$H_l$}
	$$H_l:=\{(t,z) \in H_{l-1} \mid \sigma_{l-1}z_{l-1}(t,z) \in \mathbb{C}^{-}\}$$
	and \index{$z_l$}
	$$z_l:H_l \to \mathbb{C}, (t,z) \mapsto \log(\sigma_{l-1}z_{l-1}(t,z))- \Theta_l(t).$$
	\item [(b)] Set \index{$H$}
	$$H:=\{(t,z) \in H_r \mid \sigma_rz_r(t,z) \in \mathbb{C}^-\}.$$
\end{itemize}
\end{definition}

\begin{lemma}
\label{5.12}
The following holds for $l \in \{0, \ldots ,r\}$.
\begin{itemize}
	\item [(1)] We have $C \subset D_l \subset H_l$, $H_l \subset H_{l-1}$ if $l > 0$ and $(H_l)_t$ is open for every $t \in \pi(C)$. Additionally, $z_l$ is well-defined and $(z_l)_t$ is holomorphic in $z$ for every $t \in \pi(C)$.
	\item [(2)] We have $z_l|_{D_l}=g_l$ and $z_l|_C=y_l$.
	\item [(3)] Let $E$ be a set of positive definable functions such that $\Theta_0, \ldots ,\Theta_r$ can be constructed from $E$. Then $z_l$ can be constructed from $E$.
\end{itemize}
\end{lemma}

\begin{proof}
(1): This follows by an easy induction on $l$ and the fact that the complex logarithm function is holomorphic on $\mathbb{C}^-$.

(2): This follows immediately by an easy induction on $l$ and from Definition~\ref{5.11}(a).

(3): We do an induction on $l$. We have $z_0(t,z)=z-\Theta_0(t)$ for every $(t,z) \in H_0$ and use Lemma~\ref{3.10}(2). $l-1 \to l$: We have $z_l(t,z)=\log(\sigma_{l-1}z_{l-1}(t,z))-\Theta_l(t)$ for every $(t,z) \in H_l$ and use the inductive hypothesis and Lemma~\ref{3.10}(2) (since the complex logarithm is log-analytic by Lemma~\ref{3.7}).
\end{proof}

We often consider $z_l$ as a function on $H$ for every $l \in \{0, \ldots ,r\}$. Lemma~\ref{5.12} shows that $(z_l)_t$ is a holomorphic extension of $(y_l)_t$ for every $t \in \pi(C)$ and hence, $\mathcal{Z}:=(z_0, \ldots , z_r)$ is a unary parametric global complexification of $\mathcal{Y}=(y_0, \ldots , y_r)$.

\begin{lemma}
\label{5.13}
Let $k:=k^{\textnormal{ch}}$. Then the following holds where $\mathbb{C}_{a}^-:=\mathbb{C} \setminus \mathbb{R}_{\leq a}$ and $\mathbb{C}_{a}^+:=\mathbb{C} \setminus \mathbb{R}_{\geq a}$ for $a \in \mathbb{R}$.
\begin{itemize}
\item[(1)] If $\prod_{j=0}^r\sigma_j=1$ and $k<0$ then
$$H=\{(t,z) \in \pi(C) \times \mathbb{C} \mid z \in \mathbb{C}^-_{\mu_r(t)}\}.$$
\item[(2)] If $\prod_{j=0}^r\sigma_j=-1$ and $k<0$ then
$$H=\{(t,z) \in \pi(C) \times \mathbb{C} \mid z \in \mathbb{C}^+_{\mu_r(t)}\}.$$
\item[(3)] If $\prod_{j=0}^r\sigma_j=1$ and $k \geq 0$ then
$$H=\{(t,z) \in \pi(C) \times \mathbb{C} \mid z \in \mathbb{C}^-_{\mu_r(t)} \cap \mathbb{C}^+_{\mu_k(t)}\}.$$
\item[(4)] If $\prod_{j=0}^r\sigma_j=-1$ and $k \geq 0$ then
$$H=\{(t,z) \in \pi(C) \times \mathbb{C} \mid z \in \mathbb{C}^+_{\mu_r(t)} \cap \mathbb{C}^-_{\mu_k(t)}\}.$$
\end{itemize}
\end{lemma}

\begin{proof}
We have
$$H=D \cup (\pi(C) \times \mathbb{C} \setminus \mathbb{R}),$$
because for $(t,z) \in \pi(C) \times \mathbb{C} \setminus \mathbb{R}$ we obtain $\sigma_lz_l(t,z) \in \mathbb{C} \setminus \mathbb{R}$ for every $l \in \{0, \ldots ,r\}$ by the definition of the complex logarithm function. We are done due to Lemma~\ref{5.12}(2) and Lemma~\ref{5.9}.
\end{proof}

Now we are able to compute a unary parametric global complexification of $\vert{\mathcal{Y}}\vert^{\otimes q}$ for every $q \in \mathbb{Q}^{r+1}$. 

\begin{definition}
\label{5.14}
Let $q:=(q_0, \ldots ,q_r) \in \mathbb{Q}^{r+1}$. We set \index{$(\sigma\mathcal{Z})^{\otimes{q}}$}
$$(\sigma\mathcal{Z})^{\otimes q}:H \to \mathbb{C}, (t,z) \mapsto \prod_{j=0}^r (\sigma_jz_j(t,z))^{q_j}.$$
\end{definition}

\begin{lemma}
\label{5.15}
Let $q \in \mathbb{Q}^{r+1}$. The following holds.
\begin{itemize}
	\item [(1)]
	$H_t$ is open for every $t \in \pi(C)$ and $(\sigma\mathcal{Z})^{\otimes q}$ is well-defined and holomorphic in $z$. In particular
	$$(\sigma\mathcal{Z})^{\otimes q}|_C=\vert{\mathcal{Y}}\vert^{\otimes q}.$$ 
	\item [(2)] Let $E$ be a set of positive definable functions such that $\Theta_0, \ldots ,\Theta_r$ can be constructed from $E$. Then $(\sigma\mathcal{Z})^{\otimes q}$ can be constructed from $E$.\\
\end{itemize}
\end{lemma}

\begin{proof}
Note that $\sigma_lz_l(t,z) \in \mathbb{C}^-$ for every $(t,z) \in H$ and $l \in \{0, \ldots ,r\}$ and that for $p \in \mathbb{Q}$ the function $\mathbb{C}^- \to \mathbb{C}, z \mapsto z^p,$ is holomorphic and globally subanalytic. So property (1) follows from Lemma~\ref{5.12}(1) and Lemma~\ref{5.12}(2). Property (2) follows immediately with Lemma~\ref{5.12}(3) and Lemma~\ref{3.10}(2).
\end{proof}

At the end of this subsection we introduce a further class of functions which describe $z_l$ as a function of $z-\mu_\kappa$.  Those are needed for technical issues for example to describe the asymptotic behavior of $(\sigma\mathcal{Z})^{\otimes q}$ in Section~6 or when dealing with parameterized integrals of $(\sigma\mathcal{Z})^{\otimes q}$ in Section~8.\\

For the rest of this subsection we fix $\kappa=r$ or $\kappa \in \{k^{\textnormal{ch}},r\}$ if $k^{\textnormal{ch}} \geq 0$.

\begin{definition}
\label{5.16}
Let $H_{0,\kappa}^*:=\pi(C) \times \mathbb{C}$ and $\mathcal{P}_{0,\kappa}:H_{0,\kappa}^* \to \mathbb{C}, (t,z) \mapsto z-\mu_\kappa(t),$ and inductively for $l \in \{1, \ldots ,\kappa\}$ let
$$H_{l,\kappa}^*:=\Bigl\{(t,z) \in H_{l-1,\kappa}^* \mid 1+\frac{\sigma_{l-1}\mathcal{P}_{l-1,\kappa}(t,z)}{e^{\mu_{\kappa-l,\kappa}(t)}} \in \mathbb{C}^-\Bigl\}$$
and\index{$\mathcal{P}_{l,\kappa}$}
$$\mathcal{P}_{l,\kappa}:H_{l,\kappa}^* \to \mathbb{C}, (t,z) \mapsto  \log \Bigl(1+\frac{\sigma_{l-1}\mathcal{P}_{l-1,\kappa}(t,z)}{e^{\mu_{\kappa-l,\kappa}(t)}} \Bigl).$$
Set\index{$\mathcal{P}_\kappa$} $\mathcal{P}_\kappa:=\mathcal{P}_{\kappa,\kappa}$.
\end{definition}

\begin{remark}
\label{5.17}
Let $l \in \{0, \ldots ,\kappa\}$. The function $\mathcal{P}_{l,\kappa}$ is well-defined, definable and holomorphic in $z$. Furthermore we have $(t,\mu_\kappa(t)) \in H_{l,\kappa}^*$ and $\mathcal{P}_{l,\kappa}(t,\mu_\kappa(t))=0$ for $t \in \pi(C)$.
\end{remark}

\begin{lemma}
\label{5.18}
Let $l \in \{0, \ldots ,\kappa\}$. We have $H_{l,\kappa}^*=H_l$ and if $l \neq \kappa$ then
	$$z_l(t,z)=\mathcal{P}_{l,\kappa}(t,z) + \sigma_le^{\mu_{\kappa-l-1,\kappa}(t)}$$
	for $(t,z) \in H_l$.
\end{lemma}

\begin{proof}
By definition of $H_{0,\kappa}^*$ we have $H_{0,\kappa}^*=H_0$. If $\kappa=0$ we are done with the proof since $\mu_0=\Theta_0$. So assume that $\kappa>0$. We do an induction on $l$. For $(t,z) \in H_{0,\kappa}^*$ we obtain since $\mu_\kappa=\Theta_0+\sigma_0e^{\mu_{\kappa-1,\kappa}}$ (see Definition~\ref{5.2})
\begin{align*}
	z_0(t,z)&=z-\Theta_0-\sigma_0e^{\mu_{\kappa-1,\kappa}(t)}+\sigma_0e^{\mu_{\kappa-1,\kappa}(t)}\\
	&= z-\mu_{\kappa}(t)+\sigma_0e^{\mu_{\kappa-1,\kappa}(t)}\\
	&=\mathcal{P}_{0,\kappa}(t,z) + \sigma_0e^{\mu_{\kappa-1,\kappa}(t)}.
\end{align*}

$l-1 \to l$: Note that $H_{l-1} =  H^*_{l-1,\kappa}$. Let $(t,z) \in H_{l-1}$. We get with the inductive hypothesis (since $l-1 \neq \kappa$)
\begin{align*}
	\sigma_{l-1}z_{l-1}(t,z) \in \mathbb{C}^- &\Leftrightarrow \Bigl(1+\frac{\sigma_{l-1}(z_{l-1}(t,z)-\sigma_{l-1}e^{\mu_{\kappa-l,\kappa}(t)})}{e^{\mu_{\kappa-l,\kappa}(t)}} \Bigl) \in \mathbb{C}^-\\ 
	&\Leftrightarrow \Bigl(1+\frac{\sigma_{l-1}\mathcal{P}_{l-1,\kappa}(t,z)}{e^{\mu_{\kappa-l,\kappa}(t)}} \Bigl) \in \mathbb{C}^-.
\end{align*}

So we obtain $H_l=H_{l,\kappa}^*$. Suppose $l \neq \kappa$. Let $(t,z) \in H_l$. We obtain with the inductive hypothesis and the fact that $\log(az) = \log(a)+\log(z)$ for $z \in \mathbb{C}^-$ and $a \in \mathbb{R}$
\begin{align*}
	z_l(t,z)&=\log(\sigma_{l-1}z_{l-1}(t,z))-\Theta_l(t)\\
	&=\log(\sigma_{l-1}(\mathcal{P}_{l-1,\kappa}(t,z)+\sigma_{l-1}e^{\mu_{\kappa-l,\kappa}(t)}))-\Theta_l(t)\\
	&=\log(\sigma_{l-1}\mathcal{P}_{l-1,\kappa}(t,z)+e^{\mu_{\kappa-l,\kappa}(t)}) - \Theta_l(t)\\
	&=\log \Bigl(1+\frac{\sigma_{l-1}\mathcal{P}_{l-1,\kappa}(t,z)}{e^{\mu_{\kappa-l,\kappa}(t)}} \Bigl) +\mu_{\kappa-l,\kappa}(t)- \Theta_l(t)\\
	&=\mathcal{P}_{l,\kappa}(t,z)+\mu_{\kappa-l,\kappa}(t)-\Theta_l(t)\\
	&=\mathcal{P}_{l,\kappa}(t,z)+\Theta_l(t)+\sigma_le^{\mu_{\kappa-l-1,\kappa}(t)}-\Theta_l(t)\\
	&=\mathcal{P}_{l,\kappa}(t,z)+\sigma_le^{\mu_{\kappa-l-1,\kappa}(t)}
\end{align*}
since $\mu_{\kappa-l,\kappa}(t) = \Theta_l(t)+\sigma_le^{\mu_{\kappa-l-1,\kappa}(t)}$ by Definition~\ref{5.2}. 
\end{proof}

\begin{corollary}
\label{5.19}
For $(t,z) \in H_\kappa = H_{\kappa,\kappa}^*$ we have $z_\kappa(t,z) = \mathcal{P}_\kappa(t,z)$.
\end{corollary}

\begin{proof}
Let $\kappa=0$. Then $\Theta_\kappa=\mu_\kappa$ and therefore 
$$z_\kappa(t,z)=z-\mu_\kappa(t)=\mathcal{P}_\kappa(t,z)$$
for $(t,z) \in H_\kappa$. So let $\kappa>0$. With Lemma~\ref{5.18} we have for $(t,z) \in H_\kappa$
\begin{align*}
	z_\kappa(t,z)&=\log(\sigma_{\kappa-1}z_{\kappa-1}(t,z)) - \Theta_\kappa(t)\\
	&=\log(\sigma_{\kappa-1}(\mathcal{P}_{\kappa-1,\kappa}(t,z)+\sigma_{k-1}e^{\mu_{0,\kappa}(t)})) - \Theta_\kappa(t)\\
	&=\log(\sigma_{\kappa-1}(\mathcal{P}_{\kappa-1,\kappa}(t,z)+\sigma_{k-1}e^{\Theta_\kappa(t)})) - \Theta_\kappa(t)\\
	&=\log(\sigma_{\kappa-1}\mathcal{P}_{\kappa-1,\kappa}(t,z)+e^{\Theta_\kappa(t)}) - \Theta_\kappa(t)\\
	&=\log \Bigl(1+\frac{\sigma_{\kappa-1}\mathcal{P}_{\kappa-1,\kappa}(t,z)}{e^{\Theta_\kappa(t)}}\Bigl) = \mathcal{P}_\kappa(t,z). \qedhere
\end{align*}
\end{proof}

Since $H_l \subset H$ for $l \in \{0, \ldots ,\kappa\}$ we often consider $\mathcal{P}_{l,\kappa}$ as a function on $H$. 

\begin{example}
\label{5.20}
Consider 
$$C:=\{(t,x) \in \mathbb{R} \times \mathbb{R} \mid t \in \textnormal{} ]0,1[ \textnormal{}, \textnormal{}\tfrac{1}{1+t}+e^{-t-1/t} < x < \tfrac{1}{1+t} + e^{-t/2-1/t}\}.$$
Let $\mathcal{Y}:=(y_0,y_1,y_2)$ and $\Theta:=(\Theta_0,\Theta_1,\Theta_2)$ be as in Example~\ref{5.10}. Note that $\text{sign}(\mathcal{Y}) = (1,-1,-1)$. We have 
$$H=\{(t,z) \in \pi(C) \times \mathbb{C} \mid \tfrac{1}{1+t}+e^{-1-1/t} < z < \tfrac{1}{1+t}+e^{-t/2-1/t} \textnormal{ if }z \in \mathbb{R}\}.$$
For $(t,z) \in H$ we have
$$z_0(t,z)=z-\Theta_0(t), \textnormal{} z_1(t,z)=\log(z_0(t,z))-\Theta_1(t),$$
$$z_2(t,z)=\log(-z_1(t,z))-\Theta_2(t)$$
and
$$\mathcal{P}_{0,1}(t,z)=z-\mu_1(t), \textnormal{ } \mathcal{P}_{1,1}(t,z)=\log\Bigl(1+\frac{z-\mu_1(t)}{e^{\mu_{0,1}(t)}}\Bigl),$$
$$\mathcal{P}_{0,2}(t,z)=z-\mu_2(t), \textnormal{ } \mathcal{P}_{1,2}(t,z)=\log \Bigl(1+\frac{z-\mu_2(t)}{e^{\mu_{1,2}(t)}} \Bigl),$$
$$\mathcal{P}_{2,2}(t,z)=\log \Bigl(1-\frac{\mathcal{P}_{1,2}(t,z)}{e^{\mu_{0,2}(t)}} \Bigl)$$
where $\mu_0,\mu_1,\mu_2,\mu_{0,1}, \mu_{0,2}$ and $\mu_{1,2}$ are as in Example~\ref{5.10}.
\end{example}

The following property of $\mathcal{P}_{l,\kappa}$ is also needed later when dealing with parameterized integrals in Section~8. 

\begin{lemma}
\label{5.21}
Let $l \in \{0, \ldots ,\kappa\}$. For $(t,z) \in H$ we have $(\prod_{j=l}^\kappa \sigma_j)\mathcal{P}_{l,\kappa}(t,z) \in \mathbb{C}^-$ and $(\prod_{j=l}^\kappa \sigma_j) \mathcal{P}_{l,\kappa}(t,z) \in \mathbb{R}_{>0}$ if $z \in \mathbb{R}$.
\end{lemma}

\begin{proof}
We obtain with Corollary~\ref{5.19}
$$\sigma_\kappa \mathcal{P}_{\kappa,\kappa}(t,z)=\sigma_\kappa z_\kappa(t,z) \in \mathbb{C}^-$$
and if $z \in \mathbb{R}$ then $(t,z) \in D$ (compare with Lemma~\ref{5.13} and Lemma~\ref{5.9}) and therefore $\sigma_\kappa \mathcal{P}_{\kappa,\kappa}(t,z)=\sigma_\kappa g_\kappa(t,z) \in \mathbb{R}_{>0}$ by construction of $D$ and Lemma~\ref{5.12}(2). So assume $l \neq \kappa$. With Lemma~\ref{5.18}
$$\Bigl(\prod_{j=l}^\kappa \sigma_j \Bigl)\mathcal{P}_{l,\kappa}(t,z)=\Bigl(\prod_{j=l}^\kappa \sigma_j \Bigl)z_l(t,z)-\Bigl(\prod_{j=l}^\kappa \sigma_j \Bigl) \sigma_l e^{\mu_{\kappa-l-1,\kappa}(t)}.$$
Since $z_l(t,z) \in \mathbb{C} \setminus \mathbb{R}$ for $(t,z) \in \mathbb{C} \setminus \mathbb{R}$ we may assume $z \in \mathbb{R}$. Therefore $(t,z) \in D$. We show
$$\Bigl(\prod_{j=l}^\kappa \sigma_j \Bigl)\mathcal{P}_{l,\kappa}(t,z) \in \mathbb{R}_{>0}$$
by induction on $l \in \{0, \ldots ,\kappa-1\}$. For $l=0$ we obtain
$$\Bigl(\prod_{j=0}^\kappa \sigma_j \Bigl)\mathcal{P}_{0,\kappa}(t,z)=\Bigl(\prod_{j=0}^\kappa \sigma_j \Bigl) (z-\mu_\kappa(t)).$$
We get the following with Lemma~\ref{5.9}: If $\prod_{j=0}^\kappa \sigma_j=1$ then $\mu_\kappa(t)<z$ and therefore $z-\mu_\kappa(t) \in \mathbb{R}_{>0}$. 
If $\prod_{j=0}^\kappa \sigma_j=-1$ then $z<\mu_\kappa(t)$ and therefore $\Bigl(\prod_{j=0}^{\kappa} \sigma_j\Bigl)(z-\mu_\kappa(t)) = \mu_\kappa(t)-z \in \mathbb{R}_{>0}$.

$l-1 \to l:$ With the inductive hypothesis and Lemma~\ref{5.18} we obtain
\begin{align*}
\Bigl(\prod_{j=l-1}^\kappa \sigma_j \Bigl)\mathcal{P}_{l-1,\kappa}(t,z)&=\Bigl(\prod_{j=l-1}^\kappa \sigma_j \Bigl)z_{l-1}(t,z)-\Bigl(\prod_{j=l-1}^\kappa\sigma_j \Bigl) \sigma_{l-1}e^{\mu_{\kappa-l,\kappa}(t)} \\ 
&=\sigma_{l-1}\Bigl(\prod_{j=l}^\kappa \sigma_j \Bigl)z_{l-1}(t,z)-\Bigl(\prod_{j=l}^\kappa\sigma_j \Bigl) e^{\mu_{\kappa-l,\kappa}(t)} \in \mathbb{R}_{>0}.
\end{align*}
Assume $\prod_{j=l}^\kappa \sigma_j=1$. The case ''$\prod_{j=l}^\kappa \sigma_j=-1$'' is handled completely similarly. We obtain
$$\sigma_{l-1}z_{l-1}(t,z) > e^{\mu_{\kappa-l,\kappa}(t)}$$
and by taking logarithm
$$\log(\sigma_{l-1}z_{l-1}(t,z)) > \mu_{\kappa-l,\kappa}(t).$$
This gives with Definition~\ref{5.2}
$$z_l(t,z) > \mu_{\kappa-l,\kappa}(t) - \Theta_l(t)=\sigma_le^{\mu_{\kappa-l-1,\kappa}(t)}$$
and therefore
$$\Bigl(\prod_{j=l}^\kappa \sigma_j \Bigl) \mathcal{P}_{l,\kappa}(t,z)=\Bigl(\prod_{j=l}^\kappa \sigma_j \Bigl)z_l(t,z) - \Bigl(\prod_{j=l}^\kappa \sigma_j \Bigl) \sigma_le^{\mu_{\kappa-l-1,\kappa}(t)} \in \mathbb{R}_{>0}.$$
This finishes the proof.
\end{proof}

\subsection{Log-Analytically Prepared Functions}

In this subsection we compute a unary parametric global complexification $F:\Gamma \to \mathbb{C}$ of a log-analytically prepared function $f:C \to \mathbb{R}$ as follows. 
At first we determine a unary parametric global complexification $(\sigma\mathcal{Z})^{\otimes q}:H \to \mathbb{C}$ of $\vert{\mathcal{Y}}\vert^{\otimes q}$, where $\mathcal{Y}$ is the logarithmic scale occuring in the preparation of $f$ (as described in Subsection~5.1) and afterwards we extend the \emph{unit} $u=v \circ \phi$ by expanding $v$ to a convergent complex power series $V$ in a suitable way.\\

For this subsection we fix $r \in \mathbb{N}_0$ and a function $f: C \to \mathbb{R}$ which is $r$-log-analytically prepared in $x$. Let $(r,\mathcal{Y},a,q,s,v,b,P)$ be an LA-preparing tuple for $f$ with $b:=(b_1, \ldots ,b_s)$ and $P:=(p_1, \ldots ,p_s)^t$ for $p_i \in \mathbb{Q}^{r+1}$. Fix $R>1$ such that $v$ converges absolutely on an open neighborhood of $Q^s(0,R)$. Recall that 
\[
D^s(0,R)=\{(w_1, \ldots ,w_s) \in \mathbb{C}^s \mid \vert{w_j}\vert<R \textnormal{ for all }j \in \{1, \ldots ,s\}\}.
\] 

\begin{definition}
	\label{5.76}
	For the $r$-log-analytically prepared function $f$ with LA-preparing tuple $(r,\mathcal{Y},a,q,s,v,b,P)$ we define a set $\Gamma \subset \pi(C) \times \mathbb{C}$ and a function $F:\Gamma \to \mathbb{C}$ which is \textbf{complex $r$-log-analytically prepared in $z$ with center $\Theta$}\index{complex!log-analytically prepared} as follows. Fix $R>1$ such that $v$ converges absolutely on an open neighborhood of $\overline{Q^s(0,R)}$. For $l \in \{1, \ldots ,s\}$ set
	$$\Gamma_l:=\{(t,z) \in H \mid \vert{b_l(t)(\sigma\mathcal{Z})^{\otimes p_l}(t,z)}\vert < R\}.$$
	Set 
	$$\Gamma:=\bigcap_{l=1}^s \Gamma_l$$
	and consider
	$$F:\Gamma \to \mathbb{C}, (t,z) \mapsto$$ $$a(t)(\sigma\mathcal{Z})^{\otimes q}(t,z) V(b_1(t) (\sigma\mathcal{Z})^{\otimes p_1}(t,z), \ldots ,b_s(t)(\sigma\mathcal{Z})^{\otimes p_s}(t,z)),$$
	where $V$ is a complex power series in $s$ variables which converges absolutely on an open neighborhood of $\overline{D^s(0,R)}$ with $V|_{[-1,1]^s}=v$.
	
	
	We call $a$ \textbf{coefficient}, $b=(b_1, \ldots , b_s)$ a tuple of \textbf{base functions}\index{coefficient}\index{base functions}, and 
	$$(r,\mathcal{Z},a,q,s,V,b,P)$$
	where $\mathcal{Z}:=(z_0, \ldots ,z_r)$ 
	a \textbf{complex LA-preparing tuple}\index{complex LA-preparing tuple} for $F$.
\end{definition}



\begin{lemma}
	\label{5.77}
	Let $F:\Gamma \to \mathbb{C}$ be as in Definition~\ref{5.76}. Then the following holds.
	\begin{itemize}
	\item[(1)] The function $F$ is well-defined, definable and holomorphic in $z$. 
	\item[(2)] We have $C \subset \Gamma$ and $F|_C=f$. 
	\item[(3)] Let $E$ be a set of positive definable functions on $\pi(C)$ such that the functions $a,b_1, \ldots ,b_s$ and $\Theta_0, \ldots , \Theta_r$ can be constructed from $E$. Then $F$ can be constructed from $E$.
	\end{itemize} 
\end{lemma}

\begin{proof}
	(1): Follows by construction and Lemma~\ref{5.15} since $\Gamma_t$ is open for every $t \in \pi(C)$.\\
	(2): Follows from the fact that $C \subset H$, $(\sigma\mathcal{Z})^{\otimes p}(t,x)=\vert{\mathcal{Y}}\vert^{\otimes p}(t,x)$ for every $p \in \mathbb{Q}^{r+1}$ and every $(t,x) \in C$, and $\vert{b_j(t)}\vert \vert{\mathcal{Y}}\vert^{\otimes p_j}(t,x) \leq 1$ for every $j \in \{1, \ldots , s\}$ and every $(t,x) \in C$.\\ 
	(3): With Lemma~\ref{3.10}(2) and Lemma~\ref{5.12}(3) it is enough to show the following Claim~\ref{claim-prepared-log-an}.
	
	\begin{claim}
		\label{claim-prepared-log-an}
		Let $\eta:=(a,b_1, \ldots ,b_s)$. Then there is a globally subanalytic $J:\mathbb{C}^{s+1} \times \mathbb{C}^{r+1} \to \mathbb{C}$ such that
		$$F(t,z)=J(\eta(t,z),z_0(t,z), \ldots ,z_r(t,z))$$
		for every $(t,z) \in \Gamma$. 
	\end{claim}
	
	\begin{proof}
		Let $u:=(u_0, \ldots ,u_s)$ range over $\mathbb{C}^{s+1}$ and $w:=(w_0, \ldots ,w_r)$ range over $\mathbb{C}^{r+1}$. Let $(p_{00}, \ldots ,p_{0r}):=(q_0, \ldots , q_r) = q$. For $l \in \{0, \ldots ,s\}$ let
		$$\alpha_l:\mathbb{C}^{s+1} \times \mathbb{C}^{r+1} \to \mathbb{C}, (u,w) \mapsto $$
		$$\left\{\begin{array}{ll} u_l \prod_{j=0}^r (\sigma_jw_j)^{p_{lj}}, & \sigma_jw_j \in \mathbb{C}^- \textnormal{ for every } j \in \{0, \ldots ,r\}, \\
			0, & \textnormal{else.} \end{array}\right.$$
		Set
		$$J:\mathbb{C}^{s+1} \times \mathbb{C}^{r+1} \to \mathbb{C}, (u,w) \mapsto $$
		$$\left\{\begin{array}{ll} \alpha_0(u,w)V(\alpha_1(u,w), \ldots ,\alpha_s(u,w)), & \vert{\alpha_i(u,w)}\vert \leq R \textnormal{ for all } i \in \{1, \ldots ,s\}, \\
			0, & \textnormal{else}. \end{array}\right.$$
		Then $J$ is globally subanalytic. Note that $\sigma_lz_l(t,z) \in \mathbb{C}^-$ for every $l \in \{0, \ldots ,r\}$ and $(t,z) \in H$. Therefore 
		$$F(t,z)=J(\eta(t,z),z_0(t,z), \ldots ,z_r(t,z))$$
		for every $(t,z) \in \Gamma$. 
	\end{proof}
This finishes the proof since $z_0, \ldots , z_r$ can be also constructed from $E$ by Lemma~\ref{5.12}(3).
\end{proof}

\subsection{Log-Exp-Analytically-Prepared Functions}
In this subsection we compute a unary parametric global complexification of a restricted log-exp-analytically prepared function. This procedure is more sofisticated than in the case of a log-analytical preparation, since it involves iterations of logarithmic scales, power series and \emph{exponentials} of functions that are themselves log-exp-analytically prepared. Another challenge is that one could possibly lose definability when extending a real analytic exponential function to a holomorphic one with unbounded imaginary part.\\

For this section we set the following. Fix $e \in \mathbb{N}_0 \cup \{-1\}$, $r \in \mathbb{N}_0$ and a function $f:C \to \mathbb{R}$ which is log-exp-analytically $(e,r)$-prepared in $x$ with center $\Theta$ with respect to a finite set $E$ of positive definable functions such that every $g \in \log(E)$ is log-exp-analytically $(l,r)$-prepared in $x$ with respect to $E$ for $l \in \{-1, \ldots ,e-1\}$. For $g \in \log(E) \cup \{f\}$ we fix a preparing tuple $(r,\mathcal{Y},a_g,\exp(d_{0,g}),q_g,s,v_g,b_g,\exp(d_g),P_g)$ where $b_g:=(b_{1,g}, \ldots ,b_{s,g})$ and $\exp(d_g):=(\exp(d_{1,g}), \ldots ,\exp(d_{s,g}))$ (compare with Definition~\ref{4.9}). Fix a finite set $\mathcal{E}$ of $C$-heirs such that $\Theta$ and $a_g,b_g$ can be constructed from $\mathcal{E}$ for every $g \in \log(E) \cup \{f\}$. \\
Fix $R>1$ such that $v_g$ converges absolutely on an open neighborhood of $\overline{Q^s(0,R)}$ and a complex power series $V_g$ which converges absolutely on an open neighborhood of $\overline{D^s(0,R)}$ with $V_g|_{[-1,1]^s}=v_g$ for every $g \in \log(E) \cup \{f\}$. 

\begin{definition}
	\label{5.78}
	For a log-exp-analytically $(l,r)$-prepared function $g$ with $g \in \log(E) \cup \{f\}$ we define by induction on $l \in \{-1, \ldots ,e\}$ a set $\Lambda_g \subset \pi(C) \times \mathbb{C}$ \index{$\Lambda_g$} and a function $\Phi_g:\Lambda_g \to \mathbb{C}$ \index{$\Phi_g$} which is \textbf{complex log-exp-analytically $(l,r)$-prepared in $z$}\index{complex log-exp-analytically!$(l,r)$-prepared} and a \textbf{complex preparing tuple $\mathcal{J}_{\mathbb{C}}$}\index{complex!preparing tuple} for $\Phi_g$. To this preparation we assign a set $E_{\Phi_g}$\index{$E_{\Phi_g}$} of \textbf{complex definable functions} on $\Lambda_g$. \\
	
	{\bf Base case}: $l=-1$: We set $\Lambda_g=\pi(C) \times \mathbb{C}$, $\Phi_g=0$, $E_{\Phi_g}=\emptyset$ and $\mathcal{J}_{\mathbb{C}}=(0)$.\\
	
	{\bf Inductive Step}: Let
	$$(r,\mathcal{Y},a,e^{d_0},q,s,v,b,e^d,P):=(r,\mathcal{Y},a_g,e^{d_{0,g}},q_g,s,v_g,b_g,e^{d_g},P_g)$$
	where $(b_1, \ldots ,b_s):=(b_{1,g}, \ldots ,b_{s,g})$, $(e^{d_1}, \ldots ,e^{d_s}):=(e^{d_{1,g}}, \ldots ,e^{d_{s,g}})$, $P_g:=(p_1, \ldots ,p_s)^t$ for $p_1, \ldots , p_s \in \mathbb{Q}^{r+1}$ and for $d_i$ the corresponding $\Lambda_{d_i}$, the complex log-exp-analytically $(l-1,r)$-prepared $\Phi_{d_i}$ in $z$ and the corresponding set $E_{\Phi_{d_i}}$ have already been defined for $i \in \{0, \ldots ,s\}$. Set
	$$\Lambda_{0,g}:=\{(t,z) \in \Lambda_{d_0} \cap H \mid \vert{\textnormal{Im}(\Phi_{d_0}(t,z))}\vert<\pi\}$$
	and for $i \in \{1, \ldots ,s\}$ set
	$$\Lambda_{i,g}:=\{(t,z) \in \Lambda_{d_i} \cap H \mid \vert{b_i(t)(\sigma\mathcal{Z})^{\otimes p_i}(t,z)\exp(\Phi_{d_i}(t,z))}\vert < R,$$
	$$\vert{\textnormal{Im}(\Phi_{d_i}(t,z))}\vert<\pi\}.$$
	Set $\Lambda_g=\bigcap_{i=0}^s \Lambda_{i,g}$ and
	$$\Phi_g:\Lambda_g \to \mathbb{C}, (t,z) \mapsto a(t)(\sigma\mathcal{Z})^{\otimes q}(t,z) e^{\Phi_{d_0}(t,z)}$$
	$$V_g(b_1(t)(\sigma\mathcal{Z})^{\otimes p_1}(t,z)e^{\Phi_{d_1}(t,z)}, \ldots ,b_s(t)(\sigma\mathcal{Z})^{\otimes p_s}(t,z)e^{\Phi_{d_s}(t,z)}).$$
	We set
	$$E_{\Phi_g}=(E_{\Phi_{d_0}} \cup E_{\Phi_{d_1}} \cup  \ldots  \cup E_{\Phi_{d_s}} \cup \{e^{\Phi_{d_0}}, \ldots ,e^{\Phi_{d_s}}\}) |_{\Lambda_g}$$ 
	and 
	$$\mathcal{J}_{\mathbb{C}}=(r,\mathcal{Z},a,\exp(\Phi_{d_0}),q,s,V_g,b,\exp(\Phi_d),P)$$
	for $\Phi_g$ where $\exp(\Phi_d):=(\exp(\Phi_{d_1}), \ldots ,\exp(\Phi_{d_s}))$.
	\index{$\Lambda_g$}\index{$\Phi_g$}\index{complex! log-exp-analytically $(l,r)$-prepared}\index{complex!preparing tuple}\index{$E_{\Phi_g}$}
\end{definition}

\begin{remark}
	\label{rem:prepared-log-exp-vs-log-analytic}
	Suppose that $f$ is log-exp-analytically $(l,r)$-prepared in $x$ with respect to $E:=\{1\}$. Then $f$ is log-analytically prepared. In particular, $\Phi_f$ is complex log-analytically prepared in $z$.
\end{remark}

\begin{lemma}
	\label{5.79}
	Let $g \in \log(E) \cup \{f\}$ be log-exp-analytically $(l,r)$-prepared in $x$ with respect to $E$ for $l \in \{-1, \ldots ,e\}$. Then the following holds, showing among other things that $\Phi_g$ is a unary parametric global complexification of $g$.
	\begin{itemize}
	\item[(1)] The set $(\Lambda_g)_t$ is open for every $t \in \pi(C)$ and $\Phi_g$ is holomorphic in $z$.
	\item[(2)] We have $C \subset \Lambda_g$ and $\Phi_g|_C = g$.
	\item[(3)] The extension $\Phi_g$ can be constructed from $E_{\Phi_g} \cup \mathcal{E}$.
	\end{itemize}
\end{lemma}

\begin{proof}
	(1)/(2): We do an induction on $e$. For $e=-1$ it is clear since $\Lambda_g=\pi(C) \times \mathbb{C}$, $g=0$ and $\Phi_g=0$. So suppose that $e \geq 0$. By construction we see that $(\Lambda_g)_t$ is open since $(\sigma\mathcal{Z})^{\otimes \lambda}$ is continuous in $z$ for every $\lambda \in \mathbb{Q}^{r+1}$. By composition of holomorphic functions we also see with the inductive hypothesis that $\Phi_g$ is holomorphic in $z$ since $(\sigma\mathcal{Z})^{\otimes \lambda}$ is holomorphic in $z$ for every $\lambda \in \mathbb{Q}^{r+1}$ (by Lemma~\ref{5.15}), $V_g$ is a complex power series converging absolutely on $D^s(0,R)$ and the complex exponential function is also holomorphic. 
	Since $(\sigma\mathcal{Z})^{\otimes \lambda}|_C = \vert{\mathcal{Y}}\vert^{\otimes \lambda}$ for $\lambda \in \mathbb{Q}^{r+1}$, $V_g |_{[-1,1]^s}=v_g$,
	$$\vert{b_i(t)(\sigma\mathcal{Z})^{\otimes p_i}(t,x)\exp(\Phi_{d_i}(t,x))}\vert = \vert{b_i(t)\vert{\mathcal{Y}}\vert^{\otimes p_i}(t,x)\exp(d_i(t,x))}\vert \leq 1$$
	for $(t,x) \in C$ by the log-exp-analyical preparation and
	$$\textnormal{Im}(\Phi_{d_i}(t,x)) = \textnormal{Im}(d_i(t,x)) = 0 < \pi$$
	for $(t,x) \in C$ (due to $\Phi_{d_i}|_C=d_i$ by the inductive hypothesis) we see $C \subset \Lambda_g$ and $\Phi_g|_C=g$.\\
	(3): We do an induction on $l$. For $l=-1$ there is nothing to show.	
	$l-1 \to l:$ Let  
	$$(r,\mathcal{Y},a,e^{d_0},q,s,v,b,e^d,P):=(r,\mathcal{Y},a_g,e^{d_{0,g}},q_g,s,v_g,b_g,e^{d_g},P_g)$$
	where $(b_1, \ldots ,b_s):=(b_{1,g}, \ldots , b_{s,g})$, $(e^{d_1}, \ldots ,e^{d_s}):=(e^{d_{1,g}}, \ldots ,e^{d_{s,g}})$ and $P_g:=(p_1, \ldots ,p_s)^t$ for $p_1, \ldots , p_s \in \mathbb{Q}^s$. Note that $a$ and $b$ can be constructed from $\mathcal{E}$. Let
	$$(r,\mathcal{Z},a,\exp(D_0),q,s,V,b,\exp(D),P)$$
	be a complex preparing tuple for $\Phi_g$ where $D_j:=\Phi_{d_j}$ for $j \in \{0, \ldots ,s\}$ and
	$$\exp(D):=(\exp(D_1), \ldots ,\exp(D_s)).$$
	Note that $z_0, \ldots ,z_r$ can be constructed from $\mathcal{E}$ by Lemma~\ref{5.12}(3) since $\Theta_0, \ldots ,\Theta_r$ can be constructed from $\mathcal{E}$. By the inductive hypothesis we have that the functions $D_0, \ldots ,D_s$ can be constructed from $E_{D_0} \cup  \ldots  \cup E_{D_s} \cup \mathcal{E}$. So we obtain that $D_0|_{\Lambda_g}, \ldots ,D_s|_{\Lambda_g}$ can be constructed from $(E_{D_0} \cup  \ldots  \cup E_{D_s})|_{\Lambda_g} \cup \mathcal{E}$ and we see with Lemma~\ref{3.10}(1) that the functions $e^{D_0}|_{\Lambda_g}, \ldots ,e^{D_s}|_{\Lambda_g}$ can be constructed from $E_{\Phi_g} \cup \mathcal{E}$ since the imaginary part of $D_j$ is bounded by $\pi$ and $e^{D_j}|_{\Lambda_g} \in E_{\Phi_g}$ for $j \in \{0, \ldots ,s\}$. Let 
	$$\eta:\Lambda_g \to \mathbb{C} \times \mathbb{C}^s \times \mathbb{C}^{s+1}, (t,z) \mapsto (a(t),b_1(t), \ldots ,b_s(t),e^{D_0(t,z)}, \ldots ,e^{D_s(t,z)}).$$
	Then $\eta$ can be constructed from $E_{\Phi_g} \cup \mathcal{E}$. With the following Claim~\ref{claim-10} we are done with the proof due to Lemma~\ref{3.10}(2).
	
	\begin{claim}
		\label{claim-10}
		There is a globally subanalytic $J:\mathbb{C} \times \mathbb{C}^s \times \mathbb{C}^{s+1} \times \mathbb{C}^{r+1} \to \mathbb{C}$ such that
		$$\Phi_g(t,z)=J(\eta(t,z),z_0(t,z), \ldots ,z_r(t,z))$$
		for every $(t,z) \in \Lambda_g$. 
	\end{claim}
	
	\begin{proof}
		Let $u:=(u_0, \ldots ,u_s)$ and $y:=(y_0, \ldots ,y_s)$ range over $\mathbb{C}^{s+1}$. Let $(p_{00}, \ldots ,p_{0r}):=(q_0, \ldots , q_r)$. Let $w:=(w_0, \ldots ,w_r)$ range over $\mathbb{C}^{r+1}$. For $l \in \{0, \ldots ,s\}$ let
		$$\alpha_l:\mathbb{C}^{s+1} \times \mathbb{C}^{s+1} \times \mathbb{C}^{r+1} \to \mathbb{C}, (u,y,w) \mapsto $$
		$$\left\{\begin{array}{ll} u_ly_l \prod_{j=0}^r (\sigma_jw_j)^{p_{lj}}, & \sigma_jw_j \in \mathbb{C}^- \textnormal{ for every } j \in \{0, \ldots ,r\}, \\
			0, & \textnormal{else.} \end{array}\right.$$
		Set
		$$J:\mathbb{C}^{s+1} \times \mathbb{C}^{s+1} \times \mathbb{C}^{r+1} \to \mathbb{C}, (u,y,w) \mapsto $$
		$$\left\{\begin{array}{ll} \alpha_0(u,y,w)V(\alpha_1(u,y,w), \ldots ,\alpha_s(u,y,w)), & \vert{\alpha_i(u,y,w)}\vert \leq R \\
		&\textnormal{for all } i \in \{1, \ldots ,s\}, \\
			0, & \textnormal{else}. \end{array}\right.$$
		Then $J$ is globally subanalytic. Note that $\sigma_lz_l(t,z) \in \mathbb{C}^-$ for every $l \in \{0, \ldots ,r\}$ and $(t,z) \in H$. Therefore 
		$$\Phi_g(t,z)=J(\eta(t,z),z_0(t,z), \ldots ,z_r(t,z))$$
		for every $(t,z) \in \Lambda_g$. 
	\end{proof}
This finishes the proof. 
\end{proof}

To close Section~5 we give one example which summarizes most of the concepts which have been introduced so far. 

\begin{example}
\label{5.80}
Let $X=C:= \textnormal{}]0,1[^2$. Let 
$$E:=\{C \to \mathbb{R}, (t,x) \mapsto \exp(t \sqrt{-\log(x)})\}.$$
Then $E$ is a set of one positive definable function with $\log(E)=\{g\}$ where $g:C \to \mathbb{R}, (t,x) \mapsto t \sqrt{-\log(x)},$ is locally bounded, but not bounded at zero. Consider 
$$f:C \to \mathbb{R}, (t,x) \mapsto (1+1/t^2)x(-\log(x))^2\exp(t \sqrt{-\log(x)}).$$
Then $f$ is $(2,X)$-restricted log-exp-analytically $(1,1)$-prepared in $x$ with center $0$ (since $f$ is $(1,1)$-prepared in $x$ with respect to $E$, the underlying $1$-logarithmic scale is $\mathcal{Y}=(x,\log(x))$ with $\text{sign}(\mathcal{Y}) = (1,-1)$ and $v_f=v_g=1$). Note that $\mu_0=0$ and $\mu_1=1$ and that the change index of $\mathcal{Y}$ is $k^{\text{ch}}=0$. So we obtain $H=\textnormal ]0,1[ \times (\mathbb{C}^-_0 \cap \mathbb{C}^+_1)$ by Lemma~\ref{5.13}. We see that $\Lambda_g=H$ and 
$$\Lambda_f=\{(t,z) \in H \mid \vert{\text{Im}(t\sqrt{-\log(z)})}\vert < \pi\},$$
$\Phi_g:\Lambda_g \to \mathbb{C}, (t,z) \mapsto t \sqrt{-\log(z)},$ and
$$\Phi_f:\Lambda_f \to \mathbb{C}, (t,z) \mapsto (1+1/t^2)z(-\log(z))^2\exp(t \sqrt{-\log(z)}).$$
\end{example}

\section{Asymptotic Behavior} 

One major goal of this paper is to compute a \emph{unary parametric global complexification} of a real analytic restricted log-exp-analytic function $f:X \to \mathbb{R}$ which is also \emph{restricted log-exp-analytic}. The idea is to use Fact~\ref{4.12} for preparing $f$ cellwise restricted log-exp-analytically and then to "glue" the single unary parametric complexifications together computed in Section~5. One challenge during the "gluing" step is to respect restricted log-exp-analyticity. This necessitates a closer look at the exponential functions which occur in the unary parametric global complexification of a restricted log-exp-analytically prepared function (constructed in Definition~\ref{5.78} for a set $E$ of locally bounded functions in $(u,x)$). Lemma~\ref{5.79}(3) identifies two types of exponentials. The first type are $C$-heirs, defined as functions of the form $\exp(\tilde{\Theta}_j)$ where $\tilde{\Theta}_j$ is a component of a center of a logarithmic scale $\tilde{\mathcal{Y}}$ on $C$. Note that $\tilde{\Theta}_j$ may differ to the components of the center of the logarithmic scale $\mathcal{Y}$ occuring in the preparation of $f|_C$. The second consists of exponentials of functions that are themselves complex log-exp-analytically prepared in $z$, particulary of the form $G:=\exp(a(t)(\sigma\mathcal{Z})^{\otimes q}(t,z)\exp(\Phi(t,z))V(t,z))$ where $\Phi(t,z)$ is also complex log-exp-analytically prepared in $z$ of a lower order than $G$ and $V$ has a special form (compare with Definition~\ref{5.78} above). 
We have to show that all those exponentials can be "glued" together is such a way that the logarithm of the resulting function is locally bounded in $(u,z)$. 

For this purpose, it is necessary to analyze the asymptotic behavior of complex log-exp-analytically prepared functions on the boundary of their domain. 
This happens in this section. The following \emph{curve selection lemma}, as introduced by Van den Dries~\cite{11}, provides a suitable tool to this end: For any definable set $A \subset \mathbb{R}^n$ and every $x \in \overline{A}$ there is a definable continuous map $\gamma:\text{}]0,1[ \text{} \to A$ with $\lim_{y \searrow 0} \gamma(y) = x$. We call such a map $\gamma$ \emph{definable curve}. For such a curve $\gamma$, the limit $\lim_{y \searrow 0}h(\gamma(y))$ always exists in $\mathbb{R} \cup \{-\infty,\infty\}$ where $h:A \to \mathbb{R}$ is a definable function, because with the \emph{monotonicity theorem} in o-minimal structures (see Theorem 1.2 in~\cite{11}, Section 3) we see that $\lim_{y \searrow 0} g(y) \in \mathbb{R} \cup \{-\infty,\infty\}$ for every definable function $g: \textnormal{}]0,1[\textnormal{} \to \mathbb{R}$. With definable curves we can also investigate whether definable functions are bounded.

\begin{lemma}
	\label{lem:unboundedness-vs-curve-selection}
	Let $m \in \mathbb{N}$, $D \subset \mathbb{R}^m$ be non-empty and let $f:D \to \mathbb{R}$ be definable. Let $w \in \overline{D}$. Then $f$ is unbounded at $w$ if and only if there is a definable curve $\gamma:\text{}]0,1[\text{} \to D$ with $\lim_{y \searrow 0 } \vert{f(\gamma(y))}\vert = \infty$.
\end{lemma}

\begin{proof}
	Suppose that there is such a definable curve $\gamma$. Then for $U=\mathbb{R}^m$ the restriction $f|_{U \cap D}$ is unbounded since $\gamma(]0,1[) \subset D$. 
	
	Suppose that $f$ is unbounded at $x$. Then $x \in \overline{D} \setminus D$ and $D$ cannot be a single point. Hence, there is $c>0$ such that for all $0<r<c$ we have $S_r(x) \cap D \neq \emptyset$ where $S_r(x):=\{w \in \mathbb{R}^m \mid \vert{x-w}\vert = r\}$ (compare with Definition~\ref{14} and Fact~\ref{15}). Fix such a $c$. Let 
	\[
	b:\text{}]0,c[\text{} \to \mathbb{R}, r \mapsto \min\{\sup\{\vert{f(w)}\vert \mid w \in D \cap S_r(x)\}-1,1/r\}.
	\] 
	Then $b$ is definable. Now consider the definable set
	\[
	M:=\{(w,r) \in D \times \text{}]0,c[\text{} \mid f(w) \geq b(r)\}.
	\]
	Note that $(x,0) \in \overline{M}$ since for every $0<r<c$ we find $w \in S_r(x) \cap D$ with $f(w) \geq b(r)$ by definition of $b$. Hence, by curve selection there is a definable curve $\gamma:\text{}]0,1[\text{} \to M$ with $\lim_{y \searrow 0} \gamma(y) = (x,0)$. Since $\lim_{r \searrow 0} b(r)=\infty$, we also see that $\lim_{y \searrow 0} \vert{f(\gamma_x(y))}\vert=\infty$ where $\gamma_x$ denote the first $m$ components of $\gamma$. This concludes the proof.  
\end{proof}

Hence, we use definable curves in order to analyze the asymptotic behavior of a unary parametric global complexification of a restricted log-exp-analytically prepared function. However, we will focus on a special class of curves which we call \emph{compatible} since non-compatible curves are not needed for the proof of our main results below. 
We also conjecture that one can analyze general curves in a similar, albeit more sophisticated way.


This section is divided as follows. The first subsection is about compatible curves. The second one investigates the asymptotic behavior of unary parametric global complexifications of logarithmic scales under compatible curves and the third one extends these results on restricted log-exp-analytically prepared functions and their unary parametric global complexifications from Definition~\ref{5.78}. \\   

For the whole section we introduce the following notation: Let $n \in \mathbb{N}_0$. Let $t$ range over $\mathbb{R}^n$ and $x$ over $\mathbb{R}$. Let $z$ range over $\mathbb{C}$. Let $\pi:\mathbb{R}^n \times \mathbb{C}, (t,z) \mapsto t,$ be the projection on the first $n$ real coordinates. We fix a non-empty definable cell $C \subset \mathbb{R}^n \times \mathbb{R}_{\neq 0}$ and an $r$-logarithmic scale $\mathcal{Y}:=(y_0, \ldots ,y_r)$ on $C$ with center $\Theta:=(\Theta_0, \ldots ,\Theta_r)$. We fix $\textrm{sign}(\mathcal{Y}):=\sigma:=(\sigma_0, \ldots ,\sigma_r) \in \{-1,1\}^{r+1}$ and the change index $k^{\text{ch}}$ of $\mathcal{Y}$ (see Definition~\ref{5.7}). 

Let $y$ range over $\mathbb{R}$. For a definable curve $\gamma:\textnormal{}]0,1[\textnormal{} \to \pi(C) \times \mathbb{C}$ with $\gamma:=(\gamma_1, \ldots ,\gamma_{n+1})$ we set $\gamma_t:=(\gamma_1, \ldots ,\gamma_n)$ for the first $n$ real components, $\gamma_z:=\gamma_{n+1}$ for the last component and if $\gamma_{n+1}(y) \in \mathbb{R}$ for every $y \in \textnormal{}]0,1[$ we also write $\gamma_x$ instead of $\gamma_z$. \index{$\gamma_t$}\index{$\gamma_x$}\index{$\gamma_z$}

\subsection{Compatible Curves}
In this subsection we define compatible curves and give some elementary properties.

\begin{definition}
	\label{5.22}
	For $t \in \pi(C)$ we define the \textbf{length of $C$}\index{length of a cell} with respect to $x$ as \index{$L_C(t)$}$L_C(t):=\sup(C_t)-\inf(C_t)$.
\end{definition}

\begin{remark}
	\label{5.25}
	Let $\gamma:\textnormal{}]0,1[\textnormal{} \to \pi(C) \times \mathbb{C}$ be a definable curve. Then $$\lim_{y \searrow 0} L_C(\gamma_t(y)) \in \mathbb{R}_{\geq 0} \cup \{\infty\}.$$
\end{remark}

By the definition of a definable cell (see Definition~\ref{13}) we have the following.

\begin{remark}
	\label{rem:length}
	We have either $L_C(t)>0$ for every $t \in \pi(C)$ or $L_C(t)=0$ for every $t \in \pi(C)$.
\end{remark}

\begin{definition}
	\label{5.23}
	We say that $C$ is \textbf{fat}\index{fat} with respect to $x$ if $L_C(t)>0$ for every $t \in \pi(C)$ (i.e. $C_t$ is an open interval for every $t \in \pi(C)$). Otherwise we call $C$ \textbf{thin}\index{thin} with respect to $x$.
\end{definition}

\begin{definition}
	\label{5.26}
	Let $\gamma:\textnormal{}]0,1[\textnormal{} \to \pi(C) \times \mathbb{C}$ be a definable curve. We say that $\gamma$ is \textbf{compatible with $C$}\index{compatible} if $\lim_{y \searrow 0} \gamma(y) \in \mathbb{R}^n \times \mathbb{C}$, $\lim_{y \searrow 0} L_C(\gamma_t(y))>0$ and neither $\lim_{y \searrow 0} \sup(C_{\gamma_t(y)})=-\infty$ nor $\lim_{y \searrow 0} \inf(C_{\gamma_t(y)})=+\infty$.  
\end{definition}


In other words, $\gamma:\textnormal{}]0,1[\textnormal{} \to \pi(C) \times \mathbb{C}$ is compatible with $C$ if there is $0<\varepsilon<1$ such that $\gamma(]0,\varepsilon[)$ is bounded, $L_C(\gamma_t(]0,\varepsilon[))>\delta$ for a $\delta>0$,  $\{\sup(C_{\gamma_t(y)}) \mid y \in \text{}]0,\varepsilon[\}$ is bounded from below and $\{\inf(C_{\gamma_t(y)}) \mid y \in \text{} ]0,\varepsilon[\}$ is bounded from above.

\begin{remark}
	\label{rem:compatibility-fat}
	If there is a definable curve $\gamma: \text{}]0,1[\text{} \to \pi(C) \times \mathbb{C}$ compatible with $C$ then $C$ is fat with respect to $x$. 
\end{remark}

Definition~\ref{5.26} simplifies as follows for a definable curve in $C$.

\begin{lemma}
	\label{compatibility}
	A definable curve $\gamma: \textnormal{}]0,1[ \textnormal{} \to C$ is compatible with $C$ if \\
	$\lim_{y \searrow 0} L_C(\gamma_t(y))>0$ and $\lim_{y \searrow 0} \gamma(y) \in \mathbb{R}^n \times \mathbb{R}$.
\end{lemma}

\begin{proof}
	Since $\lim_{y \searrow 0} \gamma(y) \in \mathbb{R}^n \times \mathbb{R}$ we see that there is $\varepsilon>0$ such that the curve $\gamma$ is bounded on $]0,\varepsilon[$ and hence $c_1:=\inf\{\gamma_t(y) \mid y \in \text{} ]0,\varepsilon[\} \in \mathbb{R}$ and $c_2:=\sup\{\gamma_t(y) \mid y \in \text{} ]0,\varepsilon[\} \in \mathbb{R}$. Further we have $\gamma(]0,\varepsilon[) \subset C$ and therefore $\{\inf(C_{\gamma_t(y)}) \mid y \in \text{}]0,\varepsilon[)\}$ is bounded from above by $c_2$ and $\{\sup(C_{\gamma_t(y)}) \mid y \in \text{}]0,\varepsilon[)\}$ is bounded from below by $c_1$.
\end{proof}

\begin{figure}[h]
	\begin{center}\includegraphics[width=7.5cm,height=5cm,keepaspectratio]{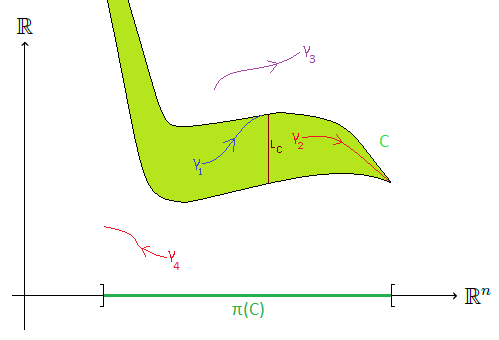}\end{center}
	\caption{The curves $\gamma_1$ and $\gamma_3$ are compatible with $C$ while the curves $\gamma_2$ and $\gamma_4$ aren't. (For $j \in \{1,2,3,4\}$ the arrows show $\gamma_j(y)$ for decreasing $y \in \textnormal{}]0,1[$.)}
\end{figure}

\begin{example}
	\label{ex:bridge-to-Tamm}
	If the cell $C$ is \emph{simple} (i.e. $C_t$ of the form $]0,d_t[$ for every $t \in \pi(C)$, see also Kaiser et al.~\cite{19} or Opris~\cite{28}) then for every $t \in \pi(C)$ the curve $\gamma:\text{}]0,1[\text{} \to C, y \mapsto d_t \cdot y$ is compatible with $C$. One could say that this curve is used in~\cite{19} and~\cite{28} to describe the asymptotic behavior of $\vert{\mathcal{Y}}\vert^{\otimes q}$ at $0$.  
\end{example}

In the following we give some elementary properties of definable curves $\gamma$ in $\pi(C) \times \mathbb{C}$ compatible with $C$ needed in several proofs below.

\begin{lemma}
	\label{5.27}
	Let $\gamma: \textnormal{} ]0,1[ \textnormal{} \to \pi(C) \times \mathbb{C}$ be a definable curve compatible with $C$. The following holds.
	
	\begin{itemize}
		\item [(1)]
		There is a definable function $\tilde{\gamma}_x: \textnormal{}]0,1[ \textnormal{} \to \mathbb{R}$ such that $\tilde{\gamma}:=(\gamma_t,\tilde{\gamma}_x)$ is a definable curve in $C$ compatible with $C$ and 
		$$\lim_{y \searrow 0} \textnormal{dist}(\tilde{\gamma}_x(y),\mathbb{R} \setminus C_{\gamma_t(y)})>0.$$
		So for every $l \in \{0, \ldots ,r\}$ we have
		$$\lim_{y \searrow 0} \vert{\tilde{\gamma}_x(y)-\mu_l(\gamma_t(y))}\vert > 0.$$ 	
		\item [(2)]
		Suppose that 
		$$\lim_{y \searrow 0} \textnormal{dist}(\gamma_z(y),C_{\gamma_t(y)})=0.$$
		Then there is a definable function $\hat{\gamma}_x: \textnormal{}]0,1[ \textnormal{} \to \mathbb{R}$ such that $\hat{\gamma}:=(\gamma_t,\hat{\gamma}_x)$ is a definable curve in $C$ compatible with $C$ and we have $$\lim_{y \searrow 0}(\gamma_t(y),\hat{\gamma}_x(y))=(t_0,x_0).$$
	\end{itemize}
\end{lemma}

\begin{proof}
	Throughout this proof let $f(t):=\inf(C_t)$ and $g(t):=\sup(C_t)$ for $t \in \pi(C)$. Note also that $C$ is fat with respect to $x$.\\
	
	(1): Since $\gamma$ is compatible with $C$ we find continuous definable functions $h_1:\pi(C) \to \mathbb{R}$ and $h_2:\pi(C) \to \mathbb{R}$ with $h_1<h_2$ such that
	$$\{(t,x) \in \pi(C) \times \mathbb{R} \mid h_1(t)<x<h_2(t)\} \subset C,$$
	$\lim_{y \searrow 0} h_j(\gamma_t(y)) \in \mathbb{R}$ for $j \in \{1,2\}$ and $\lim_{y \searrow 0} h_2(\gamma_t(y))-h_1(\gamma_t(y)) \in \mathbb{R}_{>0}$: Choose $h_1=f$ and $h_2=\min\{f+1,g\}$ if $\lim_{y \searrow 0}f(\gamma_t(y)) \neq -\infty$, $h_1=\max\{f,g-1\}$ and $h_2=g$ if $\lim_{y \searrow 0}g(\gamma_t(y)) \neq \infty$ and otherwise $h_1=\max\{f,0\}$ and $h_2=\min\{g,1\}$. Consider
	$$\tilde{\gamma}_x: \textnormal{}]0,1[ \textnormal{} \to \mathbb{R}, y \mapsto \frac{h_1(\gamma_t(y))+h_2(\gamma_t(y))}{2}.$$
	Then $]0,1[ \textnormal{} \to C, y \mapsto (\gamma_t(y),\tilde{\gamma}_x(y)),$ is a definable curve compatible with $C$ with
	$$\lim_{y \searrow 0} \textnormal{dist}(\tilde{\gamma}_x(y),\mathbb{R} \setminus C_{\gamma_t(y)}) \geq \lim_{y \searrow 0} \frac{h_2(\gamma_t(y))-h_1(\gamma_t(y))}{2} > 0.$$
	Consequently, $\lim_{y \searrow 0} \vert{\tilde{\gamma}_x(y)-\mu_j(\gamma_t(y))}\vert>0$ since $\mu_j(\gamma_t(y)) \notin C_{\gamma_t(y)}$ for every $y \in \textnormal{} ]0,1[$ and $j \in \{0, \ldots ,r\}$. \\
	
	(2): Let $\lim_{y \searrow 0} \gamma(y):=(t_0,x_0)$ where $x_0 \in \mathbb{R}$. 
	Then 
	$$\lim_{y \searrow 0} f(\gamma_t(y)) \leq x_0 \leq \lim_{y \searrow 0} g(\gamma_t(y)).$$
	Suppose that $x_0 = \lim_{y \searrow 0} g(\gamma_t(y))$. Consider
	\[
	\hat{\gamma}_x:\textnormal{}]0,1[\textnormal{} \to \mathbb{R}, y \mapsto g(\gamma_t(y))-y(g(\gamma_t(y))-\max\{f(\gamma_t(y)),g(\gamma_t(y))-1\}.
	\]
	Note that $\hat{\gamma}_x(y) \in C_{\gamma_t(y)}$ for $y \in \textnormal{}]0,1[$ since $ f(\gamma_t(y)) < \hat{\gamma}_x(y)<g(\gamma_t(y))$. Since $\lim_{y \searrow 0} \hat{\gamma}_x(y)=\lim_{y \searrow 0} g(\gamma_t(y))=x_0$ we see that $\hat{\gamma}_x$ does the job. The case "$x_0 = \lim_{y \searrow 0} f(\gamma_t(y))$" is handled completely similarly by considering
	$$\hat{\gamma}_x:\textnormal{}]0,1[\textnormal{} \to \mathbb{R}, y \mapsto f(\gamma_t(y))+y(\min\{g(\gamma_t(y)),f(\gamma_t(y))+1\}-f(\gamma_t(y))).$$
	Suppose that
	$$\lim_{y \searrow 0} f(\gamma_t(y)) < x_0 < \lim_{y \searrow 0} g(\gamma_t(y)).$$
	Then there are two definable functions $h_1:\pi(C) \to \mathbb{R}$ and $h_2:\pi(C) \to \mathbb{R}$ with $f < h_1 < h_2 < g$, $\lim_{y \searrow 0} h_1(\gamma_t(y)) \in \mathbb{R}$, $\lim_{y \searrow 0} h_2(\gamma_t(y)) \in \mathbb{R}$, and 
	$$\lim_{y \searrow 0} h_1(\gamma_t(y)) < x_0 < \lim_{y \searrow 0} h_2(\gamma_t(y)).$$
	(Choose $h_1 = \max\{f,x_0-1\}$ and $h_2 = \min\{g,x_0+1\}$.) Let $c \in \text{}]0,1[$ be such that
	$$\lim_{y \searrow 0}(h_1(\gamma_t(y)) + c(h_2(\gamma_t(y))-h_1(\gamma_t(y))) = x_0.$$
	Then 
	$$\hat{\gamma}_x:\textnormal{}]0,1[\textnormal{} \to \mathbb{R}, y \mapsto h_1(\gamma_t(y))+c(h_2(\gamma_t(y))-h_1(\gamma_t(y))),$$
	does the job since $h_1(\gamma_t(y)) \leq \hat{\gamma}_x(y) \leq h_2(\gamma_t(y))$ and hence, $\hat{\gamma}_x(y) \in C_{\gamma_t(y)}$ for every $y \in \text{}]0,1[$.
	
	Note that in all cases $\hat{\gamma}=(\gamma_t, \hat{\gamma}_x)$ is compatible with $C$ since $\gamma$ is compatible with $C$ and $\gamma_t=\hat{\gamma}_t$ (compare with Definition~\ref{5.26}).
\end{proof}

\subsection{Logarithmic Scales}

Now we analyze the asymptotic behavior of the unary parametric global complexification $(\sigma\mathcal{Z})^{\otimes q}$ of $\vert{\mathcal{Y}}\vert^{\otimes q}$ constructed in Section~5.1. 

\begin{lemma}
\label{5.28}
Let $l \in \{0, \ldots ,r\}$. Let $\gamma: \textnormal{}]0,1[ \textnormal{} \to \pi(C) \times \mathbb{C}$ be a definable curve compatible with $C$. Then the following holds.
	\begin{itemize}
		\item [$(I_l)$] $\lim_{y \searrow 0}\Theta_l(\gamma_t(y)) \in \mathbb{R}$.
		\item [$(II_l)$]  Suppose that $\gamma(]0,1[) \subset H$. Then 
		$\lim_{y \searrow 0} \vert{z_l(\gamma(y))}\vert = 0$ if and only if
		$\lim_{y \searrow 0} (\gamma_z(y)-\mu_l(\gamma_t(y)))=0$. We have $\lim_{y \searrow 0} \vert{z_l(\gamma(y))}\vert = \infty$ if and only if there is $j \in \{0, \ldots ,l-1\}$ such that $\lim_{y \searrow 0} (\gamma_z(y)-\mu_j(\gamma_t(y)))=0$.
\end{itemize}
\end{lemma}

\begin{proof}
Fix a definable curve $\tilde{\gamma}$ which fulfills the properties from Lemma~\ref{5.27}(1). By the monotonicity theorem we have that $\lim_{y \searrow 0} \Theta_l(\gamma_t(y)) \in \mathbb{R} \cup \{\pm \infty\}$. We show both statements by induction on $l$ simultaneously.

$l=0$: Assume $\lim_{y \searrow 0} \vert{\Theta_0(\gamma_t(y))}\vert=\infty$. So $\Theta_0 \neq 0$. By Definition~\ref{4.1}(d) there is $\varepsilon_0 \in \textnormal{}]0,1[$ such that
$$\vert{x-\Theta_0(t)}\vert < \varepsilon_0 \vert{x}\vert$$
for all $(t,x) \in C$ and we get
$$+\infty = \lim\limits_{y \searrow 0} \Bigl|1-\frac{\Theta_0(\gamma_t(y))}{\tilde{\gamma}_x(y)}\Bigl| \leq \varepsilon_0$$
since $\lim\limits_{y \searrow 0} \tilde{\gamma}_x(y) \in \mathbb{R}$ (due to compatibility). This is a contradiction and gives $(I_0)$. Since 
$$z_0(\gamma_z(y))=\gamma_z(y)-\Theta_0(\gamma_t(y)) = \gamma_z(y)-\mu_0(\gamma_t(y))$$
for $y \in \textnormal{} ]0,1[$ (by Definition~\ref{5.2}) we obtain $(II_0)$ (since the second part in $(II_0)$ cannot enter).

$l-1 \to l$: Assume $\lim_{y \searrow 0} \vert{\Theta_l(\gamma_t(y))}\vert=\infty$. So $\Theta_l \neq 0$. So by Definition~\ref{4.1}(d) there is $\varepsilon_l \in \textnormal{}]0,1[$ such that
$$\vert{\log(\sigma_{l-1}y_{l-1}(t,x))-\Theta_l(t)}\vert < \varepsilon_l \vert{\log(\sigma_{l-1}y_{l-1}(t,x))}\vert$$
for all $(t,x) \in C$. Therefore	
$$\lim\limits_{y \searrow 0} \Bigl|1-\frac{\Theta_l(\gamma_t(y))}{\log(\sigma_{l-1}y_{l-1}(\tilde{\gamma}(y)))}\Bigl| \leq \varepsilon_l.$$
This gives $\lim_{y \searrow 0} \sigma_{l-1}y_{l-1}(\tilde{\gamma}(y)) \in \{0,\infty\}$. With $(II_{l-1})$ we find $j \in \{0, \ldots ,l-1\}$ such that 
$$\lim_{y \searrow 0}(\tilde{\gamma}_x(y)-\mu_j(\gamma_t(y)))=0.$$
But this is a contradiction to the choice of $\tilde{\gamma}$. We obtain $(I_l)$.

$(II_l)$: We have $\lim_{y \searrow 0} \vert{z_l(\gamma(y))}\vert=\infty$ if and only if $\lim_{y \searrow 0} \vert{z_{l-1}(\gamma(y))}\vert \in \{0,\infty\}$ by $(I_l)$ if and only if there is $j \in \{0, \ldots ,l-1\}$ such that by $(II_{l-1})$
$$\lim_{y \searrow 0}(\gamma_z(y)-\mu_j(\gamma_t(y)))=0.$$
We see with $(I_a)$ for $a \in \{0, \ldots ,l\}$, Definition~\ref{5.2} and the continuity of the global exponential and logarithm function on the positive real line that 
\begin{align*}
\lim_{y \searrow 0} z_l(\gamma(y))=0  &\Leftrightarrow \lim_{y \searrow 0} \log(\sigma_{l-1}z_{l-1}(\gamma(y)))=\lim_{y \searrow 0} \Theta_l(\gamma_t(y)) \\
&\Leftrightarrow \lim_{y \searrow 0} z_{l-1}(\gamma(y))= \lim_{y \searrow 0} \sigma_{l-1}e^{\Theta_l(\gamma_t(y))} \\
&\Leftrightarrow \lim_{y \searrow 0} z_{l-1}(\gamma(y))=\lim_{y \searrow 0} \sigma_{l-1}e^{\mu_{0,l}(\gamma_t(y))}
\end{align*}
and for $j \in \{1, \ldots , l-1\}$ that
\begin{align*}
& \lim_{y \searrow 0}z_{l-j}(\gamma(y))=\lim_{y \searrow 0}\sigma_{l-j}e^{\mu_{j-1,l}(\gamma_t(y))} \\
\Leftrightarrow & \lim_{y \searrow 0}\log(\sigma_{l-j-1}z_{l-j-1}(\gamma(y)))=\lim_{y \searrow 0} (\Theta_{l-j}(\gamma_t(y)) + \sigma_{l-j}e^{\mu_{j-1,l}(\gamma_t(y))})\\
\Leftrightarrow &\lim_{y \searrow 0}\sigma_{l-j-1}z_{l-j-1}(\gamma(y))=\lim_{y \searrow 0} e^{(\Theta_{l-j}(\gamma_t(y)) + \sigma_{l-j}e^{\mu_{j-1,l}(\gamma_t(y))})}\\
\Leftrightarrow &\lim_{y \searrow 0}z_{l-j-1}(\gamma(y))=\lim_{y \searrow 0}\sigma_{l-j-1}e^{\mu_{j,l}(\gamma_t(y))}
\end{align*}
which implies by successively applying the latter equivalences for every $j \in \{1, \ldots , l-1\}$
\begin{align*}
\lim_{y \searrow 0} z_l(\gamma(y))=0 &\Leftrightarrow \lim_{y \searrow 0} z_{l-1}(\gamma(y))=\lim_{y \searrow 0} \sigma_{l-1}e^{\mu_{0,l}(\gamma_t(y))}\\
&\Leftrightarrow \lim_{y \searrow 0} z_0(\gamma(y))=\lim_{y \searrow 0} \sigma_0e^{\mu_{l-1,l}(\gamma_t(y))}\\
&\Leftrightarrow \lim_{y \searrow 0} \gamma_z(y)=\lim_{y \searrow 0} (\Theta_0(\gamma_t(y))+\sigma_0e^{\mu_{l-1,l}(\gamma_t(y))}) \\
&\Leftrightarrow \lim_{y \searrow 0} \gamma_z(y)= \lim_{y \searrow 0} \mu_l(\gamma_t(y)). \qedhere
\end{align*}
\end{proof}

An immediate consequence of Lemma~\ref{5.28} is the following.

\begin{corollary}
\label{5.29}
Let $\gamma: \textnormal{}]0,1[ \textnormal{} \to \pi(C) \times \mathbb{C}$ be a definable curve compatible with $C$. Let $l \in \{0, \ldots ,r\}$. The following holds.
	\begin{itemize}
		\item [(1)] $\lim_{y \searrow 0} \mu_{j,l}(\gamma_t(y)) \in \mathbb{R}$ for every $j \in \{0, \ldots ,l\}$.
		\item [(2)] $\lim\limits_{y \searrow 0} \vert{\mu_l(\gamma_t(y)) - \mu_j(\gamma_t(y))}\vert \in \mathbb{R}_{>0}$ for every $j \in \{0, \ldots ,r\}$ with $j \neq l$. 
	\end{itemize}
\end{corollary}

\begin{proof}
(1): This is clear with the continuity of the global exponential function, $(I_j)$ in Lemma~\ref{5.28} for $j \in \{0, \ldots ,l\}$ and Definition~\ref{5.2}.
 
(2): Let $j \in \{0, \ldots ,l-1\}$. With Definition~\ref{5.2} and $(I_a)$ in Lemma~\ref{5.28} for $a \in \{0, \ldots ,j-1\}$ we obtain that
\begin{align*}
& \lim\limits_{y \searrow 0} \vert{\mu_{l-a,l}(\gamma_t(y)) - \mu_{j-a,j}(\gamma_t(y))}\vert>0\\
\Leftrightarrow & \lim\limits_{y \searrow 0} \vert{\Theta_a + \sigma_ae^{\mu_{l-a-1,l}(\gamma_t(y))} - (\Theta_a + \sigma_ae^{\mu_{j-a-1,j}(\gamma_t(y))})}\vert>0\\
\Leftrightarrow & \lim\limits_{y \searrow 0} \vert{e^{\mu_{l-a-1,l}(\gamma_t(y))} - e^{\mu_{j-a-1,j}(\gamma_t(y))}}\vert>0\\
\Leftrightarrow &\lim\limits_{y \searrow 0} \vert{\mu_{l-a-1,l}(\gamma_t(y)) - \mu_{j-a-1,j}(\gamma_t(y))}\vert>0.
\end{align*}
So we have
$$\lim\limits_{y \searrow 0} \vert{\mu_l(\gamma_t(y)) - \mu_j(\gamma_t(y))}\vert>0$$
by applying the equivalences above for every $a \in \{0, \ldots ,j-1\}$ since 
$$\lim\limits_{y \searrow 0} \vert{\mu_{l-j,l}(\gamma_t(y)) - \mu_{0,j}(\gamma_t(y))}\vert = \lim\limits_{y \searrow 0} e^{\mu_{l-j-1,l}(\gamma_t(y))} > 0$$
with (1).
\end{proof}




The following examples show that $\lim_{y \searrow 0} \Theta_j(\gamma_t(y)) \notin \mathbb{R}$, $\lim_{y \searrow 0} \vert{\mu_{j,l}(\gamma_t(y))}\vert=\infty$ or $\lim_{y \searrow 0} \vert{\mu_l(\gamma_t(y))-\mu_j(\gamma_t(y))}\vert=0$ for $0 \leq j<l \leq r$ is indeed possible if $\gamma: \text{}]0,1[ \text{}\to \pi(C) \times \mathbb{C}$ is a definable curve not compatible with $C$ even if $C$ is fat with respect to $x$.

\begin{example}
	\label{5.30}
	The following holds.
	\begin{itemize}
		\item [(1)] Let $C$, $y_0,y_1$ and $\Theta_0,\Theta_1$ be as in Example~\ref{5.10}. Then $(y_0,y_1)$ is a $1$-logarithmic scale with center $(\Theta_0,\Theta_1)$ with $\lim_{t \searrow 0} \mu_{0,1}(t) = \lim_{t \searrow 0} \Theta_1(t)=-\infty$ and $\lim_{t \searrow 0} \mu_1(t)-\mu_0(t)=0$. However, $\lim_{t \searrow 0} L_C(t) = 0$.
		\item [(2)] Let 
		$$C:=\{(t,x) \in \textnormal{} ]0,1[ \textnormal{} \times \textnormal{} \mathbb{R} \mid 4e^{1/t} < x < (e^{1/2} + 3)e^{1/t}\}$$
		and consider $\Theta_0:=3e^{1/t}$ and $\Theta_1:=1/t$. An easy calculation shows that  there is a $1$-logarithmic scale $\mathcal{Y}$ on $C$ with center $\Theta:=(\Theta_0,\Theta_1)$ and $\text{sign}(\mathcal{Y})=(1,1)$. Although $\lim_{t \searrow 0} L_C(t)>0$  we see $\lim_{t \searrow 0} \Theta_0(t) = \infty$, $\lim_{t \searrow 0} \Theta_1(t)=\infty$, $\lim_{t \searrow 0} \mu_1(t) = \lim_{t \searrow 0} \Theta_0(t)+e^{\Theta_1(t)} = \infty$ and 
		$$\lim_{t \searrow 0} \vert{\mu_1(t)-\mu_0(t)}\vert = \lim_{t \searrow 0} \vert{\Theta_0(t)+e^{\Theta_1(t)}-\Theta_0(t)}\vert = \lim_{t \searrow 0} e^{\Theta_1(t)} = \infty.$$
		However, $\lim_{t \searrow 0} \inf(C_t) = \infty$ where $\inf(C_t)=4e^{1/t}$.\\
	\end{itemize}
\end{example}


We can slightly strengthen the above results, particularly Lemma~\ref{5.28}, if $C$ satisfies the following property.

\begin{definition}
	\label{5.45}
	We call $C$ \textbf{near with respect to $\mu_\kappa$} \index{near} if there is a definable curve $\gamma: \textnormal{} ]0,1[ \textnormal{} \to C$ compatible with $C$ such that
	$\lim_{y \searrow 0} (\mu_\kappa(\gamma_t(y))-\gamma_x(y))=0$. Otherwise we call $C$ \textbf{far with respect to $\mu_\kappa$}\index{far}.
\end{definition}

\begin{figure}[h]
	\begin{center}\includegraphics[width=9cm,height=6cm,keepaspectratio]{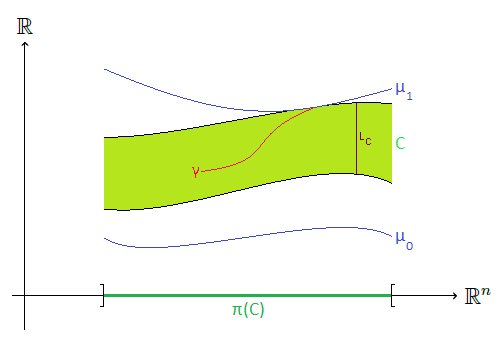}\end{center}
	\caption{The functions $\mu_j$ for a $1$-logarithmic scale $\mathcal{Y}$ with $\textnormal{sign}(\mathcal{Y}) = (1,-1)$ where $C$ is far with respect to $\mu_0$ and near with respect to $\mu_1$. (Here we have $k^{\textnormal{ch}} = 0$ and $r=1$.)}
\end{figure}


\begin{lemma}
	\label{5.47}
	Let $k:=k^{\textnormal{ch}} \geq 0$ and let $C$ be near with respect to $\mu_k$. Then $\Theta_l=0$ for every $l \in \{k+1, \ldots ,r\}$.
\end{lemma}

\begin{proof}
	Let $\gamma:\textnormal{}]0,1[\textnormal{} \to C$ be a definable curve compatible with $C$ and let
	$$\lim_{y \searrow 0}(\mu_k(\gamma_t(y))-\gamma_x(y))=0.$$
	Then $\lim_{y \searrow 0} y_k(\gamma(y))=0$ by $(II_k)$ in Lemma~\ref{5.28}. We do an induction on $l \in \{k+1, \ldots ,r\}$. 
	
	$l=k+1$: Assume $\Theta_{k+1} \neq 0$. Then there is $\varepsilon_{k+1} \in \textnormal{}]0,1[$ such that 
	$$\vert{\log(\sigma_ky_k(t,x))-\Theta_{k+1}(t)}\vert < \varepsilon_{k+1}\vert{\log(\sigma_ky_k(t,x))}\vert$$
	for all $(t,x) \in C$. We obtain
	$$\lim\limits_{y \searrow 0} \Bigl|1-\frac{\Theta_{k+1}(\gamma_t(y))}{\log(\sigma_ky_k(\gamma(y)))}\Bigl| \leq \varepsilon_{k+1}$$
	and therefore $\lim_{y \searrow 0} \vert{\Theta_{k+1}(\gamma_t(y))}\vert=\infty$, a contradiction to $(I_{k+1})$ in Lemma~\ref{5.28}.
	
	$l \to l+1$: Note that $\sigma_{k+1}=-1$ and $\sigma_{k+2}= \ldots =\sigma_l=1$. So we obtain
	$$\log(\sigma_ly_l)=\log_{l-k}(-\log(\sigma_ky_k))$$
	on $C$ by the inductive hypothesis. Assume $\Theta_{l+1} \neq 0$. Then there is $\varepsilon_{l+1} \in \textnormal{}]0,1[$ such that 
	$$\vert{\log_{l-k}(-\log(\sigma_ky_k(t,x)))-\Theta_{l+1}(t)}\vert < \varepsilon_{l+1}\vert{\log_{l-k}(-\log(\sigma_ky_k(t,x)))}\vert$$
	for all $(t,x) \in C$. We obtain
	$$\lim\limits_{y \searrow 0} \Bigl|1-\frac{\Theta_{l+1}(\gamma_t(y))}{\log_{l-k}(-\log(\sigma_ky_k(\gamma(y))))}\Bigl| \leq \varepsilon_{l+1}$$
	and therefore $\lim_{y \searrow 0} \vert{\Theta_{l+1}(\gamma_t(y))}\vert=\infty$, a contradiction to $(I_{l+1})$ in Lemma~\ref{5.28}.
\end{proof}

Now we determine the asymptotic behavior of $(\sigma\mathcal{Z})^{\otimes q}$ for $q \in \mathbb{Q}^{r+1}$ under definable curves in $H$ compatible with $C$. We start with two straightforward observations summarized in Lemma~\ref{5.31}. 

\begin{lemma}
\label{5.31}
Let $\gamma: \textnormal{}]0,1[ \textnormal{} \to H$ be a definable curve compatible with $C$. Then the following holds.
\begin{itemize}
	\item [(1)] Let $q \in \mathbb{Q}^{r+1}$ and suppose that $\lim_{y \searrow 0} \vert{(\sigma\mathcal{Z})^{\otimes q}(\gamma(y))}\vert \in \{0,\infty\}$. Then there is $l \in \{0, \ldots ,r\}$ with $\lim_{y \searrow 0} (\gamma_z(y)-\mu_l(\gamma_t(y)))=0$.

	\item [(2)] Let $\kappa=r$ or $\kappa \in \{k^{\text{ch}},r\}$ if $k^{\text{ch}} \geq 0$. Assume $\lim_{y \searrow 0} \gamma_z(y) - \mu_\kappa(\gamma_t(y))=0$. Then 
	$$\lim_{y \searrow 0} \mathcal{P}_{l,\kappa}(\gamma(y)) = 0$$ 
	for every $l \in \{0, \ldots ,\kappa\}$.
\end{itemize}
\end{lemma}

\begin{proof}
(1): This is clear with $(II_l)$ in Lemma~\ref{5.28} for $l \in \{0, \ldots ,r\}$.\\

(2): We do an induction on $l$. If $l=0$ then this is clear with Definition~\ref{5.16}.

$l-1 \to l$: By Definition~\ref{5.16} we obtain
$$\lim_{y \searrow 0}\mathcal{P}_{l,\kappa}(\gamma(y))=\lim_{y \searrow 0} \log \Bigl(1+\frac{\sigma_{l-1}\mathcal{P}_{l-1,\kappa}(\gamma(y))}{e^{\mu_{\kappa-l,\kappa}(\gamma_t(y))}} \Bigl)=0$$
by the inductive hypothesis, Corollary~\ref{5.29}(1) and the continuity of the global logarithm at $1$.
\end{proof}

With Lemma~\ref{5.31}(1) one sees immediately that $\lim_{y \searrow 0} \vert{(\sigma\mathcal{Z})^{\otimes q}(\gamma(y))}\vert \in \mathbb{R}_{>0}$ if $\gamma:\text{}]0,1[\text{} \to H$ is a definable curve compatible with $C$ with 
$$\lim_{y \searrow 0} \vert{(\gamma_z(y)-\mu_l(\gamma_t(y)))}\vert>0$$ for every $l \in \{0, \ldots , r\}$. Hence, it is left to analyze the asymptotic behavior of $(\sigma\mathcal{Z})^{\otimes q}$ when approaching the graph of a zero. Here we particularly focus on the graph of $\mu_\kappa$ where $\kappa=r$ or $\kappa \in \{k^{\text{ch}},r\}$ if $k^{\text{ch}} \geq 0$ since investigating other zeros is not needed for our purposes and can be done in a completely similar manner.

\begin{lemma}
	\label{lem:Asymptotics-neutral}
	Suppose that $(q_\kappa,\ldots,q_r) = 0$ and let $\gamma: \textnormal{}]0,1[ \textnormal{} \to H$ be a definable curve compatible with $C$ with $\lim_{y \searrow 0} \gamma_z(y) - \mu_\kappa(\gamma_t(y)) = 0$. Then $\lim_{y \searrow 0} \vert{(\sigma\mathcal{Z})^{\otimes q}(\gamma(y))}\vert \in \mathbb{R}_{>0}$.
\end{lemma}

\begin{proof}
	Suppose that $\lim_{y \searrow 0} (\sigma\mathcal{Z})^{\otimes q}(\gamma(y)) \in \{0,\infty\}$. Then by Lemma~\ref{5.31} (applied to $r=\kappa$) there is $l \in \{0, \ldots , \kappa-1\}$ such that 
	$$\lim_{y \searrow 0} \gamma_z(y)-\mu_l(\gamma_t(y)) = 0.$$ 
	But then $$\lim_{y \searrow 0} (\mu_l(\gamma_t(y))-\mu_\kappa(\gamma_t(y)))=0,$$
	a contradiction to Corollary~\ref{5.29}(2).
\end{proof}

In the following we determine the asymptotic behavior of $(\sigma\mathcal{Z})^{\otimes q}$ subject to $q$ when approaching the graph of $\mu_\kappa$ if $(q_\kappa, \ldots , q_r) \neq 0$.

\begin{definition}
\label{5.32}	
Assume $k:=k^{\text{ch}} \geq 0$. Let $q:=(q_0, \ldots , q_r) \in \mathbb{Q}^{r+1}$ with $(q_{k+1}, \ldots ,q_r) \neq 0$. Set\index{$l_{\textnormal{min}}(q)$}
$$l_{\min}(q) := \textnormal{min}\{j \in \{k+1, \ldots ,r\} \mid q_j \neq 0\}.$$  
\end{definition}

\begin{definition}
\label{5.33}	
Let $q:=(q_0, \ldots , q_r) \in \mathbb{Q}^{r+1}$ be with $(q_\kappa, \ldots ,q_r) \neq 0$. Let $k:=k^{\textnormal{ch}}$.\index{positive-$\kappa$} \index{negative-$\kappa$}
\begin{itemize}
	\item [(a)] Let $\kappa=r$. We call $q$ \textbf{$r$-positive} if $q_r>0$. We call $q$ \textbf{$r$-negative} if $q_r<0$.
	\item [(b)] Let $\kappa=k \geq 0$. We call $q$ \textbf{$k$-positive} if $q_k>0$ or if $q_k=0$ then $q_{l_{\min}}<0$. We call $q$ \textbf{$k$-negative} if $q_k<0$ or if $q_k=0$ then $q_{l_{\min}}>0$.
\end{itemize}
\end{definition}

\begin{remark}
\label{rem:kappa-positive-negative}
A vector $q:=(q_0, \ldots , q_r) \in \mathbb{Q}^{r+1}$ with $(q_\kappa, \ldots ,q_r) \neq 0$ is either $\kappa$-positive or $\kappa$-negative.
\end{remark}


\begin{lemma}
\label{5.34}
Let $q:=(q_0, \ldots ,q_r) \in \mathbb{Q}^{r+1}$ be with $(q_\kappa, \ldots , q_r) \neq 0$. Let $\gamma: \textnormal{}]0,1[ \textnormal{} \to H$ be a definable curve compatible with $C$. Suppose $\lim_{y \searrow 0} \gamma_z(y) - \mu_\kappa(\gamma_t(y)) = 0$. Then the following holds.
\begin{itemize}
	\item [(1)] If $q$ is $\kappa$-positive then $\lim_{y \searrow 0} \vert{(\sigma\mathcal{Z})^{\otimes q}(\gamma(y))}\vert = 0$.
	\item [(2)] If $q$ is $\kappa$-negative then $\lim_{y \searrow 0} \vert{(\sigma\mathcal{Z})^{\otimes q}(\gamma(y))}\vert = \infty$. 
\end{itemize}
\end{lemma}

\begin{proof}
If $\kappa=r$ we obtain with $(II_r)$ in Lemma~\ref{5.28} $\lim_{y \searrow 0} z_r(\gamma(y))=0$. With Lemma~\ref{lem:Asymptotics-neutral} we get the result. So let $0 \leq k:=k^{\text{ch}}$ and suppose $\kappa=k$. Let 
$$W:=\{(t,z) \in H \mid \vert{z_k(t,z)}\vert < 1/\exp_{r-k}(1)\},$$
$$L_1:W \to \mathbb{C}, (t,z) \mapsto \log(\vert{z_k(t,z)}\vert),$$
and for $l \in \{2, \ldots ,r-k\}$
$$L_l:W \to \mathbb{C}, (t,z) \mapsto \log_{l-1}(-\log(\vert{z_k(t,z)}\vert)).$$
Note that $L_l$ is well-defined and definable for $l \in \{1, \ldots ,r-k\}$. By considering a suitable subcurve of $\gamma$ if necessary we may assume that $\gamma(]0,1[) \subset W$ (since $\lim_{y \searrow 0} z_k(\gamma(y))=0$ by $(II_k)$ in Lemma~\ref{5.28}). Note also that $\lim_{y \searrow 0} L_1(\gamma(y))=-\infty$ and that $\lim_{y \searrow 0} L_l(\gamma(y))=\infty$ for $l \in \{2, \ldots ,r-k\}$.

\begin{claim}
\label{claim-1}
Let $l \in \{1, \ldots ,r-k\}$. There is a definable function $h_l:W \to \mathbb{C}$ with $\lim_{y \searrow 0} h_l(\gamma(y))=0$ such that 
$$z_{k+l}(t,z)=L_l(t,z)(1+h_l(t,z))$$
for every $(t,z) \in W$.
\end{claim}

\begin{proof}
We do an induction on $l$. Suppose $l=1$. We have for $(t,z) \in W$
$$z_{k+1}(t,z) = \log(\vert{z_k(t,z)}\vert) \Bigl(1-\frac{\Theta_{k+1}(t)-i\arg(\sigma_kz_k(t,z))}{\log(\vert{z_k(t,z)}\vert)} \Bigl).$$
So choose 
$$h_1:W \to \mathbb{C}, (t,z) \mapsto -\frac{\Theta_{k+1}(t)-i\arg(\sigma_kz_k(t,z))}{\log(\vert{z_k(t,z)}\vert)}.$$
We have $\lim_{y \searrow 0} h_1(\gamma(y))=0$ by $(I_{k+1})$ and $(II_k)$ in Lemma~\ref{5.28}. 

$l=2$: Note that $\sigma_{k+1}=-1$ (compare with Definition~\ref{5.7}). We have for $(t,z) \in W$ and due to $\log(az)=\log(a)+\log(z)$ for $z \in \mathbb{C}^-$ and $a \in \mathbb{R}_{>0}$
\begin{align*}
z_{k+2}(t,z) &= \log(\sigma_{k+1}z_{k+1}(t,z)) - \Theta_{k+2}(t)\\
&= \log(\sigma_{k+1}L_1(t,z)(1+h_1(t,z))) - \Theta_{k+2}(t)\\
&= \log(-L_1(t,z)) \left(1 + \frac{\log(1+h_1(t,z)) - \Theta_{k+2}(t)}{\log(-L_1(t,z))} \right).
\end{align*}
Choose
$$h_2:W \to \mathbb{C}, (t,z) \mapsto \frac{\log(1+h_1(t,z)) - \Theta_{k+2}(t)}{\log(-L_1(t,z))}.$$ 
Then $h_2$ is definable and by $(I_{k+2})$ in Lemma~\ref{5.28} we get $\lim_{y \searrow 0} h_2(\gamma(y))=0$. 

$l-1 \to l$: We may assume that $l > 2$. Note that $\sigma_{k+l-1}=1$ (again by Definition~\ref{5.7}). With the inductive hypothesis we obtain
\begin{align*}
z_{k+l}(t,z) &= \log(\sigma_{k+l-1}z_{k+l-1}(t,z)) - \Theta_{k+l}(t) \\
&= \log(\sigma_{k+l-1}L_{l-1}(t,z)(1+h_{l-1}(t,z))) - \Theta_{k+l}(t) \\
&=\log(L_{l-1}(t,z)) \left(1 + \frac{\log(1+h_{l-1}(t,z)) - \Theta_{k+l}(t)}{\log(L_{l-1}(t,z))} \right)
\end{align*}
for every $(t,z) \in W$. Choose 
$$h_l:H \to \mathbb{C}, (t,z) \mapsto \frac{\log(1+h_{l-1}(t,z)) - \Theta_{k+l}(t)}{\log(L_{l-1}(t,z))}.$$
Then $h_l$ is definable and by $(I_{k+l})$ in Lemma~\ref{5.28} we get $\lim_{y \searrow 0} h_l(\gamma(y))=0$.
\end{proof}

Let $h_l$ and $L_l$ be as in Claim~\ref{claim-1} for $l \in \{1, \ldots ,r-k\}$. Then we obtain for $y \in \textnormal{}]0,1[$
$$\vert{(\sigma\mathcal{Z})^{\otimes q}(\gamma(y))}\vert=d(\gamma(y))\vert{z_k(\gamma(y))}\vert^{q_k} \prod_{l=1}^{r-k} L_l(\gamma(y))^{q_{k+l}}(1+h_l(\gamma(y)))^{q_{k+l}}$$
where 
$$d(\gamma(y)):=\prod_{l=0}^{k-1} \vert{z_l(\gamma(y))}\vert^{q_l}.$$
Note that $\lim_{y \searrow 0} d(\gamma(y)) \in \mathbb{R}_{>0}$ by  Lemma~\ref{lem:Asymptotics-neutral}. Since $\lim_{y \searrow 0} \vert{z_k(\gamma(y))}\vert=0$ we see with the growth properties of the iterated logarithm that
$$\lim_{y \searrow 0} \vert{z_k(\gamma(y))}\vert^{q_k}\prod_{l=1}^{r-k} \vert{L_l(\gamma(y))}\vert^{q_{k+l}}=0$$
if $q$ is $k$-positive and 
$$\lim_{y \searrow 0} \vert{z_k(\gamma(y))}\vert^{q_k}\prod_{l=1}^{r-k} \vert{L_l(\gamma(y))}\vert^{q_{k+l}}=\infty$$
if $q$ is $k$-negative. 
So one sees that $\lim_{y \searrow 0} \vert{(\sigma\mathcal{Z})^{\otimes q}(\gamma(y))}\vert =0$ if $q$ is $k$-positive and $\lim_{y \searrow 0} \vert{(\sigma\mathcal{Z})^{\otimes q}(\gamma(y))}\vert = \infty$ if $q$ is $k$-negative. This proves Lemma~\ref{lem:Asymptotics-neutral}.
\end{proof}

An immediate consequence of the preceeding results is the following growth property for the derivative of $(\sigma\mathcal{Z})^{\otimes q}$ for $q \in \mathbb{Q}^{r+1}$ when approaching the graph of $\mu_\kappa$. (Compare with~\cite{19} and with~\cite{28} for corresponding results for logarithmic scales on \emph{simple} cells in the real setting.) 

\begin{corollary}
\label{5.35}
Let $q:=(q_0, \ldots , q_r) \in \mathbb{Q}^{r+1}$ with $(q_\kappa, \ldots ,q_r) \neq 0$. Let $j_\kappa(q):=\min\{j \in \{\kappa, \ldots ,r\} \mid q_j \neq 0\}$ and $$q_{\kappa,\textnormal{diff}}:=(0, \ldots ,0,q_\kappa-1, \ldots ,q_{j_\kappa(q)}-1,q_{j_\kappa(q)+1}, \ldots ,q_r) \in \mathbb{Q}^{r+1}.$$
Let $\gamma: \textnormal{}]0,1[ \textnormal{} \to H$ be a definable curve compatible with $C$ with $\lim_{y \searrow 0} \gamma_z(y)-\mu_\kappa(\gamma_t(y)) =0$. Then
$$\lim_{y \searrow 0} \Bigl|\frac{(\frac{d}{dz}(\sigma\mathcal{Z})^{\otimes q})(\gamma(y))}{(\sigma\mathcal{Z})^{\otimes q_{\kappa,\textnormal{diff}}}(\gamma(y))} \Bigl| \in \mathbb{R} \setminus \{0\}.$$
\end{corollary}

\begin{proof}
We have 
$$\frac{d}{dz}{\sigma_lz_l} = \frac{\sigma_l}{\prod_{j=0}^{l-1}z_j}$$
and hence
\begin{align*}
\frac{d}{dz}{(\sigma\mathcal{Z})^{\otimes q}} &= \sum_{j=0}^r q_j(\sigma_jz_j)^{q_j-1}\frac{\sigma_j}{\prod_{i=0}^{j-1}z_i} \prod_{i \neq j} (\sigma_iz_i)^{q_i}\\ 
&= \sum_{j=0}^r q_j \sigma_0 \cdot \ldots \cdot \sigma_j \left(\prod_{i =0}^j (\sigma_iz_i)^{q_i-1}\right) \prod_{i =j+1}^r (\sigma_iz_i)^{q_i}.
\end{align*}
Consider for $l \in \{0, \ldots ,r\}$ 
$$p_l^*:=(q_0-1, \ldots ,q_l-1,q_{l+1}, \ldots ,q_r)-q_{\kappa,\textnormal{diff}}.$$
Then 
\begin{align*}
\frac{(\frac{d}{dz}(\sigma\mathcal{Z})^{\otimes q})(t,z)}{(\sigma\mathcal{Z})^{\otimes q_{\kappa,\textnormal{diff}}}(t,z)} &= \sum_{j=0}^r q_j\sigma_0 \cdot \ldots \cdot \sigma_j (\sigma\mathcal{Z})^{\otimes p_j^*}(t,z) \\
&= q_{j_\kappa(q)} \sigma_0 \cdot \ldots \cdot \sigma_{j_\kappa(q)} (\sigma_0 z_0)^{q_0-1} \cdot \ldots \cdot (\sigma_{\kappa-1} z_{\kappa-1})^{q_{\kappa-1}-1} \\
&+ \sum_{j=0, j \neq j_\kappa(q)}^r q_j\sigma_0 \cdot \ldots \cdot \sigma_j (\sigma\mathcal{Z})^{\otimes p_j^*}(t,z).
\end{align*}
Hence we have to show for $l \neq j_\kappa(q)$ that $p_l^*$ is $\kappa$-positive if $q_l \neq 0$. 
Then we obtain the assertion with Lemma~\ref{lem:Asymptotics-neutral} and Lemma~\ref{5.34}(1). 

So suppose that $l \neq j_\kappa(q)$ and $q_l \neq 0$. If $l < \kappa$ then $(p_l^*)_\kappa=1$ and hence, $p_l^*$ is $\kappa$-positive. If $l=\kappa$ then $q_\kappa \neq 0$ and therefore $l=j_\kappa(q)$, a contradiction to $l \neq j_\kappa(q)$. So suppose that $0 \leq \kappa=k^{\textnormal{ch}}<l$. Since $q_l \neq 0$ and $j_\kappa(q) \neq l$ we have $j_\kappa(q)<l$ and therefore $(p_l^*)_\kappa = \ldots = (p_l^*)_{j_\kappa(q)} = 0$ and $(p_l^*)_{j_\kappa(q)+1} = -1$. Since $j_\kappa(q)+1>\kappa$ we see with Definition~\ref{5.33} that $p_l^*$ is $\kappa$-positive.
\end{proof}

\subsection{Restricted Log-Exp-Analytically Prepared \\
	Functions}
In this subsection we focus on the unary parametric global complexification $\Phi_f:\Lambda_f \to \mathbb{C}$ as constructed in Definition~\ref{5.78} above where $f:C \to \mathbb{R}$ is $(m+1,X)$-restricted log-exp-analytically prepared (see Definition~\ref{4.10} for the latter notion). 
Note that we already established the asymptotic behavior of $(\sigma\mathcal{Z})^{\otimes q}$ in the previous subsection. Therefore, we begin by examining the coefficient and base functions of $\Phi_f$, which also appear in $f$, and identify them with so-called \emph{$C$-consistent} functions with respect to $X$ (see Definition~\ref{5.48} below). These are functions depending only on $t$, but their asymptotic behavior can significantly depend on $X$. Surprisingly, functions such as the center of a logarithmic scale on $C$, $C$-heirs, the functions $\mu_{l}$ and $\mu_{j,l}$ from Definition~\ref{5.2} belong to this class.

Afterwards, we give some results for the asymptotic behavior of $\Phi_f$.\\

All the definitions and results are given in the parametric setting. So we fix $l,m \in \mathbb{N}_0$ such that $l+m=n$. Let $w:=(w_1, \ldots ,w_l)$ range over $\mathbb{R}^l$ and $u:=(u_1, \ldots ,u_m)$ over $\mathbb{R}^m$. Here $(u,x)$ and $(u,z)$ are serving as tuples of independent
variables of families of functions parameterized by $w$. Note that $t=(w,u)$ what we often use.

Let $\pi_l:\mathbb{R}^l \times \mathbb{R}^m \times \mathbb{R} \to \mathbb{R}^l, (w,u,x) \mapsto w$, be the projection on the first $l$ coordinates\index{$\pi_l$}. For $w \in \pi_l(C)$ and a definable curve $\gamma: \textnormal{}]0,1[ \textnormal{} \to \mathbb{R}^m \times \mathbb{C}$ we set $\gamma_u:=(\gamma_1, \ldots ,\gamma_m)$ for the first $m$ real components and $\gamma_z:=\gamma_{m+1}$ for the last complex component and if $\gamma_{m+1}(y) \in \mathbb{R}$ for every $y \in \textnormal{}]0,1[$ we also write $\gamma_x$ instead of $\gamma_z$\index{$\gamma_u$}\index{$\gamma_x$}\index{$\gamma_z$}. 

Further fix a definable $X \subset \mathbb{R}^n \times \mathbb{R}$ with $C \subset X$ and suppose that $X_w$ is open for every $w \in \mathbb{R}^l$.

The notion of compatibility (see Definition~\ref{5.26}) can be immediately transferred to the parametric setting.

\begin{definition}
	\label{5.37}
	We say that a definable curve $\gamma: \textnormal{}]0,1[ \textnormal{} \to \pi(C)_w \times \mathbb{C}$ is compatible with $C_w$\index{compatible!with $C_w$} if the curve $\hat{\gamma}: \textnormal{}]0,1[ \textnormal{} \to \pi(C) \times \mathbb{C}, y \mapsto (w,\gamma(y))$, is compatible with $C$.
\end{definition}

Now we are able to define $C$-consistent functions. 

\begin{definition}
	\label{5.48}
	Let $g:\pi(C) \to \mathbb{R}$ be definable. We call $g$ \textbf{$C$-consistent in $u$ with respect to $X$}\index{$C$-consistent} if the following holds. Let $w \in \pi_l(C)$. Then for every definable curve $\gamma: \textnormal{}]0,1[\textnormal{} \to C_w$ compatible with $C_w$ with $\lim_{y \searrow 0} \gamma(y) \in X_w$ we have
	$$\lim_{y \searrow 0} g(w,\gamma_u(y)) \in \mathbb{R}.$$
\end{definition}

In other words, a $C$-consistent function $g:\pi(C) \to \mathbb{R}$ in $u$ with respect to $X$ fulfills the following: If for $(w_0,u_0) \in \pi(X)$ there is an open neighborhood $U$ of $u_0$ in $\pi(X)_{w_0}$ and $\varepsilon,c_1,c_2 \in \mathbb{R}$ with $\delta>0$ such that $L_C(w_0,u)>\delta$, $\inf(C_{(w_0,u)})<c_1$ and $\sup(C_{(w_0,u)})>c_2$ for every $u \in \pi(C)_{w_0} \cap U$ then $g_{w_0}$ is bounded at $u_0$ (see also Lemma~\ref{lem:unboundedness-vs-curve-selection}).


\begin{figure}[h]
	\begin{center}\includegraphics[width=7.5cm,height=5cm,keepaspectratio]{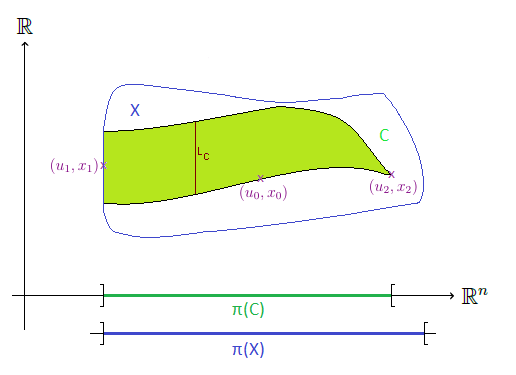}\end{center}
	\caption{Let $l=0$. A $C$-consistent function $f:\pi(C) \to \mathbb{R}$ in $u$ with respect to $X$ is bounded at $u_0$, but not necessarily bounded at $u_1$ respectively $u_2$.}
\end{figure}

The next example shows that the property of a function $g:\pi(C) \to \mathbb{R}$ to be $C$-consistent may depend on $C$ and $X$ even if $g$ depends only on $t$. 

\begin{example}
	\label{5.49}
	Let $C_1:= \textnormal{}]0,1[\textnormal{} \times \textnormal{}]0,1[$ and
	$$C_2:=\{(u,x) \in \textnormal{} ]0,1[\textnormal{} \times \textnormal{} \mathbb{R} \mid 0<x<u\}.$$
	Then $g: \text{}]0,1[ \textnormal{}\to \mathbb{R}, u \mapsto 1/u,$ is $(C_1)$-consistent in $u$ with respect to $\mathbb{R}_{>0} \times \mathbb{R}$, but not $(C_1)$-consistent in $u$ with respect to $\mathbb{R}^2$. Note that $g$ is $(C_2)$-consistent in $u$ with respect to \emph{any} open $X \subset \mathbb{R}^2$ with $C_2 \subset X$ since $\lim_{u \searrow 0} L_{C_2}(u)=0$.
\end{example}

In the following we give some elementary properties of $C$-consistent functions in $u$ with respect to $X$.

\begin{remark}
	\label{5.50}
	Let $Y \subset \mathbb{R}^n \times \mathbb{R}$ be definable with $X \subset Y$ such that $Y_w$ is open for every $w \in \mathbb{R}^l$. A function $g:\pi(C) \to \mathbb{R}$ which is $C$-consistent in $u$ with respect to $Y$ is also $C$-consistent in $u$ with respect to $X$.
\end{remark}

\begin{example}
	\label{5.51}
The following properties hold.
	\begin{itemize}
		\item [(1)] The functions $\Theta_0, \ldots ,\Theta_r$ are $C$-consistent in $u$ with respect to $\mathbb{R}^n \times \mathbb{R}$.
		\item [(2)] Let $l \in \{0, \ldots ,\kappa\}$. The function $\pi(C) \mapsto \mathbb{R}, t \mapsto \mu_{j,l}(t),$ is $C$-consistent in $u$ with respect to $\mathbb{R}^n \times \mathbb{R}$ for $j \in \{0, \ldots ,l\}$.
		\item [(3)] The logarithm of a $C$-heir is $C$-consistent in $u$ with respect to $\mathbb{R}^n \times \mathbb{R}$.
	\end{itemize}
\end{example}

\begin{proof}
	(1): Follows from $(I_j)$ in Lemma~\ref{5.28} for $j \in \{0, \ldots ,r\}$.
	
	(2): Follows from Corollary~\ref{5.29}(1). 
	
	(3): Let $g:\pi(C) \to \mathbb{R}_{>0}$ be a $C$-heir. Note that there is $l \in \mathbb{N}$, an $l$-logarithmic scale $\hat{\mathcal{Y}}$ with center $(\hat{\Theta}_0, \ldots ,\hat{\Theta}_l)$ such that $g=\exp(\hat{\Theta}_j)$ for $j \in \{1, \ldots ,l\}$. So this property follows from $(I_j)$ in Lemma~\ref{5.28} applied to $\hat{\mathcal{Y}}$.
\end{proof}

\begin{remark}
	\label{5.52}
	The following properties hold.
	\begin{itemize}
		\item [(1)] The set of $C$-consistent functions in $u$ with respect to $X$ is an $\mathbb{R}$-algebra with respect to pointwise addition and multiplication. 
		\item [(2)] The set of non-negative $C$-consistent functions in $u$ with respect to $X$ is a divisible monoid with respect to pointwise multiplication.
		\item [(3)] Let $q \in \mathbb{Q}$. Let $f:\pi(C) \to \mathbb{R}$ be $C$-consistent in $u$ with respect to $X$. Then $\exp(q \cdot f)$ is $C$-consistent in $u$ with respect to $X$. 
	\end{itemize}
\end{remark}

Now we investigate restricted log-exp-analytically prepared functions. For the rest of this subsection fix $e \in \mathbb{N}_0 \cup \{-1\}$, $r \in \mathbb{N}_0$, and a function $f:C \to \mathbb{R}, (w,u,x) \mapsto f(w,u,x),$ which is $(m+1,X)$-restricted log-exp-analytically $(e,r)$-prepared. Hence, by Definition~\ref{4.10} we can also fix a set $E$ of positive definable functions on $C$ such that every $g \in \log(E)$ is locally bounded in $(u,x)$ with reference set $X$ and $f$ is log-exp-analytically $(e,r)$-prepared in $x$ with respect to $E$. For $g \in \log(E) \cup \{f\}$ let $\Phi_g:\Lambda_f \to \mathbb{C}$ be the corresponding unary parametric global complexification of $g$ from Definition~\ref{5.78}.

Based on the properties of $E$, we can show that the base functions of the preparation of $f$ are $C$-consistent in $u$ with respect to $X$. The same holds for the coefficient if $f$ is also locally bounded in $u$ with reference set $X$.

\begin{lemma}
	\label{5.53}
	Let $(r,\mathcal{Y},a,e^{d_0},q,s,v,b,e^d,P)$ be a preparing tuple for $f$ where $b:=(b_1, \ldots , b_s)$, $\exp(d):=(\exp(d_1), \ldots , \exp(d_s))$ and $P:=(p_1, \ldots, p_s)$ for $p_i \in \mathbb{Q}^{r+1}$. Then the following properties hold. 
	\begin{itemize}
		\item [(1)] Let $j \in \{1, \ldots ,s\}$. The base function $b_j$ is $C$-consistent in $u$ with respect to $X$. 
		\item [(2)] Assume that $f$ is locally bounded in $(u,x)$ with reference set $X$. The coefficient $a$ is $C$-consistent in $u$ with respect to $X$.
	\end{itemize} 
\end{lemma}

\begin{proof}
	For $(t,x) \in C$ let
	$$W(t,x):=v(b_1(t)\vert{\mathcal{Y}}\vert^{\otimes p_1}(t,x) e^{d_1(t,x)}, \ldots ,b_s(t)\vert{\mathcal{Y}}\vert^{\otimes p_s}(t,x) e^{d_s(t,x)}).$$
	For $(t,x) \in C$ we have
	$$f(t,x)=a(t)\vert{\mathcal{Y}}\vert^{\otimes q}(t,x) \exp(c(t,x)) W(t,x).$$
	Let $w \in \pi_l(C)$ and $\gamma: \textnormal{}]0,1[ \textnormal{}\to C_w$ be a definable curve compatible with $C_w$ such that $(u_0,x_0):=\lim_{y \searrow 0} \gamma(y) \in X_w$. 
	We find $\delta>0$ and $B>0$ such that $U:=Q^m(u_0,\delta) \times \textnormal{}]x_0-\delta,x_0+\delta[ \textnormal{} \subset X_w$, $\vert{c(w,u,x)}\vert<B$ and $\vert{d_j(w,u,x)}\vert<B$ for $j \in \{1, \ldots ,s\}$ and every $(u,x) \in U \cap C_w$. 
	Similarly as in the proof of Lemma~\ref{5.27}(1) we find a definable function $\tilde{\gamma}_x \textnormal{}: \textnormal{}]0,1[ \textnormal{} \to \mathbb{R}$ such that $\tilde{\gamma}:=(\gamma_u,\tilde{\gamma}_x)$ is a definable curve in $C_w \cap U$ compatible with $C_w$ and for every $l \in \{0, \ldots ,r\}$ we have 
	$$\lim_{y \searrow 0} \vert{\tilde{\gamma}_x(y)-\mu_l(w,\gamma_u(y))}\vert > 0.$$
	
	(1): Let $j \in \{1, \ldots ,s\}$. Note that
	$$\lim_{y \searrow 0} \vert{\mathcal{Y}}\vert^{\otimes p_j}(w,\tilde{\gamma}(y)) \in \mathbb{R} \setminus \{0\}$$
	by Lemma~\ref{5.31}(1). Since
	$$b_j(t) \vert{\mathcal{Y}}\vert^{\otimes p_j}(t,x)e^{d_j(t,x)} \in [-1,1]$$
	for every $(t,x) \in C$ (and $\vert{d_j(t,x)}\vert \leq B$ for every $(t,x) \in U$), we obtain $\lim_{y \searrow 0} b_j(w,\gamma_u(y)) \in \mathbb{R}$. 
	
	
	(2): Note that there is $\rho>1$ such that $V(t,x) \in \textnormal{}]1/\rho,\rho[$ for every $(t,x) \in C$ and that
	$$\lim_{y \searrow 0} \vert{\mathcal{Y}}\vert^{\otimes q}(w,\tilde{\gamma}(y)) \in \mathbb{R} \setminus \{0\}$$
	by Lemma~\ref{5.31}(1). Since $\lim_{y \searrow 0} f(w,\tilde{\gamma}(y)) \in \mathbb{R}$ (and $\vert{c(t,x)}\vert \leq B$ for every $(t,x) \in U$), we obtain
	\begin{align*}
	&\lim_{y \searrow 0} a(w,\gamma_u(y)) \in \mathbb{R}. \qedhere
	\end{align*}
	
\end{proof}

Note that Lemma~\ref{5.53} transfers directly to $\Phi_f:\Lambda_f \to \mathbb{C}$ since its coefficient and base functions coincide with the coefficient and base functions of $f$. 

However, it may happen that $\lim_{y \searrow 0} \vert{g(\gamma_u(y))}\vert = \infty$ for a $C$-consistent $g:\pi(C) \to \mathbb{R}$ in $u$ with respect to $X$ and a definable curve $\gamma: \text{}]0,1[\text{} \to (\Lambda_f)_w$ with $(u_0,x_0):=\lim_{y \searrow 0}\gamma(y) \in X_w$ as the following example shows.

\begin{example}
	\label{ex:consistent-not-transfer}
	Let $C:=\text{} ]0,1[^2$, 
	$$X:=\{(u,x) \in \mathbb{R}_{>0} \times \mathbb{R}_{>0} \mid \text{if } x \leq 2 \text{ then }0<u<1\},$$
	and consider 
	$$f:C \to \mathbb{R}, (u,x) \mapsto g(u) \cdot e^{\sqrt{x}}$$
	where $g:\text{}]0,1[\text{} \to \mathbb{R}, u \mapsto 1/(1-u).$ Note that $X$ is open and $g$ is $C$-consistent in $u$ with respect to $X$ since 
	$$\{(1,x) \in \mathbb{R}^2 \mid 0 \leq x \leq 1\} \cap X = \emptyset$$
	and $g$ is only unbounded at $u=1$. The latter also shows that $f$ is  $(2,X)$-restricted $(1,0)$-log-exp-analytically prepared. Note also that 
	$$\Lambda_f=\{(u,z) \in \text{} ]0,1[ \text{}\times \mathbb{C}^- \mid \vert{\text{Im}(\sqrt{z})}\vert < \pi\}$$
	and $\Phi_f(u,z)=g(u) \cdot e^{\sqrt{z}}$ for $(u,z) \in \Lambda_f$. So consider $\gamma: \text{}]0,1[\text{} \to \Lambda_f, y \mapsto (1-y,3)$. Then $\gamma$ runs through $\Lambda_f$ and $\lim_{y \searrow 0}\gamma(y) = (1,3) \in X$, but $\lim_{y \searrow 0} g(\gamma_u(y))=\infty$. 
\end{example}

The following lemma shows that this feature does not occur if $x_0 \in \overline{C_w}$.

\begin{lemma}
	\label{lem:C-consistent-enlarged}
	Let $g:\pi(C) \to \mathbb{R}$ be $C$-consistent in $u$ with respect to $X$. Then for every definable curve $\gamma: \text{}]0,1[\text{} \to \pi(C)_w \times \mathbb{C}$ compatible with $C_w$ with $(u_0,x_0):=\lim_{y \searrow 0} \gamma(y) \in X_w \cap \overline{C_w}$ we have that $\lim_{y \searrow 0} g(w,\gamma_u(y)) \in \mathbb{R}$. 
\end{lemma}

\begin{proof}
	Note that
	$$\lim_{y \searrow 0} \textnormal{dist}(\gamma_z(y),C_{(w,\gamma_u(y))})=0.$$
	By Lemma~\ref{5.27}(2) there is a definable curve $\hat{\gamma}: \text{}]0,1[\text{} \to C_w$ compatible with $C_w$ with $\hat{\gamma}_u=\gamma_u$ and $\lim_{y \searrow 0} \hat{\gamma}_x(y) = x_0$. With Definition~\ref{5.48} we see
	\begin{align*}
	\lim_{y \searrow 0} g(w,\gamma_u(y)) &=  \lim_{y \searrow 0} g(w,\hat{\gamma}_u(y)) \in \mathbb{R}. \qedhere
	\end{align*}
\end{proof}

Outgoing from all these observations we finally describe the asymptotic behavior of $\Phi_f$ under a definable curve $\gamma:\text{}]0,1[\text{} \to \pi(C)_w \times \mathbb{C}$ compatible with $C_w$ where $\lim_{y \searrow 0} \gamma(y) \in \overline{C_w}$. We limit our discussions to statements needed later in Section~9 for proving the main results. This includes also considerations when $\gamma$ does not approach the graph of any zero $\mu_j$ for $j \in \{0, \ldots , r\}$. We think that even more similar and general results can be derived. 

\begin{lemma}
	\label{5.87}
	Let $g \in \log(E) \cup \{f\}$. Let $w \in \mathbb{R}^l$. Let $\gamma,\tilde{\gamma}:\textnormal{}]0,1[\textnormal{} \to (\Lambda_g)_w$ 
	be definable curves compatible with $C_w$ where $\gamma_u=\tilde{\gamma}_u$. Suppose that
	\begin{itemize}
		\item[(1)] $\lim_{y \searrow 0} (\gamma_z(y)-\tilde{\gamma}_z(y)) = 0$,
		\item[(2)] $\lim_{y \searrow 0} \vert{\gamma_x(y)-\mu_\kappa(w,\gamma_u(y))}\vert>0$
		for $\kappa=r$ and all $\kappa \in \{k^{\text{ch}},r\}$ if $k^{\text{ch}} \geq 0$,
		\item[(3)] $\lim_{y \searrow 0} \gamma(y) \in X_w \cap \overline{C_w}$.
	\end{itemize}
	Then we have for $g \in \log(E)$
	$$\lim_{y \searrow 0} \Phi_g(w,\gamma(y)) = \lim_{y \searrow 0} \Phi_g(w,\tilde{\gamma}(y)) \in \mathbb{R}.$$
	If $g=f$ then
	$$\lim_{y \searrow 0} \Phi_g(w,\gamma(y)) = \lim_{y \searrow 0} \Phi_g(w,\tilde{\gamma}(y)) \in \mathbb{R}$$
	if $\lim_{y \searrow 0} a(w,\gamma_u(y)) \in \mathbb{R}$ and
	$$\lim_{y \searrow 0} \vert{\Phi_g(w,\gamma(y))}\vert = \lim_{y \searrow 0} \vert{\Phi_g(w,\tilde{\gamma}(y))}\vert = \infty$$
	otherwise.  
\end{lemma}

\begin{proof}
	We do an induction on $l$. For $l=-1$ the statement is clear since $\Phi_g=0$. 
	
	$l-1 \to l:$ Let a preparing tuple for $g$ be
	$$(r,\mathcal{Y},a,e^{d_0},q,s,v,b,e^{d},P)$$
	where $b_g:=(b_1, \ldots ,b_s)$, $\exp(d_g):=(\exp(d_1), \ldots ,\exp(d_s))$ and $P:=(p_1, \ldots ,p_s)^t$ for $p_1, \ldots, p_s \in \mathbb{Q}^{r+1}$. Note that $\lim_{y \searrow 0} b_i(w,\gamma_u(y)) \in \mathbb{R}$ for $i \in \{1, \ldots ,s\}$ by Lemma~\ref{lem:C-consistent-enlarged} since $b_1, \ldots, b_s$ are $C$-consistent and that $\lim_{y \searrow 0} a(w,\gamma_u(y)) \in \mathbb{R}$ if $g \in \log(E)$ since $a$ is then also $C$-consistent. Let 
	$$(r,\mathcal{Z},a,\exp(D_0),q,s,V,b,\exp(D),P)$$
	be a complex preparing tuple for $\Phi_g$ where $D_j:=\Phi_{d_j}$ for $j \in \{0, \ldots ,s\}$ and $\exp(D):=(\exp(D_1), \ldots ,\exp(D_s))$. By the inductive hypothesis we have 
	$$\lim_{y \searrow 0} \exp(D_j(w,\gamma(y))) = \lim_{y \searrow 0} \exp(D_j(w,\tilde{\gamma}(y))) \in \mathbb{R}_{>0}$$
	for $j \in \{0, \ldots ,s\}$. By Lemma~\ref{5.27}(2) there is a definable curve $\hat{\gamma}:\text{}]0,1[\text{} \to C_w$ compatible with $C_w$ with $\hat{\gamma}_u=\gamma_u$ and $\lim_{y \searrow 0} \hat{\gamma}_x(y) = \lim_{y \searrow 0} \gamma_x(y)$. Note that
	$$\lim_{y \searrow 0} \vert{\hat{\gamma}_x(y)-\mu_\kappa(w,\gamma_u(y))}\vert>0$$ 
	for $\kappa=r$ and all $\kappa \in \{k^{\text{ch}},r\}$ if $k^{\text{ch}} \geq 0$ and therefore (since $\hat{\gamma}$ runs through $C_w$) 
	$$\lim_{y \searrow 0} \vert{\hat{\gamma}_x(y)-\mu_j(w,\gamma_u(y))}\vert>0$$  
	for every $j \in \{0, \ldots ,r\}$ by Lemma~\ref{5.6}(2) and Corollary~\ref{5.29}(2).
	\begin{claim}
		\label{claim-16}
		Let $\alpha \in \{0, \ldots ,r\}$. Then we have
		$$\lim_{y \searrow 0} \sigma_\alpha z_\alpha(w,\gamma(y)) = \lim_{y \searrow 0} \sigma_\alpha z_\alpha(w,\hat{\gamma}(y)) = \lim_{y \searrow 0} \sigma_\alpha y_\alpha(w,\hat{\gamma}(y)) \in \mathbb{R}_{>0}.$$
	\end{claim}
	
	\begin{proof}
		
		We show the statement by induction on $\alpha$. For $\alpha=0$ we have with $(I_0)$ and $(II_0)$ in Lemma~\ref{5.28} and by the observation that $\sigma_0z_0(t,z) = \sigma_0y_0(t,z) \in \mathbb{R}_{>0}$ for $(t,z) \in C$
		\begin{eqnarray*}
			\lim_{y \searrow 0} \sigma_0z_0(w,\gamma(y)) &=& \lim_{y \searrow 0}\sigma_0(\gamma_z(y)-\Theta_0(w,\gamma_u(y)))\\
			&=&  \lim_{y \searrow 0}\sigma_0(\hat{\gamma}_z(y)-\Theta_0(w,\gamma_u(y)))\\
			&=& \lim_{y \searrow 0} \sigma_0z_0(w,\hat{\gamma}(y)) = \lim_{y \searrow 0} \sigma_0y_0(w,\hat{\gamma}(y)) \in \mathbb{R}_{>0}.
		\end{eqnarray*}
		
		$\alpha-1 \to \alpha$: We have with the inductive hypothesis, $(I_\alpha)$ and $(II_\alpha)$ in Lemma~\ref{5.28}, by the observation that $\sigma_\alpha z_\alpha(t,z) = \sigma_\alpha y_\alpha(t,z) \in \mathbb{R}_{>0}$ for $(t,z) \in C$, and the continuity of the global logarithm
		\begin{align*}
			\lim_{y \searrow 0} \sigma_\alpha z_\alpha(w,\gamma(y)) &= \lim_{y \searrow 0}(\sigma_\alpha(\log(\sigma_{\alpha-1}z_{\alpha-1}(w,\gamma(y)))-\Theta_\alpha(w,\gamma_u(y))))\\
			&=  \lim_{y \searrow 0}(\sigma_\alpha(\log(\sigma_{\alpha-1}z_{\alpha-1}(w,\hat{\gamma}(y)))-\Theta_\alpha(w,\gamma_u(y))))\\
			&=  \lim_{y \searrow 0}\sigma_\alpha z_\alpha(w,\hat{\gamma}(y)) = \lim_{y \searrow 0}\sigma_\alpha y_\alpha(w,\hat{\gamma}(y)) \in \mathbb{R}_{>0}
		\end{align*}
	which finishes the proof of the claim.
	\end{proof}
	One can prove in a similar way 
	$$\lim_{y \searrow 0} \sigma_\alpha z_\alpha(w,\tilde{\gamma}(y)) = \lim_{y \searrow 0} \sigma_\alpha z_\alpha(w,\hat{\gamma}(y)) = \lim_{y \searrow 0} \sigma_\alpha y_\alpha(w,\hat{\gamma}(y)) \in \mathbb{R}_{>0}$$
	and therefore 
	$$\lim_{y \searrow 0} \sigma_\alpha z_\alpha(w,\gamma(y)) = \lim_{y \searrow 0} \sigma_\alpha z_\alpha(w,\tilde{\gamma}(y)) \in \mathbb{R}_{>0}.$$
	Thus we get with the continuity of the rational power function
	$$\lim_{y \searrow 0} (\sigma\mathcal{Z})^{\otimes \nu}(\gamma(y)) = \lim_{y \searrow 0}\vert{\mathcal{Y}}\vert^{\otimes \nu}(\tilde{\gamma}(y)) \in \mathbb{R}_{>0}$$
	for every $\nu \in \{q,p_1, \ldots ,p_s\}$. We are done since $V|_{[-1,1]^s}=v$,
	$$b_i(t)(\sigma\mathcal{Z})^{\otimes p_i}(t,z)\exp(D_i(t,z)) = b_i(t)\vert{\mathcal{Y}}\vert^{\otimes p_i}(t,z)\exp(d_i(t,z)) \in [-1,1]^s$$
	for $(t,z) \in C$ and $v([-1,1]^s)
	\subset \mathbb{R}_{>0}$.
\end{proof}



Finally we investigate the domain $\Lambda_f$ of the unary parametric global complexification $\Phi_f$ of $f$. 

\begin{lemma}
	\label{5.88}
	Let $g \in \log(E) \cup \{f\}$. Let $w \in \mathbb{R}^l$. Let $\gamma: \textnormal{}]0,1[ \textnormal{} \to (\Lambda_g)_w$ be a definable curve compatible with $C_w$ with
	\begin{itemize}
		\item[(1)] $\lim_{y \searrow 0} \vert{\gamma_x(y)-\mu_\kappa(w,\gamma_u(y))}\vert>0$
		for $\kappa=r$ and all $\kappa \in \{k^{\text{ch}},r\}$ if $k^{\text{ch}} \geq 0$,
		\item[(2)] $\lim_{y \searrow 0} \gamma(y) \in X_w \cap \overline{C_w}$.
	\end{itemize} 
	Then by passing to a suitable subcurve of $\gamma$ if necessary there is $\varepsilon>0$ such that for every $y \in \text{} ]0,1[$
	\[
	B(\gamma_z(y),\varepsilon) \subset (\Lambda_g)_{(w,\gamma_u(y))}
	\]
\end{lemma}

\begin{proof}
	Let $(u_0,x_0):=\lim_{y \searrow 0} \gamma(y)$. Assume that $g$ is log-exp-analytically $(l,r)$-prepared in $x$ with respect to $E$ for $l \in \{-1, \ldots ,e\}$. We do an induction on $l$. The assertion is clear for $l=-1$ (since $\Lambda_g=\pi(C) \times \mathbb{C}$).\\ 
	
	$l-1 \to l:$
	Let 
	$$(r,\mathcal{Y},a,e^{d_0},q,s,v,b,e^d,P)$$
	be a preparing tuple for $g$ where $b_g:=(b_1, \ldots ,b_s)$, $P:=(p_1, \ldots ,p_s)^t$ for $p_i \in \mathbb{Q}^{r+1}$ and $\exp(d_g):=(\exp(d_1), \ldots ,\exp(d_s))$. Let
	$$(r,\mathcal{Z},a,\exp(D_0),q,s,V,b,\exp(D),P)$$
	be a complex preparing tuple for $\Phi_g$ where $D_j:=\Phi_{d_j}$ for $j \in \{0, \ldots ,s\}$ and $\exp(D):=(\exp(D_1), \ldots ,\exp(D_s))$. We set for $j \in \{1, \ldots ,s\}$ and $(t,z) \in \Lambda_g$
	$$\phi_j(t,z):=b_j(t)(\sigma\mathcal{Z})^{\otimes p_j}(t,z)\exp(D_j(t,z)).$$
	Note that there is $R>1$ such that $V$ converges on an open neighborhood of $\overline{D^s(0,R)}$ and $\vert{\phi_j(t,x)}\vert \leq 1$ for $(t,x) \in C$. By definition of $\Lambda_g$ we also have $(w,\gamma(y)) \in \Lambda_{d_i}$ for $i \in \{0, \ldots ,s\}$ and every $y \in \textnormal{}]0,1[$. 
	
	Suppose that there is no such $\varepsilon>0$. 
	Then by curve selection there is a definable curve $\tilde{\gamma}:\textnormal{}]0,1[ \textnormal{} \to H_w$ with $\tilde{\gamma}_u=\gamma_u$, 
	$\lim_{y \searrow 0}\tilde{\gamma}_z(y)=x_0$ and $0<\delta<1$ such that $\tilde{\gamma}_z(y) \notin (\Lambda_g)_{(w,\gamma_u(y))}$ for $y \in \textnormal{} ]0,\delta[$.
	
	By shrinking $\delta$ if necessary we may assume with the induction hypothesis that $\tilde{\gamma}(y) \in (\Lambda_{d_i})_w$ for every $y \in \textnormal{}]0,\delta[$ and $i \in \{0, \ldots ,s\}$. By further shrinking $\delta$ if necessary we find $j^* \in \{0, \ldots ,s\}$ such that
	$$\vert{\textnormal{Im}(D_{j^*}(w,\tilde{\gamma}(y)))}\vert \geq \pi$$
	for every $y \in \textnormal{} ]0,\delta[$ or $j^+ \in \{1, \ldots ,s\}$ such that $\vert{\phi_{j^+}(w,\tilde{\gamma}(y))}\vert \geq R$ for every $y \in \textnormal{} ]0,\delta[$ (compare with the definition of $\Lambda_g$ in Definition~\ref{5.78}). By Lemma~\ref{5.27}(2) there is a definable curve $\hat{\gamma}: \textnormal{}]0,\delta[ \textnormal{} \to C_w$ compatible with $C_w$ with $\hat{\gamma}_u=\tilde{\gamma}_u$ and $\lim_{y \searrow 0} \hat{\gamma}(y)=(u_0,x_0)$. Since $\hat{\gamma}$ runs through $C_w$ we have $\vert{\phi_{j^+}(w,\hat{\gamma}(y))}\vert \leq 1$ for $y \in \textnormal{}]0,\delta[$ and therefore $\lim_{y \searrow 0} \vert{\phi_{j^+}(w,\hat{\gamma}(y))}\vert \leq 1$. 
	One sees with Lemma~\ref{5.87} that 
	\[
	\lim_{y \searrow 0} \vert{\phi_{j^+}(w,\tilde{\gamma}(y))}\vert = \lim_{y \searrow 0} \vert{\phi_{j^+}(w,\hat{\gamma}(y))}\vert \leq 1
	\]
	since $\phi_{j^+}$ is log-exp-analytically prepared with respect to $E$ and due to the $C$-consistency of $b$ (see Lemma~\ref{lem:C-consistent-enlarged}) $\lim_{y \searrow 0} b_{j^+}(w,\tilde{\gamma}_u(y)) \in \mathbb{R}$. So the latter cannot hold. Since $\hat{\gamma}$ runs through $C_w$ we obtain
	$$\textnormal{Im}(D_{j^*}(w,\hat{\gamma}(y))) =  \textnormal{Im}(d_{j^*}(w,\hat{\gamma}(y)))=0$$
	for $y \textnormal{} \in \textnormal{}]0,\delta[$ and so the former can also not hold by Lemma~\ref{5.87}, a contradiction.
\end{proof}

\section{A Complex Log-Analytical Preparation Theorem for Restricted Log-Exp-Analytically \\ Prepared Functions}

An essential step in the proof of our main results is constructing a unary \emph{high} parametric global complexification of a restricted log-exp-analytic function $f:X \to \mathbb{R}$. This involves finding a unary parametric global complexification of $f$ which fulfills also the properties stated in Definition~\ref{2.14}. This allows to expand unary parametric global complexifications to multivariate ones.\\   
The strategy is to use Fact~\ref{4.12} to obtain a decomposition $\mathcal{C}$ of $X$ into finitely many definable cells such that $f|_C$ is log-exp-analytically prepared for every $C \in \mathcal{C}$. Subsequently, we construct a unary parametric global complexification $\Phi_{f|_C} : \Lambda_{f_C} \to \mathbb{C}$ of $f|_C$ as stated in Definition~\ref{5.78} above. But to obtain the \emph{highness} after gluing all these together, it is necessary to extend $\Phi_{f|_C}$ to an even larger unary parametric global complexification of $f|_C$. 
For the latter step, we need definability results of parameterized integrals of the form
\[
\Psi: \Delta \to \mathbb{C}, (t,s,z) \mapsto \frac{1}{2 \pi i}\int_{\partial B(\mu_\kappa,s)} \frac{\Phi_{f|_C}(t,\xi)}{\xi-z} d \xi,
\]
where 
\[
\Delta:=\{(t,s,z) \in \pi(C) \times \mathbb{R}_{>0} \times \mathbb{C} \mid N(t) < s < M(t) \text{ and } z \in B(\mu_\kappa(t),s)\}.
\]
Here $N,M:\pi(C) \to \mathbb{R}_{\geq 0}$ are two
definable functions, and 
$\kappa=r$ or $\kappa \in \{k^{\textnormal{ch}},r\}$ if $k^{\text{ch}} \geq 0$ where $k^{\textnormal{ch}}$ denotes the change index of the $r$-logarithmic scale $\mathcal{Y}$ occuring in the preparation of $f|_C$.
(See Proposition 2.7 in Kaiser~\cite{17} for the corresponding parameterized integral in the globally subanalytic setting where $\mu_\kappa$ coincides with the center of the globally subanalytical preparation according to Definition~\ref{def:globally-subanalytical-preparation}, Definition~\ref{5.2} and Remark~\ref{rem:change:globally-subanalytical} above.) However, proving that $\Psi$ is definable is not straightforward. One has to deal with iterations of exponentials, logarithms (within $(\sigma\mathcal{Z})^{\otimes q}$) and power series occuring in $\Phi_{f|_C}$. Note also that $\mathbb{R}_{\text{an,exp}}$ is not closed under parametric integration (see Fact~\ref{21} above) and $\Phi_{f|_C}$ is not necessarily globally subanalytic. 
But on the other hand it suffices to establish the definability result for non-negative functions $M,N$ where $M$ belongs to the class of \emph{$\kappa$-persistent} functions and $N$ to the class of \emph{$\kappa$-soft} ones which we will introduce in Definition~\ref{5.38} below. Under a definable curve $\gamma$ in $C$ compatible with $C$ approaching the graph of $\mu_\kappa$, a $\kappa$-persistent function converges to a positive real number, while a $\kappa$-soft function converges to zero. 

Another crucial insight for proving the definability result of the parameterized integral above is that for $\Phi_{f|_C}$ one can find a $\kappa$-persistent $M:\pi(C) \to \mathbb{R}_{\geq 0}$ and a $\kappa$-soft $N:\pi(C) \to \mathbb{R}_{\geq 0}$ such that 
\[
D:=\{(t,z) \in H \mid N(t) < \vert{z-\mu_\kappa(t)}\vert < M(t)\} \subset \Lambda_{f|_C},
\]
and $\Phi_{f|_C}|_D(t,z)$ is complex $r$-log-analytically prepared in $z$. Hence, $\Phi_{f|_C}|_D(t,z)$ is of the form $a(t)(\sigma_0z_0(t,z))^{q_0} \cdot  \ldots  \cdot (\sigma_rz_r(t,z))^{q_r}U(t,z)$ where $z_0=z-\Theta_0(t), z_1=\log(\sigma_0z_0)-\Theta_1(t)$, \ldots , $\sigma_j \in \{-1,1\}$, the $q_j$'s are rational exponents and $U(t,z)$ is a function of a special form (compare with Definition~\ref{5.76} above). The coefficient, center and base functions of this log-analytical preparation belongs to a special class of functions which we call \emph{$C$-regular} (see Definition 5.43 below). $C$-Regular functions can be constructed from $C$-consistent functions and hence generalize $C$-nice functions. Consequently, one can assume that $\Phi_{f|_C}$ is complex log-analytically prepared in the integral above, which significantly simplifies the proof that $\Psi$ is definable since one can "neglect" the exponentials occuring in $\Phi_{f|_C}$. 

Therefore, the goals for this and the subsequent section are the following. In Subsection 7.1 we introduce and investigate $C$-regular functions, in Subsection 7.2 we look at $\kappa$-persistent and $\kappa$-soft functions and in Subsections 7.3 and 7.4 we give the desired result on complex log-analytical preparation. Section 8 is then devoted to the definability result of the parameterized integral stated above.\\

As in the section before we give all the necessary results and definitions in the parametric setting. So for the whole section we set the following. Let $n \in \mathbb{N}$, $t$ range over $\mathbb{R}^n$, $x$ over $\mathbb{R}$ and $z$ over $\mathbb{C}$. Let $\pi:\mathbb{R}^n \times \mathbb{C} \to \mathbb{R}^n, (t,z) \mapsto t,$ be the projection on the first $n$ real coordinates. Fix $l,m \in \mathbb{N}_0$ such that $l+m=n$. Let $w:=(w_1, \ldots ,w_l)$ range over $\mathbb{R}^l$ and $u:=(u_1, \ldots ,u_m)$ over $\mathbb{R}^m$. Here $(u,x)$ and $(u,z)$ are serving as tuples of independent variables of families of functions parameterized by $w$. Note that $t=(w,u)$ what we often use. 
Let $\pi_l:\mathbb{R}^l \times \mathbb{R}^m \times \mathbb{R} \to \mathbb{R}^l, (w,u,x) \mapsto w$, be the projection on the first $l$ coordinates\index{$\pi_l$}. For $w \in \pi_l(C)$ and a definable curve $\gamma: \textnormal{}]0,1[ \textnormal{} \to \mathbb{R}^m \times \mathbb{C}$ we set $\gamma_u:=(\gamma_1, \ldots ,\gamma_m)$ for the first $m$ real components and $\gamma_z:=\gamma_{m+1}$ for the last complex component and if $\gamma_{m+1}(y) \in \mathbb{R}$ for every $y \in \textnormal{}]0,1[$ we also write $\gamma_x$ instead of $\gamma_z$\index{$\gamma_u$}\index{$\gamma_x$}\index{$\gamma_z$}. 

Further let $X \subset \mathbb{R}^n \times \mathbb{R}$ be non-empty, definable and suppose that $X_w$ is open for every $w \in \mathbb{R}^l$. Let $C \subset X \cap (\mathbb{R}^n \times \mathbb{R}_{\neq 0})$ be a non-empty definable cell. Fix an $r$-logarithmic scale $\mathcal{Y}$ on $C$ with center $\Theta:=(\Theta_0, \ldots , \Theta_r)$ and fix $\kappa=r$ or $\kappa \in \{k^{\textnormal{ch}},r\}$ if $k^{\textnormal{ch}} \geq 0$ where $k^{\text{ch}}$ denotes the change index of $\mathcal{Y}$ (see Definition~\ref{5.7}). We also assume that $\Theta$ is $C$-nice.

\subsection{Regular Functions}

In this subsection we define $C$-regular functions as compositions of log-analytic functions and $C$-consistent functions 
and give some elementary properties. 

\begin{definition}
	\label{5.54}
	A function $f:\pi(C) \to \mathbb{R}$ is called \textbf{$C$-regular in $u$ with respect to $X$}\index{regular} if there is a set $E$ of positive definable functions on $\pi(C)$ such that the following holds. The set $\log(E)$ consists only of $C$-consistent functions in $u$ with respect to $X$ and $f$ can be constructed from $E$. 
\end{definition}

For the rest of the whole Section 7 we often say ''$C$-consistent'' respectively ''$C$-regular''  instead of ''$C$-consistent in $u$ with respect to $X$'' respectively ''$C$-regular in $u$ with respect to $X$'' since $u$ and $X$ are fixed.\\

\begin{lemma}
	\label{5.55}
	The following holds.
	\begin{itemize}
		\item [(1)]
		Let $f:\pi(C) \to \mathbb{R}$ be $C$-regular and $C$-consistent. Then $\exp(f)$ is $C$-regular and $C$-consistent.
		\item [(2)]
		Let $f_1, \ldots ,f_p:\pi(C) \to \mathbb{R}$ be $C$-regular. Let $G:\mathbb{R}^p \to \mathbb{R}$ be log-analytic. Then $G(f_1, \ldots ,f_p)$ is $C$-regular.
	\end{itemize}
\end{lemma}

\begin{proof}
	Let $E$ be a set of positive definable functions on $\pi(C)$ such that $\log(E)$ is a set of $C$-consistent functions and $f$ can be constructed from $E$. 
	
	(1): With Remark~\ref{5.52}(3) we see that $\exp(f)$ is $C$-consistent. With Remark~\ref{3.4}(1) we see that $f$ can be constructed from $E \cup \{\exp(f)\}$ and is therefore $C$-regular.
	
	(2): Follows with Remark~\ref{3.4}(2).
\end{proof}

\begin{example}
	\label{5.56}
	The following holds.
	\begin{itemize}
		\item [(1)] A $C$-nice function is $C$-regular.
		\item [(2)] Let $l \in \{0, \ldots ,r\}$. The function $\pi(C) \to \mathbb{R}, t \mapsto \mu_{j,l}(t),$ is $C$-regular for $j \in \{0, \ldots ,l\}$.
		\item [(3)] The function $\pi(C) \to \mathbb{R}, t \mapsto \prod_{j=0}^{\kappa-1}e^{q_j \mu_{j,\kappa}(t)},$ is $C$-regular for $q_0, \ldots ,q_{\kappa-1} \in \mathbb{Q}$.
		\item [(4)] Suppose that $X \subset \mathbb{R}^2$ is open with $0 \in X$. Let $C:= \textnormal{}]0,1[^2$ and let $h: \textnormal{}]0,1[\textnormal{} \to \mathbb{R}, u \mapsto e^{-1/u}$. Then $h$ is not $C$-regular.
	\end{itemize}
\end{example}

\begin{proof}
	(1): Note that a $C$-nice function can be constructed by a set $E$ of positive definable functions such that every $g \in \log(E)$ is a $C$-heir. By Example~\ref{5.51} $g$ is $C$-consistent in $u$ with respect to $\mathbb{R}^n \times \mathbb{R}$ and hence, also with respect to $X$ (see Remark~\ref{5.50}).
		
	(2): By Lemma~\ref{5.55}(1), Remark~\ref{5.52}(1) and an easy induction on $j \in \{0, \ldots ,l\}$ we see that $\mu_{j,l}$ is $C$-regular.
	
	(3): With (2) and Lemma~\ref{5.55}(1) we obtain that $e^{\mu_{j,l}}$ is $C$-regular. Then, Lemma~\ref{5.55}(2) provides the $C$-regularity of $\prod_{j=0}^{\kappa-1}e^{q_j \mu_{j,\kappa}}$. 
	Note that $\exp(\mu_{j,l})$ is $C$-consistent. Therefore $\prod_{j=0}^{\kappa-1}e^{q_j \mu_{j,\kappa}}$ is $C$-consistent (by Remark~\ref{5.52}(1),(3)).
	
	(4): Note that $\pi(C) = \textnormal{}]0,1[$. Suppose the contrary. Let $E$ be a set of positive definable functions such that every function from $\log(E)$ is $C$-consistent and $h$ can be constructed from $E$. Note that every $g \in \log(E)$ is bounded at $0$. Therefore $\lim_{u \searrow 0} h(u)/\beta(u)=1$ for a log-analytic $\beta: \textnormal{} ]0,1[ \textnormal{} \to \mathbb{R}_{>0}$ at zero. But $\lim_{u \searrow 0} h(u)/\beta(u)=0$ for every log-analytic $\beta: \textnormal{} ]0,1[ \textnormal{} \to \mathbb{R}_{>0}$: By Fact~\ref{4.7} and Example~\ref{ex:C-simple} there are $a \in \mathbb{R} \setminus \{0\}$, a non-negative integer $r \in \mathbb{N}_0$, and $q_0, \ldots , q_r \in \mathbb{Q}$ such that 
	\begin{align*}
	& \lim_{y \searrow 0} \frac{\beta(y)}{y^{q_0}(-\log(y))^{q_1} \cdot \ldots \cdot \log_{r-1}(-\log(y))^{q_r}}=a. \qedhere
	\end{align*}
\end{proof}

We often work with restrictions of $C$-consistent and $C$-regular functions. Hence, we introduce the following to simplify notation.

\begin{convention}
	\label{con:regular-consistent-subsets}
	Let $B \subset \pi(C)$ be definable and $g:B \to \mathbb{R}$ be a function. If $g$ is the restriction of a $C$-regular ($C$-consistent) function then we call $g$ also $C$-regular ($C$-consistent) in $u$ with respect to $X$.
\end{convention}

Note that all the properties stated above (in particular Remark~\ref{5.50}, Remark~\ref{5.52}, Lemma~\ref{5.55}) 
transfer to the corresponding functions from Convention~\ref{con:regular-consistent-subsets}.

\subsection{Persistent and Soft Functions}

In this subsection we define $\kappa$-persistent and $\kappa$-soft functions and give also some elementary properties. Similar to $C$-consistent functions, these functions depend only on $t$, but their asymptotic behavior may crucially  depend on $X$. 


\begin{definition}
\label{5.38}
Let $g:\pi(C) \to \mathbb{R}_{\geq 0}$ be a function.
\begin{itemize}
	\item [(a)]
	We say that $g$ is \textbf{$(C,\mathcal{Y},\kappa)$-persistent in $u$ with respect to $X$}\index{persistent} if $g$ is definable and the following holds. Let $w \in \pi_l(C)$. Then for every definable curve $\gamma:\textnormal{}]0,1[\textnormal{} \to C_w$ compatible with $C_w$ with
	$\lim_{y \searrow 0} \gamma(y) \in X_w$ and
	$$\lim_{y \searrow 0} (\gamma_x(y) - \mu_\kappa(w,\gamma_u(y)))=0$$
	we have $\lim_{y \searrow 0}g(w,\gamma_u(y)) \in \mathbb{R}_{>0} \cup \{\infty\}$.
	\item [(b)]
	We say that $g$ is \textbf{$(C,\mathcal{Y},\kappa)$-soft in $u$ with respect to $X$}\index{soft} if $g$ is definable and the following holds. Let $w \in \pi_l(C)$. Then for every definable curve $\gamma:\textnormal{}]0,1[\textnormal{} \to C_w$ compatible with $C_w$ with
	$\lim_{y \searrow 0} \gamma(y) \in X_w$ and
	$$\lim_{y \searrow 0} (\gamma_x(y) - \mu_\kappa(w,\gamma_u(y)))=0$$
	we have $\lim_{y \searrow 0}g(w,\gamma_u(y))=0$.
\end{itemize}
\end{definition}


The next remark and two lemmas follow straightforwardly from Definition~\ref{5.38}.

\begin{remark}
	\label{5.39}
	Let $g:\pi(C) \to \mathbb{R}_{\geq 0}$ be a function. Let $Y \subset \mathbb{R}^n \times \mathbb{R}$ be definable with $X \subset Y$ such that $Y_w$ is open for every $w \in \mathbb{R}^l$. If $g$ is $(C,\mathcal{Y},\kappa)$-persistent ($(C,\mathcal{Y},\kappa)$-soft) in $u$ with respect to $Y$ then $g$ is $(C,\mathcal{Y},\kappa)$-persistent ($(C,\mathcal{Y},\kappa)$-soft) in $u$ with respect to $X$. 
\end{remark}

Since $C$, $\mathcal{Y}$ and $X$ are fixed we often say ''$\kappa$-persistent'' respectively ''$\kappa$-soft'' instead of ''$(C,\mathcal{Y},\kappa)$-persistent in $u$ with respect to $X$'' respectively ''$(C,\mathcal{Y},\kappa)$-soft in $u$ with respect to $X$'' for the rest of the whole Section 7. 

\begin{lemma}
\label{lem:kappa-persistence-distance}
Let $g:\pi(C) \to \mathbb{R}_{\geq 0}$ be a definable function. If for every $w \in \pi_l(C)$ and every definable curve $\gamma: \textnormal{}]0,1[\textnormal{} \to C_w$ compatible with $C_w$ with $\lim_{y \searrow 0} \gamma(y) \in X_w$ and 
$$\lim_{y \searrow 0} \textnormal{dist}(\mu_\kappa(w,\gamma_u(y)),C_{(w,\gamma_u(y))})=0$$
we have $\lim_{y \searrow 0}g(w,\gamma_u(y))>0$ ($\lim_{y \searrow 0}g(w,\gamma_u(y))=0$) then $g$ is $\kappa$-persistent ($\kappa$-soft).
\end{lemma}

\begin{proof}
Suppose the former. Let $\gamma: \textnormal{}]0,1[\textnormal{} \to C_w$ be a definable curve compatible with $C_w$ with $\lim_{y \searrow 0} \gamma(y) \in X_w$ and 
$$\lim_{y \searrow 0} (\gamma_x(y) - \mu_\kappa(w,\gamma_u(y)))=0.$$
Since $\gamma(]0,1[) \subset \pi(C)$ we also obtain with Lemma~\ref{5.6}(2)	
$$\lim_{y \searrow 0} \textnormal{dist}(\mu_\kappa(w,\gamma_u(y)),C_{(w,\gamma_u(y))})=0$$
and hence $\lim_{y \searrow 0}g(w,\gamma_u(y))>0$. This shows that $g$ is $\kappa$-persistent. The $\kappa$-softness is obtained similarly if the latter is assumed.
\end{proof}


\begin{lemma}
	\label{5.40}
	Let $w \in \pi_l(X)$ and $(u_0,x_0) \in \overline{C_w} \cap X_w$. Let $M:\pi(C) \to \mathbb{R}_{\geq 0}$ be $\kappa$-persistent and $N:\pi(C) \to \mathbb{R}_{\geq 0}$ be $\kappa$-soft. Let $\gamma:\text{}]0,1[\text{} \to \pi(C)_w \times \mathbb{C}$ be a definable curve compatible with $C_w$ with $\lim_{y \searrow 0} (\gamma_u(y),\gamma_z(y))=(u_0,x_0)$ and $\lim_{y \searrow 0}\vert{\mu_\kappa(w,\gamma_u(y)) - \gamma_z(y)}\vert = 0$. Then $\lim_{y \searrow 0} M(w,\gamma_u(y))>0$ and $\lim_{y \searrow 0} N(w,\gamma_u(y))=0$.  
\end{lemma}

\begin{proof}
	We have that $\lim_{y \searrow 0}\textnormal{dist}(\gamma_z(y),C_{(w,\gamma_u(y))})=0$ and with Lemma~\ref{5.27}(2) there is a definable curve $\hat{\gamma}$ in $C_w$ compatible with $C_w$ with $\hat{\gamma}_u = \gamma_u$ and
	$$\lim_{y \searrow 0} (\hat{\gamma}_u(y),\hat{\gamma}_x(y)) = (u_0,x_0).$$ 
	With Definition~\ref{5.38} we obtain the desired results, because also 
	\begin{align*}
	\lim_{y \searrow 0} &\vert{\mu_\kappa(w,\hat{\gamma}_u(y))-\hat{\gamma}_x(y)}\vert = 0. \qedhere
	\end{align*}
\end{proof}

The following remark and example show that both notions from Definition~\ref{5.38} are not mutually exclusive. 

\begin{remark}
	\label{5.46}
	The following holds.
	\begin{itemize}
	\item[(1)] If $C$ is thin with respect to $x$ then there is no definable curve $\gamma:\text{}]0,1[\text{} \to C$ compatible with $C$ and therefore every definable function $\pi(C) \to \mathbb{R}_{\geq 0}$ is $\kappa$-persistent and $\kappa$-soft.
	\item[(2)] If $C$ is fat with respect to $x$ and far with respect to $\mu_\kappa$ then every definable curve $\gamma:\text{}]0,1[\text{} \to C$ compatible with $C$
	fulfills 
	$$\lim_{y \searrow 0} \vert{\gamma_x(y) - \mu_\kappa(\gamma_u(y))}\vert > 0$$
	and therefore every definable function $\pi(C) \to \mathbb{R}_{\geq 0}$ is $\kappa$-persistent and $\kappa$-soft.
	\end{itemize}
\end{remark}

\begin{example}
	\label{5.41}
	Let $X=\mathbb{R}^2$ and let $C$ and $\mathcal{Y}=(y_0,y_1,y_2)$ be as in Example~\ref{5.10} (i.e. $r=2$ and $k^{\textnormal{ch}}=1$). Consider the definable functions $g_1:\textnormal{}]0,1[\textnormal{} \to \mathbb{R}_{>0}, u \mapsto u,$ and $g_2:\textnormal{}]0,1[\textnormal{} \to \mathbb{R}_{>0}, u \mapsto 1-u$. The following holds.
	
	\begin{itemize}
		\item [(1)] Let $X=\mathbb{R}^2$. Then $g_1$ is $2$-persistent and $g_2$ is $2$-soft. 
		\item [(2)] Let $X=C$. Then $g_1$ and $g_2$ are $\kappa$-persistent and $\kappa$-soft for $\kappa \in \{1,2\}$.
	\end{itemize}
\end{example}

\begin{proof}
	We have 
	$$C=\{(u,x) \in \mathbb{R} \times \mathbb{R} \mid u \in \textnormal{} ]0,1[ \textnormal{}, \textnormal{}\tfrac{1}{1+u}+e^{-u-1/u} < x < \tfrac{1}{1+u} + e^{-u/2-1/u}\}.$$
	Note that $C$ is open in $\mathbb{R}^2$. For $u \in \textnormal{}]0,1[$ we have $L_C(u)=e^{-1/u}(e^{-u/2}-e^{-u})$. Note that $L_C(u)>0$ for $u \in \textnormal{} ]0,1[$, that $\lim_{u \searrow 0} L_C(u)=0$ and that $\lim_{u \nearrow 1} L_C(u)>0$. 
	For $u \in \textnormal{} ]0,1[$ we have $\mu_1(u)=\frac{1}{1+u}+e^{-1/u}$ and $\mu_2(u)=\frac{1}{1+u}+e^{-1-1/u}$. 
	
	(1): We see that 
	$$\lim_{u \nearrow 1}\textnormal{dist}(\mu_2(u),C_u)=\lim_{u \nearrow 1} \Bigl(e^{-u-1/u} - e^{-1-1/u} \Bigl)= 0,$$
	$\lim_{u \nearrow 1} g_1(u)>0$ and $\lim_{u \nearrow 1} g_2(u)=0$. Hence, due to $\textnormal{dist}(\mu_2(u),C_u)>0$ for every $u \in  \textnormal{} ]0,1[$, we have with Lemma~\ref{lem:kappa-persistence-distance}
	that $g_1$ and $g_2$ fulfill the desired properties. 
	
	(2): Since $C$ is open we see with Lemma~\ref{5.6}(2) that there is no definable curve $\gamma: \textnormal{}]0,1[ \textnormal{} \to C$ with $\lim_{y \searrow 0}\gamma(y) \in C$ and 
	$$\lim_{y \searrow 0}(\gamma_x(y)-\mu_\kappa(\gamma_u(y)))=0.$$
	So the assertion follows.
\end{proof}

The next remark shows that there are functions which are neither $(C,\mathcal{Y},\kappa)$-persistent in $u$ with respect to $X$ nor $(C,\mathcal{Y},\kappa)$-soft in $u$ with respect to $X$.

\begin{remark}
	\label{rem:closure-set-persistence}
	Suppose that $C$ is fat with respect to $x$ and that $\mu_{\kappa}(t) \in \overline{C_t}$ for $t \in \pi(C)$. Then $g_1(t)=0$ and $g_2(t)>0$ for every $\kappa$-persistent $g_1:\pi(C) \to \mathbb{R}_{\geq 0}$ and every $\kappa$-soft $g_2:\pi(C) \to \mathbb{R}_{\geq 0}$. Consequently, the function $$g:\text{}]0,1[\text{} \to \mathbb{R}, u \mapsto (u-1/2)^2,$$ is neither $\kappa$-persistent nor $\kappa$-soft if $C=\text{}]0,1[^2$, $X=\mathbb{R}^2$ and $\mathcal{Y}:=(y_0)$ with $y_0:=x$. 
\end{remark}

Now we give some more basic properties and examples of $\kappa$-persistent and $\kappa$-soft functions. 

\begin{remark}
	\label{5.42}
	The following holds.
	\begin{itemize}
		\item [(1)] The set of all $\kappa$-persistent ($\kappa$-soft) functions on $\pi(C)$ is closed under pointwise addition and multiplication.
		\item [(2)] Let $g:\pi(C) \to \mathbb{R}_{\geq 0}$ be $\kappa$-persistent ($\kappa$-soft) and let $q \in \mathbb{R}_{>0}$. Then $\pi(C) \to \mathbb{R}_{\geq 0}, t \mapsto \sqrt[q]{g(t)}$, is $\kappa$-persistent ($\kappa$-soft).
		\item [(3)] Let $f_1,f_2:\pi(C) \to \mathbb{R}_{\geq 0}$ be $\kappa$-persistent ($\kappa$-soft). Then 
		$$\pi(C) \mapsto \mathbb{R}_{\geq 0}, t \mapsto \min\{f_1(t),f_2(t)\},$$
		and
		$$\pi(C) \mapsto \mathbb{R}_{\geq 0}, t \mapsto \max\{f_1(t),f_2(t)\},$$
		are $\kappa$-persistent ($\kappa$-soft).
		\item[(4)] Let $g:\pi(C) \to \mathbb{R}_{\geq 0}$ be $\kappa$-soft and $h:\pi(C) \to \mathbb{R}_{\geq 0}$ be $C$-consistent. Then $g \cdot h$ is $\kappa$-soft.
	\end{itemize}
\end{remark}

\begin{example}
	\label{5.43}
	Let $X=\mathbb{R}^n \times \mathbb{R}$. Let $l \in \{1, \ldots ,r\}$ and $j \in \{0, \ldots ,l-1\}$. Then
	$$\pi(C) \to \mathbb{R}_{>0}, t \mapsto e^{\mu_{j,l}(t)},$$
	and
	$$\pi(C) \to \mathbb{R}_{>0}, t \mapsto \vert{\mu_l(t)-\mu_j(t)}\vert,$$
	are $\kappa$-persistent. Note that $X$ can be chosen arbitrarily such that $C \subset X$ and $X_w$ is open for every $w \in \pi_l(C)$. 
\end{example}

\begin{proof}
	Let $w \in \pi_l(C)$. Then for every definable curve $\gamma:\textnormal{}]0,1[\textnormal{} \to C_w$ compatible with $C_w$ we have $\lim_{y \searrow 0}e^{\mu_{j,l}(w,\gamma_u(y))} > 0$ by Corollary~\ref{5.29}(1) and 
	$$\lim_{y \searrow 0} \vert{\mu_l(w,\gamma_u(y))-\mu_j(w,\gamma_u(y))}\vert > 0$$
	if $j \neq l$ by Corollary~\ref{5.29}(2).
\end{proof}

The next lemma shows a simple connection between $C$-consistent and $\kappa$-persistent functions.  

\begin{lemma}
	\label{lem:consistent-vs-persistence}
The following properties hold.
\begin{itemize}
	\item [(1)] Let $g:\pi(C) \to \mathbb{R}_{> 0}$ be $C$-consistent in $u$ with respect to $X$. Then 
	$\pi(C) \to \mathbb {R}_{>0}, t \mapsto 1/g(t),$ is $(C,\mathcal{Y},\kappa)$-persistent in $u$ with respect to $X$. 
	\item [(2)] Let $h:\pi(C) \to \mathbb{R}_{\geq 0}$ be $C$-consistent in $u$ with respect to $X$ and $\kappa$-persistent in $u$ with respect to $X$. Then 
	$$h_{\text{pers}}:\pi(C) \to \mathbb {R}_{\geq 0}, t \mapsto \left\{\begin{array}{lll} 1/h(t),& h(t) \neq 0, \\
		0,&\textnormal{else,} \end{array}\right.$$
	is $(C,\mathcal{Y},\kappa)$-persistent in $u$ with respect to $X$. 
\end{itemize}
\end{lemma}

\begin{proof}
Let $w \in \pi_l(C)$ and let $\gamma: \text{}]0,1[ \text{} \to C_w$ be a definable curve compatible with $C_w$ with $\lim_{y \searrow 0} \gamma_u(y) \in X_w$ and $\lim_{y \searrow 0} (\gamma_x(y)-\mu_\kappa(w,\gamma_u(y))) = 0$. \\
(1): With Definition~\ref{5.48} we see that $\lim_{y \searrow 0} g(w,\gamma_u(y)) \in \mathbb{R}$ since $\lim_{y \searrow 0} \gamma(y) \in X_w$. This implies the $\kappa$-persistence of $1/g$. \\
(2): Since $h$ is $\kappa$-persistent we see that $\lim_{y \searrow 0} h(w, \gamma_u(y)) > 0$ and hence we find $\varepsilon>0$ such that $h(w,\gamma_u(]0,\varepsilon[)) \subset \mathbb{R}_{>0}$. Thus, 
$$\lim_{y \searrow 0} h_{\text{pers}}(w, \gamma_u(y)) = \lim_{y \searrow 0} 1/(h(w, \gamma_u(y))) > 0$$
since $h$ is $C$-consistent in $u$ with respect to $X$.
\end{proof}

In the following we consider an $(m+1,X)$-restricted log-exp-analytically $(e,r)$-prepared function $f:C \to \mathbb{R}, (w,u,x) \mapsto f(w,u,x),$ and check its coefficient and base functions for $\kappa$-persistence and $\kappa$-softness.

\begin{lemma}
\label{5.44}
Let $e \in \mathbb{N}_0 \cup \{-1\}$ and $r \in \mathbb{N}_0$. Let $f:C \to \mathbb{R}, (w,u,x) \mapsto f(w,u,x),$ be $(m+1,X)$-restricted log-exp-analytically $(e,r)$-prepared in $x$ with preparing tuple 
$$(r,\mathcal{Y},a,\exp(c),q,s,v,b,\exp(d),P)$$
with $b:=(b_1, \ldots ,b_s)$, $\exp(d):=(\exp(d_1), \ldots ,\exp(d_s))$ and $P:=(p_1, \ldots ,p_s)^t$ for $p_i \in \mathbb{Q}^{r+1}$. Then the following holds. 
\begin{itemize}
	\item [(1)] Let $j \in \{1, \ldots ,s\}$. The function $\pi(C) \to \mathbb{R}_{>0}, t \mapsto \vert{1/b_j(t)}\vert,$ is $\kappa$-persistent. Additionally, if $p_j$ is $\kappa$-negative then $\vert{b_j}\vert$ is $\kappa$-soft. 
	\item [(2)] Suppose that $a$ does not have any zero and that $f$ is locally bounded in $(u,x)$ with reference set $X$. Then the function $\pi(C) \to \mathbb {R}_{>0}, t \mapsto \vert{1/a(t)}\vert,$ is $\kappa$-persistent. Additionally, if $q$ is $\kappa$-negative then $\vert{a}\vert$ is $\kappa$-soft.
\end{itemize} 
\end{lemma}

\begin{proof}
Note that $b$ is $C$-consistent and if $f$ is locally bounded in $(u,x)$ with reference set $X$ then also $a$ is $C$-consistent. 
The $\kappa$-persistence of $\vert{1/a}\vert$ (if $f$ is locally bounded in $(u,x)$ with reference set $X$) and $\vert{1/b}\vert$ follows with Lemma~\ref{5.53} and Lemma~\ref{lem:consistent-vs-persistence}(1). Hence, it remains to deal with the $\kappa$-softness. For $(t,x) \in C$ let
$$V(t,x):=v(b_1(t)\vert{\mathcal{Y}}\vert^{\otimes p_1}(t,x) e^{d_1(t,x)}, \ldots ,b_s(t)\vert{\mathcal{Y}}\vert^{\otimes p_s}(t,x) e^{d_s(t,x)}).$$
For $(t,x) \in C$ we have 
$$f(t,x)=a(t)\vert{\mathcal{Y}}\vert^{\otimes q}(t,x) e^{c(t,x)} V(t,x).$$
Let $w \in \pi_l(C)$ and $\gamma: \textnormal{}]0,1[ \textnormal{}\to C_w$ be a definable curve compatible with $C_w$ such that $(u_0,x_0):=\lim_{y \searrow 0} \gamma(y) \in X_w$ and 
$$\lim_{y \searrow 0}\vert{\gamma_x(y)-\mu_\kappa(w,\gamma_u(y))}\vert=0.$$
We find $\delta>0$ and $B>0$ such that $U:=Q^m(u_0,\delta) \times \textnormal{}]x_0-\delta,x_0+\delta[ \textnormal{} \subset X_w$, $\vert{c(w,u,x)}\vert<B$ and $\vert{d_j(w,u,x)}\vert<B$ for $j \in \{1, \ldots ,s\}$ and every $(u,x) \in U \cap C_w$. 

(1): Let $j \in \{1, \ldots ,s\}$. Suppose that $p_j$ is $\kappa$-negative. Then by Lemma~\ref{5.34}(2) and Lemma~\ref{5.15}(1) we have
$$\lim_{y \searrow 0} \vert{\mathcal{Y}}\vert^{\otimes p_j}(w,\gamma(y)) = \infty$$
and since $b_j(t)\vert{\mathcal{Y}}\vert^{\otimes p_j}(t,x)e^{d_j(t,x)} \in [-1,1]$ for every $(t,x) \in C$
$$\lim_{y \searrow 0}b_j(w,\gamma_u(y))=0.$$

(2): Note that there is $\rho>1$ such that $V(t,x) \in \textnormal{}]1/\rho,\rho[$ for every $(t,x) \in C$. If $q$ is $\kappa$-negative we get with Lemma~\ref{5.34}(2) and Lemma~\ref{5.15}(1)
$$\lim_{y \searrow 0} \vert{\mathcal{Y}}\vert^{\otimes q}(w,\gamma(y)) = \infty.$$
Hence, $\lim_{y \searrow 0}a(w,\gamma_u(y))=0$ since $\lim_{y \searrow 0} f(w,\gamma_u(y)) \in \mathbb{R}$.
\end{proof}

Finally, we introduce notation for a class of definable sets needed in the next sections to construct for $\Phi_f:\Lambda_f \to \mathbb{C}$ a set
$$D = \{(t,z) \in H \mid N(t)<\vert{z-\mu_\kappa(t)}\vert<M(t)\}$$
such that $D \subset \Lambda_f$ and $\Phi_f|_D$ is complex log-analytically prepared. Here $\Phi_f$ is the unary parametric global complexification of a restricted log-exp-analytically prepared $f:C \to \mathbb{R}$ from Definition~\ref{5.78}. We also provide some simple properties of these sets. Also $D$ belongs to this class.

\begin{definition}
	\label{5.59}
	For definable functions $N,M:\pi(C) \to \mathbb{R}_{\geq 0}$ and a definable function $G:H \to \mathbb{C}$ we set
	$$\mathcal{B}(G,M):=\{(t,z) \in H \mid \vert{G(t,z)}\vert < M(t)\}$$\index{$\mathcal{B}(G,M)$}
	and
	$$\mathcal{A}(G,N,M):=\{(t,z) \in H \mid N(t) < \vert{G(t,z)}\vert < M(t)\}.$$\index{$\mathcal{A}(G,N,M)$}
\end{definition}

\begin{remark}
	\label{5.60}
	Let $N,M:\pi(C) \to \mathbb{R}_{\geq 0}$ be definable. For a definable $G:H \to \mathbb{C}$ the sets $\mathcal{B}(G,M) \subset \mathbb{R}^n \times \mathbb{C}$ and $\mathcal{A}(G,N,M) \subset \mathbb{R}^n \times \mathbb{C}$ are definable. We have 
	$$\mathcal{B}(\mathcal{P}_{0,\kappa},M)=\{(t,z) \in H \mid 0 < \vert{z-\mu_\kappa(t)}\vert < M(t)\}$$
	\index{$\mathcal{B}(\mathcal{P}_{0,\kappa},M)$}\index{$\mathcal{A}(\mathcal{P}_{0,\kappa},N,M)$}and
	$$\mathcal{A}(\mathcal{P}_{0,\kappa},N,M)=\{(t,z) \in H \mid N(t) < \vert{z-\mu_\kappa(t)}\vert < M(t)\}.$$
\end{remark}

To close this subsection we demonstrate with the next lemma and corollary how the sets from Remark~\ref{5.60} can be used to construct even more $\kappa$-persistent and $\kappa$-soft functions. 

\begin{lemma}
	\label{lem:persistence-soft-restriction}
	Let $M:\pi(C) \to \mathbb{R}_{\geq 0}$ be $\kappa$-persistent and $N:\pi(C) \to \mathbb{R}_{\geq 0}$ be $\kappa$-soft. Let $B:=\mathcal{A}(\mathcal{P}_{0,\kappa},N,M)$. Then the following properties hold. 
	\begin{itemize}
		\item[(1)] Let $T:\pi(C) \to \mathbb{R}_{\geq 0}$ be $\kappa$-persistent. Then the function 
		$$\hat{T}:\pi(C) \to \mathbb{R}_{\geq 0}, t \mapsto \left\{\begin{array}{lll} T(t),& t \in \pi(B), \\
			0,&\textnormal{else,} \end{array}\right.$$
		is $\kappa$-persistent. 
		\item[(2)] Let $S:\pi(C) \to \mathbb{R}_{\geq 0}$ be $\kappa$-soft. Then the function 
		$$\hat{S}:\pi(C) \to \mathbb{R}_{\geq 0}, t \mapsto \left\{\begin{array}{lll} S(t),& t \in \pi(B), \\
			1,&\textnormal{else,} \end{array}\right.$$
		is $\kappa$-soft. 
	\end{itemize}
\end{lemma}

\begin{proof}
Let $w \in \pi_l(C)$ and let $\gamma: \text{}]0,1[ \text{} \to C_w$ be a definable curve compatible with $C_w$ with $\lim_{y \searrow 0} \gamma_u(y) \in X_w$ and $\lim_{y \searrow 0} (\gamma_x(y)-\mu_\kappa(w,\gamma_u(y))) = 0$. Since $M$ is $\kappa$-persistent and $N$ is $\kappa$-soft we can assume by passing to a suitable subcurve of $\gamma$ if necessary that $(w,\gamma_u(]0,1[)) \subset \pi(B)$ (since $\lim_{y \searrow 0} M(w,\gamma_u(y)) > 0$ and $\lim_{y \searrow 0} N(w,\gamma_u(y))=0$).\\
(1): Note that
\[
\lim_{y \searrow 0} \hat{T}(w,\gamma_u(y)) = \lim_{y \searrow 0} T(w,\gamma_u(y)) > 0
\] 
since $T$ is $\kappa$-persistent.\\
(2): Note that
\[
\lim_{y \searrow 0} \hat{S}(w,\gamma_u(y)) = \lim_{y \searrow 0} S(w,\gamma_u(y)) = 0
\] 
since $S$ is $\kappa$-soft.
\end{proof}

\begin{corollary}
	\label{cor:persistence-soft-restriction}
	Let $M:\pi(C) \to \mathbb{R}_{\geq 0}$ be $\kappa$-persistent and $N:\pi(C) \to \mathbb{R}_{\geq 0}$ be $\kappa$-soft. Let $B:=\mathcal{A}(\mathcal{P}_{0,\kappa},N,M)$. Let $s \in \mathbb{N}$ and let $h_1, \ldots , h_s:\pi(B) \to \mathbb{R}_{>0}$ be $C$-consistent. Then 
	\[
	T:\pi(C) \to \mathbb{R}_{\geq 0}, t \mapsto \left\{\begin{array}{lll} \min\{1/h_1(t), \ldots , 1/h_s(t)\},& t \in \pi(B), \\
		0,&\textnormal{else,} \end{array}\right.
	\]
	is $\kappa$-persistent.  
\end{corollary}

\begin{proof}
	Note that $h_1, \ldots , h_s$ extend to $C$-consistent functions $\hat{h}_1, \ldots , \hat{h}_s:\pi(C) \to \mathbb{R}_{>0}$ and hence, by Lemma~\ref{lem:consistent-vs-persistence}, the functions $1/\hat{h} _1, \ldots , 1/\hat{h}_s$ are $\kappa$-persistent. By Remark~\ref{5.42}(3) the function $\min\{1/\hat{h}_1, \ldots , 1/\hat{h}_s\}$ is also $\kappa$-persistent. Now use Lemma~\ref{lem:persistence-soft-restriction}(1).
\end{proof}

\subsection{Logarithmic Scales}
This subsection contains technical results regarding logarithmic scales in the context of $\kappa$-persistent and $\kappa$-soft functions which are necessary for proving the desired preparation theorem in the next subsection. The procedure is as follows:

Let $U$ be a given $\kappa$-persistent function and $W$ be a $\kappa$-soft one. In the first lemma we construct a $\kappa$-persistent function $M_1:\pi(C) \to \mathbb{R}_{\geq 0}$ and a $\kappa$-soft function $N_1:\pi(C) \to \mathbb{R}_{\geq 0}$ such that 
$$\mathcal{A}(\mathcal{P}_{0,\kappa},N_1,M_1) = \{(t,z) \in H \mid N_1(t)<\vert{z-\mu_{\kappa}(t)}\vert<M_1(t)\} \subset \mathcal{A}(z_\kappa,W,U).$$ 
In the second and third lemmas we extend this result on $(\sigma\mathcal{Z})^{\otimes q}$ for $q:=(q_0, \ldots , q_r) \in \mathbb{Q}^{r+1}$ with $(q_\kappa, \ldots, q_r)\neq 0$, showing that 
$$\mathcal{A}(\mathcal{P}_{0,\kappa},N_2,M_2) \subset \mathcal{A}((\sigma\mathcal{Z})^{\otimes q},W,U)$$
for a suitable $\kappa$-persistent $M_2$ and $\kappa$-soft $N_2$. These results demonstrate that for $R>1$ one can choose $N_2$ and $M_2$ such that for such a $p_i \in \mathbb{Q}^{r+1}$,
$$\vert{b_i(t)(\sigma\mathcal{Z})^{\otimes p_i}(t,z)}\vert \leq R$$
for all $(t,z) \in \mathcal{A}(\mathcal{P}_{0,\kappa},N_2,M_2)$ where $b_i$ is a base function occuring in a restricted log-exp-analytic preparation. This is possible, because $R/\vert{b_i(t)}\vert$ is $\kappa$-persistent by Lemma~\ref{5.44}(1). Another possible choice of $M_2$ and $N_2$ even implies
$$\vert{b_i(t)(\sigma\mathcal{Z})^{\otimes p_i}(t,z)\exp(D_i(t,z))}\vert \leq R$$
for all $(t,z) \in \mathcal{A}(\mathcal{P}_{0,\kappa},N_2,M_2)$ where $b_i$ and $p_i$ are as above and $D_i$ is the unary parametric global complexification of a restricted log-exp-analytically prepared function $d_i \in \log(E)$ from Definition~\ref{5.78}. The reason is that on $\mathcal{A}(\mathcal{P}_{0,\kappa},N_2,M_2)$ one can write $D_i=D_{i,1}+D_{i,2}$ with $D_{i,1}:\pi(C) \to \mathbb{R}$ being $C$-consistent, and $D_{i,2}$ being a complex log-analytically prepared function with $\vert{\text{Im}(D_{i,2})}\vert<\pi$. According to Lemma~\ref{5.44}(1), $R/\vert{b_i(t)\exp(D_{i,1}(t))}\vert$ is also $\kappa$-persistent. (See Theorem~\ref{5.84} and Claim~\ref{claim-13} below for the details.) 

We will not prove similar results for the case that $(q_\kappa, \ldots , q_r)=0$, but instead, we construct a $\kappa$-persistent function $T:\pi(C) \to \mathbb{R}_{\geq 0}$ to re-prepare products of the form  $\prod_{j=0}^{\kappa-1}(\sigma_jz_j)^{q_j}$ on 
$$\mathcal{B}(\mathcal{P}_{0,\kappa},T) = \{(t,z) \in H \mid \vert{z-\mu_\kappa(t)}\vert<T(t)\}$$
suitably for our purposes in a complex log-analytically way (compare with Lemma~\ref{lem:complex-preparation-simplified} below).\\



All the results in this subsection are formulated for $X=\mathbb{R}^n \times \mathbb{R}$ and therefore for every possible choice of $l,m \in \mathbb{N}_0$ (the number of components of $w$ and $u$) such that $n=l+m$. This includes every possible choice of $X$ with $C \subset X$ such that $X_w$ is open.  Therefore we work with the variable $t$ instead of $w$ and $u$.

The functions $\mathcal{P}_{j,\kappa}$ from Definition~\ref{5.16} play a major role in our proofs which we consider as functions on $H$ (as also remarked after Corollary~\ref{5.19}).

Note that the center $\Theta$ of the logarithmic scale $\mathcal{Y}$ is $C$-nice and therefore $C$-regular (by Example~\ref{5.56}(1)).

\begin{lemma}
	\label{5.67}
	Let $U:\pi(C) \to \mathbb{R}_{\geq 0}$ be $\kappa$-persistent. Then the following holds.
	
	\begin{itemize}
		\item [(1)] Let $0<\delta \leq 1/2$. There is a $C$-regular $\kappa$-persistent $\zeta:\pi(C) \to \mathbb{R}_{\geq 0}$ such that $K:=\mathcal{B}(\mathcal{P}_{0,\kappa},\zeta) \subset \mathcal{B}(\mathcal{P}_{j,\kappa},U)$ and
		$$\Bigl|\frac{\mathcal{P}_{j-1,\kappa}(t,z)}{e^{\mu_{\kappa-j,\kappa}(t)}}\Bigl| < \delta$$
		for every $(t,z) \in K$ and every $j \in \{1, \ldots ,\kappa\}$. If $U$ is $C$-regular then $\zeta$ can be chosen in a $C$-regular way.
		
		\item [(2)] Let $W:\pi(C) \to \mathbb{R}_{\geq 0}$ be $\kappa$-soft. Then there is a $C$-regular $\kappa$-persistent $M :\pi(C) \to \mathbb{R}_{\geq 0}$ and a $C$-regular $\kappa$-soft $N:\pi(C) \to \mathbb{R}_{\geq 0}$ such that $\mathcal{A}(\mathcal{P}_{0,\kappa},N,M) \subset \mathcal{A}(z_\kappa,W,U)$. If $W$ and $U$ are $C$-regular then $M$ and $N$ can be chosen in a $C$-regular way.
	\end{itemize}
\end{lemma}

\begin{proof}
	To prepare the proof we need the following Claim~\ref{claim-2} to handle the iterated logarithms occuring in $\mathcal{P}_{j,\kappa}$.
	
	\begin{claim}
		\label{claim-2}
		Let $0<\delta \leq 1/2$. There are globally subanalytic functions $\lambda_1,\lambda_2: \mathbb{R}_{\geq 0} \to \textnormal{} [0,\delta]$ such that the following holds.
		\begin{itemize}
			\item [(i)] $\lambda_j$ is monotone increasing for $j \in \{1,2\}$.
			\item[(ii)] For $r \in \mathbb{R}_{\geq 0}$ and $z \in \mathbb{C}$ we have $\vert{\log(1+z)}\vert<r$ if $\vert{z}\vert < \lambda_1(r)$.
			\item[(iii)] For $r \in \mathbb{R}_{\geq 0}$ and $z \in \mathbb{C}$ we have $r<\vert{\log(1+z)}\vert$ if $\lambda_2(r) < \vert{z}\vert < \delta$ for $z \in \mathbb{C}$.
			\item[(iv)] We have $\lim_{r \searrow 0} \lambda_j(r)=0$ and $\lim_{r \nearrow r_0} \lambda_j(r)>0$ for $r_0 \in \text{}]0,\delta]$ and $j \in \{1,2\}$.
		\end{itemize}
	\end{claim}
	
	\begin{proof}
		Let 
		$$\log^*:\mathbb{C} \to \mathbb{C}, z \mapsto \left\{\begin{array}{lll} \log(z),&& \vert{z}\vert \in [1/2,3/2] \setminus \mathbb{R}_{\leq 0}, \\
		&\textnormal{if}&\\
		0,&&\textnormal{else.} \end{array}\right.$$
		Define
		$$\lambda_1:\mathbb{R}_{\geq 0} \to \textnormal{}[0,\delta], r \mapsto$$
		$$\sup\{c \in \textnormal{}[0,\delta[ \textnormal{} \mid \vert{\log^*(1+z)}\vert < r \textnormal{ for every $z \in \mathbb{C}$ with } \vert{z}\vert < c\}.$$
		
		Note that $\lambda_1$ is well-defined, globally subanalytic (since $\log^*$ is globally subanalytic) and monotone increasing. The properties~(ii) and~(iv) are fulfilled by construction of $\lambda_1$ and the continuity of $\log$. 
		For $r \in \mathbb{R}_{\geq 0}$ set
		$$s(r):= \inf\{c \in \textnormal{}[0,\delta[ \textnormal{} \mid r < \vert{\log^*(1+z)}\vert \textnormal{ for every $z \in \mathbb{C}$ with } c<\vert{z}\vert < \delta\}.$$
		
		Note that $s(r_1) \leq s(r_2)$ for $0 \leq r_1 \leq r_2$. Set
		$$\lambda_2:\mathbb{R}_{\geq 0} \to [0,\delta], r \to \left\{\begin{array}{ll} s(r),& s(r) \neq \infty, \\
			\delta,&\textnormal{else.} \end{array}\right.$$
		
		Note that $\lambda_2$ is well-defined, globally subanalytic and monotone increasing. The properties~(iii) and~(iv) are fulfilled by construction of $\lambda_1$ and the continuity of $\log$.
	\end{proof}
	
	(1): Let $\lambda_1$ be as in Claim~\ref{claim-2}. 
	Let $\zeta_\kappa:\pi(C) \to \mathbb{R}_{\geq 0}, t \mapsto U(t)$. Define by descending induction on $j \in \{1, \ldots ,\kappa\}$
	$$\zeta_{j-1}:\pi(C) \to \mathbb{R}_{\geq 0}, t \mapsto \lambda_1(\zeta_j(t)) e^{\mu_{\kappa-j,\kappa}(t)}.$$
	By an easy descending induction on $j \in \{0, \ldots ,\kappa\}$ one sees with Claim~\ref{claim-2} that $\zeta_j$ is a well-defined $C$-regular $\kappa$-persistent function (compare with Remark~\ref{5.42}(1) and Example~\ref{5.43} since $\Theta$ is $C$-nice and therefore $C$-regular). If $\vert{\mathcal{P}_{j-1,\kappa}(t,z)}\vert<\zeta_{j-1}(t)$ then $$\Bigl|\frac{\mathcal{P}_{j-1,\kappa}(t,z)}{e^{\mu_{\kappa-j,\kappa}(t)}}\Bigl|<\lambda_1(\zeta_j(t)) \leq \delta$$
	and therefore 
	$$\vert{\mathcal{P}_{j,\kappa}(t,z)}\vert = \left| \log \left(1+\frac{\mathcal{P}_{j-1,\kappa}(t,z)}{e^{\mu_{\kappa-j,\kappa}(t)}} \right) \right| < \zeta_j(t)$$
	for $j \in \{1, \ldots , \kappa\}$ and $(t,z) \in H$ by definition of $\lambda_1$. So we take 
	$\zeta:=\zeta_0$.
	
	(2): Let $0 < \delta \leq 1/2$. Let $\lambda_2$ be as in Claim~\ref{claim-2}. By Corollary~\ref{5.19} we have that $\mathcal{P}_{\kappa,\kappa} = z_\kappa$. By (1) there is a $C$-regular $\kappa$-persistent $M:\pi(C) \to \mathbb{R}_{\geq 0}$ such that
	$K:=\mathcal{B}(\mathcal{P}_{0,\kappa},M) \subset \mathcal{B}(\mathcal{P}_{\kappa,\kappa},U) = \mathcal{B}(z_\kappa,U)$ and we have
	$$\Bigl|\frac{\mathcal{P}_{j-1,\kappa}(t,z)}{e^{\mu_{\kappa-j,\kappa}(t)}}\Bigl| < \delta$$
	for every $(t,z) \in K$ and every $j \in \{1, \ldots ,\kappa\}$. 
	Let $N_\kappa:=W$. Define by descending induction on $j \in \{1, \ldots ,\kappa\}$
	$$N_{j-1}:\pi(C) \to \mathbb{R}_{\geq 0}, t \mapsto \lambda_2(N_j(t)) e^{\mu_{\kappa-j,\kappa}(t)}.$$
	By an easy descending induction on $j \in \{0, \ldots ,\kappa\}$ one sees with Remark~\ref{5.42}(4) that $N_j$ is a well-defined $C$-regular $\kappa$-soft function. (From Claim~\ref{claim-2} one obtains that $\lambda(N_j)$ is $\kappa$-soft and with Example~\ref{5.51}(2) and Remark~\ref{5.52}(3) one sees that $e^{\mu_{\kappa-j,\kappa}}$ is $C$-consistent.) 
	
	This gives the following: If $N_{j-1}(t)<\vert{\mathcal{P}_{j-1,\kappa}(t,z)}\vert$ then 
	$$\lambda_2(N_j(t)) < \Bigl|\frac{\mathcal{P}_{j-1,\kappa}(t,z)}{e^{\mu_{\kappa-j,\kappa}(t)}}\Bigl| < \delta$$
	and therefore 
	$$N_j(t) < \left| \log \left(1+\frac{\mathcal{P}_{j-1,\kappa}(t,z)}{e^{\mu_{\kappa-j,\kappa}(t)}} \right) \right| = \vert{\mathcal{P}_{j,\kappa}(t,z) \vert}$$
	for $j \in \{1, \ldots ,\kappa\}$ and $(t,z) \in \mathcal{B}(\mathcal{P}_{0,\kappa},M)$ by definition of $\lambda_2$. So set $N:=N_0$. We obtain $\mathcal{A}(\mathcal{P}_{0,\kappa},N,M) \subset \mathcal{A}(z_\kappa,W,U)$. This finishes the proof.
\end{proof}

For technical reasons we need two more small notions. For a definable function $T:\pi(C) \to \mathbb{R}_{\geq 0}$ we set 
$$T_{\textnormal{pers}}:\pi(C) \to \mathbb{R}_{\geq 0}, t \mapsto \left\{\begin{array}{lll} 1/T(t),& T(t) \neq 0, \\
0 ,&\textnormal{else,} \end{array}\right.$$
and
$$T_{\textnormal{soft}}:\pi(C) \to \mathbb{R}_{\geq 0}, t \mapsto \left\{\begin{array}{lll} 1/T(t),& T(t) \neq 0, \\
	1,&\textnormal{else.} \end{array}\right.$$
Note that $T_{\textnormal{soft}}$ is $C$-regular if $T$ is $C$-regular (by Lemma~\ref{5.55}).  

\begin{lemma}
	\label{5.81}
	Let $T:\pi(C) \to \mathbb{R}_{\geq 0}$ be a $\kappa$-persistent function. Let $q:=(q_0, \ldots, q_r) \in \mathbb{Q}^{r+1}$ be with $(q_\kappa, \ldots ,q_r) \neq 0$ (i.e. $q$ is either $\kappa$-positive or $\kappa$-negative). Let $k:=k^{\text{ch}}$. Then the following holds.
	\begin{itemize}
		\item [(1)] Let $q_\kappa>0$. Then there is a $\kappa$-persistent $M:\pi(C) \to \mathbb{R}_{\geq 0}$ such that $\mathcal{B}(\mathcal{P}_{0,\kappa},M) \subset \mathcal{B}((\sigma\mathcal{Z})^{\otimes q},T)$.
		\item [(2)] Let $k \geq 0$. Suppose that $\kappa=k$ and let $q$ be $k$-positive with $q_k=0$. Suppose that $T_{\text{pers}}$ is consistent. 
		Then there is a $\kappa$-persistent $M:\pi(C) \to \mathbb{R}_{\geq 0}$ such that $\mathcal{B}(\mathcal{P}_{0,\kappa},M) \subset \mathcal{B}((\sigma\mathcal{Z})^{\otimes q},T)$.
		\item [(3)] Suppose that $T_{\textnormal{soft}}$ is $\kappa$-soft and $q_\kappa<0$. Then there is a $\kappa$-soft $N:\pi(C) \to \mathbb{R}_{\geq 0}$ and a $\kappa$-persistent $M:\pi(C) \to \mathbb{R}_{\geq 0}$ such that $\mathcal{A}(\mathcal{P}_{0,\kappa},N,M) \subset \mathcal{B}((\sigma\mathcal{Z})^{\otimes q},T).$
	\end{itemize}
Additionally, if $T$ is $C$-regular then $M$ and $N$ can be chosen in a $C$-regular way in (1) to (3).
\end{lemma}

\begin{proof}
	We may assume that $C$ is near with respect to $\mu_\kappa$ (compare with Definition~\ref{5.45} above). Otherwise we are done with the proof with Remark~\ref{5.46}(1) by choosing $M=0$ in (1) respectively (2) and $N=M=0$ in (3). Note that if $\kappa=k$ then $\sigma_{k+1}=-1$ and $\sigma_{k+2}= \ldots =\sigma_r=1$ and $\Theta_{\kappa+1}= \ldots =\Theta_r=0$ by Lemma~\ref{5.47} and therefore 
	$$\sigma_{k+j}z_{k+j}(t,z) = \log_{j-1}(-\log(\sigma_kz_k(t,z))) =  \log_j(1/(\sigma_kz_k(t,z))) \textnormal{    } (*)$$
	for $(t,z) \in H$ and $j \in \{1, \ldots ,r-k\}$. By Lemma~\ref{5.18} we have
	$$\prod_{j=0}^r \vert{z_j(t,z)}\vert^{q_j} = \left(\prod_{j=0}^{\kappa-1} \vert{\sigma_j\mathcal{P}_{j,\kappa}(t,z)+e^{\mu_{\kappa-j-1,\kappa}(t)}}\vert^{q_j} \right) \prod_{j=\kappa}^r \vert{z_j(t,z)}\vert^{q_j}$$
	for $(t,z)\in H$. We have $(t,z) \in \mathcal{B}((\sigma\mathcal{Z})^{\otimes q},T)$ if
	$$\prod_{j=\kappa}^r \vert{z_j(t,z)}\vert^{q_j} < T(t) \prod_{j=0}^{\kappa-1} \frac{e^{-q_j\mu_{\kappa-j-1,\kappa}(t)}}{\bigl|\frac{\sigma_j\mathcal{P}_{j,\kappa}(t,z)}{e^{\mu_{\kappa-j-1,\kappa}(t)}}+1\bigl|^{q_j}}$$
	for $(t,z) \in H$. By Lemma~\ref{5.67}(1) there is a $C$-regular $\kappa$-persistent $M^+:\pi(C) \to \mathbb{R}_{ \geq 0}$ such that for every $(t,z) \in \mathcal{B}(\mathcal{P}_{0,\kappa},M^+)$ and $j \in \{0, \ldots ,\kappa-1\}$ we have
	$$\Bigl|\frac{\sigma_j\mathcal{P}_{j,\kappa}(t,z)}{e^{\mu_{\kappa-j-1,\kappa}(t)}} \Bigl| < 1/2.$$
	
	So we have $(t,z) \in \mathcal{B}((\sigma\mathcal{Z})^{\otimes q},T)$ if
	\begin{align*}
	\prod_{j=\kappa}^r \vert{z_j(t,z)}\vert^{q_j} < \rho T(t)\prod_{j=0}^{\kappa-1} e^{-q_j\mu_{\kappa-j-1,\kappa}(t)} \textnormal{ }(+)
	\end{align*}
	for $(t,z) \in \mathcal{B}(\mathcal{P}_{0,\kappa},M^+)$ where $\rho:=\prod_{j=0}^{\kappa-1} (1/3)^{\vert{q_j}\vert}$. Set
	$$\phi:\pi(C) \to \mathbb{R}_{\geq 0}, t \mapsto \rho T(t)\prod_{j=0}^{\kappa-1} e^{-q_j\mu_{\kappa-j-1,\kappa}(t)}.$$
	Note that $\phi$ is $C$-regular by Example~\ref{5.56}(2) and Lemma~\ref{5.55}(2) if $T$ is $C$-regular and $\kappa$-persistent by Remark~\ref{5.42}(1) (compare also with Example~\ref{5.43}).
	
	To be able to do the necessary estimations with the iterated logarithms in $(+)$ and to construct the desired functions $M$ and $N$ in the appropriate cases below we need the following Claim~\ref{claim-11}.
	
	\begin{claim}
		\label{claim-11}
		Let $j \in \mathbb{N}$, $\nu \in \mathbb{N}$ and $z \in \mathbb{C}^-$ with $\exp_j(\nu^2\pi)<\vert{z}\vert$. Then $\vert{\log_j(z)}\vert < \sqrt[2^j \nu]{\vert{z}\vert}$.
	\end{claim}
	
	\begin{proof}
		For $p \in \mathbb{N}$ we set $c_p:=\exp(p^2\pi)$. We do an induction on $j$.
		
		$j=1$: Assume $c_\nu<\vert{z}\vert$. We show
		$$\sqrt{\vert{z}\vert^{-1/\nu}(\log^2(\vert{z}\vert)+\arg^2(z))} < 1.$$
		One sees that
		$$[c_\nu,\infty[ \textnormal{} \to \mathbb{R}_{\geq 0}, x \mapsto \log^2(x)x^{-1/\nu},$$
		is strictly monotone decreasing and that the function
		$$\mathbb{R}_{\geq 0} \to \mathbb{R}_{\geq 0}, x \mapsto x^4\pi^2/\exp(x\pi),$$ 
		takes its values on $]0,1/2[$. Therefore
		$$\log^2(c_\nu)c_\nu^{-1/\nu}=\log^2(c_\nu)/\exp(\nu\pi) = \nu^4\pi^2/\exp(\nu \pi)<1/2$$
		and consequently
		$$\log^2(\vert{z}\vert)\vert{z}\vert^{-1/\nu}<1/2.$$
		Clearly 
		$$\textnormal{arg}^2(z) \vert{z}\vert^{-1/\nu}<1/2.$$
		
		$j-1 \to j$: Assume $\exp_j(\nu^2 \pi) < \vert{z}\vert$. We have by the definition of the complex logarithm
		$$c_1 \leq c_m < \log_{j-1}(\vert{z}\vert) \leq \vert{\log_{j-1}(z)}\vert.$$
		Therefore by the base case for $\nu=1$ applied to $\log_{j-1}(z)$
		$$\vert{\log_j(z)}\vert  = \vert{\log(\log_{j-1}(z))}\vert < \sqrt{\vert{\log_{j-1}(z)}\vert}$$ 
		and by the inductive hypothesis
		\begin{align*}
		& \sqrt{\vert{\log_{j-1}(z)}\vert} < \sqrt{\sqrt[2^{j-1} \nu]{\vert{z}\vert}} = \sqrt[2^j \nu]{\vert{z}\vert}
		\end{align*}
	which proves the claim.
	\end{proof}
	
	If $k \geq 0$ we set
	$$\nu:=\max\{\lceil{\vert{q_j/q_b}\vert}\rceil \mid j \in \{b, \ldots ,r\}\}$$
	where $b=k$ if $q_k \neq 0$ and $b=l_{\min}(q)$ if $q_k=0$ (see Definition~\ref{5.32} for the notion of $l_{\min}(q)$). Note that $\nu \geq 1$.\\
	
	(1): Suppose that $q_\kappa>0$. We do two cases.
	
	\textbf{Case 1}:
	Assume $\kappa=r$ or if $\kappa=k$ (i.e. $k \geq 0$) then $q_{k+1}= \ldots =q_r=0$. Let
	$$\hat{M}:\pi(C) \to \mathbb{R}_{\geq 0}, t \mapsto \sqrt[q_\kappa]{\phi(t)}.$$
	Then $\hat{M}$ is $\kappa$-persistent by Remark~\ref{5.42}(2) and $C$-regular by Lemma \ref{5.55}(2) if $\phi$ and hence, if $T$ is $C$-regular. If $\vert{z_\kappa(t,z)}\vert<\hat{M}(t)$ then $\vert{z_\kappa(t,z)}\vert^{q_\kappa} < \phi(t)$ for $(t,z) \in H$ (so the inequality in $(+)$ is fulfilled for such a $(t,z)$). By Lemma~\ref{5.67}(1) there is a $\kappa$-persistent $M^*:\pi(C) \to \mathbb{R}_{\geq 0}$ such that $\mathcal{B}(\mathcal{P}_{0,\kappa},M^*) \subset \mathcal{B}(z_\kappa,\hat{M})$ which is also $C$-regular if $\hat{M}$ (and hence if $\phi$ and thus $T$) is $C$-regular. So take $M:\pi(C) \to \mathbb{R}_{\geq 0}, t \mapsto \min\{M^*(t),M^+(t)\}$. Clearly $M$ is $\kappa$-persistent by Remark~\ref{5.42}(3) and also $C$-regular by Lemma~\ref{5.55}(2) if $\phi$ and hence $T$ is $C$-regular. With $(+)$ we have $\mathcal{B}(\mathcal{P}_{0,\kappa},M) \subset \mathcal{B}((\sigma\mathcal{Z})^{\otimes q},T)$.
	
	\textbf{Case 2}: Assume $\kappa=k$ (i.e. $k \geq 0$ and $q_k>0$), and $(q_{k+1}, \ldots ,q_r) \neq 0$. Set
	$$\hat{M}:\pi(C) \to \mathbb{R}_{\geq 0}, t \mapsto \min\{\sqrt[s]{\phi(t)},1/\exp_{r-k}(\nu^2\pi)\},$$
	where $s:=\frac{q_k}{2^{r-k}}$. Then $\hat{M}$ $k$-persistent and $C$-regular if $T$ is $C$-regular.
	Let $(t,z) \in \mathcal{B}(z_k,\hat{M})$. We show that the estimation in $(+)$ holds. We have $\vert{z_k(t,z)}\vert < \sqrt[q_k]{\phi(t)^{2^{r-k}}}$ implying
	$$\vert{z_k(t,z)}\vert \prod_{j=1}^{r-k} \vert{z_k(t,z)}\vert^{-1/2^j} = \vert{z_k(t,z)}\vert^{1/2^{r-k}} < \sqrt[q_k]{\phi(t)}$$
	and therefore
	$$\vert{z_k(t,z)}\vert \prod_{j=1}^{r-k} (1/\vert{z_k(t,z)}\vert)^{1/2^j} < \sqrt[q_k]{\phi(t)}.$$
	We obtain
	$$\vert{z_k(t,z)}\vert \prod_{j=1}^{r-k} (1/\vert{z_k(t,z)}\vert)^{\Bigl|\tfrac{q_{k+j}}{2^j \nu q_k}\Bigl|} < \sqrt[q_k]{\phi(t)}$$
	since $1/\vert{z_k(t,z)}\vert = 1/\vert{\sigma_kz_k(t,z)}\vert > \exp_{r-k}(\nu^2\pi)$ and $\vert{\tfrac{q_{k+j}}{\nu q_k}}\vert \leq 1$. With Claim~\ref{claim-11} applied to $1/(\sigma_kz_k(t,z))$ we obtain
	$$\vert{\log_j(1/(\sigma_kz_k(t,z)))}\vert < \sqrt[2^j \nu]{1/\vert{z_k(t,z)}\vert}$$
	for $j \in \{1, \ldots ,r-k\}$ and therefore
	$$\vert{z_k(t,z)}\vert \prod_{j=1}^{r-k} \vert{\log_j(1/(\sigma_kz_k(t,z)))}\vert^{\vert{q_{k+j}/q_k}\vert} < \sqrt[q_k]{\phi(t)}.$$
	Because $\vert{\log_j(1/(\sigma_kz_k(t,z)))}\vert>1$ for every $j \in \{1, \ldots ,r-k\}$ (due to the inequality $1/\vert{z_k(t,z)}\vert>\exp_{r-k}(\nu^2 \pi)$) we obtain
	$$\vert{z_k(t,z)}\vert \prod_{j=1}^{r-k} \vert{\log_j(1/(\sigma_kz_k(t,z)))}\vert^{q_{k+j}/q_k} < \sqrt[q_k]{\phi(t)}$$
	which yields 
	$$\vert{z_k(t,z)}\vert^{q_k} \prod_{j=1}^{r-k} \vert{\log_j(1/(\sigma_kz_k(t,z)))}\vert^{q_{k+j}} < \phi(t).$$
	With $(*)$ we get
	$$\prod_{j=k}^r \vert{z_j(t,z)}\vert^{q_j} < \phi(t)$$
	which shows that the inequality in $(+)$ is satisfied. 
	Now define $M$ outgoing from $\hat{M}$ in the same way as in Case 1. With $(+)$ we have $\mathcal{B}(\mathcal{P}_{0,k},M) \subset \mathcal{B}((\sigma\mathcal{Z})^{\otimes q},T)$.
	
	(2): Suppose that $\kappa=k$, that $q$ is $k$-positive with $q_k=0$ and that $T_{\textnormal{pers}}$ is $C$-consistent.
	We have that $\phi_{\textnormal{pers}}$ is $C$-consistent by Example~\ref{5.51}(2), Remark~\ref{5.52}(3) (since $T_{\textnormal{pers}}$ is $C$-consistent). 
	Let $\mathfrak{l}:=l_{\min}(q)$. Since $q$ is $\kappa$-positive we have $q_\mathfrak{l}<0$. Set
	$$\hat{M}:\pi(C) \to \mathbb{R}_{\geq 0}, t \mapsto $$
	$$\left\{\begin{array}{lll} \min \{1/\exp_{\mathfrak{l}-k}(\sqrt[s]{\phi_{\textnormal{pers}}(t)}),1/\exp_{r-k}(\nu^2\pi)\} ,& \phi_{\textnormal{pers}}(t) \neq 0, \\
		0,&\textnormal{else,} \end{array}\right.$$
	where $s:=-(1/2)^{r-\mathfrak{l}}q_\mathfrak{l}$ and $\phi_{\textnormal{pers}}(t)=1/\phi(t)$ if $\phi(t)\neq 0$ and $0$ otherwise. Since $\phi_{\text{pers}}$ is $C$-consistent we obtain that $\sqrt[s]{\phi_{\textnormal{pers}}}$ is $C$-consistent by Remark~\ref{5.52}(2). By Remark~\ref{5.52}(3) we see that $\exp_j(\sqrt[s]{\phi_{\textnormal{pers}}})$ is also $C$-consistent for $j \in \{1, \ldots ,\mathfrak{l}-k\}$. So we see that $\hat{M}$ is $C$-regular by Lemma~\ref{5.55}(2) if $\phi_{\text{pers}}$ (and hence, if $T$) is $C$-regular. By Corollary~\ref{cor:persistence-soft-restriction} we see that $\hat{M}$ is $\kappa$-persistent.
	Let $(t,z) \in \mathcal{B}(z_k,\hat{M})$. Then $\phi_{\textnormal{pers}}(t) \neq 0$. We show that the inequality in $(+)$ is satisfied for this $(t,z)$. Let $y:=\log_{\mathfrak{l}-k}(1/(\sigma_kz_k(t,z)))$. Note that
	$$\max\bigl\{\exp_{\mathfrak{l}-k}\bigl(\sqrt[s]{\phi_{\textnormal{pers}}(t)}\bigl),\exp_{r-k}(\nu^2\pi)\bigl\} < 1/\vert{z_k(t,z)}\vert$$
	and hence
	$$\max\bigl\{\sqrt[s]{\phi_{\textnormal{pers}}(t)},\exp_{r-\mathfrak{l}}(\nu^2\pi)\bigl\} < \log_{\mathfrak{l}-k}(1/\vert{z_k(t,z)}\vert)\leq \vert{y}\vert.$$
	The last inequality is due to $\vert{\log_a(z)}\vert \geq  \log_a(\vert{z}\vert)$ 
	for $z \in \mathbb{C}$ with $\vert{z}\vert>\exp_a(0)$. With $\phi_{\textnormal{pers}}(t)=1/\phi(t)$ we see $ \vert{y}\vert > \sqrt[q_{\mathfrak{l}}]{\phi(t)^{2^{r-\mathfrak{l}}}}$ which yields $\vert{y}\vert^{1/2^{r-\mathfrak{l}}} > \sqrt[q_{\mathfrak{l}}]{\phi(t)}$. The latter implies
	$$\vert{y}\vert \prod_{j=\mathfrak{l}+1}^r \vert{y}\vert^{-1/2^{j-\mathfrak{l}}} > \sqrt[q_\mathfrak{l}]{\phi(t)}$$
	and therefore (due to $q_{\mathfrak{l}}<0$)
	$$\vert{y}\vert \prod_{j=\mathfrak{l}+1}^r (1/\vert{y}\vert)^{1/2^{j-\mathfrak{l}}} > \sqrt[q_\mathfrak{l}]{\phi(t)}.$$
	This gives
	$$\vert{y}\vert \prod_{j=\mathfrak{l}+1}^r (1/\vert{y}\vert)^{\Bigl|\tfrac{q_j}{2^{j-\mathfrak{l}} \nu q_\mathfrak{l}}\Bigl|} > \sqrt[q_\mathfrak{l}]{\phi(t)}$$
	since $\vert{y}\vert>1$ and $\vert{\tfrac{q_j}{\nu q_\mathfrak{l}}}\vert \leq 1$. Since $\exp_{r-\mathfrak{l}}(\nu^2\pi)<\vert{y}\vert$ we obtain with Claim~\ref{claim-11} applied to $y$
	$$\vert{\log_{j-\mathfrak{l}}(y)}\vert < \sqrt[2^{j-\mathfrak{l}}\nu]{\vert{y}\vert}$$
	for $j \in \{\mathfrak{l}+1, \ldots ,r\}$ and therefore
	$$\vert{y}\vert \prod_{j=\mathfrak{l}+1}^r \vert{\log_{j-\mathfrak{l}}(y)}\vert^{-\vert{q_j/q_\mathfrak{l}}\vert} > \sqrt[q_\mathfrak{l}]{\phi(t)}.$$
	Because $\vert{\log_{j-\mathfrak{l}}(y)}\vert>1$ for every $j \in \{\mathfrak{l}+1, \ldots ,r\}$ we obtain by pluggin in $\log_{\mathfrak{l}-k}(1/(\sigma_kz_k(t,z)))$ for $y$ and the previous inequality
	$$\prod_{j=\mathfrak{l}}^r \vert{\log_{j-k}(1/(\sigma_kz_k(t,z)))}\vert^{q_j/q_\mathfrak{l}} > \sqrt[q_\mathfrak{l}]{\phi(t)}$$
	and therefore (due to $q_{\mathfrak{l}}<0$)
	$$\prod_{j=\mathfrak{l}}^r \vert{\log_{j-k}(1/(\sigma_kz_k(t,z)))}\vert^{q_j} < \phi(t).$$
	Using $(*)$ we see
	$$\prod_{j=\mathfrak{l}}^r \vert{z_j(t,z)}\vert^{q_j} < \phi(t)$$
	which shows that the inequality in $(+)$ is satisfied.  
	Now define $M$ outgoing from $\hat{M}$ in the same way as in $(1)$ in Case 1. With $(+)$ we have $\mathcal{B}(\mathcal{P}_{0,k},M) \subset \mathcal{B}((\sigma\mathcal{Z})^{\otimes q},T)$.
	
	(3): Suppose that $T_{\textnormal{soft}}$ is $\kappa$-soft and $q_\kappa<0$. We see with Remark~\ref{5.52}(3) and Remark~\ref{5.42}(4) that $\phi_{\textnormal{soft}}$ is also $\kappa$-soft. We have to consider two more cases.
	
	\textbf{Case 1}: Assume $\kappa=r$ or if $\kappa=k$ then $q_{k+1}= \ldots =q_r=0$. Let
	$$\hat{N}:\pi(C) \to \mathbb{R}_{\geq 0}, t \mapsto \sqrt[-q_\kappa]{\phi_{\textnormal{soft}}(t)}.$$
	Then $\hat{N}$ is $\kappa$-soft and also $C$-regular if $T$ is $C$-regular. For $(t,z) \in H$ with $\hat{N}(t)<\vert{z_\kappa(t,z)}\vert<1$ we see that the inequality in $(+)$ is satisfied (since $\phi(t) \neq 0$ and hence $\phi_{\text{soft}}(t) = 1/\phi(t)$). By Lemma~\ref{5.67}(2) there is a $C$-regular $\kappa$-soft $N:\pi(C) \to \mathbb{R}_{\geq 0}$ and a $C$-regular $\kappa$-persistent $M^*:\pi(C) \to \mathbb{R}_{\geq 0}$ such that $\mathcal{A}(\mathcal{P}_{0,\kappa},N,M^*) \subset \mathcal{A}(z_\kappa,\hat{N},1)$. So take $M:\pi(C) \to \mathbb{R}_{\geq 0}, t \mapsto \min\{M^*(t),M^+(t)\}$. With (+) we see that $\mathcal{A}(\mathcal{P}_{0,\kappa},N,M) \subset \mathcal{B}((\sigma\mathcal{Z})^{\otimes q},T)$. 
	
	\textbf{Case 2:} Assume $\kappa=k$ and $(q_{k+1}, \ldots ,q_r) \neq 0$. Let $s:=-(2-(1/2)^{r-k})q_k$. Set
	$$\hat{N}:\pi(C) \to \mathbb{R}_{>0}, t \mapsto \sqrt[s]{\phi_{\textnormal{soft}}(t)},$$
	and
	$$\hat{M}:\pi(C) \to \mathbb{R}_{>0}, t \mapsto 1/\exp_{r-k}(\nu^2 \pi).$$
	Then $\hat{N}$ is $k$-soft and $C$-regular if $T$ is $C$-regular. 
	Let $(t,z) \in \mathcal{A}(z_\kappa,\hat{N},\hat{M})$. Then $\phi(t) \neq 0$. We show that the inequality in $(+)$ is satisfied for this $(t,z)$. We obtain $\vert{z_k(t,z)}\vert^{2-(1/2)^{r-k}}>\sqrt[-q_k]{\phi_{\text{soft}}(t)}$ and hence with $\phi(t)=1/\phi_{\text{soft}}(t)$ and the geometric sum $\sum_{j=0}^{r-k} (1/2)^j = 2-(1/2)^{r-k}$
	$$\vert{z_k(t,z)}\vert \prod_{j=1}^{r-k} \vert{z_k(t,z)}\vert^{\tfrac{1}{2^j}} > \sqrt[q_k]{\phi(t)}$$
	and therefore
	$$\vert{z_k(t,z)}\vert \prod_{j=1}^{r-k} (1/\vert{z_k(t,z)}\vert)^{-\tfrac{1}{2^j}} > \sqrt[q_k]{\phi(t)}.$$
	We obtain
	$$\vert{z_k(t,z)}\vert \prod_{j=1}^{r-k} (1/\vert{z_k(t,z)}\vert)^{-\Bigl|\tfrac{q_{k+j}}{2^j \nu q_k}\Bigl|} > \sqrt[q_k]{\phi(t)}$$
	since $1/\vert{z_k(t,z)}\vert>\exp_{r-k}(\nu^2 \pi)$ and $\vert{\tfrac{q_{k+j}}{\nu q_k}}\vert \leq 1$. With Claim~\ref{claim-11} applied to $1/(\sigma_kz_k(t,z))$ we obtain
	$$\vert{\log_j(1/(\sigma_kz_k(t,z)))}\vert < \sqrt[2^j \nu]{1/\vert{z_k(t,z)}\vert}$$
	for $j \in \{1, \ldots ,r-k\}$ and therefore
	$$\vert{z_k(t,z)}\vert \prod_{j=1}^{r-k} \vert{\log_j(1/(\sigma_kz_k(t,z)))}\vert^{-\vert{q_{k+j}/q_k}\vert} > \sqrt[q_k]{\phi(t)}.$$
	Because $\vert{\log_j(1/(\sigma_kz_k(t,z)))}\vert>1$ for every $j \in \{1, \ldots ,r-k\}$ we obtain
	$$\vert{z_k(t,z)}\vert \prod_{j=1}^{r-k} \vert{\log_j(1/(\sigma_kz_k(t,z)))}\vert^{q_{k+j}/q_k} > \sqrt[q_k]{\phi(t)}$$
	and due to $q_k<0$
	$$\vert{z_k(t,z)}\vert^{q_k} \prod_{j=1}^{r-k} \vert{\log_j(1/(\sigma_kz_k(t,z)))}\vert^{q_{k+j}} < \phi(t).$$
	This implies with $(*)$
	$$\prod_{j=k}^{r} \vert{z_j(t,z)}\vert^{q_j} < \phi(t),$$
	i.e. the inequality in $(+)$ is satisfied. Now define $N$ and $M$ analoguously as in Case 1 (outgoing from $\hat{N}$ and $\hat{M}$ instead from $\hat{N}$ and $1$). With $(+)$ we have $\mathcal{A}(\mathcal{P}_{0,\kappa},N,M) \subset \mathcal{B}((\sigma\mathcal{Z})^{\otimes q},T)$. In total we are done with the proof of Lemma~\ref{5.81}.
\end{proof}

\begin{lemma}
	\label{5.82}
	Let $T:\pi(C) \to \mathbb{R}$ be $\kappa$-persistent. Suppose that $T_{\textnormal{soft}}$ is $\kappa$-soft. Let $q:=(q_0, \ldots, q_r) \in \mathbb{Q}^{r+1}$ be $\kappa$-negative. Then there is a $\kappa$-soft $N:\pi(C) \to \mathbb{R}_{\geq 0}$ and a $\kappa$-persistent $M:\pi(C) \to \mathbb{R}_{\geq 0}$ such that
	$$\mathcal{A}(\mathcal{P}_{0,\kappa},N,M) \subset \mathcal{B}((\sigma\mathcal{Z})^{\otimes q},T).$$
	If $T$ is $C$-regular then $N$ and $M$ can be chosen in a $C$-regular way.
\end{lemma}

\begin{proof}
	Note that $\sqrt{T}$ is $\kappa$-persistent and that $\sqrt{T_{\textnormal{soft}}}$ is $\kappa$-soft which are both $C$-regular if $T$ is $C$-regular. Let $q_\kappa^*:=\vert{q_\kappa}\vert+1$. Let 
	$$\hat{\mathcal{Z}}:=(\sigma_\kappa z_\kappa)^{q_\kappa-q_\kappa^*}\prod_{j \neq \kappa}(\sigma_jz_j)^{q_j}.$$
	Then we have 
	$$\mathcal{B}((\sigma_\kappa z_\kappa)^{q_\kappa^*},\sqrt{T}) \cap \mathcal{B}(\hat{\mathcal{Z}},\sqrt{T}) \subset \mathcal{B}((\sigma\mathcal{Z})^{\otimes q},T).$$
	With Lemma~\ref{5.81}(1) respectively Lemma~\ref{5.81}(3) we find $\kappa$-persistent functions $M_1,M_2:\pi(C) \to \mathbb{R}_{\geq 0}$ and a $\kappa$-soft function $N:\pi(C) \to \mathbb{R}_{\geq 0}$ such that
	$\mathcal{B}(\mathcal{P}_{0,\kappa},M_1) \subset \mathcal{B}((\sigma_\kappa z_\kappa)^{q^*_\kappa},\sqrt{T})$ and $\mathcal{A}(\mathcal{P}_{0,\kappa},N,M_2) \subset \mathcal{B}(\hat{\mathcal{Z}},\sqrt{T})$. Note that $N,M_1$ and $M_2$ can be chosen in a $C$-regular way if $T$ is $C$-regular. Set 
	\[
	M:\pi(C) \to \mathbb{R}_{\geq 0}, t \mapsto \min\{M_1(t),M_2(t)\}.
	\]
	We see that $\mathcal{A}(\mathcal{P}_{0,\kappa},N,M) \subset \mathcal{B}((\sigma\mathcal{Z})^{\otimes q},T)$.
\end{proof}

The next result is about "re-preparing" products of the form $\prod_{j=0}^{\kappa-1} (\sigma_\kappa z_\kappa(t,z))^{\nu_j}$ for $\nu:=(\nu_0, \ldots , \nu_{\kappa-1}) \in \mathbb{Q}^\kappa$: We construct a $\kappa$-persistent $T:\pi(C) \to \mathbb{R}_{\geq 0}$ such that 
the function $$\mathcal{B}(\mathcal{P}_{0,\kappa},T) \to \mathbb{C}, (t,z) \mapsto \prod_{j=0}^{\kappa-1} (\sigma_\kappa z_\kappa(t,z))^{\nu_j},$$
is log-analytically prepared with preparing tuple $(r,\mathcal{Z},a_\nu,q_\nu,s_\nu,V_\nu,b_\nu,P_\nu)$ such that $q_\nu=0$, and the first $\kappa$ and last $r-\kappa$ columns of $P_\nu$ are zero.\\ 

But at first we have to define complex log-analytically prepared functions in a more general framework than in Definition~\ref{5.76}. 

\begin{definition}
\label{def:complex-log-prepared-general}
Let $Y \subset H$ be definable. We call a function $F:Y \to \mathbb{C}$ \textbf{complex $r$-log-analytically prepared in $z$ with center $\Theta$}\index{complex!log-analytically prepared} 
If for $(t,z) \in Y$
\[
F(t,z) = a(t)(\sigma\mathcal{Z})^{\otimes q}(t,z)V(b_1(t)(\sigma\mathcal{Z})^{\otimes p_1}(t,z), \ldots ,b_s(t)(\sigma\mathcal{Z})^{\otimes p_1}(t,z))
\]
where $a,b_1, \ldots, b_s:\pi(Y) \to \mathbb{R}$ are definable functions where $a=0$ or $a$ does not have a zero, $b_1, \ldots, b_s$ do not have a zero, $q,p_1, \ldots, p_s \in \mathbb{Q}^{r+1}$, and $V$ is a complex power series which converges absolutely on an open neighborhood of $\overline{D^s(0,1)}$ and for $(t,z) \in Y$ 
\[
\vert{b_j(t)(\sigma\mathcal{Z})^{\otimes p_j}(t,z)}\vert \leq 1.
\]
We call $a$ \textbf{coefficient} and $b:=(b_1, \ldots, b_s)$ a tuple of \textbf{base functions} for $F$ \index{coefficient}\index{base functions}. Further
$(r,\mathcal{Z},a,q,s,V,b,P)$ is a \textbf{complex LA-preparing tuple}\index{complex LA-preparing tuple} for $F$ where $p_j:=(p_{j0}, \ldots , p_{jr})$ and
$$P:=\left(\begin{array}{cccc}
	p_{10}&\cdot&\cdot&p_{1r}\\
	\cdot&& &\cdot\\
	\cdot&& &\cdot\\
	p_{s0}&\cdot&\cdot&p_{sr}\\
\end{array}\right)\in M\big(s\times (r+1),\mathbb{Q}).$$
\end{definition}

\begin{remark}
	\label{rem:complex-log-prepared-general-vs-special}
	A complex $r$-log-analytically prepared function $F:\Gamma \to \mathbb{C}$ according to Definition~\ref{5.76} is also complex $r$-log-analytically prepared according to Definition~\ref{def:complex-log-prepared-general}. 
\end{remark}

\begin{lemma}
	\label{lem:complex-preparation-simplified}
	Let $1<\delta$ and $\nu:=(\nu_0, \ldots , \nu_{\kappa-1}) \in \mathbb{Q}^\kappa$. There is a $C$-regular $\kappa$-persistent $T:\pi(C) \to \mathbb{R}_{\geq 0}$ such that $\prod_{j=0}^{\kappa-1} (\sigma_jz_j)^{\nu_j}|_{\mathcal{B}(\mathcal{P}_{0,\kappa},T)}$ is complex $r$-log-analytically prepared in $z$ with complex LA-preparing tuple 
	$$(r,\mathcal{Z},a,q,s,V,b,P)$$
	such that $q=0$, $\vert{V(\overline{D^s(0,1)})}\vert \subset [1/\delta,\delta]$, $a,b$ are $C$-regular and $C$-consistent without a zero, and
	\[
	P:=\left(\begin{array}{cccc}
		p_{10}&\cdot&\cdot&p_{1r}\\
		\cdot&& &\cdot\\
		\cdot&& &\cdot\\
		p_{s0}&\cdot&\cdot&p_{sr}\\
	\end{array}\right)\in M\big(s \times (r+1),\mathbb{Q})
	\]
	with 
	$p_{l \kappa} = 1$ and $p_{lj} =0$ for $j \in \{0, \ldots ,r\} \setminus \{\kappa\}$ and $l \in \{1, \ldots ,s\}$.
\end{lemma}

\begin{proof}
	Let $0 < \beta \leq 1/2$ be such that $1/\delta<\prod_{j=0}^{\kappa-1}\vert{1+z}\vert^{\nu_j}<\delta$ for $z \in B(0,\beta)$. 
	We have by Lemma~\ref{5.18} for $(t,z) \in H$
	\begin{align*}
	\prod_{j=0}^{\kappa-1}(\sigma_jz_j(t,z))^{\nu_j}=\prod_{j=0}^{\kappa-1} e^{\nu_j\mu_{\kappa-j-1,\kappa}(t)} \Bigl(1+\frac{\sigma_j\mathcal{P}_{j,\kappa}(t,z)}{e^{\mu_{\kappa-j-1,\kappa}(t)}} \Bigl)^{\nu_j} \textnormal{ }(+)
	\end{align*}
	The procedure is to prepare all the $\mathcal{P}_{j,\kappa}$ according to the following property $(*)_{\textnormal{prep}}$ at first and then to conclude on the final preparation according to $(+)$. \\
	We say that a function $J:Z \to \mathbb{C}$ where $Z \subset H$ is definable fulfills property $(*)_{\textnormal{prep}}$ if the following holds:
	$J$ is complex $r$-log-analytically prepared in $z$ with complex LA-preparing tuple
	$$(r,\hat{\mathcal{Z}},\hat{a},\hat{q},\hat{s},\hat{V},\hat{b},\hat{P})$$
	such that $\hat{a}$ and $\hat{b}$ are $C$-regular and $C$-consistent without a zero, $\hat{q}:=(\hat{q}_0, \ldots \hat{q}_r)$ is such that $\hat{q}_\kappa=1$ and $\hat{q}_j=0$ for $j \in \{0, \ldots ,r\} \setminus \{\kappa\}$, 
	\[
	\hat{P}:=\left(\begin{array}{cccc}
		\hat{p}_{10}&\cdot&\cdot&\hat{p}_{1r}\\
		\cdot&& &\cdot\\
		\cdot&& &\cdot\\
		\hat{p}_{\hat{s}0}&\cdot&\cdot&\hat{p}_{\hat{s} r}\\
	\end{array}\right)\in M\big(\hat{s} \times (r+1),\mathbb{Q})
	\]
	is such that $\hat{p}_{i \kappa}=1$ and $\hat{p}_{ij}=0$ for $j \in \{0, \ldots ,r\} \setminus \{\kappa\}$ and $i \in \{1, \ldots ,\hat{s}\}$ and $\hat{V}(\overline{D^{\hat{s}}(0,1)}) \subset \overline{B(0,\beta)}$.\\
	\begin{claim}
	\label{claim:property-prep}
	There is a $C$-regular $\kappa$-persistent $M:\pi(C) \to \mathbb{R}_{\geq 0}$ such that $\mathcal{P}_{j,\kappa}|_{\mathcal{B}(\mathcal{P}_{0,\kappa},M)}$ fulfills property $(*)_{\textnormal{prep}}$ for every $j \in \{0, \ldots, \kappa\}$.
	\end{claim}
	
	\begin{proof}
		We show by descending induction on $j \in \{0, \ldots, \kappa\}$ that there is a $C$-regular $\kappa$-persistent $T_j:\pi(C) \to \mathbb{R}_{\geq 0}$ with $T_j \leq T_{j+1}$ for $j>0$ such that $\mathcal{P}_{j,\kappa}|_{\mathcal{B}(\mathcal{P}_{0,\kappa},T_j)}$ fulfills property $(*)_{\textnormal{prep}}$. Then we are done with the proof of this claim by taking $M:=T_1$. For $j=\kappa$ we obtain the statement by choosing $T_\kappa:=1$ since $\mathcal{P}_{\kappa,\kappa}|_{\mathcal{B}(\mathcal{P}_{0,\kappa},1)}=z_\kappa|_{\mathcal{B}(\mathcal{P}_{0,\kappa},1)}$ (by Definition~\ref{5.16}) and $z_\kappa$ fulfills $(*)_{\textnormal{prep}}$ since we have $z_\kappa = 1/\beta \cdot z_\kappa \cdot \beta$ and consider $\tilde{V}:=\beta$ and $\hat{a}:=1/\beta$.\\
		$j+1 \to j$: For $(t,z) \in H$ we have 
		\[
		\mathcal{P}_{j,\kappa}(t,z) = \sigma_je^{\mu_{\kappa-j-1,\kappa}(t)}(e^{\mathcal{P}_{j+1,\kappa}(t,z)}-1)
		\]
		(see Definition~\ref{5.16} above). Let 
		$$(r,\mathcal{Z},\tilde{a},\tilde{q},\tilde{s},\tilde{V},\tilde{b},\tilde{P})$$
		be a corresponding complex LA-preparing tuple for $\mathcal{P}_{j+1,\kappa}|_{\mathcal{B}(\mathcal{P}_{0,\kappa},T_{j+1})}$. We have for $(t,z) \in \mathcal{B}(\mathcal{P}_{0,\kappa},T_{j+1})$ 
		\[
		\mathcal{P}_{j,\kappa}(t,z) = \sigma_je^{\mu_{\kappa-j-1,\kappa}(t)}(\exp(\tilde{a}(t)\sigma_\kappa z_\kappa(t,z)Q(t,z))-1)
		\]
		where
		\[
		Q(t,z) := \tilde{V}(\tilde{b}_1(t)\sigma_\kappa z_\kappa(t,z), \ldots ,\tilde{b}_{\tilde{s}}(t)\sigma_\kappa z_\kappa(t,z)).
		\]
		We see
		$$\mathcal{P}_{j,\kappa}(t,z)=\sigma_je^{\mu_{\kappa-j-1,\kappa}(t)}\tilde{a}(t)\sigma_\kappa z_\kappa(t,z) \frac{\exp(\tilde{a}(t)\sigma_\kappa z_\kappa(t,z)Q(t,z))-1}{\tilde{a}(t)\sigma_\kappa z_\kappa(t,z)}. \textnormal{  }(++)$$
		Note that $\tilde{a}$ does not have any zero. Because $\tilde{a}$ is $C$-consistent and $C$-regular (by the inductive hypothesis) we see that $1/\vert{\tilde{a}}\vert$ is $C$-regular and by Corollary~\ref{cor:persistence-soft-restriction} there is a $C$-regular $\kappa$-persistent $\tilde{a}_{\text{pers}}:\pi(C) \to \mathbb{R}_{\geq 0}$ such that $\tilde{a}_{\text{pers}}|_{\pi(\mathcal{B}(\mathcal{P}_{0,\kappa},T_{j+1}))} = 1/\vert{\tilde{a}}\vert$ (i.e. $\tilde{a}_{\text{pers}}(t)=1/\vert{\tilde{a}(t)}\vert$ if $t \in \pi(\mathcal{B}(\mathcal{P}_{0,\kappa},T_{j+1}))$ and $\tilde{a}_{\text{pers}}(t)=0$ otherwise). 
		Now let $\alpha>0$ be such that 
		\[
		\vert{(\exp(\alpha w_1 w_2)-1)/(\alpha w_1)}\vert \leq \beta \textnormal{  } (\#)
		\]
		for $w_1 \in \overline{B(0,1)} \setminus \{0\}$ and $w_2 \in \overline{B(0,\beta)}$. 
		By Lemma~\ref{5.81}(1) there is a $C$-regular $\kappa$-persistent $\hat{T}:\pi(C) \to \mathbb{R}_{\geq 0}$ such that $\mathcal{B}(\mathcal{P}_{0,\kappa},\hat{T}) \subset \mathcal{B}(z_\kappa,\alpha \cdot \tilde{a}_{\text{pers}})$. Note that $\vert{\tilde{a}(t)\sigma_\kappa z_\kappa(t,z)/\alpha}\vert<1$ for $(t,z) \in \mathcal{B}(\mathcal{P}_{0,\kappa},\hat{T})$. Note that $\sigma_je^{\mu_{\kappa-j-1,\kappa}}\tilde{a}$ is $C$-consistent and $C$-regular (since $e^{\mu_{\kappa-j-1,\kappa}}$ is $C$-consistent and $C$-regular). Take $T_j:\pi(C) \to \mathbb{R}_{\geq 0}, t \mapsto \min\{T_{j+1}(t),\hat{T}(t)\}$ which is also $C$-regular and $\kappa$-persistent. Then by using the exponential series and composition of power series we see with  the inductive hypothesis, $(++)$ and $(\#)$ that $\mathcal{P}_{j,\kappa}|_{\mathcal{B}(\mathcal{P}_{0,\kappa},T_j)}$ fulfills $(*)_{\textnormal{prep}}$ (where we plug in $\tilde{a}(t)\sigma_\kappa z_\kappa(t,z)/\alpha$ for $w_1$ and $Q(t,z)$ for $w_2$) which proves the claim (since $\tilde{V}(\overline{D^{\tilde{s}}(0,1)}) \subset \overline{B(0,\beta)}$). 
	\end{proof}

	Fix $M$ from Claim~\ref{claim:property-prep} and let $(r,\mathcal{Z},a^*_j,q^*_j,s^*,V^*_j,b^*,P^*_j)$ be a corresponding complex LA-preparing tuple for $\mathcal{P}_{j,\kappa}|_{\mathcal{B}(\mathcal{P}_{0,\kappa},M)}$ for every $j \in \{0, \ldots ,\kappa-1\}$. (Note that $b^*$ and $s^*$ can be chosen independently of $j$.) Then with $(+)$
	\begin{align*}
	\prod_{j=0}^{\kappa-1}(\sigma_j z_j)^{\nu_j}|_{\mathcal{B}(\mathcal{P}_{0,\kappa},M)} = \prod_{j=0}^{\kappa-1} e^{\nu_j\mu_{\kappa-j-1,\kappa}(t)} \Bigl(1+\frac{\sigma_ja_j^*\sigma_\kappa z_\kappa Q_j^*}{e^{\mu_{\kappa-j-1,\kappa}(t)}} \Bigl)^{\nu_j} \textnormal{ }
	\end{align*}
	where 
	\[
	Q_j^*:=V_j^*(b^*_1\sigma_\kappa z_\kappa, \ldots , b^*_{s^*}\sigma_\kappa z_\kappa)
	\]
	for $j \in \{0, \ldots, \kappa-1\}$.   
	 Since $1/e^{\mu_{\kappa-j-1,\kappa}}$ is $C$-regular and $C$-consistent (see Example~\ref{5.52}(2), Remark~\ref{5.53}(3) and Example~\ref{5.56}) and $a^*_j$ is also $C$-regular and $C$-consistent, we see with Corollary~\ref{cor:persistence-soft-restriction} that for every $j \in \{0, \ldots , \kappa-1\}$ there is a $C$-regular $\kappa$-persistent $T_j^*:\pi(C) \to \mathbb{R }_{\geq 0}$ with $T_j^*(t)=e^{\mu_{\kappa-j-1,\kappa}(t)}/\vert{a_j^*(t)}\vert$ for $t \in \pi(\mathcal{B}(\mathcal{P}_{0,\kappa},M))$. 
	 Hence, we find by Lemma~\ref{5.81}(1) a $C$-regular $\kappa$-persistent $\tilde{M}_j:\pi(C) \to \mathbb{R}_{\geq 0}$ such that $\mathcal{B}(\mathcal{P}_{0,\kappa},\tilde{M}_j) \subset \mathcal{B}(z_\kappa,T_j^*)$ for $j \in \{0, \ldots , \kappa-1\}$. Note that $$\Bigl\vert{\frac{\sigma_ja_j^*(t)\sigma_\kappa z_\kappa(t,z)}{e^{\mu_{\kappa-j-1,\kappa}(t)}}}\Bigl\vert = \Bigl\vert{\frac{a_j^*(t) z_\kappa(t,z)}{e^{\mu_{\kappa-j-1,\kappa}(t)}}}\Bigl\vert < 1$$
	 for $(t,z) \in \mathcal{B}(\mathcal{P}_{0,\kappa},\tilde{M}_j)$. So choose 
	 $$T:\pi(C) \to \mathbb{R}_{\geq 0}, t \mapsto \min\{\tilde{M}_1(t), \ldots , \tilde{M}_{\kappa-1}(t),M(t)\}.$$
	 Then $T$ is $C$-regular and $\kappa$-persistent. Also for $(w_1, \ldots , w_{s^*+1}) \in \overline{D^{s^*+1}(0,1)}$
	$$\vert{w_1Q^*_j(w_2, \ldots , w_{s^*+1})}\vert \leq \vert{Q^*_j(w_2, \ldots , w_{s^*+1})}\vert \leq \beta$$
	(by the property $(*)_{\text{prep}}$) and therefore
	\[
	W:=\prod_{j=0}^{\kappa-1} |1+w_1Q_j^*(w_2, \ldots , w_{s^*+1})|^{\nu_j}
	\]
	fulfills $1/\delta \leq W \leq \delta$ by the choice of $\beta$ above. Since $\beta \leq 1/2$ we see by composition of power series that $W$ can be written as a power series in $(w_1, \ldots ,w_{s^*+1})$ which converges absolutely on an open neighborhood of $\overline{D^{s^*+1}(0,1)}$. Thus, with the equation in (+) above and composition of power series we obtain the desired preparation for $\prod_{j=0}^{\kappa-1}(\sigma_jz_j)^{\nu_j}|_{\mathcal{B}(\mathcal{P}_{0,\kappa},T)}$ since
	\[
	\Bigl|\frac{\sigma_ja_j^*(t) \sigma_\kappa z_\kappa(t,z)}{e^{\mu_{\kappa-j-1,\kappa}(t)}} \Bigl| \leq 1 
	\]
	for $(t,z) \in \mathcal{B}(\mathcal{P}_{0,\kappa},T)$ and every $j \in \{0, \ldots, \kappa-1\}$, 
	and $\prod_{j=0}^{\kappa-1}e^{\nu_j \mu_{\kappa-j-1,\kappa}}$ is $C$-regular and $C$-consistent.
\end{proof}

\subsection{Restricted Log-Exp-Analytically Prepared \\ Functions}

In this subsection we prove the desired preparation theorem for the unary parametric global complexification $\Phi_f:\Lambda_f \to \mathbb{C}$ of a restricted log-exp-analytically prepared function $f:C \to \mathbb{R}$. Specifically, we construct a $\kappa$-persistent $M:\pi(C) \to \mathbb{R}_{\geq 0}$ and a $\kappa$-soft $N:\pi(C) \to \mathbb{R}_{\geq 0}$ such that $D:=\mathcal{A}(\mathcal{P}_{0,\kappa},N,M) \subset \Lambda_f$ and $\Phi_f|_D$ is complex log-analytically prepared. This constitutes the main theorem of this section.\\

Recall that $l,m \in \mathbb{N}_0$ with $l+m=n$, that $w:=(w_1, \ldots ,w_l)$ range over $\mathbb{R}^l$ and $u:=(u_1, \ldots , u_m)$ over $\mathbb{R}^m$. Here $(u,x)$ or $(u,z)$ are serving as tuples of independent variables of families of functions parameterized by $w$. Fix an $(m+1,X)$-restricted log-exp-analytically prepared function $f:C \to \mathbb{R}, (w,u,x) \mapsto f(w,u,x),$ in $x$ with center $\Theta$, the corresponding unary parametric global complexification $\Phi_f:\Lambda_f \to \mathbb{C}$ from Definition~\ref{5.78}. Further let $E$ be a finite set of positive definable functions on $C$ such that every $g \in \log(E)$ is locally bounded in $(u,x)$ with reference set $X$ such that $f$ is log-exp-analytically prepared in $x$ with respect to $E$. For $0 \leq \rho \leq r+1$ we denote by $M_\rho(s \times (r+1),\mathbb{Q})$ the set of all $(s \times (r+1))$-matrices whose first $\rho$ columns are zero.

We also need the following small notions from the previous Section~7.3. For a definable function $T:\pi(C) \to \mathbb{R}_{\geq 0}$ we set
$$T_{\textnormal{pers}}:\pi(C) \to \mathbb{R}_{\geq 0}, t \mapsto \left\{\begin{array}{lll} 1/T(t),& T(t) \neq 0, \\
	0 ,&\textnormal{else,} \end{array}\right.$$
and
$$T_{\textnormal{soft}}:\pi(C) \to \mathbb{R}_{\geq 0}, t \mapsto \left\{\begin{array}{lll} 1/T(t),& T(t) \neq 0, \\
1,&\textnormal{else.} \end{array}\right.$$

\begin{theorem}
	\label{5.84}
	Let $g \in \log(E) \cup \{f\}$. There is a $C$-regular $\kappa$-persistent $M:\pi(C) \to \mathbb{R}_{\geq 0}$ and a $C$-regular $\kappa$-soft $N:\pi(C) \to \mathbb{R}_{\geq 0}$ such that $\mathcal{A}(\mathcal{P}_{0,\kappa},N,M) \subset \Lambda_g$ and $G:=\Phi_g|_{\mathcal{A}(\mathcal{P}_{0,\kappa},N,M)}$ is complex $r$-log-analytically prepared in $z$ with a complex LA-preparing tuple
	$$(r,\mathcal{Z},a_G,q_G,s_G,V_G,b_G,P_G)$$
	where $q_G:=((q_G)_0, \ldots ,(q_G)_r)$, $b_G:=(b_{1,G}, \ldots ,b_{s_G,G})$ and the following holds.
	\begin{itemize}
		\item [(1)] We have $(q_G)_0, \ldots ,(q_G)_{\kappa-1}=0$ and $P_G \in M_\kappa(s \times (r+1), \mathbb{Q})$.
		\item [(2)] $a_G$ is $C$-regular. If $g \in \log(E)$ the following holds in addition: $a_G$ is $C$-consistent and if $q_G$ is $\kappa$-negative then 
		$a_G$ can also be extended to a $\kappa$-soft function $\hat{a}_G:\pi(C) \to \mathbb{R}_{\geq 0}$.
		\item [(3)] $b_G$ is $C$-regular, $C$-consistent and if $p_{i,G}$ is $\kappa$-negative for $i \in \{1, \ldots ,s_G\}$ then $b_{i,G}$ can also be extended to a $\kappa$-soft function $\hat{b}_{i,G}:\pi(C) \to \mathbb{R}_{\geq 0}$ where $P_G:=(p_{1,G}, \ldots ,p_{s_G,G})^t$. 
	\end{itemize}
\end{theorem}

\begin{proof}
	Let $g$ be log-exp-analytically $(l,r)$-prepared in $x$ with respect to $E$ for $l \in \{-1, \ldots ,e\}$. We do an induction on $l$. For $l=-1$ we have that $f=0$, and $f=g$ for every $g \in \log(E) \cup \{f\}$ (since $\log(E)=\emptyset$) and hence, there is nothing to show (since one can choose $a_G=0$). 
	
	$l-1 \to l:$ Let 
	$$(r,\mathcal{Y},a,e^{d_0},q,s,v,b,e^d,P)$$
	be a preparing tuple for $g$ where $b:=(b_1, \ldots ,b_s)$, $P:=(p_1, \ldots , p_s)^t$ for $p_i \in \mathbb{Q}^{r+1}$, and $\exp(d):=(\exp(d_1), \ldots , \exp(d_s))$.  
	By redefining $a$ and $v$ suitably we may assume that $\vert{v}\vert<\pi$.
	
	We may assume that $a$ does not have any zero. Otherwise $a=0$ and the statement is clear. Note that $d_i$ is locally bounded in $(u,x)$ with reference set $X$ for every $i \in \{0, \ldots ,s\}$. For $i \in \{0, \ldots ,s\}$ let $\Lambda_i:=\Lambda_{d_i}$. Let
	$$(r,\mathcal{Z},a,\exp(\tilde{D}_0),q,s,V,b,\exp(\tilde{D}),P)$$
	be a complex preparing tuple for $\Phi_g$ where $\tilde{D}_j:=\Phi_{d_j}$ for $j \in \{0, \ldots ,s\}$ and $\exp(\tilde{D}):=(\exp(\tilde{D}_1), \ldots ,\exp(\tilde{D}_s))$ (according to Definition~\ref{5.78}).
	
	By the inductive hypothesis there is a $C$-regular $\kappa$-soft $N^*:\pi(C) \to \mathbb{R}_{\geq 0}$ and a $C$-regular $\kappa$-persistent $M^*:\pi(C) \to \mathbb{R}_{\geq 0}$ such that for every $i \in \{0, \ldots ,s\}$
	$$B^*:= \mathcal{A}(\mathcal{P}_{0,\kappa},N^*,M^*) \subset \Lambda_i$$
	and $D_i:=\tilde{D}_i|_{B^*}$ is complex $r$-log-analytically prepared in $z$ with complex LA-preparing tuple $$(r,\mathcal{Z},a_{D_i},q_{D_i},s_{D_i},V_{D_i},b_{D_i},P_{D_i})$$
	where $q_{D_i}:=((q_{D_i})_0, \ldots ,(q_{D_i})_r) \in \mathbb{Q}^{r+1}$, $P_{D_i}:=(p_{1,D_i}, \ldots ,p_{s_{D_i},D_i})^t$ where $p_{j,D_i} \in \mathbb{Q}^{r+1}$, $b_{D_i}:=(b_{1,D_i}, \ldots ,b_{s_{D_i},D_i})$ and the following properties (1) - (3) hold:
	\begin{itemize}
		\item [(1)] $(q_{D_i})_0= \ldots =(q_{D_i})_{\kappa-1}=0$ and $P_{D_i} \in M_{\kappa}(s_{D_i} \times (r+1), \mathbb{Q})$.
		\item [(2)] $a_{D_i}$ is $C$-regular. If $g \in \log(E)$ the following holds in addition: $a_{D_i}$ is $C$-consistent and if $q_{D_i}$ is $\kappa$-negative then $a_{D_i}$ can also be extended to a $\kappa$-soft function $\hat{a}_{D_i}:\pi(C) \to \mathbb{R}_{\geq 0}$. 
		\item [(3)] $b_{D_i}$ is $C$-regular, $C$-consistent and if $p_{j,D_i}$ is $\kappa$-negative then 
		$b_{j,D_i}$ can also be extended to a $\kappa$-soft function $\hat{b}_{j,D_j}:\pi(C) \to \mathbb{R}_{\geq 0}$ for $j \in \{1, \ldots ,s_{D_i}\}$.
	\end{itemize}
	
	At first we deal with $\exp(D_0), \ldots ,\exp(D_s)$ and prepare those functions log-analytically such that properties (1)-(3) hold. This happens in Claim~\ref{claim-13}. Due to technical reasons we give even a few additional properties.
	
	\begin{claim}
		\label{claim-13}
		Let $0<\delta<1/2$. Then there is a $C$-regular $\kappa$-soft $S_\delta:\pi(C) \to \mathbb{R}_{\geq 0}$  and a $C$-regular $\kappa$-persistent $T_\delta:\pi(C) \to \mathbb{R}_{\geq 0}$ with $Y_\delta:=\mathcal{A}(\mathcal{P}_{0,\kappa},S_\delta,T_\delta) \subset B^*$ and for $i \in \{0, \ldots ,s\}$ definable functions $D_{i,1}:\pi(Y_\delta) \to \mathbb{R}, t \mapsto D_{i,1}(t),$ and $D_{i,2}:Y_\delta \to \mathbb{C}, (t,z) \mapsto D_{i,2}(t,z),$ with the following properties.
		\begin{itemize}
			\item [(1)*] 
			$D_i|_{Y_\delta}=D_{i,1}+D_{i,2}$.
			\item [(2)*] $D_{i,1}$ is $C$-regular and $C$-consistent.
			\item [(3)*] $D_{i,2}$ is bounded by $\delta$. Additionally, $\exp(D_{i,2})$ is complex $r$-log-analytically prepared in $z$ with complex LA-preparing tuple $(r,\mathcal{Z},a^*,q^*,s^*,V^*,b^*,P^*)$ such that properties (1) - (3) above are fulfilled where $a^*=1$ and $q^*=0$. 
		\end{itemize}
	
	Hence, $\exp(D_i)|_{Y_\delta}$ is complex log-analytically prepared such that properties (1)-(3) above are satisfied.
	\end{claim}
	
	\begin{proof}
		Let $K \in \{D_0, \ldots ,D_s\}$. By Remark~\ref{5.42}(3) and Lemma~\ref{5.55}(2) it suffices to find a $C$-regular $\kappa$-soft $\hat{S}_\delta:\pi(C) \to \mathbb{R}_{\geq 0}$ and a $C$-regular $\kappa$-persistent $\hat{T}_\delta:\pi(C) \to \mathbb{R}_{\geq 0}$ with $W:= \mathcal{A}(\mathcal{P}_{0,\kappa},\hat{S}_{\delta},\hat{T}_{\delta}) \subset B^*$ and corresponding functions $K_1:\pi(W) \to \mathbb{R}$ and $K_2:W \to \mathbb{C}$ which fulfill properties (1)* - (3)*. Let $$(r,\mathcal{Z},\mathfrak{a},\mathfrak{q},\mathfrak{s},\mathfrak{V},\mathfrak{b},\mathfrak{P}):=(r,\mathcal{Z},a_K,q_K,s_K,V_K,b_K,P_K).$$
		Suppose $\mathfrak{a}$ does not have any zero. Otherwise the claim follows. Fix $R>1$ such that $\mathfrak{V}$ converges absolutely on an open neighborhood of $\overline{D^\mathfrak{s}(0,R)}$. Let $\sum_{\alpha \in \mathbb{N}_0^{\mathfrak{s}}} c_\alpha Z^{\alpha}$ be the Taylor series expansion of $\mathfrak{V}$ where $Z:=(Z_1, \ldots ,Z_s)$. Fix $L \in \mathbb{R}_{>0}$ with $\delta<L$ such that $\sum_{\alpha \in \mathbb{N}^{\mathfrak{s}}}\vert{c_\alpha w^{\alpha}}\vert<L$ for $w \in D^\mathfrak{s}(0,R)$. We look at two cases.\\ 
		
		\textbf{1. Case:} $(\mathfrak{q}_\kappa, \ldots ,\mathfrak{q}_r) \neq 0$.\\ 
		Consider the definable function
		\[
		T_\delta^*: \pi(C) \to \mathbb{R}_{\geq 0}, t \mapsto \left\{\begin{array}{ll} \frac{\delta}{L\vert{\mathfrak{a}(t)}\vert} , & t \in \pi(B^*),\\
		0, & \textnormal{else}. \end{array}\right.
		\]
		Note that if $\vert{(\sigma\mathcal{Z})^{\otimes \mathfrak{q}}(t,z)}\vert < T_\delta^*(t)$ then $\vert{K(t,z)}\vert < \delta$ for $(t,z) \in B^*$. Note that $T_\delta^*$ is $C$-regular by Lemma~\ref{5.55}(2) and by the fact that $\mathfrak{a}$ is also $C$-regular. By Corollary~\ref{cor:persistence-soft-restriction} we see that $T_\delta^*$ is also $\kappa$-persistent.
		
		
		
		\textbf{1.1 Case:} $\mathfrak{q}$ is $\kappa$-positive.
		By Lemma~\ref{5.81}(1) if $q_\kappa > 0$ respectively Lemma~\ref{5.81}(2) if $\kappa=k \geq 0$ and $q_k=0$ there is a $C$-regular $\kappa$-persistent $M_\delta^+:\pi(C) \to \mathbb{R}_{ \geq 0}$ such that $\mathcal{B}(\mathcal{P}_{0,\kappa},M_\delta^+) \subset \mathcal{B}((\sigma\mathcal{Z})^{\otimes \mathfrak{q}},T_\delta^*)$ (since $(T_\delta^*)_{\text{pers}}$ is $C$-consistent since $\mathfrak{a}$ is $C$-consistent). Take $\hat{S}_\delta:=N^*$,
		$$\hat{T}_\delta:\pi(C) \to \mathbb{R}_{\geq 0}, t \mapsto \min\{M_\delta^+(t),M^*(t)\},$$
		$K_1:=0$ and $K_2:=K|_W$ where $W:=\mathcal{A}(\mathcal{P}_{0,\kappa},\hat{S}_\delta,\hat{T}_\delta)$. Note that $\exp(K_2)=\exp^*(K_2)$ where
		$$\exp^*:\mathbb{C} \to \mathbb{C}, z \mapsto \left\{\begin{array}{ll} \exp(z), & \vert{z}\vert \leq \delta,\\
			0, & \textnormal{ else,}\end{array}\right.$$ 
		is globally subanalytic. Additionally, we obtain $\vert{\mathfrak{a}(t)(\sigma\mathcal{Z})^{\otimes \mathfrak{q}}(t,z)}\vert < 1$ for $(t,z) \in W$. Hence, by composition of power series we get the desired preparation for $\exp(K_2)$.
		
		\textbf{1.2 Case:} $\mathfrak{q}$ is $\kappa$-negative. 
		By Lemma~\ref{lem:persistence-soft-restriction} we see that $(T_\delta^*)_{\textnormal{soft}}$ is $\kappa$-soft (since $\mathfrak{a}$ is $\kappa$-soft). 
		By Lemma~\ref{5.82} there is a $C$-regular $\kappa$-soft $N_\delta^+:\pi(C) \to \mathbb{R}_{\geq 0}$ and a $C$-regular $\kappa$-persistent $M_\delta^+:\pi(C) \to \mathbb{R}_{\geq 0}$ such that $\mathcal{A}(\mathcal{P}_{0,\kappa},N_\delta^+,M_\delta^+) \subset \mathcal{B}((\sigma\mathcal{Z})^{\otimes \mathfrak{q}},T_\delta^*)$. Take $\hat{S}_\delta:\pi(C) \to \mathbb{R}_{\geq 0}, t \mapsto \max\{N_\delta^+(t),N^*(t)\}$, $\hat{T}_\delta:\pi(C) \to \mathbb{R}_{\geq 0}, t \mapsto \min\{M_\delta^+(t),M^*(t)\}$, $K_1:=0$ and $K_2:=K|_W$. We are done with Case~1.2 in a completely similar way as in Case 1.1.\\
		
		\textbf{2. Case:} $(\mathfrak{q}_\kappa, \ldots ,\mathfrak{q}_r) = 0$. Then $$K=\mathfrak{a}\mathfrak{V}(\mathfrak{b}_1(\sigma\mathcal{Z})^{\otimes \mathfrak{p}_1}, \ldots ,\mathfrak{b}_\mathfrak{s}(\sigma\mathcal{Z})^{\otimes \mathfrak{p}_{\mathfrak{s}}})$$
		where $\mathfrak{p}_1, \ldots , \mathfrak{p}_{\mathfrak{s}} \in \{0\}^\kappa \times \mathbb{Q}^{r-\kappa+1}$. For $i \in \{1, \ldots ,\mathfrak{s}\}$ let $\mathfrak{p}_i=(\mathfrak{p}_{i0}, \ldots , \mathfrak{p}_{ir})$. Further we set 
		$$\mathcal{S}:=\{i \in \{1, \ldots ,\mathfrak{s}\} \mid (\mathfrak{p}_{i\kappa}, \ldots ,\mathfrak{p}_{ir}) \neq 0\}.$$
		If $\mathcal{S}=\emptyset$ then $K=\mathfrak{a} \mathfrak{V}(\mathfrak{b}_1, \ldots ,\mathfrak{b}_\mathfrak{s})$ is a $C$-regular and $C$-consistent function considered as function on $\pi(B^*)$. So we are done with Case~2 by taking $\hat{S}_\delta:=N^*$, $\hat{T}_\delta:=M^*$, $K_1:=K$ (where $K$ is considered as a function on $\pi(B^*)$) and $K_2:=0$. So suppose $\mathcal{S} \neq \emptyset$. Let
		$$\Gamma_1:=\{\alpha \in \mathbb{N}_0^{\mathfrak{s}} \mid {^t}P \alpha = 0\}$$
		and $\Gamma_2:= \mathbb{N}_{0}^\mathfrak{s} \setminus \Gamma_1.$ Note that $\Gamma_2 \neq \emptyset$. Set $\mathfrak{V}_1:=\sum_{\alpha \in \Gamma_1}c_{\alpha}Z^{\alpha}$ and $\mathfrak{V}_2:=\sum_{\alpha \in \Gamma_2}c_{\alpha}Z^{\alpha}$. For $l \in \{1,2\}$ let
		$$J_l :=  \mathfrak{a} \mathfrak{V}_l(\mathfrak{b}_1(\sigma\mathcal{Z})^{\otimes \mathfrak{p}_1}, \ldots ,\mathfrak{b}_\mathfrak{s}(\sigma\mathcal{Z})^{\otimes \mathfrak{p}_\mathfrak{s}}).$$
		Note that $K=J_1+J_2$. Let $\phi_i:=\mathfrak{b}_i(\sigma\mathcal{Z})^{\otimes \mathfrak{p}_i}$ for $i \in \{1, \ldots ,\mathfrak{s}\}$. For $j \in \mathcal{S}$ we consider the definable function
		$$U_{j,\delta}:\pi(C) \to \mathbb{R}_{\geq 0}, t \mapsto \left\{\begin{array}{ll} \min\Bigl\{\Bigl|\frac{\delta}{L \mathfrak{a}(t)\mathfrak{b}_j(t)} \Bigl|,\Bigl|\frac{\delta}{L \mathfrak{b}_j(t)} \Bigl|\}, & t \in \pi(B^*),\\
		0, & \textnormal{else}. \end{array}\right. $$
		Note that $U_{j,\delta}$ is $C$-regular (by Lemma~\ref{5.55}(2) since $\mathfrak{a}$ and $\mathfrak{b}_j$ are $C$-regular), and that $\vert{a(t)b_j(t) (\sigma\mathcal{Z})^{\otimes \mathfrak{p}_j}(t,z)}\vert \leq 1$ and $\vert{b_j(t) (\sigma\mathcal{Z})^{\otimes \mathfrak{p}_j}(t,z)}\vert \leq 1$ for $(t,z) \in \mathcal{B}((\sigma\mathcal{Z})^{\otimes \mathfrak{p}_j},U_{j,\delta})$ respectively. With Corollary~\ref{cor:persistence-soft-restriction} one sees for $j \in \mathcal{S}$ that $U_{j,\delta}$ is $\kappa$-persistent (since $\mathfrak{b}_j$ and $\mathfrak{a}$ are $C$-consistent by properties (2) and (3) above) and that if $\mathfrak{p}_j$ is $\kappa$-negative then $(U_{j,\delta})_{\textnormal{soft}}$ is $\kappa$-soft (by Remark~\ref{5.42}(4) since $\mathfrak{a}$ is $C$-consistent and $b_j$ is then $\kappa$-soft by properties (2) and (3) above). Additionally, $(U_{j,\delta})_{\textnormal{pers}}$ is $C$-consistent for $j \in \mathcal{S}$.
		
		Consequently, for $j \in \mathcal{S}$ there is a $C$-regular $\kappa$-soft $\hat{S}_{j,\delta}:\pi(C) \to \mathbb{R}_{\geq 0}$ and a $C$-regular $\kappa$-persistent $\hat{T}_{j,\delta}:\pi(C) \to \mathbb{R}_{\geq 0}$ such that
		$$\mathcal{A}(\mathcal{P}_{0,\kappa},\hat{S}_{j,\delta},\hat{T}_{j,\delta}) \subset \mathcal{B}((\sigma\mathcal{Z})^{\otimes \mathfrak{p}_j},U_{j,\delta}).$$
		(By Lemma~\ref{5.81}(1) if $\mathfrak{p}_j$ is $\kappa$-positive and $\mathfrak{p}_{j \kappa} \neq 0$ respectively Lemma~\ref{5.81}(2) if $\mathfrak{p}_j$ is $\kappa$-positive, $\kappa=k \geq 0$ and $\mathfrak{p}_{jk} = 0$. If $\mathfrak{p}_j$ is $\kappa$-negative we use Lemma~\ref{5.82}. If $\mathfrak{p}_j$ is $\kappa$-positive take $\hat{S}_{j,\delta}=0$.)
		Consider $\hat{S}_\delta :\pi(C) \to \mathbb{R}_{\geq 0}, t \mapsto \max(\{\hat{S}_{j,\delta}(t) \mid j \in \mathcal{S}\} \cup \{N^*(t)\})$ and $\hat{T}_\delta :\pi(C) \to \mathbb{R}_{\geq 0}, t \mapsto \min(\{ \hat{T}_{j,\delta}(t) \mid j \in \mathcal{S}\} \cup \{M^*(t)\})$. Then $\hat{S}_\delta$ is $C$-regular and $\kappa$-soft and $\hat{T}_\delta$ is $C$-regular and $\kappa$-persistent. Consider $W:=\mathcal{A}(\mathcal{P}_{0,\kappa},\hat{S}_\delta,\hat{T}_\delta)$. 
		
		To obtain Claim~\ref{claim-13} it remains to show that $\vert{J_2(t,z)}\vert<\delta$ for $(t,z) \in W$, $\exp(J_2)|_W$ is complex $r$-log-analytically prepared in $z$ as desired (which is done in Subclaim~1) and that $J_1|_W$ coincides with a $C$-regular and $C$-consistent function on $\pi(W)$.\\
		
		\textbf{Subclaim~1:} {\it For $(t,z) \in W$ we have $\vert{J_2(t,z)}\vert<\delta$. The function $\exp(J_2)|_W$ is complex $r$-log-analytically prepared in $z$ as required in property (3)*.}
		
		\begin{proof}
			Set $\phi:=(\phi_1, \ldots ,\phi_\mathfrak{s})$. For every $\alpha:=(\alpha_1, \ldots ,\alpha_\mathfrak{s}) \in \Gamma_2$ fix $i_\alpha \in \{1, \ldots ,\mathfrak{s}\}$ with $i_\alpha \in \mathcal{S}$ and $\alpha_{i_{\alpha}} \neq 0$. Therefore we have for $(t,z) \in W$
			$$\vert{\mathfrak{a}(t)\phi_{i_\alpha}(t,z)}\vert<\delta/L  < 1.$$
			We obtain
			\begin{eqnarray*}
				\vert{J_2(t,z)}\vert &=& \vert{\mathfrak{a}(t)\mathfrak{V}_2(\phi(t,z))}\vert\\
				&=& \vert{\mathfrak{a}(t) \sum_{\alpha \in \Gamma_1}c_{\alpha} \phi_{i_{\alpha}}(t,z)^{\alpha_{i_\alpha}} \prod_{j \neq i_{\alpha}} \phi_j(t,z)^{\alpha_j}}\vert\\
				&=&\vert{\sum_{\alpha \in \Gamma_1}c_{\alpha}\mathfrak{a}(t)\phi_{i_{\alpha}}(t,z) \phi_{i_{\alpha}}(t,z)^{\alpha_{i_\alpha}-1} \prod_{j \neq i_{\alpha}} \phi_j(t,z)^{\alpha_j}}\vert\\
				&<& \delta/L \sum_{\alpha \in \Gamma_1}\vert{c_{\alpha} \phi_{i_{\alpha}}(t,z)^{\alpha_{i_\alpha}-1} \prod_{j \neq i_{\alpha}} \phi_j(t,z)^{\alpha_j}}\vert\\
				&<& \delta.
			\end{eqnarray*}
			
			Let $\mathcal{S}=\{j_1, \ldots ,j_\beta\}$ where $\beta \in \mathbb{N}$. 
			Let $(w_{j_1}, \ldots ,w_{j_\beta})$ be a new tuple of complex variables. Consider for $\gamma \in \{1, \ldots ,\beta\}$
			$$\tilde{\phi}_\gamma:W \to \overline{B(0,1)}, (t,z) \mapsto a(t)b_{j_\gamma}(t)(\sigma\mathcal{Z})^{\otimes \mathfrak{p}_{j_\gamma}}(t,z).$$
			Consider 
			$$\hat{V}:\overline{D^{\mathfrak{s}+\beta}(0,1)} \to \mathbb{C}, (z_1, \ldots ,z_s,w_{j_1}, \ldots ,w_{j_\beta}) \mapsto \sum_{\alpha \in \Gamma_2} c_{\alpha}w_{i_\alpha}z_{i_\alpha}^{\alpha_{i_\alpha}-1} \prod_{j \neq i_\alpha} z_j^{\alpha_j}.$$
			Note that $\hat{V}$ defines a power series which converges absolutely on an open neighborhood of $\overline{D^{\mathfrak{s}+\beta}(0,1)}$ (and is therefore globally subanalytic). We have for $(t,z) \in W$
			$$J_2(t,z)=\hat{V}(\phi_1(t,z), \ldots ,\phi_s(t,z),\tilde{\phi}_1(t,z), \ldots ,\tilde{\phi}_\beta(t,z))$$
			and therefore
			$$\exp(J_2(t,z))=\exp^*(\hat{V}(\phi_1(t,z), \ldots ,\phi_\mathfrak{s}(t,z),\tilde{\phi}_1(t,z), \ldots ,\tilde{\phi}_\beta(t,z)))$$
			where
			$$\exp^*:\mathbb{C} \to \mathbb{C}, z \mapsto \left\{\begin{array}{ll} \exp(z), & \vert{z}\vert \leq \delta, \\
				0, & \textnormal{else}. \end{array}\right.$$
			By using the exponential series and composition of power series we see that $\exp(J_2)$ has the desired properties since $\vert{\phi_j(t,z)}\vert \leq 1$ for $(t,z) \in W$ and $j \in \{1, \ldots , \mathfrak{s}\}$.
		\end{proof}
		
		Now we show that $J_1|_W$ coincides with a $C$-regular and $C$-consistent function on $\pi(W)$ (i.e. property $(2)^*$). Set
		$$\Xi:\pi(W) \to \mathbb{R}, t \mapsto \sum_{\alpha \in \Gamma_1}\mathfrak{a}(t) c_\alpha \prod_{j=1}^{\mathfrak{s}}\mathfrak{b}_j(t)^{\alpha_j}.$$
		Then we obtain
		$J_1(t,z)=\Xi(t)$ for $(t,z) \in W$. We show that $\Xi$ is $C$-consistent and $C$-regular and are done with the proof of Claim~\ref{claim-13} since we can set $K_1:=\Xi$ and $K_2:=J_2|_W$. Note that
		$$\Xi(t)=\mathfrak{a}(t)\mathfrak{V}_1(\mathfrak{b}_1(t)(\sigma\mathcal{Z})^{\otimes \mathfrak{p}_1}(t,\lambda(t)), \ldots ,\mathfrak{b}_\mathfrak{s}(t)(\sigma\mathcal{Z})^{\otimes \mathfrak{p}_\mathfrak{s}}(t,\lambda(t)))$$
		for $t \in \pi(W)$ where
		$$\lambda:\pi(W) \to \mathbb{R}, t \mapsto \mu_\kappa(t)+(\prod_{j=0}^\kappa \sigma_j)\frac{\hat{S}_\delta(t)+\hat{T}_\delta(t)}{2},$$
		is $C$-regular (since $(t,\lambda(t)) \in W$ for $t \in \pi(W)$: If $\prod_{j=0}^\kappa \sigma_j = 1$ then $\hat{S}_\delta(t)+\mu_\kappa(t) < \lambda(t) < T_{\delta}(t)+\mu_\kappa(t)$ and hence, $\hat{S}_\delta(t) < \vert{\lambda(t) - \mu_\kappa(t)}\vert < T_{\delta}(t)$. If $\prod_{j=0}^\kappa \sigma_j = -1$ then $-\hat{T}_\delta(t)+\mu_\kappa(t) < \lambda(t) < -\hat{S}_{\delta}(t)+\mu_\kappa(t)$ and hence, $\hat{S}_\delta(t) < \vert{\lambda(t) - \mu_\kappa(t)}\vert < T_{\delta}(t)$ since $\lambda(t) < \mu_\kappa(t)$.) So we see that $\Xi$ is $C$-consistent since $\mathfrak{a}$ is $C$-consistent.\\
		
		{\bf Subclaim~2:} {\it The function $\Xi$ is $C$-regular.}
		
		\begin{proof}
			Similarly as in the proof of the claim in Lemma~\ref{5.79} we find a globally subanalytic function $G:\mathbb{C}^{s+1} \times \mathbb{C}^{r-\kappa+1} \to \mathbb{C}$ such that for $t \in \pi(W)$ 
			$$\Xi(t)=G(\eta(t),\mathcal{Z}_\kappa(t,\lambda(t)))$$
			where $\mathcal{Z}_\kappa:=(z_\kappa, \ldots ,z_r)$ and $\eta:\pi(W) \to \mathbb{R}^{s+1}, t \mapsto (\mathfrak{a}(t),\mathfrak{b}_1(t), \ldots ,\mathfrak{b}_s(t))$. Note that $\eta$ and $\Theta$ are $C$-regular (since $\Theta$ is $C$-nice). So with Lemma~\ref{5.55}(2) we see that $\Xi$ is $C$-regular if $\chi_l:\pi(W) \to \mathbb{C}, t \mapsto z_l(t,\lambda(t)),$ is $C$-regular which we show briefly by induction on $l \in \{0, \ldots ,r\}$.
			
			$l=0$: We have $z_0(t,\lambda(t))=\lambda(t)-\Theta_0(t)$ for $t \in \pi(W)$ and hence $\chi_0$ is a $C$-regular function on $\pi(C)$ by Lemma~\ref{5.55}(2) since $\Theta_0$ is $C$-regular.
			
			$l-1 \to l$:  We have $z_l(t,\lambda(t))=\log(\sigma_{l-1}z_{l-1}(t,\lambda(t)))-\Theta_l(t)$ for $t \in \pi(W)$. Hence, with Lemma~\ref{5.55}(2) and the inductive hypothesis $\chi_l$ is $C$-regular since $\Theta_l$ is $C$-regular.
		\end{proof}
		
		Hence, the proof of Claim~\ref{claim-13} is accomplished.
	\end{proof}
	Fix $0<\delta<\min\{\log(\sqrt[4]{R}),\pi\}$. To this $\delta$ fix the corresponding $C$-regular $\kappa$-soft $S_\delta$, the $C$-regular $\kappa$-persistent $T_\delta$ from Claim~\ref{claim-13} and the corresponding $D_{i,1}$ and $D_{i,2}$. Then $1/\sqrt[4]{R}<\vert{\exp(D_{i,2}(t,z))}\vert<\sqrt[4]{R}$ for $(t,z) \in Y_\delta$ and $i \in \{0, \ldots ,s\}$. (+)\\
	
	We obtain for $(t,z) \in \Lambda_g \cap Y_\delta$ 
	$$\Phi_g(t,z)=a(t)\exp(D_{0,1}(t))(\sigma\mathcal{Z})^{\otimes q}(t,z)\exp(D_{0,2}(t,z))$$
	$$V(b_1(t)\exp(D_{1,1}(t))(\sigma\mathcal{Z})^{\otimes p_1}(t,z)\exp(D_{1,2}(t,z)),\ldots ,$$
	$$b_s(t)\exp(D_{s,1}(t))(\sigma\mathcal{Z})^{\otimes p_s}(t,z)\exp(D_{s,2}(t,z))).$$

	Now let $p_0:=q$ with $p_0=(p_{00}, \ldots , p_{0r})$ for $p_{0j} \in \mathbb{Q}$ where $j \in \{0, \ldots , r\}$. By Lemma~\ref{lem:complex-preparation-simplified} we find a $C$-regular $\kappa$-persistent function $T_R:\pi(C) \to \mathbb{R}_{\geq 0}$ such that for every $i \in \{0, \ldots , s\}$
	$$\prod_{j=0}^{\kappa-1}(\sigma_jz_j)^{p_{ij}}|_{\mathcal{B}(\mathcal{P}_{0,\kappa},T_R)} = a_{w_i} V_{w_i}((b_{w_i})_1 \sigma_\kappa z_\kappa, \ldots ,(b_{w_i})_{s_{w_i}} \sigma_\kappa z_\kappa):=a_{w_i}Q_i$$
	where the functions $a_{w_i}$ and $b_{w_i}$ are $C$-regular and $C$-consistent without a zero, and $\vert{(b_{w_i}(t))_j\sigma_\kappa z_\kappa(t,z)}\vert < 1$ for $(t,z) \in \mathcal{B}(\mathcal{P}_{0,\kappa},T_R)$ and $j \in \{1, \ldots , s_{w_i}\}$. Additionally, $V_{w_i}$ is a power series in $s_{w_i}$ variables which converges absolutely on an open neighborhood of $\overline{D^{s_{w_i}}(0,1)}$ with $1/\sqrt[4]{R} \leq V_{w_i}(\overline{D^{s_{w_i}}(0,1)}) \leq \sqrt[4]{R}$. (++)\\
	
	Let $Y:=\mathcal{A}(\mathcal{P}_{0,\kappa},S_\delta,\min\{T_\delta,T_R\})$. (Note that $Y \subset \mathcal{A}(\mathcal{P}_{0,\kappa},N^*,M^*)$.) For $\lambda:=(\lambda_0, \ldots ,\lambda_r) \in \mathbb{Q}^{r+1}$ and $(t,z) \in \Lambda_g \cap Y$ let
	$$(\sigma\mathcal{Z})_o^{\otimes \lambda}(t,z):=\prod_{j=\kappa}^r(\sigma_jz_j(t,z))^{\lambda_j}.$$
	Then for $(t,z) \in \Lambda_g \cap Y$
	$$\Phi_g(t,z)=a(t)a_{\omega_0}(t)\exp(D_{0,1}(t))(\sigma\mathcal{Z})_o^{\otimes q}(t,z)Q_0(t,z)\exp(D_{0,2}(t,z))$$
	$$V(b_1(t)a_{\omega_1}(t)\exp(D_{1,1}(t))(\sigma\mathcal{Z})_o^{\otimes p_1}(t,z)Q_1(t,z)\exp(D_{1,2}(t,z)), \ldots ,$$
	$$b_s(t)a_{\omega_s}(t)\exp(D_{s,1}(t))(\sigma\mathcal{Z})_o^{\otimes p_s}(t,z)Q_s(t,z)\exp(D_{s,2}(t,z))).$$
	Let $b_0:=a$ and for $i \in \{0, \ldots ,s\}$
	$$A_i: \pi(\Lambda_g \cap Y) \to \mathbb{R}, t \mapsto b_i(t)a_{\omega_i}(t)\exp(D_{i,1}(t)).$$
	Note that for $i \in \{1, \ldots ,s\}$ the function $A_i$ is $C$-consistent (by Remark~\ref{5.52}(3) $\exp(D_{i,1})$ is $C$-consistent, by Lemma~\ref{5.53}(1) $b_i$ is also $C$-consistent) and $C$-regular (by Lemma~\ref{5.55}(2) since $b_i$ is $C$-nice and therefore $C$-regular). Note that $A_0$ is $C$-regular (since $a$ is $C$-nice and therefore $C$-regular), and if $g \in \log(E)$ then $A_0$ is $C$-consistent (by Lemma~\ref{5.53}(2) $b_0$ is $C$-consistent if $g \in \log(E)$, by Remark~\ref{5.52}(3) $\exp(D_{0,1})$ is $C$-consistent). 
	
	Note that $\vert{\textnormal{Im}(D_{i,2}(t,z))}\vert<\delta < \pi$ for $(t,z) \in Y$ and $i \in \{0, \ldots ,s\}$. Therefore 
	$$\Lambda_g \cap Y = \{(t,z) \in Y \mid \vert{A_i(t)(\sigma\mathcal{Z})_o^{\otimes p_i}(t,z) Q_i(t,z) \exp(D_{i,2}(t,z))}\vert < R$$ $$\textnormal{ for every }i \in \{1, \ldots ,s\}\}.$$
	
	With the following claim we are almost done with the proof of this theorem.
	
	\begin{claim}
		\label{claim-15}	
		There is a $C$-regular $\kappa$-persistent $M^+: \pi(C) \to \mathbb{R}_{\geq 0}$ and a $C$-regular $\kappa$-soft $N^+: \pi(C) \to \mathbb{R}_{\geq 0}$ such that $\mathcal{A}(\mathcal{P}_{0,\kappa},N^+,M^+) \subset \Lambda_g \cap Y$.
	\end{claim}
	
	\begin{proof}
		For $i \in \{1, \ldots ,s\}$ let 
		$$\phi_i: Y \to \mathbb{C}, (t,z) \mapsto A_i(t)(\sigma\mathcal{Z})_o^{\otimes p_i}(t,z)Q_i(t,z)\exp(D_{i,2}(t,z)).$$
		It is enough to show that there is a $C$-regular $\kappa$-persistent $M_i^+: \pi(C) \to \mathbb{R}_{\geq 0}$ and a $C$-regular $\kappa$-soft $N_i^+: \pi(C) \to \mathbb{R}_{\geq 0}$ such that $\mathcal{A}(\mathcal{P}_{0,\kappa},N_i^+,M_i^+) \subset \mathcal{B}(\phi_i,R) \cap Y$. As in the proof of Claim~\ref{claim-13} we consider two cases.\\ 
		
		{\bf Case 1:} Let $(p_{i\kappa}, \ldots ,p_{ir}) \neq 0$. Consider the $C$-regular function
		$$U:\pi(C) \to \mathbb{R}_{\geq 0}, t \mapsto \left\{\begin{array}{ll} \frac{\sqrt{R}}{\vert{A_i(t)}\vert}, & t \in \pi(Y), \\
			0, & \textnormal{ else.} \end{array}\right.$$
		Since $A_i$ is $C$-consistent we see with Corollary~\ref{cor:persistence-soft-restriction} that $U$ is $\kappa$-persistent. For $(t,z) \in \mathcal{B}((\sigma\mathcal{Z})^{\otimes p_i}_o,U) \cap Y$ we have $\vert{\phi_i(t,z)}\vert<R$ since 
		$$\vert{Q_i(t,z)\exp(D_{i,2}(t,z))}\vert < \sqrt{R}.$$
		Suppose that $p_i$ is $\kappa$-positive. By Lemma~\ref{5.81}(1) if $p_{i \kappa} \neq 0$ respectively Lemma~\ref{5.81}(2) if $\kappa=k \geq 0$ and $p_{i \kappa} = 0$ (since $U_{\text{pers}}$ is $C$-consistent since $A_i$ is also $C$-consistent) there is a $C$-regular $\kappa$-persistent $\hat{M}_i^+:\pi(C) \to \mathbb{R}_{ \geq 0}$ such that $\mathcal{B}(\mathcal{P}_{0,\kappa},\hat{M}_i^+) \subset \mathcal{B}((\sigma\mathcal{Z})^{\otimes p_i},U)$. Choose $N_i^+:=S_\delta$ and 
		$$M_i^+:\pi(C) \to \mathbb{R}_{\geq 0}, t \mapsto \min\{\hat{M}_i^+(t),T_R(t),T_\delta(t)\}.$$
		Then $\mathcal{A}(\mathcal{P}_{0,\kappa},N_i^+,M_i^+) \subset \mathcal{B}((\sigma\mathcal{Z})^{\otimes p_i},U) \cap Y$. Suppose that $p_i$ is $\kappa$-negative. Note that $U_{\textnormal{soft}}$ is $\kappa$-soft (since $\vert{b_i}\vert$ is $\kappa$-soft by Lemma~\ref{5.44}(1)). With Lemma~\ref{5.82} we find a $C$-regular $\kappa$-soft $\hat{N}_i^+:\pi(C) \to \mathbb{R}_{\geq 0}$ and a $C$-regular $\kappa$-persistent $\hat{M}_i^+:\pi(C) \to \mathbb{R}_{\geq 0}$ such that $\mathcal{A}(\mathcal{P}_{0,\kappa},\hat{N}_i^+,\hat{M}_i^+) \subset \mathcal{B}((\sigma\mathcal{Z})^{\otimes p_i},U)$. Choose $M_i^+$ as in the $\kappa$-positive case and $N_i^+:\pi(C) \to \mathbb{R}_{\geq 0}, t \mapsto \max\{S_\delta(t),\hat{N}_i^+(t)\}$.\\
				
		{\bf Case 2:} Let $(p_{i\kappa}, \ldots ,p_{ir}) = 0$. Take $N_i^+:=S_\delta$. 
		If $\vert{A_i(t)}\vert<\sqrt{R}$ then 
		$$\vert{\phi_i(t,z)}\vert = \vert{A_i(t)Q_i(t,z)\exp(D_{i,2}(t,z))\vert}<\sqrt{R} \cdot \sqrt[4]{R} \cdot \sqrt[4]{R}=R$$
		for $(t,z) \in Y$ by (+) and (++) and therefore $(t,z) \in \mathcal{B}(\phi_i(t,z),R)$.  Let
		$$M_i^+:\pi(C) \to \mathbb{R}_{\geq 0}, t \mapsto \left\{\begin{array}{ll} \min\{T_R(t),T_\delta(t)\}, & t \in \pi(Y) \textnormal{ and } \vert{A_i(t)}\vert<\sqrt{R}, \\
			0, & \textnormal{ else.} \end{array}\right.$$
		Note that $M_i^+$ is $C$-regular and by construction that 
		$$\mathcal{A}(\mathcal{P}_{0,\kappa},N_i^+,M_i^+) \subset \mathcal{B}(\phi_i,R) \cap Y.$$
		We show that $M_i^+$ is $\kappa$-persistent and are done with the proof of Claim~\ref{claim-15}. Note that
		$$\vert{A_i(t)Q_i(t,x)\exp(D_{i,2}(t,x))}\vert \leq 1$$
		for every $(t,x) \in Y \cap C$ since the left is the same as $b_i(t)\vert{\mathcal{Y}}\vert^{\otimes p_i}(t,x)\exp(d_{i}(t,x))$ which is bounded by $1$ in its absolute value. Let $w \in \pi_l(C)$. Let $\gamma:\textnormal{}]0,1[\textnormal{} \to C_w$ be a definable curve compatible with $C_w$ with $\lim_{y \searrow 0} (\gamma(y)) \in X_w$ and 
		$$\lim_{y \searrow 0} (\gamma_x(y)-\mu_\kappa(w,\gamma_u(y)))=0.$$
		By passing to a suitable subcurve of $\gamma$ if necessary we may assume that $(w,\gamma_u(y)) \in \pi(Y)$ for $y \in \textnormal{}]0,1[$ (compare with the proof of Lemma~\ref{lem:persistence-soft-restriction}) and that there is a definable continuous $\hat{\gamma}_x: \textnormal{}]0,1[\textnormal{} \to \mathbb{R}$ with $(\gamma_u(y),\hat{\gamma}_x(y)) \in C_w \cap Y_w$ for $y \in \textnormal{}]0,1[$: 
		If $\mu_\kappa<C$ consider for $t \in \pi(C)$ 
		$$I_t:=\text{}]\inf(C_t)+N(t),\min\{\mu_\kappa(t)+M(t),\sup(C_t)\}[,$$
		and if 
		$\mu_\kappa>C$ consider for $t \in \pi(C)$
		$$I_t:=\text{}]\max\{\mu_\kappa(t)-M(t),\inf(C_t)\}, \sup(C_t)-N(t)[.$$ 
		Note that $I_t \subset C_t \cap Y_t$ for every $t \in \pi(C)$, and that 
		$$\lim_{y \searrow 0} \sup(I_{(w,\gamma_u(y))})-\inf(I_{(w,\gamma_u(y))})>0.$$
		We obtain
		$$\vert{A_i(w,\gamma_u(y))Q_i(w,\gamma_u(y),\hat{\gamma}_x(y))\exp(D_{i,2}(w,\gamma_u(y),\hat{\gamma}_x(y)))}\vert \leq 1$$
		and therefore
		\begin{eqnarray*}
			\vert{A_i(w,\gamma_u(y))}\vert  \leq [Q_i(w,\gamma_u(y),\hat{\gamma}_x(y))\exp(D_{i,2}(w,\gamma_u(y),\hat{\gamma}_x(y)))]^{-1} \leq \sqrt{R}
		\end{eqnarray*}
		for $y \in \textnormal{}]0,1[$. So we have 
		$$\lim_{y \searrow 0} M_i^+(w,\gamma_u(y)) = \lim_{y \searrow 0} \min\{T_R(w,\gamma_u(y)),T_\delta(w,\gamma_u(y))\}>0$$
		which gives Claim~\ref{claim-15}.
	\end{proof}
	
	Take $N:=N^+$ and $M:=M^+$ where $N^+$ and $M^+$ are from Claim~\ref{claim-15}. 
	We obtain the following: If $p_i$ is $\kappa$-negative for $i \in \{1, \ldots ,s\}$ then $A_i|_{\pi(\mathcal{A}(\mathcal{P}_{0,\kappa},N,M))}$ can be extended to a $\kappa$-soft function on $\pi(C)$ (since by Lemma~\ref{5.44}(1) $\vert{b_i}\vert$ is $\kappa$-soft) and if $g \in \log(E)$ and $q$ is $\kappa$-negative then $A_0|_{\pi(\mathcal{A}(\mathcal{P}_{0,\kappa},N,M))}$ can be extended to a $\kappa$-soft function on $\pi(C)$ (since by Lemma~\ref{5.44}(2) $\vert{a}\vert$ is $\kappa$-soft for $g \in \log(E)$) and we are done with the whole proof of Theorem~\ref{5.84} by composition of power series.
\end{proof}

The preparation theorem described above also applies for a complex log-analytically prepared function with $C$-nice coefficient, center and base functions (according to Definition~\ref{5.76}), because such a function is the unary parametric global complexification of a log-analytically prepared function with $C$-nice coefficient, center and base functions which is clearly log-exp-analytically prepared with respect to $E:=\{1\}$.


Opris proved in~\cite{28} a similar, albeit a simpler preparation theorem of restricted log-exp-analytically prepared functions in the real setting, specifically on \emph{simple} cells (which means that $C_t$ of the form $]0,d_t[$ for every $t \in \pi(C)$). 

Using this theorem from~\cite{28} differentiability properties of restricted log-exp-analytic functions have been derived in~\cite{28} in the univariate case, which were subsequently generalized to multiple variables. 

Here we give also a consequence of our preparation result stated in Theorem~\ref{5.84} concerning the holomorphic extension of $\Phi_f|_{\mathcal{A}(\mathcal{P}_{0,\kappa},N,M)}$ in $z$ at $\mu_\kappa(t)$ when $N(t)=0$. Recall that $B(c,r)=\{z \in \mathbb{C} \mid \vert{z-c}\vert<r\}$ for $c \in \mathbb{C}$ and $r \in \mathbb{R}_{\geq 0}$.

\begin{lemma}
	\label{5.86}
	There is a $C$-regular $\kappa$-persistent $M:\pi(C) \to \mathbb{R}_{\geq 0}$ and a $C$-regular $\kappa$-soft
	$N:\pi(C) \to \mathbb{R}_{\geq 0}$ such that the following holds where $W:=\mathcal{A}(\mathcal{P}_{0,\kappa},N,M)$ 
	\begin{itemize}
		\item [(1)] $F|_W$ is complex $r$-log-analytically prepared as in Theorem~\ref{5.84}.
		\item [(2)] Let $t \in \pi(C)$ with $N(t)=0$, $M(t)>0$ and suppose that $F(t,-)$ has a holomorphic extension at $\mu_\kappa(t)$ (i.e. there is an open neighborhood $U$ of $\mu_\kappa(t)$ in $B(\mu_\kappa(t),M(t))$ and a holomorphic $J:U \to \mathbb{C}$ such that $J|_{U \cap W_t} = F(t,-)$). Then there is a holomorphic $h: B(\mu_\kappa(t),M(t)) \to \mathbb{C}$ such that $h|_{W_t}=F(t,-)$.
	\end{itemize}
\end{lemma}

\begin{proof}
	For $\nu \in \mathbb{Q}^{r+1}$ let $\nu:=(\nu_0, \ldots , \nu_r)$. With Theorem~\ref{5.84} choose a $\kappa$-persistent $M:\pi(C) \to \mathbb{R}_{\geq 0}$ and a $\kappa$-soft $N:\pi(C) \to \mathbb{R}_{\geq 0}$ such that $F|_{\mathcal{A}(\mathcal{P}_{0,\kappa}N,M)}$ is complex $r$-log-analytically prepared in $z$ with complex LA-preparing tuple
	$$(r,\mathcal{Z},a,q,s,V,b,P)$$
	such that $q_0= \ldots =q_{\kappa-1}=0$ and 
	$$P :=\left(\begin{array}{cccc}
		p_{10}&\cdot&\cdot&p_{1r}\\
		\cdot&& &\cdot\\
		\cdot&& &\cdot\\
		p_{s0}&\cdot&\cdot&p_{sr}\\
	\end{array}\right)\in M_\kappa(s \times (r+1), \mathbb{Q}).$$
	For $i \in \{1, \ldots , s\}$ let $p_i:=(p_{i0}, \ldots , p_{ir})$.
	We can further assume for every $j \in \{0, \ldots ,\kappa-1\}$
	$$\Bigl|\frac{\mathcal{P}_{j,\kappa}(t,z)}{e^{\mu_{\kappa-j-1,\kappa}(t)}} \Bigl| < 1/2$$
	for $(t,z) \in \Omega:=\mathcal{A}(\mathcal{P}_{0,\kappa},N,M)$ (compare with Lemma~\ref{5.67}(1)) $(+)$\\ 
	
	For $(t,z)\in \Omega$ we have
	$$F(t,z)=a(t)(\sigma\mathcal{Z})^{\otimes q}(t,z)V(b_1(t)(\sigma\mathcal{Z})^{\otimes p_1}(t,z), \ldots ,b_s(t) (\sigma\mathcal{Z})^{\otimes p_s}(t,z)).$$
	Since for $(t,z) \in \Omega$
	$$z_\kappa(t,z)=\mathcal{P}_{\kappa,\kappa}(t,z) = \log\Bigl(1+\frac{\mathcal{P}_{\kappa-1,\kappa}(t,z)}{e^{\mu_{0,\kappa}(t)}}\Bigl),$$
	$$\mathcal{P}_{j+1,\kappa}(t,z) = \log\Bigl(1+\frac{\mathcal{P}_{j,\kappa}(t,z)}{e^{\mu_{\kappa-j-1,\kappa}(t)}} \Bigl)$$
	for $j \in \{0, \ldots , \kappa-2\}$, and $\mathcal{P}_{0,\kappa}(t,z) = z-\mu_\kappa(t)$ we obtain with $(+)$ by successively using the logarithmic series
	$$z_\kappa(t,z)=L(z-\mu_\kappa(t)) \textnormal{ } (*)$$
	for $(t,z) \in \Omega$ where $L$ is a power series around zero which converges absolutely on $B(0,M(t))$ (compare with Definition~\ref{5.16} where we consider the functions $\mathcal{P}_{j,\kappa}$ as functions on $H$). 
	We set $Z:=(Z_1, \ldots ,Z_s)$. Let $\sum_{\alpha \in \mathbb{N}_0^s}c_\alpha Z^{\alpha}$ be the Taylor series expansion of $V$. 
	
	We may assume that $a$ does not have any zero. Otherwise $a=0$ and the statement follows. Fix $t \in \pi(C)$ with $N(t)=0$. Then $(p_{j\kappa}, \ldots ,p_{jr}) = 0$ or $p_j$ is $\kappa$-positive for $j \in \{1, \ldots ,s\}$: If $p_j$ is $\kappa$-negative for a $j \in \{1, \ldots ,s\}$ we have that $\lim_{z \to \mu_\kappa(t)} \vert{b(t)(\sigma\mathcal{Z})^{\otimes q}(t,z)}\vert = \infty$ by Lemma~\ref{5.34} since $b(t) \neq 0$ which is not possible since $\vert{b(t) (\sigma\mathcal{Z})^{\otimes q}(t,z)}\vert \leq 1$ for $z \in \Omega_t$ (since  $F|_\Omega$ is complex $r$-log-analytically prepared, see Definition~\ref{def:complex-log-prepared-general}) and $\mu_\kappa(t) \in \overline{\Omega_t}$ (due to $N(t)=0$). Let
	$$\Gamma_1:=\{\alpha \in \mathbb{N}_0^s \textnormal{ } \big\vert\;^tP\alpha+q \in \{0\}^\kappa \times \mathbb{N}_0 \times \{0\}^{r-\kappa}\}$$
	and $\Gamma_2:=\mathbb{N}_0^s \setminus \Gamma_1$. Set $V_1:=\sum_{\alpha \in \Gamma_1}c_\alpha Z^\alpha$ and $V_2:=\sum_{\alpha \in \Gamma_2} c_\alpha Z^{\alpha}$. Since $P \in M_\kappa(s \times (r+1),\mathbb{Q})$ (i.e. the first $\kappa$ columns of $P$ are zero) we have that $V=V_1+V_2$. For $l \in \{1,2\}$ let
	$$G_l:\Omega_t \to \mathbb{C}, z \mapsto a(t)(\sigma\mathcal{Z})^{\otimes q}(t,z)V_l(b_1(t)(\sigma\mathcal{Z})^{\otimes p_1}(t,z), \ldots ,b_s(t)(\sigma\mathcal{Z})^{\otimes p_s}(t,z)).$$
	So we see by composition of power series and $(*)$ that 
	$$G_1(z)=\sum_{j=0}^{\infty} d_j(z-\mu_\kappa(t))^j$$
	for suitable $d_j \in \mathbb{C}$ and $z \in D_t$. Note that the series to the right converges absolutely on $B(\mu_\kappa(t),M(t))$. So $G_1$ has a holomorphic extension $\hat{G}_1$ on $B(\mu_\kappa(t),M(t))$. Let $0<\varepsilon<M(t)$ be such that $F|_{B(\mu_\kappa(t),\varepsilon) \cap \Omega_t}$ has a holomorphic extension $\hat{F}$ on $Y:=B(\mu_\kappa(t),\varepsilon)$. Set 
	$$\hat{G}_2:Y \to \mathbb{C}, z \mapsto \hat{F}(z)-\hat{G}_1(z).$$ 
	Note that $\hat{G}_2$ is holomorphic and that $\hat{G}_2|_{\Omega_t \cap Y}=G_2$. We show that $\hat{G}_2=0$ and we are done with the complete proof by the identity theorem (since then $F_t$ coincides with $G_1$ on $Y \cap \Omega_t$). We define for $j \in \{0, \ldots , r\}$
	$$(^tP\alpha+q)_j:=q_j + \sum_{i=1}^s p_{ij}\alpha_i$$
	and let
	$$K:=\{(0, \ldots, 0,\omega_\kappa, \ldots ,\omega_r) \in \{0\}^\kappa \times \mathbb{Q}^{r-\kappa+1} \mid \textnormal{ there is }\alpha \in \Gamma_2 \textnormal{ such that }$$ 
	$$(^t P \alpha + q)_j=\omega_j \textnormal{ for } j \in \{\kappa, \ldots ,r\}\}.$$
	Then $K \subset \{0\}^\kappa \times (\mathbb{Q}^{r-\kappa+1} \setminus (\mathbb{N}_0 \times \{0\}^{r-\kappa}))$. For $\omega \in K$ let
	$$\Gamma_{2,\omega}:=\{\alpha \in \mathbb{N}_0^s \mid (^tP\alpha + q)_j = \omega_j \textnormal{ for } j \in \{\kappa, \ldots ,r\}\}.$$
	We have for $z \in \Omega_t$
	$$\hat{G}_2(z)=\sum_{\omega \in K}e_\omega \prod_{j=\kappa}^r (\sigma_jz_j(t,z))^{\omega_j}$$
	where
	$$e_\omega=a(t)\sum_{\alpha \in \Gamma_{2,\omega}}c_\alpha \prod_{j=1}^s b_j(t)^{\alpha_j}.$$
	Let 
	$$K^*:=\{\omega \in K \mid e_{\omega} \neq 0 \}.$$
	If $K^* = \emptyset$ then obviously $\hat{G}_2=0$. So assume $K^* \neq \emptyset$.
	
	\begin{claim}
	\label{claim:asymptotic-behavior}
	There is $\xi:=(0, \ldots, 0,\xi_\kappa, \ldots ,\xi_r) \in K^*$ such that 
	$$\lim_{z \to \mu_\kappa(t)}\prod_{j=\kappa}^r \vert{z_j(t,z)}\vert^{\omega_j-\xi_j} = 0$$
	for all $\omega \in K^*$ with $\omega \neq \xi$.
	\end{claim}

	\begin{proof}
	Suppose that for every $i \in \{1, \ldots ,s\}$ and $j \in \{0, \ldots , r\}$ the enumerator of $p_{ij}$ is coprime to its denominator and that all the denominators occuring in $P$ are positive (i.e. enumerator and denominator of $p_{ij}$ are uniquely determined). 
	Then there is a maximum denominator $\nu$ of all entries from $P$. Let 
	\[
	\Xi_\kappa:=\{\lambda \in \mathbb{Q} \mid \textnormal{ there is }\alpha \in \Gamma_2 \textnormal{ such that } (^t P \alpha + q)_{\kappa}=\lambda\}.
	\]
	Then $\Xi_\kappa \subset 1/\nu \cdot \mathbb{N}_0 +q:=\{a/\nu + q \mid a \in \mathbb{N}_0\}$ since $p_j=0$ or $p_j$ is $\kappa$-positive for every $j \in \{1, \ldots , s\}$ (and hence, $p_{j\kappa} \geq 0$). Since $1/\nu \cdot \mathbb{N}_0+q$ is well-ordered there is a smallest element $\xi_\kappa \in \Xi_\kappa$. 
	Suppose for an induction that $\Xi_\kappa, \ldots , \Xi_{j-1}$ and that $\xi_{j-1}$ have already been defined for $j \in \{\kappa+1, \ldots , r\}$. We set 
	\[
	\Xi_j:=\{\lambda \in \mathbb{Q} \mid \textnormal{ there is }\alpha \in \Gamma_2 \textnormal{ such that } (^t P \alpha + q)_l = \xi_l \text{ for $l \in \{\kappa, \ldots, j-1\}$}\\
	\]
	\[
	\text{ and }(^t P \alpha + q)_j=\lambda\}.
	\]
	Then $\Xi_j \subset 1/\nu \cdot \mathbb{Z} + q := \{a/\nu + q \mid a \in \mathbb{Z}\}$. Let 
	\[
	S:=\{i \in \{1, \ldots, s\} \mid (p_{i \kappa}, \ldots , p_{i(j-1)}) \neq 0\} = \{i_1, \ldots, i_\eta\}
	\]
	for $\eta \leq s$ and pairwise distinct $i_1, \ldots  ,i_\eta \in \{1, \ldots, s\}$.  
	Then $p_{ij}>0$ for $i \in \{1, \ldots ,s\}$ is only possible if $i \in S$ (since $p_i=0$ or $p_i$ is $\kappa$-positive). We show that $\Xi_j$ is bounded from above which implies that $\Xi_j$ has a largest element $\xi_j$ (again by the well-ordering of $1/\nu \cdot \mathbb{N}_0 + a$ for $a \in \mathbb{R}$). Suppose not. Then
	\begin{align*} 
	 \Gamma^*:&=\{(^t P \alpha)_j \mid \text{there is } \alpha \in \Gamma_2 \text{ with }(^t P \alpha+q)_l = \xi_l \text{ for } l \in \{\kappa, \ldots, j-1\}\} \\
	&= \{\sum_{i=1}^s p_{ij}\alpha_i \mid \text{there is } \alpha \in \Gamma_2 \text{ with }(^t P \alpha+q)_l = \xi_l \text{ for } l \in \{\kappa, \ldots, j-1\}\}
	\end{align*}
	is also not bounded from above and hence,
	\[
	A:=\{(\alpha_{i_1}, \ldots, \alpha_{i_\eta}) \in \mathbb{N}_0^\eta \mid \alpha \in \Gamma_2, (^t P \alpha+q)_l = \xi_l \text{ for } l \in \{\kappa, \ldots, j-1\}\}
	\]
	is not finite (
	since $p_{ij}>0$ for $i \in \{1, \ldots ,s\}$ implies $i \in S=\{i_1, \ldots , i_\eta\}$). We show that the latter is false which implies the existence of a largest element $\xi_j \in \Xi_j$. Let 
	$S_\kappa:=\{i \in S \mid p_{i\kappa} > 0\}$ and for $\kappa < l \leq j-1$ we set $S_l:=\{i \in S \mid p_{i\kappa}= \ldots = p_{i(j-1)} = 0 \text{ and } p_{ij}<0\}$. Then $S$ is the disjoint union of $S_\kappa, \ldots , S_{j-1}$ (since $p_i$ is $\kappa$-positive for $i \in \{1, \ldots , s\}$) and for $i^* \in S_\kappa$ we see that $\{\alpha_{i^*} \mid (\alpha_{i_1}, \ldots , \alpha_{i_\eta}) \in A\}$ is finite (since $(^t P \alpha+q)_{\kappa} = q_\kappa+ \sum_{i=1}^s p_{i\kappa} \alpha_i = \xi_\kappa$ for $\alpha \in \Gamma_2$ with $(\alpha_{i_1}, \ldots , \alpha_{i_{\eta}}) \in A$; note that $p_{i\kappa} \geq 0$ for $i \in \{1, \ldots ,s\}$). Inductively for $l \in \{\kappa+1, \ldots , j-1\}$ we also argue that $\{\alpha_{i^*} \mid (\alpha_{i_1}, \ldots , \alpha_{i_\eta}) \in A\}$ is finite where $i^* \in S_l$: \\
	Suppose not. Since $(^tP\alpha + q)_l = q_l+\sum_{i=1}^s p_{il}\alpha_i=\xi_l$ for $\alpha \in \Gamma_2$ with $(\alpha_{i_1}, \ldots , \alpha_{i_\eta}) \in A$ there is $\hat{i} \in S$ with $p_{\hat{i}l}>0$ such that 
	$$\{\alpha_{\hat{i}} \mid \alpha \in \Gamma_2, (^t P \alpha+q)_l = \xi_l \text{ for } l \in \{\kappa, \ldots, j-1\}\}$$
	is not bounded. But $\hat{i} \in S_\kappa \cup \ldots \cup S_{l-1}$ since $p_{\hat{i}}$ is $\kappa$-positive (and $l>\kappa$) and therefore 
	$$\{\alpha_{\hat{i}} \mid \alpha \in \Gamma_2, (^t P \alpha+q)_l = \xi_l \text{ for } l \in \{\kappa, \ldots, j-1\}\} = \{\alpha_{\hat{i}} \mid (\alpha_{i_1}, \ldots , \alpha_{i_\eta}) \in A\}$$
	is finite by the inductive hypothesis, a contradiction. Hence, we see that the set $\{\alpha_{i^*} \mid (\alpha_{i_1}, \ldots , \alpha_{i_\eta}) \in A\}$ is finite and by induction that $A$ is finite in each component. This implies that $A$ is finite and hence, $\Xi_j$ is bounded from above. \\
	Now let $\xi:=(0, \ldots , 0,\xi_\kappa, \ldots , \xi_r) \in \{0\}^\kappa \times \mathbb{Q}^{r-\kappa+1}$. It fulfills the desired properties by Lemma~\ref{5.34} (since for every $\omega \in K$ with $\omega \neq \xi$ we see by construction of $\xi$ that $\omega-\xi$ is $\kappa$-positive). 
	\end{proof}

	\begin{claim}
		\label{claim:asymptotic-behavior-derivative}
		If $\xi:=(0, \ldots, 0,\xi_\kappa, \ldots ,\xi_r),\omega:=(0, \ldots , 0,\omega_\kappa, \ldots , \omega_r) \in K$ are such that $$\lim_{z \to \mu_\kappa(t)} \prod_{j=\kappa}^r \vert{z_j(t,z)}\vert^{\omega_j-\xi_j}=0$$
		then 
		\[
		\lim_{z \to \mu_\kappa(t)} \frac{\frac{d}{dz} \prod_{j=\kappa}^r \vert{z_j(t,z)}\vert^{\omega_j}}{\frac{d}{dz} \prod_{j=\kappa}^r \vert{z_j(t,z)}\vert^{\xi_j}} = 0.
		\]
	\end{claim}
	
	\begin{proof}
	Note that $\omega-\xi$ is $\kappa$-positive (otherwise $\lim_{z \to \mu_\kappa(t)} \prod_{j=\kappa}^r \vert{z_j(t,z)}\vert^{\omega_j-\xi_j} \neq 0$ by Lemma~\ref{5.34}(1)). 
	For $\lambda:=(\lambda_0, \ldots , \lambda_r) \in \mathbb{Q}^{r+1}$ with $(\lambda_\kappa, \ldots ,\lambda_r) \neq 0$ let $j_\kappa(\lambda):=\min\{j \in \{\kappa, \ldots, r\} \mid \lambda_j \neq 0\}$ and 
	$$\lambda_{\kappa, \text{diff}} := (0, \ldots, 0,\lambda_\kappa-1, \ldots ,\lambda_{j_\kappa(\lambda)}-1, \lambda_{j_\kappa(\lambda)+1}, \ldots , \lambda_r) \in \mathbb{Q}^{r+1}$$
	(compare with Corollary~\ref{5.35}). 
	We show that $\omega_{\kappa,\text{diff}}-\xi_{\kappa,\text{diff}}$ is $\kappa$-positive and are done with the proof of this claim by Lemma~\ref{5.34}(1) and Corollary~\ref{5.35}. If $j_\kappa(\omega) = j_\kappa(\xi)$ we have $\omega_{\kappa,\text{diff}}-\xi_{\kappa,\text{diff}} = \omega-\xi$ which is $\kappa$-positive. If $j_\kappa(\xi) > j_\kappa(\omega)$ then $\omega_{j_\kappa(\omega)}>0$ if $j_\kappa(\omega)=\kappa$ and $\omega_{j_\kappa(\omega)}<0$ if $j_\kappa(\omega)>\kappa$ since $\omega-\xi$ is $\kappa$-positive. Since
	$(\omega_{\kappa,\text{diff}}-\xi_{\kappa,\text{diff}})_j = 0$ for $j<j_\kappa(\omega)$ and $(\omega_{\kappa,\text{diff}}-\xi_{\kappa,\text{diff}})_{j_\kappa(\omega)} = \omega_{j_\kappa(\omega)}$ one sees that $\omega_{\kappa,\text{diff}}-\xi_{\kappa,\text{diff}}$ is $\kappa$-positive. Analoguously, if $j_\kappa(\omega) > j_\kappa(\xi)$ then $\xi_{j_\kappa(\xi)}<0$ if $j_\kappa(\xi)=\kappa$ and $\xi_{j_\kappa(\xi)}>0$ if $j_\kappa(\xi)>\kappa$ since $\omega-\xi$ is $\kappa$-positive. Since
	$(\omega_{\kappa,\text{diff}}-\xi_{\kappa,\text{diff}})_j = 0$ for $j<j_\kappa(\xi)$ and $(\omega_{\kappa,\text{diff}}-\xi_{\kappa,\text{diff}})_{j_\kappa(\xi)} = -\xi_{j_\kappa(\xi)}$ one sees that $\omega_{\kappa,\text{diff}}-\xi_{\kappa,\text{diff}}$ is $\kappa$-positive.
	\end{proof}

	Hence, with Claim~\ref{claim:asymptotic-behavior} and Claim~\ref{claim:asymptotic-behavior-derivative} we are done with the whole proof of Lemma~\ref{5.86} by investigating the following two cases.
	
	{\bf Case 1:} If $\xi_\kappa \in \mathbb{N}_0$ (note that $\xi \neq 0$) we differentiate $\hat{G}_2$ $m$-times with respect to $z$ where $m:=\xi_\kappa+1$ and we see with Corollary~\ref{5.35} and Claim~\ref{claim:asymptotic-behavior-derivative} (note that $\xi \neq 0$) that there is $\beta=(-1,\beta_{\kappa+1}, \ldots ,\beta_r) \in \mathbb{Q}^{r-\kappa+1}$ such that    
	$$\lim_{z \to \mu_\kappa(t)} \Bigl|\frac{d^m \hat{G}_2/dz^m(z)}{\prod_{j=\kappa}^r z_j(t,z)^{\beta_j}} \Bigl| \in \mathbb{R}^*.$$
	But we obtain
	$$\lim_{z \to \mu_\kappa(t)} \frac{d^m \hat{G}_2}{dz^m}(z) =  \frac{d^m \hat{G}_2}{dz^m}(\mu_\kappa(t)) \in \mathbb{C},$$
	a contradiction.
	
	{\bf Case 2:} If $\xi_\kappa \notin \mathbb{N}_0$ (note that $\xi \neq 0$) we differentiate $\hat{G}_2$ $m$-times with respect to $z$ where $m:=\lceil{\xi_\kappa}\rceil$ and we see with Corollary~\ref{5.35} and Claim~\ref{claim:asymptotic-behavior-derivative} that there is $\gamma=(\gamma_\kappa, \ldots ,\gamma_r) \in \mathbb{Q}^{r-k+1}$ with $\gamma_\kappa<0$ such that 
	$$\lim_{z \to \mu_\kappa(t)} \Bigl|\frac{d^m \hat{G}_2/dz^m(z)}{\prod_{j=\kappa}^r z_j(t,z)^{\gamma_j}} \Bigl| \in \mathbb{R}^*.$$
	We get the same contradiction as in Case 1.
\end{proof}

\section{Definability Results of Parameterized Integrals}

In this section we construct a $\kappa$-persistent $M:\pi(C) \to \mathbb{R}_{\geq 0}$ and a $\kappa$-soft $N:\pi(C) \to \mathbb{R}_{\geq 0}$ such that the parameterized integral
\[
\Psi: \Delta \to \mathbb{C}, (t,s,z) \mapsto \frac{1}{2 \pi i}\int_{\partial B(\mu_\kappa,s)} \frac{\Phi_f(t,\xi)}{\xi-z} d \xi,
\]
is definable where
\[
\Delta:=\{(t,s,z) \in \pi(C) \times \mathbb{R}_{>0} \times \mathbb{C} \mid N(t) < s < M(t) \text{ and } z \in B(\mu_\kappa(t),s)\}
\] 
and $\Phi_f:\Lambda \to \mathbb{C}$ is the unary parametric global complexification of a restricted log-exp-analytically prepared function $f:C \to \mathbb{R}$ from Definition~\ref{5.76}.
Note that for fixed $t \in \pi(C)$ and $N(t)<s<M(t)$ the integral $\Psi$ coincides with $\Phi_f$ by Cauchy's integral formula if $(\Phi_f)_t$ can be extended to a holomorphic function on $B(\mu_\kappa(t),\tilde{s})$ for a $0<s<\tilde{s}<M(t)$. 

This section is structured as follows. At first we give a general method how definability results of parameterized integrals can be obtained with Fact~\ref{19} and apply this method on $\Psi$. In this context  \emph{free-regular} functions play a crucial role which we introduce in the second subsection. 
In the third subsection we look closer at logarithmic scales and in the last one we give the desired result on integration for $\Phi_f$.

For Section~8 we introduce the following notation. Let $n \in \mathbb{N}_0$. Let $t$ range over $\mathbb{R}^n$, $z$ over $\mathbb{C}$ and $\tau$ over $\mathbb{R}$. Let $\pi:\mathbb{R}^n \times \mathbb{C} \to \mathbb{R}^n, (t,z) \mapsto t,$ be the projection on the first $n$ real coordinates.
Let $s$ range over $\mathbb{R}$. Let $\pi^+:\mathbb{R}^n \times \mathbb{R} \times \mathbb{C} \to \mathbb{R}^n \times \mathbb{R}, (t,s,z) \mapsto (t,s),$ be the projection on the first $n+1$ real coordinates\index{$\pi^+$}. 
For a definable function $g:Y \to \mathbb{C}, (t,s,\tau) \mapsto g(t,s,\tau),$ where $Y \subset \mathbb{R}^n \times \mathbb{R}$ is definable we say that $g$ is bounded in $\tau$\index{bounded in $\tau$} if $g_{(t,s)}$ is bounded on $Y_{(t,s)}=\{\tau \in \mathbb{R} \mid (t,s,\tau) \in Y\}$ for every $(t,s) \in \pi^+(Y)$. 

For the whole section fix a non-empty definable cell $C \subset \mathbb{R}^n \times \mathbb{R}_{\neq 0}$ and a logarithmic scale $\mathcal{Y}:=(y_0, \ldots , y_r)$ with $C$-nice center $\Theta$. Fix also the change index $k^{\text{ch}}$ of $\mathcal{Y}$ and let $\kappa=r$ or $\kappa \in \{k^{\text{ch}},r\}$ if $k^{\text{ch}} \geq 0$.

\subsection{A General Method}

Fact~\ref{19} above states that the parameterized integral of a globally subanalytic function is constructible and hence, definable. In this subsection we show how this fact can be used to obtain definability results of parameterized integrals. These integrals involve integrands which are not necessarily constructible, but depend on the variables of integration in a globally subanalytic (or constructible) manner. For example the parameterized integral 
\[
\int_{-\pi}^\pi \frac{e^{1/t^2}}{e^{2/t^2}+x^2}dx 
\]
belongs to this class and is definable. Indeed, we obtain 
\[
\int_{-\pi}^\pi \frac{e^{1/t^2}}{e^{2/t^2}+x^2}dx = \Bigl[\arctan(e^{-1/t^2}x)\Bigl]_{x=-\pi}^{x=\pi} = 2\arctan (e^{-1/t^2}\pi).
\]
We formulate all the results for integrals with one variable of integration because this is sufficient for our purposes. \\

For this subsection we fix a field $\mathbb{K} \in \{\mathbb{R},\mathbb{C}\}$ and let $Z \subset \mathbb{R}^n \times \mathbb{K}$ be definable.

\begin{definition}
	\label{24}
	Let $p \in \mathbb{N}_0$ and $g:\pi(Z) \to \mathbb{R}^p$ be definable. We call a function $f:Z \to \mathbb{C}$ \textbf{globally subanalytic in $z$ (with support function $g$)}\index{globally subanalytic!with support function} respectively \textbf{constructible in $z$ (with support function $g$)}\index{constructible!with support function} if there is a globally subanalytic respectively constructible function $F:\mathbb{R}^p \times \mathbb{K} \to \mathbb{C}$ such that 
	$$f(t,z)=F(g(t),z)$$
	for every $(t,z) \in Z$.
\end{definition}

\begin{lemma}
	\label{25}
	Let $p \in \mathbb{N}_0$ and $g:\pi(Z) \to \mathbb{R}^p$ be definable. Let $q \in \mathbb{N}$ and let $G:\mathbb{C}^q \to \mathbb{C}$ be a globally subanalytic function. Let $J_1, \ldots ,J_q:Z \to \mathbb{C}$ be globally subanalytic in $z$ with support function $g$. Then 
	$$H:Z \to \mathbb{C}, (t,z) \mapsto G(J_1(t,z), \ldots ,J_q(t,z)),$$
	is globally subanalytic in $z$ with support function $g$.
\end{lemma}

\begin{proof}
	(1): Clear.
	(2): Let $K_1, \ldots ,K_q:\mathbb{R}^p \times \mathbb{K} \to \mathbb{C}$ be globally subanalytic such that $J_j(t,z)=K_j(g(t),z)$ for every $(t,z) \in Z$ and $j \in \{1, \ldots ,q\}$. Let $w$ range over $\mathbb{R}^p$. Then 
	$$F:\mathbb{R}^p \times \mathbb{K} \to \mathbb{C}, (w,z) \mapsto G(K_1(w,z), \ldots ,K_q(w,z)),$$
	is globally subanalytic. We have
	$$F(g(t),z)=G(J_1(t,z), \ldots ,J_q(t,z))$$
	for every $(t,z) \in Z$.
\end{proof}

The following lemma presents the general method we are using throughout this section. It states that the parameterized integral of a globally subanalytic function in $\tau$ with support function $g$ remains definable. In particular, it is the composition of a constructible function and $g$.

\begin{lemma}
	\label{26}
	Let $\mathbb{K}=\mathbb{R}$ and let $Y \subset \mathbb{R}^n$ be definable. Let $p \in \mathbb{N}_0$ and let $a,b \in \mathbb{R}$ be with $a \leq b$. Let $Z=Y \times [a,b]$. Let $f:Z \to \mathbb{C}, (t,\tau) \mapsto f(t,\tau),$ be globally subanalytic in $\tau$ with support function $g:Y \to \mathbb{R}^p$ and bounded in $\tau$. Consider
	$$F:Y \to \mathbb{C}, t \mapsto \int_a^b f(t,\tau) d\tau.$$
	Then $F$ is well-defined and there is a constructible $H:\mathbb{R}^p \to \mathbb{C}$ such that $F(t)=H(g(t))$ for every $t \in Y$.
\end{lemma}

\begin{proof}
	Let $w$ range over $\mathbb{R}^p$. Let $J:\mathbb{R}^p \times \mathbb{R} \to \mathbb{C}, (w,\tau) \mapsto J(w,\tau),$ be globally subanalytic such that 
	$$f(t,\tau)=J(g(t),\tau)$$
	for every $(t,\tau) \in Z$. Let $J=J_1+iJ_2$ where $J_1$ is the real part and $J_2$ the imaginary part of $J$. For $j \in \{1,2\}$ we consider 
	$$G_j:\mathbb{R}^p \times \mathbb{R} \to \mathbb{R}, (w,\tau) \mapsto \left\{\begin{array}{lll} J_j(w,\tau),&& \tau \in [a,b], \\
		&\textnormal{if}&\\
		0,&& \textnormal{else.} \end{array}\right.$$
	Then $G_j$ is globally subanalytic. Let $G:=G_1+iG_2$. By Fact~\ref{19}(1) the set
	$$\textnormal{Fin}(G) := \Bigl\{w \in \mathbb{R}^p \mid \int_\mathbb{R} \vert{G_j(w,\tau)}\vert d\tau < \infty \textnormal{ for $j \in \{1,2\}$} \Bigl\}$$
	is globally subanalytic. We obtain $g(Y) \subset \textnormal{Fin}(G)$, because $f_t$ is bounded for every $t \in Y$. With Fact~\ref{19}(2) we obtain that the function $H:\mathbb{R}^p \to \mathbb{C}$ defined by 
	$$H(w) = \int_\mathbb{R}  G(w,\tau)d\tau = \int_a^b J(w,\tau)d\tau$$
	if $w \in \textnormal{Fin}(G)$ and $H(w)=0$ otherwise is constructible which implies $F(t)=H(g(t))$.
\end{proof}

Now we apply Lemma~\ref{26} on the parameterized integral $\Psi$. First we parameterize the circle $\partial B(\mu_\kappa(t),s)$ suitably and then describe $\mathcal{A}(\mathcal{P}_{0,\kappa},N,M)$ using these parameterizations.

\begin{convention}
	\label{5.58}
	For $(t,s) \in \pi(C) \times \mathbb{R}_{>0}$ we parameterize the circle \\ 
	$\partial B(\mu_{\kappa}(t),s) \subset \mathbb{C}$ with the counterclockwise oriented curve
	$$\rho:[-\pi,\pi] \to \mathbb{R}, \tau \mapsto \mu_\kappa(t)+(\prod_{j=0}^{\kappa}\sigma_j)s e^{i \tau}.$$
\end{convention}

\begin{definition}
	\label{5.61}
	Let $D \subset \pi(C) \times \mathbb{C}$ be definable. We set\index{$D(\mu_\kappa)$}
	$$D(\mu_\kappa):=\Bigl\{(t,s,\tau) \in \pi(C) \times \mathbb{R}_{>0} \times [-\pi,\pi] \mid \mu_{\kappa}(t)+(\prod_{j=0}^{\kappa}\sigma_j)se^{i \tau} \in D_t \Bigl\}.$$
\end{definition}


\begin{lemma}
	\label{5.62}
	Let $M,N:\pi(C) \to \mathbb{R}_{\geq 0}$ be definable functions. Suppose that $M<\vert{\mu_r-\mu_k}\vert$ if $k:=k^{\text{ch}} \geq 0$. Let
	$$D:=\mathcal{A}(\mathcal{P}_{0,\kappa},N,M) = \{(t,z) \in H \mid N(t)<\vert{z-\mu_\kappa(t)}\vert<M(t)\}.$$
	Then
	$$D(\mu_\kappa) =\{(t,s) \in \pi(C) \times \mathbb{R}_{>0} \mid N(t)<s<M(t)\} \textnormal{}\times \textnormal{} ]-\pi,\pi[.$$
\end{lemma}

\begin{proof}
	Let
	$$B:=\{(t,s) \in \pi(C) \times \mathbb{R}_{>0} \mid N(t)<s<M(t)\} \textnormal{}\times \textnormal{} ]-\pi,\pi[.$$
	Suppose $\kappa=r$. Let $\prod_{j=0}^r \sigma_j=1$. The case ''$\prod_{j=0}^r \sigma_j=-1$'' is handled completely similarly. For  $t \in \pi(C)$, $s \in \mathbb{R}_{>0}$ and $\tau \in \textnormal{} [-\pi,\pi]$ we set 
	$$\delta(t,s,\tau):=\mu_r(t)+ se^{i \tau}.$$
	Note for $(t,s,\tau) \in B$ that we have 
	$\delta(t,s,\tau) \in \mathbb{C}^-_{\mu_r(t)}$ if $k < 0$ and $\delta(t,s,\tau) \in \mathbb{C}^-_{\mu_r(t)} \cap \mathbb{C}^+_{\mu_k(t)}$ if $k \geq 0$ respectively. Lemma~\ref{5.13} gives $B \subset D(\mu_\kappa)$. Let $(t,s,\tau) \in D(\mu_\kappa)$. Then $t \in \pi(C)$, $N(t)<s<M(t)$ and since $\delta(t,s,\tau) \in D_t$ we have $\tau \in \textnormal ]-\pi,\pi[$. (If $\tau \in \{-\pi,\pi\}$ we obtain $\delta(t,s,\tau)=\mu_r(t)-s<\mu_r(t)$, i.e. $\delta(t,s,\tau) \notin \mathbb{C}^-_{\mu_r}$ and therefore $\delta(t,s,\tau) \notin D_t$ by Lemma~\ref{5.13}.) We obtain $(t,s,\tau) \in B$. The case "$\kappa=k$" is handled completely similarly. 
\end{proof}

\begin{lemma}
	\label{5.63}
	Let $M,N:\pi(C) \to \mathbb{R}_{\geq 0}$ be definable functions. Suppose that $M<\vert{\mu_r-\mu_k}\vert$ if $k:=k^{\text{ch}} \geq 0$. Let $D:=\mathcal{A}(\mathcal{P}_{0,\kappa},N,M)$. Let $F:D \to \mathbb{C}$ be definable such that 
	$$F^*:D(\mu_\kappa) \to \mathbb{C}, (t,s,\tau) \mapsto F(t,\mu_\kappa(t)+\Bigl(\prod_{j=0}^\kappa \sigma_j \Bigl) se^{i\tau}),$$
	is bounded in $\tau$. Suppose that $F^*$ is globally subanalytic in $\tau$ with support function $g:\pi^+(D(\mu_\kappa)) \to \mathbb{R}^p$ for an $p \in \mathbb{N}_0$. Let
	$$\Delta:=\{(t,s,z) \in \pi(C) \times \mathbb{R}_{>0} \times \mathbb{C} \mid N(t) < s < M(t), z \in B(\mu_{\kappa}(t),s)\}.$$ 
	Then the function
	$$\Psi:\Delta \to \mathbb{C}, (t,s,z) \mapsto \frac{1}{2 \pi i}\int_{\partial B(\mu_{\kappa}(t),s)}\frac{F(t,\xi)}{\xi-z} d\xi,$$
	is well-defined and constructible in $z$ with support function 
	$$\pi^+(\Delta) \to \mathbb{R}^p \times \mathbb{R} \times \mathbb{R}, (t,s) \mapsto (g(t,s),\mu_\kappa(t),s).$$
	This means that there is a $C$-regular $\beta:\pi(C) \to \mathbb{R}^p$, a log-analytic $h:\mathbb{R}^p \times \mathbb{R}^2 \times \mathbb{C} \to \mathbb{C}$ such that for $(t,s,z) \in \Delta$ 
	$$\Psi(t,s,z) =h(g(t,s),\mu_\kappa(t),s,z).$$ 
\end{lemma}

\begin{proof}
Let $G:\mathbb{R}^p \times \mathbb{R} \to \mathbb{R}$ be globally subanalytic such that for $(t,s,\tau) \in D(\mu_\kappa)$
$$F^*(t,s,\tau) = G(g(t,s),\tau).$$
For $(t,s,z) \in \Delta$ we have
\begin{eqnarray*}
	\int_{\partial B(\mu_{\kappa}(t),s)}\frac{F(t,\xi)}{\xi-z} d\xi&=&\int_{-\pi}^{\pi}\frac{F^*(t,s,\tau)}{\mu_{\kappa}(t)+(\prod_{j=0}^{\kappa}\sigma_j)se^{i \tau}-z} i(\prod_{j=0}^{\kappa}\sigma_j)se^{i\tau} d\tau\\
	&=&\int_{-\pi}^{\pi}\frac{G(g(t,s),\tau)}{\mu_{\kappa}(t)+(\prod_{j=0}^{\kappa}\sigma_j)se^{i \tau}-z} i(\prod_{j=0}^{\kappa}\sigma_j)se^{i\tau} d\tau.
\end{eqnarray*}

Note that 
$$\Delta \times \textnormal{}]-\pi,\pi[ \textnormal{} \to \mathbb{C}, (t,s,z,\tau) \mapsto \frac{(\prod_{j=0}^{\kappa}\sigma_j)se^{i\tau}}{\mu_{\kappa}(t)+(\prod_{j=0}^{\kappa}\sigma_j)se^{i \tau}-z},$$
is bounded in $\tau$. So $\Psi$ is well-defined. Let $v:=(v_1, \ldots ,v_p)$ range over $\mathbb{R}^p$ and $y$ over $\mathbb{R}$. Set
$$Q:=\{(v,x,y,z,\tau) \in \mathbb{R}^p \times \mathbb{R} \times \mathbb{R}_{>0} \times \mathbb{C} \times  \textnormal{}]-\pi,\pi[ \textnormal{} \mid \vert{x-z}\vert<y\}$$
and
$$\Omega:\mathbb{R}^p \times \mathbb{R} \times \mathbb{R} \times \mathbb{C} \times \mathbb{R} \to \mathbb{C}, (v,x,y,z,\tau) \mapsto $$
$$\left\{\begin{array}{ll} \frac{G(v,\tau)}{x+(\prod_{j=0}^{\kappa}\sigma_j)ye^{i \tau}-z} i (\prod_{j=0}^{\kappa} \sigma_j)ye^{i \tau}, & (v,x,y,z,\tau) \in Q,\\
	0, & \textnormal{else.} \end{array}\right.$$
Then $\Omega$ is globally subanalytic since we have that 
$$]-\pi,\pi[ \textnormal{} \to \mathbb{C}, \tau \mapsto e^{i \tau}=\cos(\tau)+i\sin(\tau),$$ 
is globally subanalytic. 
For $(t,s,z,\tau) \in \Delta \times \textnormal{} ]-\pi,\pi[$ we have $\vert{\mu_\kappa(t)-z}\vert<s$ and therefore $(g(t,s),\mu_\kappa(t),s,z,\tau) \in Q$. So we obtain
$$\frac{G(g(t,s),\tau)}{\mu_{\kappa}(t)+(\prod_{j=0}^{\kappa}\sigma_j)se^{i \tau}-z} i(\prod_{j=0}^\kappa \sigma_j) se^{i\tau}=\Omega(g(t,s),\mu_\kappa(t),s,z,\tau)$$
for every $(t,s,z,\tau) \in \Delta \times \textnormal{} ]-\pi,\pi[$. By Lemma~\ref{26} there is a constructible $H:\mathbb{R}^p \times \mathbb{R} \times \mathbb{R} \times \mathbb{C} \to \mathbb{C}$ such that $\Psi(t,s,z)=H(g(t,s),\mu_\kappa(t),s,z)$ for $(t,s,z) \in \Delta$. So we see that $\Psi$ is constructible in $z$ with the desired support function. By noting that every constructible function is log-analytic we are done with the proof of Lemma~\ref{5.63}.
\end{proof}

Lemma~\ref{5.63} suggests that we need to find a $\kappa$-persistent function $M:\pi(C) \to \mathbb{R}_{\geq 0}$ and a $\kappa$-soft function $N:\pi(C) \to \mathbb{R}_{\geq 0}$ such that $F^*|_{\mathcal{A}(\mathcal{P}_{0,\kappa},N,M)(\mu_\kappa)}$ (where $F^*$ is defined as in Lemma~\ref{5.63})
is bounded and globally subanalytic in $\tau$. However, we need some more information about $F^*$ to be able to prove that the unary high parametric global complexification, which we construct in Section~9 by gluing the single extensions together, is also restricted log-exp-analytic.
 
In the proof in Section~9 the parameter $s$ will be "replaced" by a function $\alpha:\pi(C) \to \mathbb{R}_{\geq 0}$ which naturally arises from $N$ and $M$. This ensures restricted log-exp-analyticity if the function $\Psi(t,\alpha(t),z)$ 
is the composition of a log-analytic function and $C$-regular functions. This matches also with the preparation result in Section~7, where the coefficient and base functions of the log-analytical preparation are $C$-regular. 

We achieve the desired property for $\Psi(t,\alpha(t),z)$ if $M,N$ are $C$-regular, and $F^*$ is \emph{free-regular} in $s$, a concept introduced in the next subsection.

\subsection{Free-Regular Functions}
For the remaining Section~8 we set the following. Fix $l,m \in \mathbb{N}_0$ such that $l+m=n$. Let $w:=(w_1, \ldots ,w_l)$ range over $\mathbb{R}^l$ and $u:=(u_1, \ldots ,u_m)$ over $\mathbb{R}^m$. Here $u$ is serving as tuple of independent variables of families of functions parameterized by $w$. 
Further let $X \subset \mathbb{R}^n \times \mathbb{R}$ be definable with $C \subset X$ and suppose that $X_w$ is open for every $w \in \mathbb{R}^l$. We say ''$C$-consistent'' respectively ''$C$-regular'' instead of ''$C$-consistent in $u$ with respect to $X$'' respectively ''$C$-regular in $u$ with respect to $X$'' and also ''$\kappa$-persistent'' respectively ''$\kappa$-soft'' instead of ''$(C,\mathcal{Y},\kappa)$-persistent in $u$ with respect to $X$'' respectively ''$(C,\mathcal{Y},\kappa)$-soft in $u$ with respect to $X$''.

\begin{definition}
	\label{5.64}
	Let $D \subset \pi(C) \times \mathbb{C}$ be definable. We call a function $f: D(\mu_\kappa) \to \mathbb{C}, (t,s,\tau) \mapsto f(t,s,\tau),$ \textbf{free-regular in $s$}\index{free-regular} if there is $p \in \mathbb{N}_0$ such that $f$ is globally subanalytic in $\tau$ with support function $g:\pi^+(D(\mu_\kappa)) \to \mathbb{R}^p$ and the following holds for $g$.
	There is $q \in \mathbb{N}_0$, a $C$-regular function $\beta: \pi(C) \to \mathbb{R}^q$ and a log-analytic $h:\mathbb{R}^{q+1} \to \mathbb{R}^p$ such that $g(t,s)=h(\beta(t),s)$ for every $(t,s) \in \pi^+(D(\mu_\kappa))$.
\end{definition}

\begin{remark}
	\label{5.65}
	Let $D \subset \pi(C) \times \mathbb{C}$ be definable. 
	\begin{itemize}
		\item [(1)] Let $f: D(\mu_\kappa) \to \mathbb{C}$ be the restriction of a globally subanalytic function. Then $f$ is free-regular in $s$. 
		\item [(2)] Let $f: D(\mu_\kappa) \to \mathbb{C}$ be globally subanalytic in $\tau$ with log-analytic support function. Then $f$ is free-regular in $s$.
	\end{itemize}
\end{remark}

\begin{lemma}
	\label{5.66}
	Let $D \subset \pi(C) \times \mathbb{C}$ be definable. Let $f: D(\mu_\kappa) \to \mathbb{C}$ be free-regular in $s$. The following holds.
	\begin{itemize}
		\item [(1)] $f$ is globally subanalytic in $\tau$ with support function $g:\pi^+(D(\mu_\kappa)) \to \mathbb{R}^q$ which can be constructed from a set $E$ of positive definable functions on $\pi(C)$ such that every $h \in \log(E)$ is $C$-consistent.
		\item [(2)] Let $\alpha \in \mathbb{N}$ and let $G:\mathbb{C}^\alpha \to \mathbb{C}$ be globally subanalytic. Let $J_1, \ldots ,J_{\alpha}:D(\mu_\kappa) \to \mathbb{C}$ be functions which are free-regular in $s$. Then $$H:=G(J_1, \ldots ,J_{\alpha}):D(\mu_\kappa) \to \mathbb{C}$$
		is free-regular in $s$. 
	\end{itemize}
\end{lemma}

\begin{proof}
	(1): Definition~\ref{5.64}, Definition~\ref{5.54} and Lemma~\ref{5.55}(2) show immediately that $g$ can be constructed from a set $E$ of positive definable functions on $\pi(C)$ such that every $h \in \log(E)$ is $C$-consistent for $j \in \{1, \ldots ,q\}$.
	
	(2): Follows directly with Lemma~\ref{25}.
\end{proof}

Since every constructible function is log-analytic, Lemma~\ref{5.63} applied to a free-regular function $F^*$ provides the following.

\begin{lemma}
	\label{lem:int-parameterized-free-regular}
	Let $M,N:\pi(C) \to \mathbb{R}_{\geq 0}$ be definable functions. Suppose that $M<\vert{\mu_r-\mu_k}\vert$ if $k:=k^{\text{ch}} \geq 0$. Let $D:=\mathcal{A}(\mathcal{P}_{0,\kappa},N,M)$. Let $F:D \to \mathbb{C}$ be definable such that 
	$$F^*:D(\mu_\kappa) \to \mathbb{C}, (t,s,\tau) \mapsto F(t,\mu_\kappa(t)+\Bigl(\prod_{j=0}^\kappa \sigma_j \Bigl) se^{i\tau}),$$
	is bounded in $\tau$. Suppose that $F^*$ is free-regular in $s$. Let
	$$\Delta:=\{(t,s,z) \in \pi(C) \times \mathbb{R}_{>0} \times \mathbb{C} \mid N(t) < s < M(t), z \in B(\mu_{\kappa}(t),s)\}.$$ 
	
	Then for
	$$\Psi:\Delta \to \mathbb{C}, (t,s,z) \mapsto \frac{1}{2 \pi i}\int_{\partial B(\mu_{\kappa}(t),s)}\frac{F(t,\xi)}{\xi-z} d\xi,$$
	there is $c \in \mathbb{N}_0$, a $C$-regular $\beta:\pi(C) \to \mathbb{R}^c$, a log-analytic $h:\mathbb{R}^{c} \times \mathbb{R} \times \mathbb{C} \to \mathbb{C}$ such that for $(t,s,z) \in \Delta$ 
	$$\Psi(t,s,z) =h(\beta(t),s,z).$$ 
\end{lemma}

\begin{proof}
	Let $g:\pi^+(D(\mu_\kappa)) \to \mathbb{R}^p$ be a corresponding support function for $F^*$ where $p \in \mathbb{N}_0$. Fix a $C$-regular $\beta^*: \pi(C) \to \mathbb{R}^q$ and a log-analytic $h^*:\mathbb{R}^{q+1} \to \mathbb{R}^p$ such that $g(t,s) = h^*(\beta^*(t),s)$ for every $(t,s) \in \pi^+(D(\mu_\kappa))$. By Lemma~\ref{5.63} we find a log-analytic $\hat{h}:\mathbb{R}^p \times \mathbb{R}^2 \times \mathbb{C} \to \mathbb{C}$ such that for $(t,s,z) \in \Delta$
	$$\Psi(t,s,z) = \hat{h}(g(t,s),\mu_\kappa(t),s,z) = \hat{h}(h^*(\beta^*(t),s),\mu_\kappa(t),s,z).$$
	So take $c:=q+1$ and $\beta:\pi(C) \to \mathbb{R}^c, t \mapsto (\beta^*(t),\mu_\kappa(t))$. Then $\beta$ is regular. Let $w_1, \ldots, w_c$ be a new tuple of real variables. Let 
	$$h:\mathbb{R}^c \times \mathbb{R} \times \mathbb{C} \to \mathbb{C}, (w_1, \ldots, w_c,s,z) \mapsto \hat{h}(h^*(w_1, \ldots , w_{c-1},s),w_c,s,z).$$ 
	Note that $h$ is log-analytic with $\Psi(t,s,z) = h(\beta(t),s,z)$ for $(t,s,z) \in \Delta$.
\end{proof}

\subsection{Logarithmic Scales}

In this subsection we investigate logarithmic scales and construct a suitable $\kappa$-persistent $M:\pi(C) \to \mathbb{R}_{\geq 0}$ such that for $D:=\mathcal{B}(\mathcal{P}_{0,\kappa},M)$ the function
$$F^*:D(\mu_\kappa) \to \mathbb{C}, (t,s,\tau) \mapsto F(t,\mu_\kappa(t)+\Bigl(\prod_{j=0}^\kappa \sigma_j \Bigl) se^{i\tau}),$$
is free-regular in $s$ where 
$$F:D \to \mathbb{C}, (t,z) \mapsto (\sigma\mathcal{Z})^{\otimes q}(t,z).$$ 
The choice for the $\kappa$-soft $N$ happens in Section~7.4 when we apply the preparation Theorem~\ref{5.84} to $\Phi_f$.

All the results in this subsection are formulated for $X=\mathbb{R}^n \times \mathbb{R}$ and therefore for every possible choice of $l,m \in \mathbb{N}_0$ (the number of components of $w$ and $u$) such that $n=l+m$. This includes every possible choice of $X$ with $C \subset X$ such that $X_w$ is open. Therefore we work with the variable $t$ instead of $w$ and $u$.

Further we set for the rest of this subsection
$$\log^*:\mathbb{C} \to \mathbb{C}, z \mapsto \left\{\begin{array}{lll} \log(z),&& \vert{z}\vert \in [1/2,3/2] \setminus \mathbb{R}_{\leq 0}, \\
	&\textnormal{if}&\\
	0,&&\textnormal{else.} \end{array}\right.$$

We also consider the functions $\mathcal{P}_{j,\kappa}$ from Definition~\ref{5.16} as functions on $H$.

\begin{lemma}
\label{5.68}
	There is a $C$-regular $\kappa$-persistent function $M:\pi(C) \to \mathbb{R}_{\geq 0}$ with $M<\vert{\mu_r-\mu_k}\vert$ if $k \geq 0$ such that for every $l \in \{0, \ldots ,\kappa\}$ the following holds where $D:=\mathcal{B}(\mathcal{P}_{0,\kappa},M)$.
	\begin{itemize}
		\item [(1)]
		The function
		$$\mathcal{P}_{l,\kappa}^*: D(\mu_\kappa) \to \mathbb{C}, (t,s,\tau) \mapsto \mathcal{P}_{l,\kappa}(t,\mu_\kappa(t)+(\prod_{j=0}^\kappa \sigma_j)se^{i\tau}),$$
		is free-regular in $s$.
		\item [(2)]
		The function
		$$z_l^*:D(\mu_\kappa) \to \mathbb{C}, (t,s,\tau) \mapsto z_l(t,\mu_\kappa(t)+(\prod_{j=0}^\kappa \sigma_j) se^{i\tau}),$$
		is free-regular in $s$.
\end{itemize}
\end{lemma}

\begin{proof}
By Lemma~\ref{5.67}(1) there is a $C$-regular $\kappa$-persistent $T:\pi(C) \to \mathbb{R}_{\geq 0}$ such that 
$$\Bigl|\frac{\mathcal{P}_{j-1,\kappa}(t,z)}{e^{\mu_{\kappa-j,\kappa}(t)}} \Bigl|<1/2$$
for $(t,z) \in \mathcal{B}(\mathcal{P}_{0,\kappa},T)$ and every $j \in \{1, \ldots ,\kappa\}$. Let $M:\pi(C) \to \mathbb{R}_{\geq 0}, t \mapsto \min\{\vert{\mu_r(t)-\mu_k(t)}\vert,T(t)\}$ and $D:=\mathcal{B}(\mathcal{P}_{0,\kappa},M)$.
Note that a $C$-regular function $h:\pi(C) \to \mathbb{R}$ considered as a function on $D(\mu_\kappa)$ is free-regular in $s$.

(1): We do an induction on $l$. Note that $\mathcal{P}_{0,\kappa}^*$ is the restriction of a globally subanalytic function, because $$\mathcal{P}_{0,\kappa}^*(t,s,\tau)=(\prod_{j=0}^\kappa \sigma_j) se^{i \tau}=(\prod_{j=0}^\kappa \sigma_j) s(\cos(\tau) +i\sin(\tau))$$
for $(t,s,\tau) \in D(\mu_\kappa)$. So $\mathcal{P}_{0,\kappa}^*$ is free-regular in $s$ by Remark~\ref{5.65}(1).

$l-1 \to l:$ For $(t,s,\tau) \in D(\mu_\kappa)$ we obtain
$$\mathcal{P}_{l,\kappa}^*(t,s,\tau)=\log^* \Bigl(1+\frac{\mathcal{P}_{l-1,\kappa}^*(t,s,\tau)}{e^{\mu_{\kappa-l,\kappa}(t)}} \Bigl).$$
Because $\log^*$ is globally subanalytic we see with Lemma~\ref{5.66}(2) that $\mathcal{P}_{l,\kappa}^*$ is free-regular in $s$ since $e^{\mu_{\kappa-l,\kappa}}$ is $C$-regular (and therefore $e^{\mu_{\kappa-l,\kappa}}$ is also free-regular in $s$ considered as a function on $D(\mu_\kappa)$).

(2): For $(t,s,\tau) \in D(\mu_\kappa)$ we obtain with Lemma~\ref{5.18}
$$z_l^*(t,s,\tau) = \mathcal{P}_{l,\kappa}^*(t,s,\tau)+\sigma_le^{\mu_{\kappa-l-1,\kappa}(t)}$$
for $l \in \{0, \ldots ,\kappa-1\}$ and with Corollary~\ref{5.19}
$$z_\kappa^*(t,s,\tau) = \mathcal{P}_{\kappa,\kappa}^*(t,s,\tau).$$
So $z_\kappa^*$ is $C$-regular in $s$ with (1). For every $l \in \{0, \ldots ,\kappa-1\}$ we obtain also $C$-regurality of $z_l$ in $s$ with (1) and Lemma~\ref{5.66}(2) since $e^{\mu_{\kappa-l-1,\kappa}}$ is $C$-regular.
\end{proof}

Now we are able to show that $F^*$ is free-regular in $s$ for a suitable choice of $M$ if $\kappa=r$.

\begin{corollary}
\label{5.69}
	Let $\kappa=r$. Let $q:=(q_0, \ldots ,q_r) \in \mathbb{Q}^{r+1}$. There is a $C$-regular $r$-persistent function $M:\pi(C) \to \mathbb{R}_{\geq 0}$ with $M<\vert{\mu_r-\mu_k}\vert$ if $k \geq 0$ such that for $D:=\mathcal{B}(\mathcal{P}_{0,r},M)$ the function
	$$F^*: D(\mu_r) \mapsto \mathbb{C}, (t,s,\tau) \mapsto (\sigma\mathcal{Z})^{\otimes q}(t,\mu_r(t)+(\prod_{j=0}^r\sigma_j)se^{i \tau}),$$
	is bounded in $\tau$ and free-regular in $s$.
\end{corollary}

\begin{proof}
By Lemma~\ref{5.68}(2) we find a $C$-regular $r$-persistent $T:\pi(C) \to \mathbb{R}_{\geq 0}$ such that 
$$z_l^*:\hat{D}(\mu_r) \to \mathbb{C}, (t,s,\tau) \mapsto z_l(t,\mu_r(t)+(\prod_{j=0}^r \sigma_j) se^{i \tau}),$$
is free-regular in $s$ for $l \in \{0, \ldots ,r\}$ where $\hat{D}:=\mathcal{B}(\mathcal{P}_{0,r},T)$. Then $\sigma_lz_l^*(t,s,\tau) \in \mathbb{C}^-$ for every $(t,s,\tau) \in \hat{D}(\mu_r)$ and $l \in \{0, \ldots ,r\}$ by Definition~\ref{5.11} since $\mu_r(t)+(\prod_{j=0}^r \sigma_j) se^{i \tau} \in H_t$ for $(t,s,\tau) \in \hat{D}(\mu_r)$. Let
$$\hat{F}^*: \hat{D}(\mu_r) \mapsto \mathbb{C}, (t,s,\tau) \mapsto (\sigma\mathcal{Z})^{\otimes q}(t,\mu_r(t)+(\prod_{j=0}^r\sigma_j)se^{i \tau}),$$
and
$$G:(\mathbb{C}^-)^{r+1} \to \mathbb{C}, (w_0, \ldots ,w_r) \mapsto \prod_{j=0}^r (w_j)^{q_j}.$$
Note that $G$ is globally subanalytic. We obtain 
$$\hat{F}^*(t,s,\tau)=G(\sigma_0z_0^*(t,s,\tau), \ldots ,\sigma_rz_r^*(t,s,\tau))$$
for every $(t,s,\tau) \in \hat{D}(\mu_r)$. Therefore $\hat{F}^*$ is free-regular in $s$ by Lemma~\ref{5.66}(2). Now we construct $M:\pi(C) \to \mathbb{R}_{\geq 0}$ such that $F^*:=\hat{F}^*|_{D(\mu_r)}$ is bounded in $\tau$ where $D:=\mathcal{B}(\mathcal{P}_{0,\kappa},M)$ and are done. Let
$$d:\pi(C) \to \mathbb{R}, t \mapsto \min \{\vert{\mu_r(t)-\mu_j(t)}\vert \mid j \in \{0, \ldots ,r-1\} \Bigl\}.$$
Note that $d$ is $C$-regular (by Example~\ref{5.56}(2) and Lemma~\ref{5.55}(2)) and $r$-persistent (by Example~\ref{5.43} and Remark~\ref{5.42}(3)). Choose $M:\pi(C) \to \mathbb{R}_{\geq 0}, t \mapsto \min\{T(t),d(t)\}$. Let $D:=\mathcal{B}(\mathcal{P}_{0,r},M)$.

\begin{claim}
\label{claim-3}
The function $F^*$ is bounded in $\tau$.
\end{claim}

\begin{proof}
Assume the contrary. Then there is $(t,s) \in \pi^+(D(\mu_\kappa))$ and a definable curve $\gamma_\tau: \textnormal{}]0,1[\textnormal{} \to \textnormal{}]-\pi,\pi[$ such that $(t,s,\gamma_\tau(y)) \in D(\mu_\kappa)$ for $y \in \textnormal{}]0,1[$ and 
$$\lim_{y \searrow 0} \vert{F^*(t,s,\gamma_\tau(y))}\vert = \lim_{y \searrow 0} \vert{(\sigma\mathcal{Z})^{\otimes q}(t,\mu_r(t)+(\prod_{j=0}^r \sigma_j)se^{i \gamma_\tau(y)})}\vert = \infty.$$
Note that $0<s<d$. Consider $\gamma: \textnormal{}]0,1[ \textnormal{} \to H, y \mapsto (\gamma_t(y),\gamma_z(y))$ with $\gamma_t(y)=t$ and $\gamma_z(y)=\mu_r(t)+(\prod_{j=0}^r \sigma_j)se^{i \gamma_\tau(y)}$ for $y \in \textnormal{}]0,1[$. Note that this curve is definable. Since $0<s<d$ we have
$$\lim_{y \searrow 0} (\gamma_z(y) - \mu_l(t)) = \lim_{y \searrow 0} \Bigl(\mu_r(t) + (\prod_{j=0}^r \sigma_j)se^{i \gamma_\tau(y)} - \mu_l(t)\Bigl) \neq 0$$
for $l \in \{0, \ldots ,r\}$ which is a contradiction to Lemma~\ref{5.31}(1) (which we can apply on $\gamma$ since $(II_l)$ in Lemma~\ref{5.28} holds for $l \in \{0, \ldots ,r\}$ for this curve $\gamma$ due to $\Theta_0(t), \ldots, \Theta_r(t) \in \mathbb{R}$).
\end{proof}
Hence, the proof of Corollary~\ref{5.69} is accomplished.
\end{proof}

Now we investigate the case ''$\kappa=k^{\text{ch}}=:k$'': For the rest of Subsection~8.3 we suppose that $k \geq 0$. It is left to show that we can find a $C$-regular $k$-persistent $M:\pi(C) \to \mathbb{R}_{\geq 0}$ such that $$z_l^*:\mathcal{B}(\mathcal{P}_{0,k},M)(\mu_k) \to \mathbb{C}, (t,s,\tau) \mapsto z_l(t,\mu_k(t)+(\prod_{j=0}^k \sigma_j)se^{i\tau}),$$
is free-regular in $s$ for $l \in \{k+1, \ldots ,r\}$.

We start our considerations with the following technical result which gives some general conditions on functions $f,g:D \to \mathbb{C}$ such that $\log(f \cdot g) = \log(f) + \log(g)$ where $D \subset H$. Note that $\log(z_1z_2) \neq \log(z_1)+\log(z_2)$ in general for $z_1,z_2 \in \mathbb{C}^-$ with $z_1z_2 \in \mathbb{C}^-$.  

\begin{lemma}
\label{5.70}
Let $D \subset H$ be definable such that $D_t$ is connected and $\mathbb{R} \cap D_t \neq \emptyset$ for $t \in \pi(D)$. Let $f,g:D \to \mathbb{C}^-$ be functions with the following properties.
	\begin{itemize}
		\item [(1)]
		$f$ and $g$ are continuous in $z$.
		\item [(2)]
		We have $f(t,z) \in \mathbb{R}_{>0}$ if $(t,z) \in D \cap (\pi(C) \times \mathbb{R})$. 
		\item [(3)]
		$f(t,z)g(t,z) \in \mathbb{C}^-$ for every $(t,z) \in D$.
	\end{itemize}
	Then for every $(t,z) \in D$
	$$\log(f(t,z)g(t,z))=\log(f(t,z))+\log(g(t,z)).$$
\end{lemma}

\begin{proof}
Suppose that there is $(t,z^*) \in D$ such that 
$$\log(f(t,z^*)g(t,z^*)) = \log(f(t,z^*))+\log(g(t,z^*)) \pm 2\pi i$$
(and therefore 
$$\vert{\textnormal{arg}(f(t,z^*))+\textnormal{arg}(g(t,z^*))}\vert \geq \pi.)$$
Note that 
$$D_t \to \mathbb{R}, z \mapsto \textnormal{arg}(f(t,z)) + \textnormal{arg}(g(t,z)),$$
is continuous and that $\vert{\textnormal{arg}(f(t,z)) + \textnormal{arg}(g(t,z))}\vert < \pi$ if $z \in \mathbb{R}$ for $t \in \pi(C)$ (since $\textnormal{arg}(f(t,z))=0$ and $\vert{\textnormal{arg}(g(t,z))}\vert<\pi$). With the intermediate value theorem (by taking a definable curve which connects $z^*$ and a point $a \in D_t \cap \mathbb{R}$) there is $z' \in D_t$ such that
$$\vert{\textnormal{arg}(f(t,z'))+\textnormal{arg}(g(t,z'))}\vert = \pi.$$ 
But then for $c:=\text{arg}(f(t,z'))+\text{arg}(g(t,z')) \in \{-\pi,\pi\}$
\begin{align*}
f(t,z')g(t,z') &= \vert{f(t,z')g(t,z')}\vert e^{i c} = \vert{f(t,z')g(t,z')}\vert (\cos(c) + i \sin(c))\\
&= -\vert{f(t,z')g(t,z')}\vert \in \mathbb{R}_{<0}, 
\end{align*}

a contradiction to (3).
\end{proof}

Now we investigate the case $k=0$. 
We split our considerations into three technical lemmas. In the first one we give some general properties of iterated logarithms of the form $\log_m(-\log(z))$ for $m \in \mathbb{N}_0$ which we need to show that $z_l^*$ is globally subanalytic in $\tau$ with log-analytic support function which implies that it is also free-regular in $s$.  The latter happens in the second and third lemma. 

We handle the general case for $k$ in Lemma~\ref{5.74} below. 

\begin{lemma}
\label{5.71}                                           
	Let $m \in \mathbb{N}_0$. Let $D:=\{z \in \mathbb{C}^- \mid \vert{z}\vert<1/\exp_m(1)\}$ and
	$$L_m:D \to \mathbb{C}, z \mapsto \log_m(-\log(z)).$$
	Then $L_m$ is well-defined , definable and holomorphic. We have $\log_m(-\log(z)) \in \mathbb{C}^-$ and $\log_m(-\log(z)) \in \mathbb{R}_{>0}$ if $z \in \mathbb{R}_{>0}$ for $z \in D$.
\end{lemma}

\begin{proof}
Let $z \in D$. If $z \notin \mathbb{R}_{>0}$ then
$$\log(z)=\log(\vert{z}\vert)+i\arg(z) \in \mathbb{C} \setminus \mathbb{R}$$
since $\arg(z) \neq 0$ and therefore $\log_m(-\log(z)) \in \mathbb{C}^-$. So assume $z \in \mathbb{R}_{>0}$. We obtain
$\exp_{m-1}(1) < -\log(z)$ and therefore
$$\exp_{m-1-l}(1) < \log_l(-\log(z))$$
for every $l \in \{1, \ldots ,m-1\}$ (where $\exp_0:=1$). We obtain $\log_m(-\log(z)) > 0$.
\end{proof}

\begin{lemma}
\label{5.72}
	Let $m \in \mathbb{N}$. Let
	$$c:=\frac{1}{\exp_m\Bigl(2\sqrt{\log^2(3/2)+\pi^2} \Bigl)}$$ 
	and $D:= \textnormal{}]0,c[\textnormal{} \times \textnormal{}]-\pi,\pi[$. The function
	$$T_l:D \to \mathbb{C}^-, (s,\tau) \mapsto \log_l(-\log(se^{i \tau})),$$
	is a well-defined globally subanalytic function in $\tau$ with log-analytic support function for every $l \in \{0, \ldots ,m\}$.
\end{lemma}

\begin{proof}
Note that $se^{i\tau} \in \mathbb{C}^-$ for every $(s,\tau) \in \mathbb{R}_{>0} \textnormal{} \times \textnormal{}]-\pi,\pi[$. So we obtain well-definability with Lemma~\ref{5.71}. Set 
$$S_0:D \to \mathbb{C}, (s,\tau) \mapsto \frac{i \tau}{\log(s)},$$
and inductively for $l \in  \{1, \ldots ,m-1\}$
$$S_l:D \to \mathbb{C}, (s,\tau) \mapsto \frac{\log(1+S_{l-1}(s,\tau))}{\log_l(-\log(s))}.$$

\begin{claim}
\label{claim-4}	
Let $l \in \{0, \ldots ,m-1\}$. Then $S_l$  is a well-defined globally subanalytic function in $\tau$ with log-analytic support function. For every $(s,\tau) \in D$ we have $\vert{S_l(s,\tau)}\vert<1/2$ and
$$T_l=\log_l(-\log(s))(1+S_l(s,\tau)).$$
\end{claim}

\begin{proof}
Let $(s,\tau) \in D$. We obtain
$$2 \sqrt{\log^2(3/2)+\pi^2}<\log_l(-\log(s)) \textnormal{ }(*)$$
for every $l \in \{0, \ldots ,m-1\}$. We do an induction on $l$. The case $l=0$ is clear.

$l-1 \to l:$ With the inductive hypothesis, $(*)$, and the definition of the complex logarithm we see that $S_l$ is well-defined. Note also that for $(s,\tau) \in D$
\begin{align*}
\vert{S_l(s,\tau)}\vert &= \left|\frac{\log^*(\vert{1+S_{l-1}(s,\tau)}\vert) + i \text{arg}(\log^*(1+S_{l-1}(s,\tau)))}{\log_l(-\log(s))}\right| \\
&\leq \frac{\sqrt{\log^2(3/2) + \pi^2 }}{2\sqrt{\log^2(3/2)+\pi^2}} = \frac{1}{2}
\end{align*}
since the argument of a complex number is bounded by $\pi$ and $\log^*(\vert{z}\vert) \leq \log(3/2)$ for $z \in \mathbb{C}$. With the inductive hypothesis and Lemma~\ref{25} it follows that $S_l$ is a globally subanalytic function in $\tau$ with log-analytic support function. Additionally, we obtain for $(s,\tau) \in D$ with the inductive hypothesis
\begin{align*}
	\log_l(-\log(s)-i\tau) &=\log(\log_{l-1}(-\log(s)-i\tau))\\
	&=\log(\log_{l-1}(-\log(s))(1+S_{l-1}(s,\tau)))\\
	&=\log_l(-\log(s))+\log(1+S_{l-1}(s,\tau))\\
	&=\log_l(-\log(s)) \Bigl(1+\frac{\log(1+S_{l-1}(s,\tau))}{\log_l(-\log(s))} \Bigl)\\
	&=\log_l(-\log(s))(1+S_l(s,\tau)). \qedhere
\end{align*}
\end{proof}

Note that by Claim~\ref{claim-4}
\begin{align*}
	T_m(s,\tau)&=\log(T_{m-1}(s,\tau))\\
	&=\log(\log_{m-1}(-\log(s))(1+S_{m-1}(s,\tau)))\\
	&=\log(\log_{m-1}(-\log(s))) + \log(1+S_{m-1}(s,\tau))\\
	&=\log_m(-\log(s)) + \log^*(1+S_{m-1}(s,\tau))
\end{align*}
which shows immediately that $T_m$ is globally subanalytic in $\tau$ with log-analytic support function. 
\end{proof}

\begin{lemma}
\label{5.73}
Let $k=k^{\text{ch}}=0$. There is a $C$-regular $0$-persistent function $M:\pi(C) \to \mathbb{R}_{\geq 0}$ with $M<\vert{\mu_r-\mu_0}\vert$ such that for $D:=\mathcal{B}(\mathcal{P}_{0,0},M)$ and every $l \in \{0, \ldots ,r\}$ the function
$$z_l^*:D(\mu_0) \to \mathbb{C}, (t,s,\tau) \mapsto z_l(\mu_0(t)+\sigma_0se^{i\tau}),$$
is globally subanalytic in $\tau$ with log-analytic support function.
\end{lemma}

\begin{proof}
We may assume that $C$ is near with respect to $\mu_0$. Otherwise we are done with the proof by choosing $M=0$ (see Remark~\ref{5.46}(1)). Note that $\sigma_1=-1$, $\sigma_2= \ldots =\sigma_r=1$ and that $\Theta_1= \ldots =\Theta_r=0$ by Lemma~\ref{5.47}. By Definition~\ref{5.2} we have $\mu_0=\Theta_0$, $\mu_1=\Theta_0+\sigma_0$ and for $j \in \{2, \ldots ,r\}$ we have $\mu_j=\Theta_0+\sigma_0e^{-e_{j-1}}$ where $e_0:=0$ and $e_m:=\exp(e_{m-1})$ for $m \in \mathbb{N}$. Note that for $(t,z) \in H$
$$z_1(t,z)=\log(\sigma_0z_0(t,z))$$
and 
$$z_l(t,z)=\log_{l-1}(-\log(\sigma_0z_0(t,z)))$$
for $l \in \{2, \ldots ,r\}$. Let 
$$c:=\frac{1}{\exp_{r-1}(2 \sqrt{\log^2(3/2)+\pi^2})}.$$
Then $c<1/\exp_{r-1}(1)=1/e_r=\vert{\mu_r(t)-\mu_0(t)}\vert$ for $t \in \pi(C)$. Consider $M:=c$. With Lemma~\ref{5.62} we see that
$$D(\mu_0)=\mathcal{B}(\mathcal{P}_{0,0},c)(\mu_0)=\pi(C) \textnormal{} \times \textnormal{} ]0,c[ \textnormal{} \times \textnormal{} ]-\pi,\pi[.$$
Note that 
$z_0^*(t,s,\tau)=\sigma_0se^{i\tau} = \sigma_0s(\cos(\tau) + i\sin(\tau))$ and $z_1^*(t,s,\tau)=\log(s)+i\tau$ for $(t,s,\tau) \in D(\mu_0)$. 
So $z_0^*$ and $z_1^*$ fulfill the desired properties since the restriction of the global sine and cosine function on $]-\pi,\pi[$ are globally subanalytic. We obtain for $l \in \{1, \ldots ,r-1\}$
\begin{eqnarray*}
	z_{l+1}^*(t,s,\tau)&=&\log_l(-\log(\sigma_0z_0(t,\Theta_0(t)+\sigma_0se^{i\tau})))\\
	&=&\log_l(-\log(se^{i \tau}))
\end{eqnarray*}
for every $(t,s,\tau) \in D(\mu_0)$. With Lemma~\ref{5.72} we see that 
$$]0,c[ \textnormal{} \times \textnormal{}]-\pi,\pi[\textnormal{} \to \mathbb{C}, (s,\tau) \mapsto \log_l(-\log(se^{i \tau})),$$
is a well-defined globally subanalytic function in $\tau$ with log-analytic support function for $l \in \{0, \ldots ,r-1\}$ and we are done with the proof of Lemma~\ref{5.73}.
\end{proof}

Now we consider the general case for $k$.

\begin{lemma}
\label{5.74}
	There is a $C$-regular $k$-persistent function $M:\pi(C) \to \mathbb{R}_{\geq 0}$ with $M<\vert{\mu_r-\mu_k}\vert$ such that for $D:=\mathcal{B}(\mathcal{P}_{0,k},M)$ and every $l \in \{0, \ldots ,r\}$ the function
	$$z_l^*:D(\mu_k) \to \mathbb{C}, (t,s,\tau) \mapsto z_l(\mu_k(t)+(\prod_{j=0}^k\sigma_j)se^{i\tau}),$$
	is free-regular in $s$.
\end{lemma}

\begin{proof}
We may assume that $C$ is near with respect to $\mu_k$. Otherwise we are done with the proof by choosing $M=0$ (see Remark~\ref{5.46}(1)). Note that $\sigma_{k+1}=-1$, $\sigma_{k+2}= \ldots =\sigma_r=1$ and that $\Theta_{k+1}= \ldots =\Theta_r=0$ by Lemma~\ref{5.47}. So for $(t,z) \in H$ we have with Corollary~\ref{5.19}
$$z_{k+1}(t,z)=\log(\sigma_kz_k(t,z))=\log(\sigma_k\mathcal{P}_k(t,z))$$
and 
$$z_l(t,z)=\log_{l-k-1}(-\log(\sigma_kz_k(t,z)))=\log_{l-k-1}(-\log(\sigma_k\mathcal{P}_k(t,z)))$$
for $l \in \{k+2, \ldots ,r\}$ (since $z_\kappa(t,z) = \mathcal{P}_k(t,z)$ for $(t,z) \in H$). If $k=0$ we are done with the proof by Lemma~\ref{5.73}. So assume $k > 0$.

Fix $\lambda \in \textnormal{}]0,1/2[$ such that 
$$\Bigl|{\frac{\log(1+z)-z}{z}}\Bigl| < 1/2$$
for $z \in \mathbb{C} \setminus \{0\}$ with $\vert{z}\vert<\lambda$. For $l \in \{1, \ldots ,k\}$ set 
$$\zeta_l:\pi(C) \to \mathbb{R}_{>0}, t \mapsto 
\min \Bigl\{\tfrac{1}{e^{2 \vert{\mu_{k-l,k}(t)}\vert}},\tfrac{1}{\exp_{r-k-1}(2\sqrt{\log^2(3/2)+\pi^2})},$$
$$\lambda e^{\mu_{k-l,k}(t)},\tfrac{e^{\mu_{k-l,k}(t)}}{\exp_{r-k-1}(2 \sqrt{\log^2(3/2)+\pi^2})} \Bigl\}.$$
We see that $\zeta_l$ is $C$-regular (by Lemma~\ref{5.55} and Example~\ref{5.56}(2)) and $k$-persistent (by Example~\ref{5.43} and Remark~\ref{5.42}(3)) for every $l \in \{1, \ldots ,k\}$. For $l \in \{1, \ldots ,k\}$ set
$$\Omega_l:=\{(t,z) \in H \mid \vert{\mathcal{P}_{l-1,k}(t,z)}\vert<\zeta_l(t)\}.$$
By Lemma~\ref{5.67}(1) there is a $C$-regular $k$-persistent $M_l:\pi(C) \to \mathbb{R}_{\geq 0}$ such that $\mathcal{B}(\mathcal{P}_{0,k},M_l) \subset \Omega_l$ for $l \in \{1, \ldots ,k\}$. Consider 
$$M^*:\pi(C) \to \mathbb{R}_{\geq 0}, t \mapsto \min\{\vert{\mu_r(t)-\mu_k(t)}\vert,M_1(t), \ldots ,M_k(t)\}.$$
Then $M^*$ is $k$-persistent and $C$-regular. By Lemma~\ref{5.68}(1) There is a $C$-regular $\kappa$-persistent $M \leq M^*$ such that for $D:=\mathcal{B}(\mathcal{P}_{0,\kappa},M)$
$$\mathcal{V}_l^*:D(\mu_k) \to \mathbb{C}, (t,s,\tau) \mapsto \mathcal{V}_l(t,\mu_k(t) + (\prod_{j=0}^k\sigma_j) se^{i \tau}),$$
is free-regular in $s$ for every $l \in \{0, \ldots ,k\}$. (Therefore 
$$z_l^*:D(\mu_k) \to \mathbb{C}, (t,s,\tau) \mapsto z_l(t,\mu_k(t)+(\prod_{j=0}^k \sigma_j)se^{i\tau}),$$
is also free-regular in $s$.) $(+++_\mathcal{V})$

We have $D \subset \bigcap_{l=1}^k \Omega_l.$
Clearly $D_t$ is connected and $D_t \cap \mathbb{R} \neq \emptyset$ for $t \in \pi(D)$. For $l \in \{1, \ldots ,k\}$ let
$$c_l:D \to \mathbb{C}, (t,z) \mapsto \frac{\sigma_{l-1}\mathcal{P}_{l-1,\kappa}(t,z)}{e^{\mu_{k-l,k}(t)}}.$$
Then for every $(t,z) \in D$ and $l \in \{1, \ldots ,k\}$ we have
$$\mathcal{P}_{l,\kappa}(t,z) = \log(1 + c_l(t,z)),$$
and $c_l$ and $\mathcal{P}_{l,k}$ are continuous in $z$ since $\mathcal{P}_{l,k}(t,z)=z_l(t,z)-\sigma_le^{\mu_{k-l-1,k}}(t)$ by Lemma~\ref{5.18} and $z_l$ is continuous in $z$ by Lemma~\ref{5.12}(1). 

Let $\mathcal{V} \in \{c,\mathcal{P}\}$. By Lemma~\ref{5.21} we have 
$$(\prod_{j=l}^k \sigma_j) \mathcal{V}_l(t,z) \in \mathbb{C}^-,$$
and
$$(\prod_{j=l}^k \sigma_j) \mathcal{V}_l(t,z) \in \mathbb{R}_{>0}$$
if $z \in \mathbb{R}$ for every $(t,z) \in D$ and $l \in \{1, \ldots ,k\}$. $\textnormal{ } (+_{\mathcal{V} })$

Note also that 
$$\vert{\mathcal{V}_l(t,z)}\vert<\frac{1}{\exp_{r-k-1} \Bigl(2\sqrt{\log^2(3/2)+\pi^2} \Bigl)}$$
for every $(t,z) \in D$ and $l \in \{1, \ldots ,k\}$. $(++_\mathcal{V})$

For $l \in \{0, \ldots ,k\}$ and $m \in \{0, \ldots ,r-k-1\}$ we define
$$R_{m,l}:D \to \mathbb{C}, (t,z) \mapsto \log_m(-\log((\prod_{j=l}^k \sigma_j)\mathcal{P}_{l,k}(t,z))).$$

By Lemma~\ref{5.71}, $(+_\mathcal{P})$ and $(++_\mathcal{P})$ (i.e. $(+_\mathcal{V})$ and $(++_\mathcal{V})$ for $\mathcal{V}=\mathcal{P}$) we obtain that the function $R_{m,l}$ is well-defined, continuous in $z$, and for $(t,z) \in D$ we have $R_{m,l}(t,z) \in \mathbb{C}^-$ and if $z \in \mathbb{R}$ then $R_{m,l}(t,z) \in \mathbb{R}_{>0}$. $(+_R)$ 

Note that $\sigma_{k+m+1}z_{k+m+1}=R_{m,k}$ for $m \in \{0, \ldots ,r-k-1\}$. Hence the goal is to show that $R_{m,l}^*:D(\mu_k) \to \mathbb{C}$ is free-regular in $s$ by an induction on the natural numbers $l$ and $m$. 
For this we have to introduce some more auxiliary functions. 

For $l \in \{1, \ldots ,k\}$ we set
$$S_{0,l}:D \to \mathbb{C}, (t,z) \mapsto \mu_{k-l,k}(t),$$
and inductively for $m \in \{1, \ldots ,r-k-1\}$
$$S_{m,l}:D \to \mathbb{C}, (t,z) \mapsto \log \Bigl(1 + \frac{S_{m-1,l}(t,z)}{R_{m-1,l-1}(t,z)} \Bigl).$$

For $l \in \{1, \ldots ,k\}$ and $m \in \{0, \ldots ,r-k-1\}$ let
$$T_{m,l}:D \to \mathbb{C}, (t,z) \mapsto \log_m(-\log((\prod_{j=l}^k \sigma_j)c_l(t,z))).$$
By Lemma~\ref{5.71}, $(+_c)$ and $(++_c)$ (i.e. $(+_\mathcal{V})$ and $(++_\mathcal{V})$ for $\mathcal{V}=c$) we obtain that the function $T_{m,l}$ is well-defined, continuous in $z$, and for $(t,z) \in D$ we have $T_{m,l}(t,z) \in \mathbb{C}^-$ and if $z \in \mathbb{R}$ then $T_{m,l}(t,z) \in \mathbb{R}_{>0}$. $(+_T)$

Let for $l \in \{1, \ldots ,k\}$
$$b_l:D \to \mathbb{C}, (t,z) \mapsto \frac{\log(1+c_l(t,z))-c_l(t,z)}{c_l(t,z)},$$
$$U_{0,l}:D \to \mathbb{C}, (t,z) \mapsto -\log(1+b_l(t,z)),$$
and inductively for $m \in \{1, \ldots ,r-k-1\}$
$$U_{m,l}:D \to \mathbb{C}, (t,z) \mapsto \log \Bigl(1+\frac{U_{m-1,l}(t,z)}{T_{m-1,l}(t,z)} \Bigl).$$
Note that $\vert{c_l(t,z)}\vert<\lambda$ since $\vert{\mathcal{P}_{l-1,k}(t,z)}\vert<\lambda e^{\mu_{k-l,k}(t)}$ for $(t,z) \in D$. So we see that $b_l$ is well-defined, continuous in $z$ (by using the logarithmic series), $\vert{b_l(t,z)}\vert < 1/2$ by the choice of $\lambda$ and $1+b_l(t,z) \in \mathbb{C}^-$ for $(t,z) \in D$. $(+_b)$

\begin{claim}
\label{claim-5}
Let $l \in \{1, \ldots ,k\}$. The function $S_{m-1,l}$ is well-defined, continuous in $z$ and for $(t,z) \in D$ we have 
$$\Bigl|\frac{S_{m-1,l}(t,z)}{R_{m-1,l-1}(t,z)}\Bigl|<1/2$$
for every $m \in \{1, \ldots ,r-k-1\}$.
\end{claim}

\begin{proof}
Note that for $(t,z) \in D$
$$\vert{\mathcal{P}_{l-1,k}(t,z)}\vert<\min \Bigl\{\tfrac{1}{e^{2\vert{\mu_{k-l,k}(t)}\vert}}, \tfrac{1}{\exp_{r-k-1}(2\sqrt{\log(3/2)^2+\pi^2})} \Bigl\}. \textnormal{ }(*)$$
We show Claim~\ref{claim-5} by induction on $m \in \{1, \ldots ,r-k-1\}$. Fix $(t,z) \in D$.

$m=1$: That $S_{0,l}$ is well-defined and continuous in $z$ is clear. We obtain with $(*)$
$$\vert{\mu_{k-l,k}(t)}\vert<-1/2 \log(\vert{\mathcal{P}_{l-1,k}(t,z)}\vert)$$
and with the definition of the complex logarithm
$$\vert{\mu_{k-l,k}(t)}\vert<1/2 \vert{\log((\prod_{j=l-1}^k \sigma_j)\mathcal{P}_{l-1,k}(t,z))}\vert.$$
This gives 
$$\Bigl|\frac{S_{0,l}(t,z)}{R_{0,l-1}(t,z)}\Bigl| = \Bigl|\frac{\mu_{k-l,k}(t)}{\log((\prod_{j=l-1}^k \sigma_j)\mathcal{P}_{l-1,k}(t,z))}\Bigl| <1/2.$$

$m-1 \to m:$ 
We obtain well-definability of $S_{m-1,l}$ by the inductive hypothesis. We get
\begin{eqnarray*}
	\vert{S_{m-1,l}(t,z)}\vert&=&\Bigl|{\log \Bigl(1+\frac{S_{m-2,l}(t,z)}{R_{m-2,l-1}(t,z)} \Bigl)} \Bigl|\\
	&<& \sqrt{\log^2(3/2)+\pi^2}
\end{eqnarray*}
by the inductive hypothesis and definition of the complex logarithm. So we obtain by $(*)$
$$\vert{\mathcal{P}_{l-1,k}(t,z)}\vert<\frac{1}{\exp_m(2\vert{S_{m-1,l}(t,z)}\vert)},$$
i.e.
$$2\vert{S_{m-1,l}(t,z)}\vert<\log_{m-1}(-\log(\vert{\mathcal{P}_{l-1,k}(t,z)}\vert)).$$
By the definition of the complex logarithm we obtain that
$$\Bigl|\frac{S_{m-1,l}(t,z)}{R_{m-1,l-1}(t,z)}\Bigl|<1/2$$
and we are done with the proof of this claim since $R_{m-1,l-1}$ is continuous in $z$.
\end{proof}

By Claim~\ref{claim-5} and $(+_R)$ we have that $S_{m,l}$ is well-defined, continuous in $z$ and
$$1+\frac{S_{m-1,l}(t,z)}{R_{m-1,l-1}(t,z)} \in \mathbb{C}^-$$
for $(t,z) \in D$, $l \in \{1, \ldots ,k\}$ and $m \in \{1, \ldots ,r-k-1\}$. $(+_S)$

\begin{claim}
\label{claim-6}
Let $l \in \{1, \ldots ,k\}$. The function $U_{m-1,l}$ is well-defined, continuous in $z$ and for $(t,z) \in D$ we have
$$\Bigl|\frac{U_{m-1,l}(t,z)}{T_{m-1,l}(t,z)}\Bigl|<1/2$$
for every $m \in \{1, \ldots ,r-k-1\}$.
\end{claim}

\begin{proof}
Fix $(t,z) \in D$. We have
$$\vert{\mathcal{P}_{l-1,k}(t,z)}\vert<\frac{e^{\mu_{k-l,k}(t)}}{\exp_{r-k-1}(2\sqrt{\log^2(3/2)+\pi^2})}$$
and by $(+_b)$
$$\vert{U_{0,l}(t,z)}\vert < \sqrt{\log^2(3/2)+\pi^2}. \textnormal{ }(**)$$
We do an induction on $m$.

$m=1$: With $(+_b)$ we obtain that $U_{0,l}$ is well-defined and continuous in $z$. With $(**)$ we obtain
$$\vert{\mathcal{P}_{l-1,k}(t,z)}\vert<\frac{e^{\mu_{k-l,k}(t)}}{\exp(2\sqrt{\log^2(3/2)+\pi^2})}$$
and therefore
$$\vert{\mathcal{P}_{l-1,k}(t,z)}\vert<\frac{e^{\mu_{k-l,k}(t)}}{\exp(2\vert{U_{0,l}(t,z)}\vert)}.$$
This gives
$$\vert{U_{0,l}(t,z)}\vert < -1/2\log(\frac{\vert{\mathcal{P}_{l-1,k}(t,z)}\vert}{e^{\mu_{k-l,k}(t)}}) = -1/2\log(\vert{c_l(t,z)}\vert)$$
and with the definition of the complex logarithm
$$\vert{U_{0,l}(t,z)}\vert<1/2 \vert{\log((\prod_{j=l}^k \sigma_j)c_l(t,z))}\vert.$$
We obtain
$$\Bigl|{\frac{U_{0,l}(t,z)}{T_{0,l}(t,z)}}\Bigl|<1/2.$$

$m-1 \to m:$
We obtain well-definability of $U_{m-1,l}$ by the inductive hypothesis. We get
\begin{eqnarray*}
	\vert{U_{m-1,l}(t,z)}\vert&=&\Bigl|{\log \Bigl(1+\frac{U_{m-2,l}(t,z)}{T_{m-2,l}(t,z)} \Bigl)}\Bigl|\\
	& < &\sqrt{\log^2(3/2)+\pi^2}
\end{eqnarray*} 
by the inductive hypothesis. This gives with $(**)$
$$\vert{\mathcal{P}_{l-1,k}(t,z)}\vert<\frac{e^{\mu_{k-l,k}(t)}}{\exp_m(2\vert{U_{m-1,l}(t,z)}\vert)}.$$
This implies
$$2\vert{U_{m-1,l}(t,z)}\vert<\log_{m-1}(-\log(\vert{c_l(t,z)}\vert))$$
and with the definition of the complex logarithm we obtain
$$\Bigl|{\frac{U_{m-1,l}(t,z)}{T_{m-1,l}(t,z)}}\Bigl|<1/2$$
and are done with the proof of this claim. (We obtain the continuity of $U_{m-1,l}$ in $z$ since $T_{m-2,l}$ is continuous in $z$.)
\end{proof}

By Claim~\ref{claim-6} we have that $U_{m,l}$ is well-defined, continuous in $z$, and we have 
$$1+\frac{U_{m-1,l}(t,z)}{T_{m-1,l}(t,z)} \in \mathbb{C}^-$$
for $(t,z) \in D$, $l \in \{1, \ldots ,k\}$ and $m \in \{1, \ldots ,r-k-1\}$. $(+_U)$

\begin{claim}
\label{claim-7}
Let $(t,z) \in D$. Let $l \in \{1, \ldots ,k\}$. It holds
$$T_{m,l}(t,z)=R_{m,l-1}(t,z)+S_{m,l}(t,z)$$
for every $m \in \{0, \ldots ,r-k-1\}$.
\end{claim}

\begin{proof}
We do an induction on $m$.

$m=0$: We obtain
\begin{eqnarray*}
	T_{0,l}(t,z)&=&-\log((\prod_{j=l}^k \sigma_j)c_l(t,z))\\
	&=&-\log((\prod_{j=l-1}^k \sigma_j)\mathcal{P}_{l-1,k}(t,z))+\mu_{k-l,k}(t)\\
	&=& R_{0,l-1}(t,z)+S_{0,l}(t,z).
\end{eqnarray*}

$m-1 \to m$:
By Lemma~\ref{5.70}, $(+_R)$, $(+_S)$, $(+_T)$ we obtain with the inductive hypothesis
\begin{eqnarray*}
	\log(T_{m-1,l}(t,z))&=&\log(R_{m-1,l-1}(t,z)+S_{m-1,l}(t,z))\\
	&=&\log(R_{m-1,l-1}(t,z))+\log \Bigl(1+\frac{S_{m-1,l}(t,z)}{R_{m-1,l-1}(t,z)} \Bigl) \\
	&=&R_{m,l-1}(t,z)+\log \Bigl(1+\frac{S_{m-1,l}(t,z)}{R_{m-1,l-1}(t,z)} \Bigl)
\end{eqnarray*}
and therefore
\begin{align*}
T_{m,l}(t,z)&=\log(T_{m-1,l}(t,z))=R_{m,l-1}(t,z)+S_{m,l}(t,z). \qedhere
\end{align*}
\end{proof}

\begin{claim}
\label{claim-8}
Let $(t,z) \in D$. Let $l \in \{1, \ldots ,k\}$. It holds
$$R_{m,l}(t,z)=T_{m,l}(t,z)+U_{m,l}(t,z)$$
for every $m \in \{0, \ldots ,r-k-1\}$.
\end{claim}

\begin{proof}
We do an induction on $m$. 

$m=0$: We obtain
\begin{eqnarray*}
	R_{0,l}(t,z) &=& -\log((\prod_{j=l}^k \sigma_j) \mathcal{P}_{l,k}(t,z))=-\log((\prod_{j=l}^k \sigma_j) \log(1+c_l(t,z)))\\
	&=&-\log((\prod_{j=l}^k \sigma_j)(c_l(t,z) + \log(1+c_l(t,z)) -c_l(t,z)))\\
	&=&-\log((\prod_{j=l}^k \sigma_j)c_l(t,z)(1 + b_l(t,z))).
\end{eqnarray*}

With $(+_b)$, $(+_c)$, $(+_\mathcal{P})$ and Lemma~\ref{5.70} we obtain
$$-\log((\prod_{j=l}^k \sigma_j)c_l(t,z)(1 + b_l(t,z))) =-\log((\prod_{j=l}^k \sigma_j)c_l(t,z)) -\log(1 + b_l(t,z)).$$
Therefore
\begin{align*}
R_{0,l}(t,z) &= T_{0,l}(t,z) + U_{0,l}(t,z).
\end{align*}

$m-1 \to m$: By Lemma~\ref{5.70}, $(+_R)$, $(+_T)$, $(+_U)$ and the inductive hypothesis we obtain
\begin{eqnarray*}
	\log(R_{m-1,l}(t,z)) &=& \log(T_{m-1,l}(t,z)+U_{m-1,l}(t,z))\\
	&=&\log(T_{m-1,l}(t,z))+\log \Bigl(1+\frac{U_{m-1,l}(t,z)}{T_{m-1,l}(t,z)} \Bigl)
\end{eqnarray*}
and therefore
\begin{align*}
	R_{m,l}(t,z) &= \log(R_{m-1,l}(t,z)) = T_{m,l}(t,z)+U_{m,l}(t,z). \qedhere
\end{align*} 
\end{proof}

In the next claim we show that the corresponding function $R_{m,l}^*$ is free-regular in $s$ by showing with a simultaneous induction on $l$ that $T_{m,l}^*$ and $R_{m,l}^*$ are free-regular in $s$ for $m \in \{0, \ldots ,r-k-1\}$.

\begin{claim}
\label{claim-9}

The following holds.
\begin{itemize}
	\item [$(I_l)$] for $l \in \{1, \ldots ,k\}$: The function
	$$T_{m,l}^*:D(\mu_k) \to \mathbb{C}, (t,s,\tau) \mapsto T_{m,l}(t,\mu_k(t)+(\prod_{j=0}^k \sigma_j)se^{i\tau}),$$
	is free-regular in $s$ for every $m \in \{0, \ldots ,r-k-1\}$.
	
	\item [$(II_l)$] for $l \in \{0, \ldots ,k\}$: The function 
	$$R_{m,l}^*:D(\mu_k) \to \mathbb{C}, (t,s,\tau) \mapsto R_{m,l}(t,\mu_k(t)+(\prod_{j=0}^k \sigma_j)se^{i\tau}),$$
	is free-regular in $s$ for every $m \in \{0, \ldots ,r-k-1\}$.
\end{itemize}
\end{claim}

\begin{proof}
We start with $(II_0)$. Let $(t,s,\tau) \in D(\mu_k)$. We have
$$0<s<\tfrac{1}{\exp_{r-k-1}(2\sqrt{\log(3/2)^2+\pi^2})}$$
since 
$$\vert{z-\mu_k(t)}\vert = \vert{\mathcal{P}_{0,k}(t,z)}\vert < \tfrac{1}{\exp_{r-k-1}(2\sqrt{\log(3/2)^2+\pi^2})}$$
for all $(t,z) \in D$. Note that
$$R_{m,0}^*(t,s,\tau)=\log_m(-\log(se^{i\tau}))$$
for $m \in \{0, \ldots ,r-k-1\}$. Therefore with Lemma~\ref{5.72} and Remark~\ref{5.65}(2) we are done with the proof of $(II_0)$.

Assume that $(I_{l-1})$ and $(II_{l-1})$ have already been shown for $l \in \{1, \ldots ,k\}$ (where $(I_0)$ ius identified with any true statement).
We show $(I_l)$. By Claim~\ref{claim-7} we have $T_{m,l}=R_{m,l-1}+S_{m,l}$ on $D$ for every $m \in \{0, \ldots ,r-k-1\}$. With $(II_{l-1})$ and Lemma~\ref{5.66}(2) it suffices to show that 
$$S_{m,l}^*:D(\mu_k) \to \mathbb{C}, (t,s,\tau) \mapsto S_{m,l}(t,\mu_k(t)+(\prod_{j=0}^k \sigma_j)se^{i\tau}),$$
is free-regular in $s$ for every $m \in \{0, \ldots ,r-k-1\}$. We do an induction on $m$. For $m=0$ this is clear since $S_{0,m}^*$ considered as a function on $\pi(D)$ is the restriction of the $C$-regular function $\pi(C) \to \mathbb{R}, t \mapsto \mu_{k-l,k}(t)$. So with Definition~\ref{5.64} one sees immediately that $S_{0,m}^*$ is free-regular in $s$. Assume $m>0$. By Claim~\ref{claim-5} we obtain that 
$$\Bigl|{\frac{S_{m-1,l}(t,z)}{R_{m-1,l-1}(t,z)}}\Bigl|<1/2$$
for every $(t,z) \in D$. So
$$S_{m,l}=\log^* \Bigl(1+\frac{S_{m-1,l}}{R_{m-1,l-1}} \Bigl).$$
With $(II_{l-1})$, the inductive hypothesis and Lemma~\ref{5.66}(2) we are done with the proof of $(I_l)$. 

We show $(II_l)$. By Claim~\ref{claim-8} we have $R_{m,l}=T_{m,l}+U_{m,l}$ on $D$ for every $m \in \{0, \ldots ,r-k-1\}$. With $(I_l)$ and Lemma~\ref{5.66}(2) it suffices to show that 
$$U_{m,l}^*:D(\mu_k) \to \mathbb{R}, (t,s,\tau) \mapsto U_{m,l}(t,\mu_k(t)+(\prod_{j=0}^k \sigma_j)se^{i\tau}),$$
is free-regular in $s$ for every $m \in \{0, \ldots ,r-k-1\}$. We do an induction on $m$.

$m=0$: Note that $U_{0,l}(t,z)=\log^*(1+b_l(t,z))$  and 
$$b_l(t,z)=\frac{\log^*(1+c_l(t,z))-c_l(t,z)}{c_l(t,z)}$$
for every $(t,z) \in D$ (compare with $(+_b)$ and $(++_c)$). With $(+++_c)$ and Lemma~\ref{5.66}(2) we are done with the case "$m=0$".

$m-1 \to m$:
By Claim~\ref{claim-6} we obtain that 
$$\Bigl|\frac{U_{m-1,l}(t,z)}{T_{m-1,l}(t,z)} \Bigl|<1/2$$
for every $(t,z) \in D$. So
$$U_{m,l}=\log^* \Bigl(1+\frac{U_{m-1,l}}{T_{m-1,l}} \Bigl).$$
With $(I_l)$, the inductive hypothesis and Lemma~\ref{5.66}(2) we are done with the complete proof of Claim~\ref{claim-9}.  
\end{proof}

We have
$$\sigma_{k+m+1}z_{k+m+1}|_D=R_{m,k}$$
for every $m \in \{0, \ldots ,r-k-1\}$. Therefore we obtain with $(II_k)$ in Claim 5 that $z_{k+m+1}^*$ is free-regular in $s$ for every $m \in \{0, \ldots ,r-k-1\}$. By $(+++_\mathcal{P})$ we also have that $z_l^*$ is free-regular in $s$ for $l \in \{0, \ldots ,k\}$. This finishes the proof of Lemma~\ref{5.74}.
\end{proof}

Now we are able to show that $F^*$ is free-regular in $s$ for a suitable choice of $M$ also if $\kappa=k$.

\begin{corollary}
	\label{5.75}
	Let $q:=(q_0, \ldots ,q_r) \in \mathbb{Q}^{r+1}$. There is a $C$-regular $k$-persistent function $M:\pi(C) \to \mathbb{R}_{\geq 0}$ with $M<\vert{\mu_r-\mu_k}\vert$ such that for $D:=\mathcal{B}(\mathcal{P}_{0,k},M)$ the function
	$$F^*: D(\mu_k) \mapsto \mathbb{C}, (t,s,\tau) \mapsto (\sigma\mathcal{Z})^{\otimes q}(t,\mu_k(t)+(\prod_{j=0}^k\sigma_j)se^{i \tau}),$$
	is bounded in $\tau$ and free-regular in $s$.
\end{corollary}

\begin{proof}
The proof for this corollary is the same as the proof of Corollary~\ref{5.69} with the minor difference that $\kappa=k$ instead of $\kappa=r$.
\end{proof}

\subsection{Restricted Log-Exp-Analytically Prepared \\ Functions}

Outgoing from the preparation result stated in Theorem~\ref{5.84} in Section~7, and Corollaries~\ref{5.75} and~\ref{5.69} from the previous subsection we are able to prove the desired result on integration for the unary parametric global complexification of a restricted log-exp-analytically prepared function from Definition~\ref{5.78}.

For the whole subsection 
fix a function $f:C \to \mathbb{R}, (w,u,x) \mapsto f(w,u,x),$ which is $(m+1,X)$-restricted log-exp-analytically $(e,r)$-prepared in $x$. Fix the unary parametric global complexification $\Phi_f:\Lambda_f \to \mathbb{C}$ of $f$ from Definition~\ref{5.78}. 

\begin{lemma}
\label{5.85}
There is a $C$-regular $\kappa$-persistent $M:\pi(C) \to \mathbb{R}_{\geq 0}$ and a $C$-regular $\kappa$-soft $N:\pi(C) \to \mathbb{R}_{\geq 0}$ such that the following holds. Let
$$\Delta:=\{(t,s,z) \in \pi(C) \times \mathbb{R}_{>0} \times\mathbb{C} \textnormal{ }\mid N(t) < s < M(t),z \in B(\mu_{\kappa}(t),s)\}.$$
Then for the function
$$\Psi:\Delta \to \mathbb{C}, (t,s,z) \mapsto \frac{1}{2\pi i} \int_{\partial B(\mu_\kappa(t),s)} \frac{\Phi_f(t,\xi)}{\xi-z}d\xi,$$
there is $c \in \mathbb{N}_0$, a $C$-regular $\beta:\pi(C) \to \mathbb{R}^c$, a log-analytic $h:\mathbb{R}^{c+1} \times \mathbb{C} \to \mathbb{C}$ such that for $(t,s,z) \in \Delta$ 
$$\Psi(t,s,z) =h(\beta(t),s,z).$$ 
\end{lemma}

\begin{proof}
By Theorem~\ref{5.84} there is a $C$-regular $\kappa$-persistent $\hat{M}:\pi(C) \to \mathbb{R}_{\geq 0}$ and a $C$-regular $\kappa$-soft $N:\pi(C) \to \mathbb{R}_{\geq 0}$ such that $F|_{\hat{W}}$ is complex $r$-log-analytically prepared in $z$ with center $\Theta$ with $C$-regular coefficient and base functions where $\hat{W}:=\mathcal{A}(\mathcal{P}_{0,\kappa},N,\hat{M})$. Let 
$$(r,\mathcal{Z},a,q,s,V,b,P)$$
be a complex LA-preparing tuple for $F|_{\hat{W}}$ where $a$ and $b$ are $C$-regular.

By Lemma~\ref{5.68}(2) if $\kappa=r$ respectively Lemma~\ref{5.74} if $\kappa=k$ we find a $C$-regular $\kappa$-persistent $\tilde{M}:\pi(C) \to \mathbb{R}_{\geq 0}$ (with $\tilde{M}<\vert{\mu_r-\mu_k}\vert$ if $k \geq 0$) such that for $W:=\mathcal{A}(\mathcal{P}_{0,\kappa},N,\tilde{M})$ we have $W \subset \hat{W}$ and
$$z_\alpha^*:W(\mu_{\kappa}) \to \mathbb{C}, (t,s,\tau) \mapsto z_\alpha(t,\mu_\kappa(t)+(\prod_{j=0}^\kappa \sigma_j) se^{i \tau}),$$
is free-regular in $s$ for $\alpha \in \{0, \ldots ,r\}$. (Note that
$$W(\mu_\kappa)=\{(t,s) \in \pi(C) \times \mathbb{R}_{>0} \textnormal{ }\mid N(t) < s < \hat{M}(t)\} \times \textnormal{}]-\pi,\pi[.)$$

Let $\eta:=(a,b_1, \ldots ,b_s)$. Note that $\eta$ is $C$-regular. Similarly as in the proof of Claim~\ref{claim-10} there is a globally subanalytic function $G:\mathbb{C}^{s+1} \times \mathbb{C}^{r+1} \to \mathbb{C}$ such that 
$$F(t,z)=G(\eta(t),z_0(t,z), \ldots ,z_r(t,z))$$
for every $(t,z) \in W$. With Lemma~\ref{5.66}(2) we have that
$$F^*:W(\mu_\kappa) \to \mathbb{C}, (t,s,\tau) \mapsto F(t,\mu_\kappa(t)+(\prod_{j=0}^\kappa \sigma_j) se^{i \tau}),$$
is free-regular in $s$ (since the components of $\eta$ considered as functions on $W(\mu_\kappa)$ are free-regular in $s$). Further we find a $C$-regular $M$ with $M \leq \hat{M}$ such that $F^*|_{\mathcal{A}(\mathcal{P}_{0,\kappa},N,M)(\mu_\kappa)}$ is bounded in $\tau$. (By Corollary~\ref{5.69} if $\kappa=r$ respectively Corollary~\ref{5.75} if $\kappa=k$ we obtain that 
$$\mathcal{A}(\mathcal{P}_{0,\kappa},N,M)(\mu_\kappa) \to \mathbb{C}, (t,s,\tau) \mapsto (\sigma\mathcal{Z})^{\otimes q}(t,\mu_\kappa(t)+(\prod_{j=0}^\kappa \sigma_j) se^{i \tau}),$$
is bounded in $\tau$.) So we obtain the result with Lemma~\ref{lem:int-parameterized-free-regular}. 
\end{proof}

Note that Lemma~\ref{5.85} covers also the scenario when $\Phi_f$ is complex log-analytic-ally prepared according to Definition~\ref{5.76} with $C$-nice coefficient, center and base functions. Namely, the latter is the unary parametric global complexification of a log-exp-analytically prepared function $f$ with respect to $E=\{1\}$, implying that $f$ is also restricted log-exp-analytically prepared.

\section{Global Complexification of Restricted Log-Exp-Analytic Functions}

In this section we will establish the main result of this paper and prove that a real analytic restricted log-exp-analytic function $f:X \to \mathbb{R}$ where $X \subset \mathbb{R}^n$ is open has a global compexification $F:Z \to \mathbb{C}$. Specifically there exists an open set $Z \subset \mathbb{C}^n$ with $X \subset Z$ and a holomorphic function $F: Z \to \mathbb{C}$ with $F|_X=f$, and $F$ is again restricted log-exp-analytic. We adapt Kaiser's arguments from~\cite{17} to our context as follows. Firstly we construct a unary \emph{high} parametric global complexification of $f$ which is again restricted log-exp-analytic (see Subsection~9.1). Secondly we do a technical induction on the number of variables (see Subsection~9.2). The \emph{highness} is needed to enable this induction.


\subsection{The Univariate Case}
We start with the following result from analysis.

\begin{lemma}
\label{5.89}
Let $Y$ be a non-empty topological space. Let $g:Y \to \mathbb{R}$ be a function that is bounded from below. Let 
$$f:Y \to \mathbb{R},\textnormal{ } x \mapsto \liminf\limits_{y \to x} g(y).$$
Then $f$ is lower semi-continuous, i.e. for every $x \in Y$ and every $a<f(x)$ there is an open neighborhood $V$ of $x$ such that $a<f(y)$ for all $y \in V$.
\end{lemma}

\begin{proof}
See for example Bourbaki (~\cite{4}, Chapter IV, §6, Section 2, Proposition 4).
\end{proof}

For the whole Subsection~9.1 we set the following.

Let $n \in \mathbb{N}_0$. Let $t$ range over $\mathbb{R}^n$ and $x$ over $\mathbb{R}$.
Let $l,m \in \mathbb{N}_0$ be with $n=l+m$. Let $w$ range over $\mathbb{R}^l$ and $u$ over $\mathbb{R}^m$. Let $\pi_l:\mathbb{R}^l \times \mathbb{R}^m \times \mathbb{C} \to \mathbb{R}^l, (w,u,z) \mapsto w,$ be the projection on the first $l$ coordinates. Here $(u,x)$ and $(u,z)$ are serving as tuples of independent variables of families of functions parameterized by $w$. For a definable curve $\gamma:\textnormal{}]0,1[\textnormal{} \to \mathbb{R}^m \times \mathbb{C}$ with $\gamma:=(\gamma_1, \ldots ,\gamma_{m+1})$ we set $\gamma_u:=(\gamma_1, \ldots ,\gamma_m)$ for the first $m$ real components, $\gamma_z:=\gamma_{n+1}$ for the last component and if $\gamma_{m+1}(y) \in \mathbb{R}$ for every $y \in \textnormal{}]0,1[$ we also write $\gamma_x$ instead of $\gamma_z$.

For $p \in \mathbb{N}$, definable sets $Y,Z \subset \mathbb{C}^p$ with $Y \subset Z$, a function $F:Y \to \mathbb{C}$ and a set $E$ of functions on $Z$ with values in $\mathbb{C}^-$ we say that $F$ can be constructed from $E$ if $F$ can be constructed from $E|_Y$ according to Definition~\ref{3.8}.\\

For the proof of our main result of Section~9.1 we need the following technical notion.

\begin{definition}
\label{def:zeta}
For a definable set $Z \subset \mathbb{R}^n \times \mathbb{C}$ and $(t,x) \in \mathbb{R}^n \times \mathbb{R}$ we set \index{$\zeta_Z$}
$$\zeta_Z:=\sup\{r \geq 0 \mid x + iQ(0,r) \subset Z_t\}.$$
\end{definition}


Now we establish the main result of Section~9.1. Remember for $r \in \mathbb{R}_{>0}$ and $c \in \mathbb{C}$ that $B(c,r):=\{z \in \mathbb{C} \mid \vert{z-c}\vert<r\}$ and $Q(c,r):=\{x \in \mathbb{R} \mid \vert{x-c}\vert<r\}$ if $c \in \mathbb{R}$.

\begin{theorem}
\label{5.90}
	Let $X \subset \mathbb{R}^{n+1}$ be definable such that $X_w$ is open for every $w \in \pi_l(X)$. Let $f:X \to \mathbb{R}, (w,u,x) \mapsto f(w,u,x),$ be restricted log-exp-analytic in $(u,x)$ such that $f_w$ is real analytic for all $w \in \mathbb{R}^l$. Then $f$ has a unary high parametric global complexification $F:Z \to \mathbb{C}, (w,u,z) \mapsto F(w,u,z),$ with respect to $(u,x)$ which is restricted log-exp-analytic in $(u,z)$.
\end{theorem}

\begin{proof}
For a definable cell $C \subset X$ we say the following: The expression ''$C$-consistent'' means always ''$C$-consistent in $u$ with respect to $X$'' and ''$C$-regular'' means always ''$C$-regular in $u$ with respect to $X$''.

Let $e \in \mathbb{N}_0$ be such that $f$ is restricted log-exp-analytic in $(u,x)$ of order at most $e$. By Fact~\ref{4.12} there is $r \in \mathbb{N}_0$ and a definable cell decomposition $\mathcal{C}'$ of $X_{\neq 0}$ such that $f|_{C'}$ is $(m+1,X)$-restricted log-exp-analytically $(e,r)$-prepared in $x$ for every $C' \in \mathcal{C}'$. Let $\mathcal{C} \subset \mathcal{C}'$ be the set of fat cells with respect to $x$. 

For every $C \in \mathcal{C}$ we consider the following. Fix a finite set $E_C$ of locally bounded functions in $(u,x)$ with reference set $X$ such that every $g \in \log(E_C) \cup \{f|_C\}$ is log-exp-analytically $(l,r)$-prepared in $x$ with respect to $E_C$ for $l \in \{-1, \ldots ,e\}$. For $g \in \log(E_C) \cup \{f|_C\}$ fix a corresponding preparing tuple 
$$(r,\mathcal{Y}_C,a_{C,g},\exp(d_{C,0,g}),q_{C,g},s,v_{C,g},b_{C,g},\exp(d_{C,g}),P_{C,g})$$
where $b_{C,g}:=(b_{C,1,g}, \ldots ,b_{C,s,g})$, $\exp(d_{C,g}):=(\exp(d_{C,1,g}), \ldots ,\exp(d_{C,s,g}))$. By redefining $a_{C,g}$ and $v_{C,g}$ suitably we can assume $\vert{v_{C,g}}\vert<\pi$. (Since $\mathcal{C}$ and $E_C$ are finite, $s$ and $r$ can be chosen independently of the fat cell $C$ and of the function $g \in \log(E_C)$). Note that the center $\Theta_C:=(\Theta_{C,0}, \ldots ,\Theta_{C,r})$ of $\mathcal{Y}_C$ and for every $g \in \log(E_C) \cup \{f|_C\}$ the function $b_{C,g}$ is $C$-consistent (compare with Example~\ref{5.51}(1) and Lemma~\ref{5.53}(1)) and $C$-regular (since they are $C$-nice and every $C$-nice function is $C$-regular by Example~\ref{5.56}(1)). For every $g \in \log(E_C) \cup \{f|_C\}$ we have that $a_g$ is $C$-regular since it is $C$-nice and if $g \in \log(E_C)$ then $a_g$ is also $C$-consistent (by Lemma~\ref{5.53}(2)).  

For every $C \in \mathcal{C}$ fix the sign $\sigma_C:=(\sigma_{C,0}, \ldots ,\sigma_{C,r})$ of $\mathcal{Y}_C$, its change index 
$$k_C:=(k^{\textnormal{ch}})_C:=\max\{l \in \{0, \ldots ,r\} \mid \sigma_{C,l}=-1\}-1,$$
its corresponding functions $\mu_{C,1}, \ldots, \mu_{C,r}$ from Definition~\ref{5.2} which are $C$-regular by Example~\ref{5.56}(2). 

The expression ''$(C,\kappa)$-pesistent'' and ''$(C,\kappa)$-soft'' means always
''$(C,\mathcal{Y}_C,\kappa)$-persistent in $u$ with respect to $X$'' and ''$(C,\mathcal{Y}_C,\kappa)$-soft in $u$ with respect to $X$'' respectively for $\kappa=r$ and all $\kappa \in \{k_C,r\}$ if $k_C \geq 0$.\\

Let $H_C$ be from Definition~\ref{5.11}(b). For every $g \in \log(E_C)$ fix $\Phi_{C,g}$ from Definition~\ref{5.78} and for $\Phi_{C,g}$ a complex preparing tuple 
$$(r,\mathcal{Z}_C,a_{C,g},\exp(D_{C,0,g}),q_{C,g},s,V_{C,g},b_{C,g},\exp(D_{C,g}),P_{C,g})$$
where
$$\exp(D_{C,g}):=(\exp(D_{C,1,g}), \ldots ,\exp(D_{C,s,g})).$$
Note that the coefficient, center and base functions of $\Phi_{C,g}$ coincides with the corresponding coefficient, center and base functions of $g$. 
Let $\Lambda_C:=\Lambda_{C,f}$, $\Phi_C:=\Phi_{C,f}$ and $E_{\Phi_C}:=E_{\Phi_{C,f}}$. 

By Lemma~\ref{5.86} and Lemma~\ref{5.85} there is for $\kappa=r$ and all $\kappa \in \{k_C,r\}$ if $k_C \geq 0$ a $C$-regular $(C,\kappa)$-persistent $M_{C,\kappa}:\pi(C) \to \mathbb{R}_{\geq 0}$ with $M_{C,\kappa} < 1/2\vert{\mu_{C,r}-\mu_{C,k_C}}\vert$ if $k_C \geq 0$  and a $C$-regular $(C,\kappa)$-soft $N_{C,\kappa}:\pi(C) \to \mathbb{R}_{\geq 0}$ such that 
$$K_{C,\kappa}:=\{(t,z) \in H_C \mid N_{C,\kappa}(t)<\vert{z-\mu_{C,\kappa}(t)}\vert<M_{C,\kappa}(t)\} \subset \Lambda_C,$$
$\Phi_C|_{K_{C,\kappa}}$ is complex $r$-log-analytically prepared in $z$ with center $\Theta_C$ as in Theorem~\ref{5.84}, and $(\Phi_C|_{K_{C,\kappa}})_t$ has a holomorphic extension on $B(\mu_{C,\kappa}(t),M_{C,\kappa}(t))$ for $t \in \pi(C)$ if $N_{C,\kappa}(t)=0$, $M_{C,\kappa}(t)>0$, and $(\Phi_C|_{K_{C,\kappa}})_t$ has a holomorphic extension at $\mu_{C,\kappa}(t)$ (see Lemma~\ref{5.86}(2)). $(+_1)$\\
Let 
$$\alpha_{C,\kappa}: \pi(C) \to \mathbb{R}_{\geq 0}, t \mapsto \left\{\begin{array}{ll} 2N_{C,\kappa}(t), & N_{C,\kappa}(t) \neq 0,\\
	M_{C,\kappa}(t)/2, & N_{C,\kappa}(t)=0. \end{array}\right.$$
Note that $N_{C,\kappa}(t) < \alpha_{C,\kappa}(t)<M_{C,\kappa}(t)$ for $t \in \pi(C)$ with $M_{C,\kappa}(t) \neq 0$ and that $\alpha$ is $C$-regular. Therefore, for
$$\Delta_{C,\kappa}:=\{(t,z) \in \pi(C) \times \mathbb{R}_{>0} \times \mathbb{C} \mid N_{C,\kappa}(t) < M_{C,\kappa}(t), z \in B(\mu_{C,\kappa}(t),\alpha_{C,\kappa}(t))\}$$
the function
$$\Psi_{C,\kappa}:\Delta_{C,\kappa} \to \mathbb{C}, (t,z) \mapsto \frac{1}{2\pi i} \int_{\partial B(\mu_{C,\kappa}(t),\alpha(t))} \frac{\Phi_C(t,\xi)}{\xi-z}d\xi,$$
is well-defined (by Lemma~\ref{5.85} where $s$ is replaced by $\alpha(t)$). 

We can choose a finite set $\mathcal{E}_C$ of positive definable functions on $\pi(C)$ such that every $g \in \log(\mathcal{E}_C)$ is $C$-consistent and the following holds.
\begin{itemize}
	\item [(1)] 
	We have $\exp(\Theta_{C,0}), \ldots ,\exp(\Theta_{C,r}) \in \mathcal{E}_C$ and $\exp(a_{C,g}) \in \mathcal{E}_C$ if $g \in \log(E_C)$ (since the exponential of a $C$-consistent function is again $C$-consistent by Remark~\ref{5.52}(3)).
	\item [(2)]
	The function $\Phi_C$ can be constructed from $\mathcal{E}_C \cup E_{\Phi_C}$ (by Lemma~\ref{5.79}). $(ß)$ 
	\item [(3)]
	The functions $\Phi_C|_{K_{C,\kappa}}$ and $\Psi_{C,\kappa}$ can be constructed from $\mathcal{E}_C$ for $\kappa=r$ and all $\kappa \in \{k_C,r\}$ if $k_C \geq 0$ (by Claim~\ref{claim-10} applied on the complex log-analytical prepared $\Phi|_{K_{C,\kappa}}$ for the former, and Lemma~\ref{5.85} in combination with Lemma~\ref{lem:construction-complex} for the latter). $(ßß)$
\end{itemize}

Fix $l_C^* \in \mathbb{N}$ and positive definable functions $e_{C,1}, \ldots ,e_{C,l_C^*}$ on $\pi(C)$ such that $\mathcal{E}_C=\{e_{C,1}, \ldots ,e_{C,l_C^*}\}$.\\


The aim is to construct a unary parametric global complexification of $f|_C$ for every $C \in \mathcal{C}$ such that, after gluing the single extensions together, we obtain a unary \emph{high} parametric global complexification of $f$ with respect to $(u,x)$ which is \emph{restricted log-exp-analytic} in $(u,z)$.\\
By Lemma~\ref{5.79} $\Phi_C$ is a unary parametric global complexification of $f|_C$, but to be able to prove the \emph{highness} (which happens in Step~3 below) we have to extend $\Phi_C$ to an even larger unary parametric global complexification of $f|_C$ which happens in Step~1.\\

{\bf Step 1: Construct a suitable unary parametric complexification $G_C:\hat{\Xi}_C \to \mathbb{C}, (t,z) \mapsto F(t,z),$ of $f|_C:C \to \mathbb{R}, (t,x) \mapsto f(t,x),$ for every $C \in \mathcal{C}$.}\\

Let $\kappa = r$ or $\kappa \in \{k_C,r\}$ if $k_C \geq 0$. We set
$$B_{C,\kappa}:=\{t \in \pi(C) \mid \alpha_{C,\kappa}(t) < M_{C,\kappa}(t) \textnormal{ and } B(\mu_{C,\kappa}(t),\alpha_{C,\kappa}(t)) \to \mathbb{C},$$
$$ z \mapsto \Psi_{C,\kappa}(t,z), \textnormal { is holomorphic}\}.$$

Then $B_{C,\kappa}$ is definable since $\Psi_{C,\kappa}$ is definable, $\alpha_{C,\kappa}$ is definable and since the property of being holomorphic is definable by the Cauchy-Riemann equations. We set
$$\mathcal{G}_{C,\kappa}:=\{(t,z) \in H_C \mid t \in B_{C,\kappa} \textnormal{ and } z \in B(\mu_{C,\kappa}(t),\alpha_{C,\kappa}(t))\}.$$
Note that $\mathcal{G}_{C,\kappa} \subset \Delta_{C,\kappa}$. Hence, the function
$$\chi_{C,\kappa}:\mathcal{G}_{C,\kappa} \to \mathbb{C}, (t,z) \mapsto \Psi_{C,\kappa}(t,z),$$
can also be constructed from $\mathcal{E}_C$. $(ßßß)$

Let
$$\hat{B}_{C,\kappa}:=\{t \in B_{C,\kappa} \mid (\Phi_C)_t=(\chi_{C,\kappa})_t \textnormal{ on } (K_C)_t \cap (\mathcal{G}_{C,\kappa})_t \}$$
and
$$\mathcal{T}_{C,\kappa}:=\{(t,z) \in H_C \mid t \in \hat{B}_{C,\kappa} \textnormal{ and } z \in B(\mu_{C,\kappa}(t),M_{C,\kappa}(t)) \}.$$

We have that 
$\mathcal{T}_{C,\kappa} \subset K_{C,\kappa} \cup \mathcal{G}_{C,\kappa} \subset \Lambda_C \cup \mathcal{G}_{C,\kappa}$ (in particular $(\mathcal{T}_{C,\kappa})_t = (K_{C,\kappa})_t \cup (\mathcal{G}_{C,\kappa})_t$ if $t \in \hat{B}_{C,\kappa}$ since $N_{C,\kappa}(t) < \alpha_{C,\kappa}(t)<M_{C,\kappa}(t)$ holds) and that $\mathcal{T}_{C,k_C} \cap \mathcal{T}_{C,r} = \emptyset$ since $M_{C,\kappa}<1/2\vert{\mu_{C,r}-\mu_{C,k_C}}\vert$ if $k_C \geq 0$ (which also implies $\mathcal{G}_{C,k_C} \cap \mathcal{G}_{C,r} = \emptyset$ if $k_C \geq 0$).\\

For the proof that the unary parametric global complexification $F: (t,u,z) \mapsto F(t,u,z)$ of $f$ (which we construct in Step~2 and Step~3 below) is restricted log-exp-analytic in $(u,z)$ (which happens in Step~4 below) we have to shrink $\Lambda_C$ using the following definable $H_{C,\textnormal{small}} \subset H_C$ where $\mathfrak{d}_{C,\kappa}(t):=\textnormal{dist}(\mu_{C,\kappa}(t), C_t)/2$ for $t \in \pi(C)$ and $\kappa=r$ or $\kappa \in \{k^{\text{ch}},r\}$ if $k^{\text{ch}} \geq 0$.\\

If $\prod_{j=0}^r \sigma_{C,j}=1$ (i.e. $\mu_{C,r}<C$) and $k_C < 0$ let
\begin{align*}
H_{C,\textnormal{small}}:=&\{(t,z) \in \pi(C) \times \mathbb{C} \mid \textnormal{Re}(z) \in \textnormal{} ]\mu_{C,r}(t)+\mathfrak{d}_{C,r}(t),\infty[\}\\
=&\{(t,z) \in \pi(C) \times \mathbb{C} \mid \textnormal{Re}(z) \in \textnormal{} ]\inf(C_t) - \mathfrak{d}_{C,r}(t),\infty[\}.
\end{align*}
If $\prod_{j=0}^r \sigma_{C,j}=-1$ (i.e. $\mu_{C,r}>C$) and $k_C < 0$ let
\begin{align*}
H_{C,\textnormal{small}}:=&\{(t,z) \in \pi(C) \times \mathbb{C} \mid \textnormal{Re}(z) \in \textnormal{} ]-\infty,\mu_{C,r}(t)-\mathfrak{d}_{C,r}(t)[\}\\
=&\{(t,z) \in \pi(C) \times \mathbb{C} \mid \textnormal{Re}(z) \in \textnormal{} ]-\infty,\sup(C_t)+\mathfrak{d}_{C,r}(t)[\}.
\end{align*}
If $\prod_{j=0}^r \sigma_{C,j}=1$ and $k_C \geq 0$ (i.e. $\mu_{C,r}<C<\mu_{C,k_C}$) let
\begin{align*}
H_{C,\textnormal{small}}:=&\{(t,z) \in \pi(C) \times \mathbb{C} \mid \textnormal{Re}(z) \in \textnormal{} ]\mu_{C,r}(t)+\mathfrak{d}_{C,r}(t),\mu_{C,k_C}(t)-\mathfrak{d}_{C,k_C}(t)[\}\\
=&\{(t,z) \in \pi(C) \times \mathbb{C} \mid \textnormal{Re}(z) \in \textnormal{} ]\inf(C_t)-\mathfrak{d}_{C,r}(t),\sup(C_t)+\mathfrak{d}_{C,k_C}(t)[\}.
\end{align*}
If $\prod_{j=0}^r \sigma_{C,j}=-1$ and $k_C \geq 0$ (i.e. $\mu_{C,k_C}<C<\mu_{C,r}$) let
\begin{align*}
H_{C,\textnormal{small}}:=&\{(t,z) \in \pi(C) \times \mathbb{C} \mid \textnormal{Re}(z) \in \textnormal{} ]\mu_{C,k_C}(t)+\mathfrak{d}_{C,k_C}(t),\mu_{C,r}(t) - \mathfrak{d}_{C,r}(t)[\}\\
=&\{(t,z) \in \pi(C) \times \mathbb{C} \mid \textnormal{Re}(z) \in \textnormal{} ]\inf(C_t)-\mathfrak{d}_{C,k_C}(t),\sup(C_t) + \mathfrak{d}_{C,r}(t)[\}.
\end{align*}

Note that $C \subset H_{C,\textnormal{small}} \subset H_C$ (compare with Lemma~\ref{5.13}). 
Let $\Gamma_C := \Lambda_C \cap H_{C,\textnormal{small}}$. Then $C \subset \Gamma_C$ 
and $(\Gamma_C)_t$ is open for $t \in \pi(C)$ (since $(H_{C,\textnormal{small}})_t$ is open for $t \in \pi(C)$). Let
$$(\hat{\Xi}_C)_t:=\left\{\begin{array}{ll} (\Gamma_C)_t \cup (\mathcal{T}_{C,r})_t, & t \in \hat{B}_{C,r}, \\
(\Gamma_C)_t, & t \in \pi(C) \setminus \hat{B}_{C,r} , \\
\emptyset, & t \notin \pi(C). \end{array}\right. $$
We define $G_C:\hat{\Xi}_C \to \mathbb{C}$ by letting
$$(G_C)_t(z):= \left\{\begin{array}{ll} (\Phi_C)_t(z) \textnormal{ if } z \in ((\Gamma_C)_t \setminus (\mathcal{T}_{C,r})_t) \cup (K_{C,r})_t 
\textnormal{ and } \\
(\chi_{C,r})_t(z) \textnormal{ if } z \in (\mathcal{G}_{C,r})_t, & t \in \hat{B}_{C,r}, \\
(\Phi_C)_t(z), & t \in \pi(C) \setminus \hat{B}_{C,r}. \end{array}\right.$$
Note that $(\hat{\Xi}_C)_t$ is open for every $t \in \pi(C)$ and that $G_C$ is a well-defined unary parametric global complexification of $f|_C:C \to \mathbb{R}$ (due to $(\mathcal{T}_{C,r})_t = (K_{C,r})_t \cup (\mathcal{G}_{C,r})_t$ for $t \in \hat{B}_{C,r}$).

Suppose that $k:=k_C \geq 0$. Let for $t \in \mathbb{R}^n$
$$
(\hat{\Xi}_C)_t:=\left\{\begin{array}{ll} (\Gamma_C)_t \cup (\mathcal{T}_{C,k})_t \cup (\mathcal{T}_{C,r})_t, & t \in \hat{B}_{C,k} \cap \hat{B}_{C,r}, \\
(\Gamma_C)_t \cup (\mathcal{T}_{C,k})_t, & t \in \hat{B}_{C,k} \setminus \hat{B}_{C,r} , \\
(\Gamma_C)_t \cup (\mathcal{T}_{C,r})_t, & t \in \hat{B}_{C,r} \setminus \hat{B}_{C,k}, \\
(\Gamma_C)_t, & t \in \pi(C) \setminus (\hat{B}_{C,k} \cup \hat{B}_{C,r}) , \\
\emptyset, & t \notin \pi(C). \end{array}\right. $$
Finally define $G_C:\hat{\Xi}_C \to \mathbb{C}$ by letting
$$(G_C)_t(z):=$$
$$
\left\{\begin{array}{ll} (\Phi_C)_t(z) \textnormal{ if } z \in ((\Gamma_C)_t \setminus (\mathcal{T}_{C,r} \cup \mathcal{T}_{C,k})_t) \cup (K_{C,r})_t \cup (K_{C,k})_t \\ 
\textnormal{ and } (\chi_{C,r})_t(z) \textnormal{ if } z \in (\mathcal{G}_{C,r})_t \textnormal{ and } 
(\chi_{C,k})_t(z) \textnormal{ if } z \in (\mathcal{G}_{C,k})_t, & t \in \hat{B}_{C,k} \cap \hat{B}_{C,r} , \\
(\Phi_C)_t(z) \textnormal{ if } z \in  ((\Gamma_C)_t \setminus (\mathcal{T}_{C,k})_t) \cup (K_{C,k})_t \textnormal{ and } \\ 
(\chi_{C,k})_t(z) \textnormal{ if } z \in (\mathcal{G}_{C,k})_t, & t \in \hat{B}_{C,k} \setminus \hat{B}_{C,r}, \\
(\Phi_C)_t(z) \textnormal{ if } z \in  ((\Gamma_C)_t \setminus (\mathcal{T}_{C,r})_t) \cup (K_{C,r})_t 
\textnormal{ and }\\
(\chi_{C,r})_t(z) \textnormal{ if } z \in (\mathcal{G}_{C,r})_t, & t \in \hat{B}_{C,r} \setminus \hat{B}_{C,k}, \\
(\Phi_C)_t(z), & \textnormal{else.} \end{array}\right.$$

Note that $(\hat{\Xi}_C)_t$ is open, and by construction $G_C$ is a well-defined unary parametric global complexification of $f|_C:C \to \mathbb{R}$ (due to $(\mathcal{T}_{C,\kappa})_t = (K_{C,\kappa})_t \cup (\mathcal{G}_{C,\kappa})_t$ for $t \in \hat{B}_{C,\kappa}$ and $\kappa \in \{k,r\}$). \\

{\bf Step 2: Construct a suitable unary parametric global complexification $G:Z \to \mathbb{R}, (w,u,z) \mapsto G(w,u,z),$ of $f:X \to \mathbb{R}, (w,u,x) \mapsto f(w,u,x),$ by gluing the single extensions $G_C$ of $f|_C$ together.}\\



Note that $t=(w,u)$. We say that $(t,z) \in \hat{\Xi}_C$ fulfills
\begin{itemize}
	\item[-] $(\omega_1)$ if $(t,\textnormal{Re}(z)) \in X$,
	\item[-] $(\omega_2)$ if  
	for every $z:=z_1+iz_2 \in (\hat{\Xi}_C)_t$ we have that $z_1+i]0,z_2[ \textnormal{ } \subset (\hat{\Xi}_C)_t$ if $\textnormal{Re}(z_2)>0$ or $z_1+i]z_2,0[ \textnormal{ } \subset (\hat{\Xi}_C)_t$ if $\textnormal{Re}(z_2)<0$ where $z_1:=\text{Re}(z)$ and $z_2:=\text{Im}(z)$, 
	\item[-] $(\omega_3)$ if $\textnormal{dist}(z,C_t) < \mathcal{L}_C(t)$ for 
	\[
	\mathcal{L}_C(t):=\min\{L_C(t),\exp(-\vert{\log(e_{C,1}(t))}\vert), \ldots ,\exp(-\vert{\log(e_{C,l_C^*}(t))}\vert)\}.
	\]
\end{itemize}

In the construction of $G$ properties $(\omega_1)$ and $(\omega_2)$ are needed to obtain well-definability of $G$ and property $(\omega_3)$ to get restricted log-exp-analyticity of $G$ in $(u,z)$. 

Define
$$\Xi_C^*:=\{(t,z) \in \hat{\Xi}_C \mid (t,z) \textnormal{ fulfills properties } (\omega_1), (\omega_2) \text{ and } (\omega_3)\}.$$
Then $\Xi_C^*$ is definable with $C \subset \Xi_C^*$ and $(\Xi^*_C)_t$ is open for $t \in \pi(C)$. Further let
$$\Xi_C:=\{(t,z) \in \Xi_C^* \mid z \textnormal{ is in the same connected component of $(\Xi_C^*)_t$ as $C_t$}\}.$$
Then $\Xi_C$ is definable with $C \subset \Xi_C$ and $(\Xi_C)_t$ is open since $X_t$ is open for $t \in \pi(C)$. Let
$$\Omega:=X \cup \bigcup_{C \in \mathcal{C}} \Xi_C.$$

Then $\Omega$ is definable. We define
$$G:\Omega \to \mathbb{C}, (t,z) \mapsto \left\{\begin{array}{ll} G_C(t,z),& (t,z) \in \Xi_C,\\
f(t,z),& (t,z) \in X. \end{array}\right.$$
We establish that $G$ is a well-defined unary parametric global complexification of $f:X \to \mathbb{R}$ with the following three claims. Claim~\ref{claim-17} is a technical preparation for the proof that $\Omega_t$ is open for $t \in \mathbb{R}^n$ which happens in Claim~\ref{claim-open}. Claim~\ref{claim-well-defined} is devoted to well-definability.

\begin{claim}
	\label{claim-17}
	Let $C \in \mathcal{C}$. Let $\kappa=r$ or $\kappa \in \{k_C,r\}$ if $k_C \geq 0$. Let $t \in \pi(C)$ be with $N_{C,\kappa}(t)=0$ and $0<M_{C,\kappa}(t)$. Suppose there is $0<\varepsilon<\min\{M_{C,\kappa}(t),\mathcal{L}_C(t)/2\}$ such that $Q(\mu_{C,\kappa}(t),\varepsilon) \subset X_t$ and $f_t(x)=p(x-\mu_{C,\kappa}(t))$ for $x \in Q(\mu_{C,\kappa}(t),\varepsilon)$ where $p$ is a convergent power series on $]-\varepsilon,\varepsilon[$. 
	Suppose that $C_t \cap B(\mu_{C,\kappa}(t),\varepsilon) \neq \emptyset$.  Then $t \in \hat{B}_{C,\kappa}$ and therefore 
	$$\mathfrak{B}:=B(\mu_{C,\kappa}(t),\varepsilon) \cap (H_C)_t \subset (\Xi_C)_t \cap (\mathcal{T}_{C,\kappa})_t.$$
\end{claim}

\begin{proof}
	By Taylor series expansion we see that $f_t$ has a holomorphic extension at $\mu_\kappa(t)$ (i.e. there is an open neighborhood $U$ of $\mu_\kappa(t)$ and a holomorphic $J:U \to \mathbb{C}$ such that $J|_{U \cap X_t} = f_t$) and hence by $(+_1)$ that $(\Phi_C|_{K_C})_t$ has a holomorphic extension $h$ on $B(\mu_{C,\kappa}(t),M_{C,\kappa}(t))$. Since $\alpha_{C,\kappa}(t)=M_{C,\kappa}(t)/2$ we get by Cauchy's integral formula that
	$$h(z)=\frac{1}{2 \pi i} \int_{\partial B(\mu_{C,\kappa}(t),\alpha_{C,\kappa}(t))} \frac{(\Phi_C)_t(\xi)}{\xi - z} d \xi$$
	for all $z \in B(\mu_{C,\kappa}(t),\alpha_{C,\kappa}(t))$. So we get that 
	$$B(\mu_{C,\kappa}(t),\alpha_{C,\kappa}(t)) \to \mathbb{C}, z \mapsto \Psi_{C,\kappa}(t,z),$$ 
	coincides with $h$ and is therefore holomorphic. Hence $t \in \hat{B}_{C,\kappa}$ (due to $\chi_{C,\kappa}(t,z) = h(z)=\Phi_C(t,z)$ for every $z \in (K_C)_t \cap (\mathcal{G}_{C,\kappa})_t$) which implies
	\begin{align*}
	B(\mu_{C,\kappa}(t),M_{C,\kappa}(t)) &\cap (H_C)_t = (\mathcal{T}_{C,\kappa})_t \subset (\hat{\Xi}_C)_t
	\end{align*}
	and hence $\mathfrak{B} \subset (\mathcal{T}_{C,\kappa})_t \cap (\hat{\Xi}_C)_t$. So to obtain the claim we have to show that $\mathfrak{B} \subset (\Xi_C)_t$. Let $z \in \mathfrak{B}$.
	Since $z \in B(\mu_{C,\kappa}(t),\varepsilon)$ 
	we have that $(t,\text{Re}(z)) \in X$ (due to $Q(\mu_{C,\kappa}(t),\varepsilon) \subset X_t$). It is also obvious that $(t,z)$ fulfills property $(\omega_2)$. Further we see with the triangle inequality and $\mathfrak{B} \cap C_t \neq \emptyset$
	\begin{align*}
		\text{dist}(z,C_t) &\leq \text{dist}(z,\mu_{C,\kappa}(t)) + \text{dist}(\mu_{C,\kappa}(t),C_t) \leq 2\varepsilon  \\
		&\leq \min\{\mathcal{L}_C(t),2M_{C,\kappa}(t)\} \leq \mathcal{L}_C(t)
	\end{align*}
	which yields property $(\omega_3)$. Hence $z \in (\Xi_C^*)_t$. Since $\mathfrak{B} \cap C_t \neq \emptyset$ we see that $z$ is in the same connected component of $(\Xi_C^*)_t$ as $C_t$ which implies that $z \in (\Xi_C)_t$.
\end{proof}

\begin{claim}
	\label{claim-open}
	$\Omega_t$ is open for every $t \in \mathbb{R}^n$. 
\end{claim}

\begin{proof}
	Let $(t,z) \in \Omega$. If $z \in (\Xi_C)_t$ for a $C \in \mathcal{C}$ then we find an open neighborhood $U$ of $z$ in $\mathbb{C}$ with $U \subset \Omega$ since $(\Xi_C)_t$ is open. So suppose that $z \notin (\Xi_C)_t$ for every $C \in \mathcal{C}$. Then $z \in X$ and there is a thin $C' \in \mathcal{C}'$ with respect to $x$ with $(t,z) \in C'$. Since $X_t$ is open we find a fat $C^* \in \mathcal{C}$ with $C^*<C'$ and $z \in (\overline{C^*})_t$. If $z \neq \mu_{C^*,\kappa}(t)$ for $\kappa=r$ and all $\kappa \in \{k_{C^*},r\}$ if $k_{C^*} \geq 0$ we see that $z \in (\Gamma_{C^*})_t$ and also $z \in (\Xi_{C^*})_t$ by Lemma~\ref{5.88} and definition of $H_{C,\text{small}}$, a contradiction. So suppose that $z=\mu_{C^*,\kappa}(t)$ for a $\kappa=r$ or $\kappa \in \{k_{C^*},r\}$ if $k_{C^*} \geq 0$. Then $N_{C^*,\kappa}(t)=0$ and $M_{C^*,\kappa}(t)>0$ by Remark~\ref{rem:closure-set-persistence} since $N_{C^*,\kappa}$ is $\kappa$-soft, $M_{C^*,\kappa}$ is $\kappa$-persistent and $L_{C^*}(t)>0$. Fix $0<\varepsilon<\min\{\mathcal{L}_{C^*}(t)/2,M_{C^*,\kappa}(t)\}$ such that $Q(\mu_{C^*,\kappa}(t),\varepsilon) \subset X_t$ and $f_t$ can be written as convergent power series on $]-\varepsilon,\varepsilon[$. Since $(C^*)_t \cap B(\mu_{C,\kappa}(t),\varepsilon) \neq \emptyset$ we obtain with Claim~\ref{claim-17} that 
	$$(B(\mu_{C^*,\kappa}(t),\varepsilon) \cap (H_{C^*})_t) \cup X_t \subset (\Xi_{C^*})_t \cup X_t \subset \Omega_t.$$ Since $Q(\mu_{C^*,\kappa}(t),\varepsilon) \subset X_t$ we also obtain 
	$$B(\mu_{C^*,\kappa}(t),\varepsilon) \subset (B(\mu_{C^*,\kappa}(t),\varepsilon) \cap H_{C^*}) \cup X_t \subset \Omega_t.$$
	Since $z=\mu_{C^*,\kappa}(t)$ the set to the left is an open neighborhood of $z$ in $\Omega_t$.     
\end{proof}

\begin{claim}
\label{claim-well-defined}
$G:\Omega \to \mathbb{C}$ is a well-defined unary parametric global complexification of $f:X \to \mathbb{R}$.
\end{claim}	
	
\begin{proof}
Let $(t,z) \in \mathbb{R}^n \times \mathbb{C}$ and let $C_1,C_2 \in \mathcal{C}$ be with $(t,z) \in \Xi_{C_1} \cap \Xi_{C_2}$. At first we show that $G_{C_1}(t,z)=G_{C_2}(t,z)$. Note that $t \in \pi(C_1) \cap \pi(C_2)$ and hence $\pi(C_1)=\pi(C_2)$ (by Definition~\ref{13} and~\ref{14} since $\mathcal{C} \subset \mathcal{C'}$ and $\mathcal{C'}$ is a definable cell decomposition of $X$). Further we have that $(C_1)_t$ is in the same connected component of $X_t$ as $(C_2)_t$ by property $(\omega_1)$ (and since every $z \in \Xi_{C_j}$ is in the same connected component as $C_j$ for $j \in \{1,2\}$). 
Denote this component by $X_t^{\text{comp}}$. By property $(\omega_2)$ and Taylor series expansion of $f$ there is an open $Y \subset \mathbb{C}$ with $X_t^{\text{comp}} \subset Y$ such that $(\Xi_{C_1})_t \cap Y$ is connected (since the union of open balls pairwise intersecting each other is again connected) and a holomorphic $h:Y \to \mathbb{C}$ with $h|_{X_t^{\text{comp}}} = f_t|_{X_t^{\text{comp}}}$. Note that neither $Y$ nor $h$ are necessarily definable. By the identity theorem we obtain $G_{C_1}|_{Y \cap (\Xi_{C_1})_t}=h|_{Y \cap (\Xi_{C_1})_t}$ since 
$(C_1)_t$ is fat with respect to $x$ and
$$(C_1)_t \subset \{z \in (\Xi_{C_1})_t \cap Y \mid G_{C_1}(t,z) = h(z)\}.$$
(For $x \in (C_1)_t$ we have that $G_{C_1}(t,x) = h(x)=f(t,x)$.)
In a similar way we see $G_{C_2}|_{Y \cap (\Xi_{C_2})_t}=h|_{Y \cap (\Xi_{C_2})_t}$ and hence, 
$$G_{C_1}|_{Y \cap (\Xi_{C_1})_t \cap (\Xi_{C_2})_t} = G_{C_2}|_{Y \cap (\Xi_{C_1})_t \cap (\Xi_{C_2})_t}.$$
Finally with property $(\omega_2)$ and the identity theorem we obtain
$$G_{C_1}|_{(\Xi_{C_1})_t \cap (\Xi_{C_2})_t} = G_{C_2}|_{(\Xi_{C_1})_t \cap (\Xi_{C_2})_t}$$
which implies $G_{C_1}(t,z)=G_{C_2}(t,z)$.\\
Note that $G_t|_{(\Xi_C)_t}$ is holomorphic and therefore also $G_t|_{\bigcup_{C \in \mathcal{C}} (\Xi_C)_t}$. Since the complement $\Omega_t \setminus (\bigcup_{C \in \mathcal{C}} (\Xi_C)_t)$ is finite (since $\mathcal{C}' \setminus \mathcal{C}$ consists only of finitely many thin cells) and $f_t$ is real analytic the holomorphy of $G_t$ follows again by Taylor series expansion on every $x \in \Omega_t \setminus (\bigcup_{C \in \mathcal{C}} (\Xi_C)_t)$ and then applying the identity theorem.
\end{proof}

{\bf Step 3: Construct a definable $Z \subset \Omega$ such that $F:Z \to \mathbb{C}, (w,u,z) \mapsto G(w,u,z),$ is a unary high parametric global complexification of $f$ with respect to $(u,x)$.}\\

We need two further technical claims before we are able to construct $Z$. Note that $t=(w,u)$.


\begin{claim}
	\label{claim-18}
	Let $C \in \mathcal{C}$. Let $\kappa=r$ or $\kappa \in \{k_C,r\}$ if $k_C \geq 0$. Let $t \in \pi(C)$ be with $N_{C,\kappa}(t) \neq 0$. Suppose there is $0<\varepsilon<\min\{\mathcal{L}_C(t)/2,M_{C,\kappa}(t)\}$ such that $Q(\mu_{C,\kappa}(t),\varepsilon) \subset X_t$ and $f_t(x)=p(x-\mu_{C,\kappa}(t))$ for $x \in Q(\mu_{C,\kappa}(t),\varepsilon)$ where $p$ is a convergent power series on $]-\varepsilon,\varepsilon[$. Suppose that $\alpha_{C,\kappa}(t) <\varepsilon$ and that $C_t \cap B(\mu_{C,\kappa}(t),\varepsilon) \neq \emptyset$. Then $t \in \hat{B}_{C,\kappa}$ and therefore 
	$$\mathfrak{B}:=B(\mu_{C,\kappa}(t),\varepsilon) \cap (H_C)_t \subset (\Xi_C)_t \cap (\mathcal{T}_{C,\kappa})_t.$$
\end{claim}

\begin{proof}
	Note that $\alpha_{C,\kappa}(t) = 2N_{C,\kappa}(t)$. By Taylor series expansion there is a holomorphic extension $h:B(\mu_{C,\kappa}(t),\varepsilon) \to \mathbb{C}$ of $f_t: Q(\mu_{C,\kappa}(t),\varepsilon) \textnormal{} \to \mathbb{R}$. Let 
	$$A:=\{z \in (H_C)_t \mid N_{C,\kappa}(t) < \vert{z-\mu_{C,\kappa}(t)}\vert < \varepsilon\}.$$ 
	Note that $h|_A = \Phi_C|_A$ by the identity theorem since $A$ is connected and $A \cap C_t \neq \emptyset$ (due to $C_t \cap B(\mu_{C,\kappa}(t),\varepsilon) \neq \emptyset$ and $\varepsilon<L_C(t)$). Since $N_{C,\kappa}(t)<\alpha_{C,\kappa}(t)<\varepsilon$ we obtain by Cauchy's integral formula for all $z \in B(\mu_{C,\kappa}(t),\alpha_{C,\kappa}(t))$
	$$h(z)=\frac{1}{2 \pi i} \int_{\partial B(\mu_{C,\kappa}(t),\alpha_{C,\kappa}(t))} \frac{h(\xi)}{\xi-z} d\xi=\frac{1}{2 \pi i} \int_{\partial B(\mu_{C,\kappa}(t),\alpha_{C,\kappa}(t))} \frac{(\Phi_C)_t(\xi)}{\xi-z} d\xi.$$
	So we get that
	$$B(\mu_{C,\kappa}(t),\alpha_{C,\kappa}(t)) \to \mathbb{C}, z \mapsto \Psi_{C,\kappa}(t,z),$$
	coincides with $h$ and is therefore holomorphic. Hence $t \in \hat{B}_{C,\kappa}$ (since $\Phi_{C}(t,z)=h(z)=\chi_{C,\kappa}(t,z)$ for $z \in A$ and $(K_{C,\kappa})_t \cap (\mathcal{G}_{C,\kappa})_t \subset A$). We obtain
	\begin{align*}
		B(\mu_{C,\kappa}(t),M_{C,\kappa}(t)) &\cap (H_C)_t = (\mathcal{T}_{C,\kappa})_t \subset (\hat{\Xi}_C)_t
	\end{align*}
	and hence $\mathfrak{B} \subset (\hat{\Xi}_C)_t$. The rest of the proof is the same as for Claim~\ref{claim-17}. 
\end{proof}

For the next Claim~\ref{claim-19} we need Definition~\ref{def:zeta}.

\begin{claim}
\label{claim-19}
Let $w \in \mathbb{R}^l$ and $(u_0,x_0) \in X_w$. Then
$$\liminf_{(s,v) \to (u_0,x_0)} \zeta_\Omega(w,s,v) >0.$$
\end{claim}

\begin{proof}
We find $\delta>0$ such that $W^*:=Q^m(u_0,\delta) \textnormal{} \times ]x_0 - \delta, x_0 + \delta[ \textnormal{} \subset X_w$ and $f_w(u,x)=p(u-u_0,x-x_0)$ for $(u,x) \in W^*$ for a convergent power series $p$ on $Q^{m+1}(0,\delta)$ (since $X_w$ is open and $f_w$ is real analytic for every $w \in \pi_l(C)$). We assume
$$\liminf_{(s,v) \to (u_0,x_0)} \zeta_\Omega(w,s,v) = 0.$$
Then there is a definable curve $\gamma: \textnormal{} ]0,1[ \textnormal{} \to Q^m(u_0,\delta) \times ]x_0-\delta/2,x_0+\delta/2[, y \mapsto (\gamma_u(y),\gamma_x(y)),$ with $\lim_{y \searrow 0} \gamma(y)=(u_0,x_0)$ such that $\lim_{y \searrow 0}\zeta_\Omega(w,\gamma(y))=0$. For $y \in \textnormal{}]0,1[$ let $t(y):=(w,\gamma_u(y))$ and
$$I(y):=\{C \in \mathcal{C} \mid C_{t(y)} \cap \textrm{} ]x_0-\delta/2,x_0+\delta/2[\textrm{} \neq \emptyset\}.$$
By passing to a suitable subcurve of $\gamma$ if necessary we may assume that $I:=I(y)$ is independent of $y$. We set
$$I_{\textrm{big}}:=\{C \in I: \lim_{y \searrow 0} L_C(t(y))>0\}$$
and
$$I_{\textrm{small}}:=\{C \in I: \lim_{y \searrow 0} L_C(t(y))=0\}.$$
Because $\sum_{C \in I}L_C(t(y)) \geq \delta$ for all $y \in \textnormal{} ]0,1[$ we see that $I_{\textrm{big}} \neq \emptyset$. For $y \in \textnormal{} ]0,1[$ we choose $C_y \in \mathcal{C}$ such that $(t(y)) \in \overline{C_y}$. We may assume that $C^*:=C_y$ is independent of $y$ by passing to a suitable subcurve of $\gamma$ if necessary. We also have $(u_0,x_0) \in \overline{(C^*)_w}$. Let $k^*:=k_{C^*}$. We consider two cases.\\

\textbf{Case 1:} $C:=C^* \in I_{\textrm{big}}$. Note that $\gamma$ is compatible with $C_w$. Further we obtain with Lemma~\ref{lem:C-consistent-enlarged} that $\lim_{y \searrow 0} \mathcal{L}_C(t(y))>0$.\\

{\bf Subcase 1.1:} Assume $\lim_{y \searrow 0} \vert{\gamma_x(y) - \mu_{C,\kappa}(t(y))}\vert>0$ for $\kappa=r$ and all $\kappa \in \{k_C,r\}$ if $k_C \geq 0$. 
Since $\lim_{y \searrow 0} \text{dist}(\gamma_x(y),C_{t(y)})=0$ we can pass to a suitable subcurve if necessary and assume for $y \in \text{}]0,1[$ that there is $0<\varepsilon<\mathcal{L}_C(t(y))/2$ such that $B(\gamma_x(y),\varepsilon) \subset (H_{C,\textnormal{small}})_w$ and $B(\gamma_x(y),\varepsilon) \cap C_{t(y)} \neq \emptyset$ (due to $(u_0,x_0) \in \overline{C_w}$). 
Then with Lemma~\ref{5.88} we can shrink $\varepsilon$ and further refine $\gamma$ if necessary such that $B(\gamma_x(y),\varepsilon) \subset (\Lambda_C)_{t(y)}$ for $y \in \text{}]0,1[$. 
So we obtain for every $y \in \text{}]0,1[$ and $z \in B(\gamma_x(y),\varepsilon)$
$$\text{dist}(z, C_{t(y)}) \leq \text{dist}(z, \gamma_x(y)) + \text{dist}(\gamma_x(y), C_{t(y)}) \leq 2\varepsilon < \mathcal{L}_C(t(y))$$
and thus $(t(y),z)$ fulfills property $(\omega_1), (\omega_2)$ and $(\omega_3)$ (due to $\gamma(]0,1[) \subset X_w$) which implies $(t(y),z) \in \Xi_C$ since it is in the same connected component of $(\Xi^*_C)_t$ as $C_t$ (due to $B(\gamma_x(y),\varepsilon) \cap C_{t(y)} \neq \emptyset$).
This yields
$$\lim_{y \searrow 0} \zeta_{\Omega}(t(y)) \geq \lim_{y \searrow 0} \zeta_{\Xi_C}(t(y)) \geq \varepsilon>0,$$ 
a contradiction. 

{\bf Subcase 1.2:} Suppose that
$$\lim_{y \searrow 0} \vert{\gamma_x(y) - \mu_{C,\kappa}(t(y))}\vert=0$$
for $\kappa=r$ or a $\kappa \in \{k_C,r\}$ if $k_C \geq 0$. Fix this $\kappa$. Since $M_{C,\kappa}$ is $(C,\kappa)$-persistent respectively $N_{C,\kappa}$ is $(C,\kappa)$-soft we have $\lim_{y \searrow 0} M_{C,\kappa}(t(y)) > 0$ and $\lim_{y \searrow 0} N_{C,\kappa}(t(y)) = 0$  with Lemma~\ref{5.40}. By passing to a suitable subcurve of $\gamma$ if necessary we find $0<\varepsilon<\min\{\mathcal{L}_C(t(y))/2,M_{C,\kappa}(t(y)),\delta/2\}$ for $y \in \text{}]0,1[$ such that the following properties hold for every such $y$.
\begin{itemize}
\item[$(1^*)$] $Q^{m+1}((\gamma_u(y),\mu_{C,\kappa}(t(y))),\varepsilon) \subset Q^{m+1}((u_0,x_0),\delta) \subset X_w$ (which implies that for all $y \in \text{]0,1[}$ there is a convergent power series $p_{t(y)}$ on $]-\varepsilon,\varepsilon[$ such that $f_t(x)=p_{t(y)}(x-\mu_\kappa(t))$),
\item[$(2^*)$] $B(\mu_{C,\kappa}(t(y)),\varepsilon) \cap C_{(t(y))} \neq \emptyset$, 
\item[$(3^*)$] $2N_{C,\kappa}(t(y))<\varepsilon$,
\item[$(4^*)$] $\vert{\gamma_x(y)-\mu_{C,\kappa}(t(y))}\vert < \varepsilon/2$.
\end{itemize}
We may also assume that either $N_{C,\kappa}(t(y))>0$ for all $y \in \textnormal{} ]0,1[$ or $N_{C,\kappa}(t(y))=0$ for all $y \in \textnormal{} ]0,1[$.

Assume the former. Then we have $\alpha_{C,\kappa}(t(y))=2N_{C,\kappa}(t(y))$ and therefore $\alpha_{C,\kappa}(t(y))<\varepsilon$ for every $y \in \textnormal{} ]0,1[$. By Claim~\ref{claim-18} 
we have 
$$B(\mu_{C,\kappa}(t(y)),\varepsilon) \cap (H_{C})_{t(y)} \subset (\Xi_{C})_{t(y)}$$
and therefore (due to $Q(\mu_{C,\kappa}(t(y)),\varepsilon) \subset X_t$)
$$B(\mu_{C,\kappa}(t(y)),\varepsilon) \subset \Omega_{t(y)}$$
for all $y \in \textnormal{} ]0,1[$. This gives with $(4^*)$ and Pythagoras 
$$\Bigl\{\gamma_x(y)+iv \in \mathbb{C} \mid v \in \textnormal{} \Bigl]0,\frac{\sqrt{3}}{2} \varepsilon \Bigl[\Bigl\} \subset \Omega_{t(y)}$$
for all $y \in \textnormal{} ]0,1[$. So we have 
$$\lim_{y \searrow 0}\zeta_\Omega(t(y)) \geq \lim\limits_{y \searrow 0}\frac{\sqrt{3}}{2} \varepsilon > 0,$$
but this is a contradiction.\\ 
Assume the latter. Then by the choice of $\varepsilon$ above we can apply 
Claim~\ref{claim-17} and obtain
$$B(\mu_{C,\kappa}(t(y)),\varepsilon) \cap (H_{C})_{t(y)} \subset (\Xi_{C})_{t(y)}$$
for all $y \in \textnormal{} ]0,1[$. The rest of the argumentation is the same as in the former assumption above.

\textbf{Case 2:} $C^* \in I_{\textrm{small}}$.\\
Note that there is $C \in I_{\textnormal{big}}$ with $C<C^*$ such that
$$\lim_{y \searrow 0} \textnormal{dist}(\gamma_x(y),C_{t(y)})=0.$$  
So $(u_0,x_0) \in \overline{C_w}$. Note that $\gamma$ is compatible with $C_w$ and that $\lim_{y \searrow 0} \mathcal{L}_C(t(y))>0$ with Lemma~\ref{lem:C-consistent-enlarged}.

{\bf Subcase 2.1:} Suppose that
$$\lim_{y \searrow 0} \vert{\gamma_x(y)-\mu_{C,\kappa}(t(y))}\vert > 0$$
for $\kappa=r$ and all $\kappa \in \{k_C,r\}$ if $k_C \geq 0$. By a completely similar argument as in Subcase~1.1 applied to $C$ instead of $C^*$ we obtain the desired contradiction.

{\bf Subcase 2.2:} Suppose that
$$\lim_{y \searrow 0} \vert{\gamma_x(y)-\mu_{C,\kappa}(t(y))}\vert = 0$$
for $\kappa=r$ or a $\kappa \in \{k_C,r\}$ if $k_C \geq 0$. Fix such a $\kappa$. Then 
$$\lim_{y \searrow 0}N_{C,\kappa}(t(y))=0$$
and $\lim_{y \searrow 0} M_{C,\kappa}(t(y))>0$ by Lemma~\ref{5.40}.
A completely similar argument as in Subcase~1.2 applied to $C$ instead of $C^*$ provides the desired contradiction. 
\end{proof}

We set
$$Z:=\{(w,u,z) \in \mathbb{R}^l \times \mathbb{R}^m \times \mathbb{C} \mid (w,u,\textnormal{Re}(z)) \in X \textnormal{ and }$$ $$\textnormal{Im}(z) \in Q(0,\rho_\Omega(w,u,\textnormal{Re}(z)))\}$$
where for $(w,u,x) \in \mathbb{R}^l \times \mathbb{R}^m \times \mathbb{R}$ 
$$\rho_\Omega(w,u,x):=\liminf_{(s,v) \to (u,x)} \zeta_{\Omega}(w,s,v).$$
Note that $Z$ is definable with $X \subset Z$ by Claim~\ref{claim-19}.

\begin{claim}
\label{claim-20}
$Z_w$ is open for every $w \in \mathbb{R}^l$. 
\end{claim}
	
\begin{proof}
Let $w \in \mathbb{R}^l$ and $(u_0,z_0) \in Z_w$. Let $z_0:=x_0+iy_0$ for $x_0,y_0 \in \mathbb{R}$. Then there is $r>0$ such that $y_0 < r < \rho_\Omega(w,u_0,x_0)$. Since the function 
$$\rho_\Omega:\mathbb{R}^{m} \times \mathbb{R} \to \mathbb{R}, (u,x) \mapsto \liminf_{(s,v) \to (u,x)} \zeta_\Omega(w,s,v),$$
is lower semi-continuous by Lemma~\ref{5.89} we find an open neighborhood $U$ of $(u_0,x_0)$ contained in the open set $X_w$ such that $\rho_\Omega(w,u,x)>r$ for every $(u,x) \in U$. This shows that $V:=U + iQ(0,r) \subset Z_w$ and we have $(u_0,z_0) \in V$.
\end{proof}

Claim~\ref{claim-20} shows that $F:=G|_Z$ is a unary high parametric global complexification of $f$ with respect to $(u,x)$.\\

{\bf Step 4: Prove that $F:Z \to \mathbb{R}, (w,u,z) \mapsto F(w,u,z),$ is restricted log-exp-analytic in $(u,z)$.}\\

The strategy is to construct a partition of $Z$ into finitely many definable sets such that for every set of this partition the restriction of $F$ on this set is restricted log-exp-analytic in $(u,z)$ with reference set $Z$. Then by Definition~\ref{3.20}, Lemma~\ref{lem:construction-complex} and Remark~\ref{3.13} we see that $F$ is restricted log-exp-analytic in $(u,z)$.\\  

Let $\mathcal{B}:=\{\pi(C) \mid C \in \mathcal{C}\}$. Then $\mathcal{B}$ is a definable cell decomposition of $\pi(X)$ since $X_t$ is open for every $t \in \pi(C)$ (by Definition~\ref{12} and~\ref{13}). Fix $B \in \mathcal{B}$. Set 
$$\mathcal{C}_B:=\{C \in \mathcal{C} \mid \pi(C)=B\}$$
and let $X_B:=\{(t,x) \in X \mid t \in B\}$ and $Z_B:=\{(t,z) \in Z \mid t \in B\}$. Then $X=\dot\bigcup_{B \in \mathcal{B}} X_B$ and $Z=\dot\bigcup_{B \in \mathcal{B}} Z_B$. Let $u_B$ be the cardinality of $\mathcal{C}_B$. Let $C_{B,1}, \ldots ,C_{B,u_B} \in \mathcal{C}_B$ be the unique cells with $C_{B,1} <  \ldots < C_{B,u_B}$. Note that $C_{B,j}$ is fat with respect to $x$ for $j \in \{1, \ldots ,u_B\}$. Further for $j \in \{1, \ldots ,u_B\}$ we set $L_{B,j}:=L_{C_{B,j}}$ (the length of $C_{B,j}$), $\Gamma_{B,j}:=\Gamma_{C_{B,j}}$, $\hat{\Xi}_{B,j}:=\hat{\Xi}_{C_{B,j}}$ $\Xi_{B,j}:=\Xi_{C_{B,j}}$, $\mathcal{E}_{B,j}:=\mathcal{E}_{C_{B,j}}$, $\Phi_{B,j}:=\Phi_{C_{B,j}}$, $k_{B,j}:=k_{C_{B,j}}$, $G_{B,j}:=G_{C_{B,j}}$ and for $\kappa=r$ and all $\kappa \in \{k_{B,j},r\}$ if $k_{B,j} \geq 0$ we set $K_{B,j,\kappa}:=K_{C_{B,j},\kappa}$, $\mathcal{G}_{B,j,\kappa}:=\mathcal{G}_{C_{B,j},\kappa}$, $\mathcal{T}_{B,j,\kappa}:=\mathcal{T}_{C_{B,j},\kappa}$, and $\chi_{B,j,\kappa}:=\chi_{C_{B,j},\kappa}$. From now on we consider every $g \in \mathcal{E}_{B,j}$ (which is a function on $\pi(C_{B,j})$) as a function on $\Xi_{B,j}$. Let $Z_{B,j}:=(\Xi_{B,j} \setminus X_B) \cap Z_B$ for $j \in \{1, \ldots ,u_B\}$. Note that 
that $Z_B \setminus X_B=\cup_{j=1}^{u_B} Z_{B,j}$ which is not necessarily a disjoint union.


If $k_{B,1} < 0$ we set $Z_{B,1,\mathcal{T}}:=\mathcal{T}_{B,1,r} \cap Z_{B,1}$ and if $k_{B,1} \geq 0$ we set $Z_{B,1,\mathcal{T}}:=(\mathcal{T}_{B,1,k_{B,1}} \dot{\cup} \textnormal{ } \mathcal{T}_{B,1,r}) \cap Z_{B,1}$. Set $Z_{B,1,\Gamma}:=(\Gamma_{B,1} \cap Z_{B,1}) \setminus Z_{B,j,\mathcal{T}}$. Suppose inductively that $Z_{B,j-1,\mathcal{T}}$ and $Z_{B,j-1,\Gamma}$ have already been constructed for $j \in \{2, \ldots , u_B\}$ (i.e. for the cells $1, \ldots, j-1$ with base $B$ counted from below). Let 
$$D_{j-1}:=\bigcup_{i=1}^{j-1} (Z_{B,i,\mathcal{T}} \cup Z_{B,i,\Gamma}).$$
Then if $k_{B,j} < 0$ we set $Z_{B,j,\mathcal{T}}:=(\mathcal{T}_{B,j,r} \cap Z_{B,j}) \setminus D_{j-1}$ and if $k_{B,j} \geq 0$ we set $Z_{B,j,\mathcal{T}}:=((\mathcal{T}_{B,j,k_{B,j}} \dot{\cup} \textnormal{ } \mathcal{T}_{B,j,r}) \cap Z_{B,j}) \setminus D_{j-1}$. Set $Z_{B,j,\Gamma}:=(\Gamma_{B,j} \cap Z_{B,j}) \setminus (Z_{B,j,\mathcal{T}} \cup D_{j-1})$. 
Since $\Xi_{B,j} \subset \hat{\Xi}_{B,j}$ and $\hat{\Xi}_{B,j}= \mathcal{T}_{B,j,r} \cup \Gamma_{B,j}$ if $k_{B,j}<0$ and $\hat{\Xi}_{B,j}= \mathcal{T}_{B,j,r} \cup \mathcal{T}_{B,j,k_{B,j}} \cup \Gamma_{B,j}$ otherwise for $j \in\{1, \ldots , u_B\}$,
we see $\bigcup_{j=1}^\nu Z_{B,j} = \overset{.}{\bigcup}_{j=1}^\nu (Z_{B,j,\mathcal{T}} \text{ } \dot{\cup} \text{ } Z_{B,j,\Gamma})$ for $\nu \in \{1, \ldots , u_B\}$ and $Z_B \setminus X_B=\dot{\bigcup}_{j=1}^{u_B}(Z_{B,j,\mathcal{T}} \text{ } \dot{\cup} \text{ } Z_{B,j,\Gamma})$. 
For $B \in \mathcal{B}$ and $j \in \{1, \ldots ,u_B\}$ let
$$E_{B,j,\mathcal{T}}:=\mathcal{E}_{B,j}|_{Z_{B,j,\mathcal{T}}}$$
and
$$E_{B,j,\Gamma}:=(\mathcal{E}_{B,j} \cup E_{\Phi_{B,j}})|_{Z_{B,j,\Gamma}}$$
where every $g \in \mathcal{E}_{B,j}$ is considered as function on $Z_{B,j,\Gamma}$.
Additionally, consider
$$F_{B,j,\mathcal{T}}:Z_{B,j,\mathcal{T}} \to \mathbb{C}, (t,z) \mapsto G_{B,j}(t,z),$$
and
$$F_{B,j,\Gamma}:Z_{B,j,\Gamma} \to \mathbb{C}, (t,z) \mapsto \Phi_{B,j}(t,z).$$

Note that $F_{B,j,\Gamma}$ can be constructed from $E_{B,j,\Gamma}$ by $(ß)$.

\begin{claim}
\label{claim-21}
$F_{B,j,\mathcal{T}}$ can be constructed from $E_{B,j,\mathcal{T}}$.
\end{claim}

\begin{proof}
For $\kappa=r$ and all $\kappa \in \{k_{B,j},r\}$ if $k_{B,j} \geq 0$ we set $\mathcal{W}_\kappa:=G_{B,j}|_{\mathcal{T}_{B,j,\kappa}}$. Then we have for $(t,z) \in \mathcal{T}_{B,j,\kappa}$
(since then $(\mathcal{T}_{B,j,\kappa})_t = (K_{B,j,\kappa})_t \cup (\mathcal{G}_{B,j,\kappa})_t$)
$$\mathcal{W}_\kappa(t,z)=\left\{\begin{array}{ll} \Phi_{B,j}(t,z), & (t,z) \in K_{B,j,\kappa}, \\
\chi_{B,j,\kappa}(t,z), & (t,z) \in \mathcal{G}_{B,j,\kappa},
\end{array}\right.$$
which is well-defined since $\Phi_{B,j}(t,z)= \chi_{B,j,\kappa}(t,z)$ for $(t,z) \in K_{B,j,\kappa} \cap \mathcal{G}_{B,j,\kappa}$ by the definition of $\mathcal{T}_{B,j,\kappa}$ from above. We have that $\Phi_{B,j}|_{K_{B,j,\kappa}}$ can be constructed from $\mathcal{E}_{B,j}$ by $(ßß)$ and that $\chi_{B,j,\kappa}$ can be constructed from $\mathcal{E}_{B,j}$ by $(ßßß)$ for $\kappa=r$ and all $\kappa \in \{k_{B,j},r\}$ if $k_{B,j} \geq 0$. So with Lemma~\ref{lem:construction-complex} (applied on the disjoint sets $\mathcal{T}_{B,j,\kappa} \cap K_{B,j,\kappa}$ and $(\mathcal{T}_{B,j,\kappa} \cap \mathcal{G}_{B,j,\kappa}) \setminus K_{B,j,\kappa}$) we see that $\mathcal{W}_\kappa$ can be constructed from $\mathcal{E}_{B,j}$ since every $g \in \mathcal{E}_{B,j}$ is a positive definable function which depends only on $t$.

If $k_{B,j}<0$ we have that $F_{B,j,\mathcal{T}}=\mathcal{W}_r|_{Z_{B,j,\mathcal{T}}}$. So $F_{B,j,\mathcal{T}}$ can be constructed from $E_{B,j,\mathcal{T}}$. So suppose $k_{B,j} \geq 0$. Since for all $\kappa \in \{k_{B,j},r\}$ $$F_{B,j,\mathcal{T}}|_{\mathcal{T}_{B,j,\kappa} \cap Z_{B,j}} = \mathcal{W}_\kappa|_{\mathcal{T}_{B,j,\kappa} \cap Z_{B,j}}$$
we have that $F_{B,j,\mathcal{T}}|_{\mathcal{T}_{B,j,\kappa} \cap Z_{B,j}}$ can be constructed from $E_{B,j,\mathcal{T}}|_{\mathcal{T}_{B,j,\kappa} \cap Z_{B,j}}$. Since $Z_{B,j,\mathcal{T}} \subset (\mathcal{T}_{B,j,k_{B,1}} \cap Z_{B,j}) \textnormal{ } \dot{\cup} \textnormal{ } (\mathcal{T}_{B,j,r} \cap Z_{B,j})$ we see that $F_{B,j,\mathcal{T}}$ can be constructed from $E_{B,j,\mathcal{T}}$ by Lemma~\ref{lem:construction-complex} since every $g \in E_{B,j,\mathcal{T}}$ is a positive definable function which depends only on $t$.
\end{proof}

Since
\begin{align*}
Z&= \bigcup_{B \in \mathcal{B}} Z_B = X \text{ }\dot{\cup}\text{ } \left(\dot{\bigcup}_{B \in \mathcal{B}} (Z_B \setminus X_B)\right) \\
&= X \text{ }\dot{\cup}\text{ } \left(\dot{\bigcup}_{B \in \mathcal{B}} \dot{\bigcup}_{j=1}^{u_B}Z_{B,j,\mathcal{T}} \text{ } \dot{\cup} \text{ } Z_{B,j,\Gamma}\right)
\end{align*}
we can finally write $F:Z \to \mathbb{C}$ in the following way. If $(t,z) \in Z_{B,j,\mathcal{T}}$ for $B \in \mathcal{B}$ and $j \in \{1, \ldots ,u_B\}$ let $F(t,z) = F_{B,j,\mathcal{T}}(t,z)$. If $(t,z) \in Z_{B,j,\Gamma}$ for $B \in \mathcal{B}$ and $j \in \{1, \ldots ,u_B\}$ let $F(t,z) = F_{B,j,\Gamma}(t,z)$. If $(t,z) \in X$ let $F(t,z)=f(t,z)$.\\

Fix a set $E_f$ of positive definable functions on $X$ such that every $g \in \log(E_f)$ is locally bounded in $(u,x)$ and $f$ can be constructed from $E_f$ (which exists since $f$ is restricted log-exp-analytic in $(u,x)$). Let
$$E:= \{g \mid g: Z \to \mathbb{C} \textnormal{ is function with } g|_X \in E_f,\textnormal{ } g|_{Z_{B,j,\mathcal{T}}} \in E_{B,j,\mathcal{T}}, \textnormal{ }$$ 
$$g|_{Z_{B,j,\Gamma}} \in E_{B,j,\Gamma} \textnormal{ for }B \in \mathcal{B}, j \in \{1, \ldots ,u_B\}\}.$$
By Lemma~\ref{lem:construction-complex} we see that $F$ can be constructed from $E$. 
So it remains to show with Definition~\ref{3.20} and Remark~\ref{3.13} that $(E_f)^{\text{Re}}$, $(E_{B,j,\mathcal{T}})^{\text{Re}}$ and $(E_{B,j,\Gamma})^{\text{Re}}$ consist of positive definable functions whose logarithms are locally bounded with reference set $Z$ where $Z$, $Z_{B,j,\mathcal{T}}$, $Z_{B,j,\Gamma}$ are considered as subsets of $\mathbb{R}^n \times \mathbb{R}^2$ (and $X$ as subset of $\mathbb{R}^{n+1} \times \{0\} \subset \mathbb{R}^n \times \mathbb{R}^2$). Note that $(E_f)^{\text{Re}} = E_f$ and $(E_{B,j,\mathcal{T}})^{\text{Re}} = E_{B,j,\mathcal{T}}$. Further we consider $z \in \mathbb{C}$ as vector in $\mathbb{R}^2$.

\begin{claim}
\label{claim-22}
Let $h \in E_f$. Then the function $\log(h)$ is locally bounded in $(u,z)$ with reference set $Z$.
\end{claim}

\begin{proof}
Let $w \in \mathbb{R}^l$. Let $(u_0,z_0) \in Z_w$. If $\textnormal{Im}(z_0) \neq 0$ we are clearly done by taking an open neighborhood of $(u_0,z_0)$ in $Z_w$ disjoint from $X$ (by Definition~\ref{3.11}). 
So assume $\textnormal{Im}(z_0) = 0$. Then $(u_0,z_0) \in X_w$ and we find an open neighborhood $U$ of $(u_0,z_0)$ in $X_w$ such that $\log(h)|_U$ is bounded. So for an open neighborhood $V$ of $(u_0,z_0)$ in $Z_w$ with $V \cap X \subset U$ we see that $\log(h)|_{V \cap X_w}$ is bounded and are done with the proof of this claim by Definition~\ref{3.11}.
\end{proof} 


\begin{claim}
\label{claim-23}
Let $B \in \mathcal{B}$ and $j \in \{1, \ldots ,u_B\}$. Let $h \in (E_{B,j,\Gamma})^{\text{Re}}$. The function $\log(h)$ is locally bounded in $(u,z)$ with reference set $Z$.
\end{claim}

\begin{proof}
Let $C:=C_{B,j}$, $\mathcal{K}:=Z_{B,j,\Gamma}$, $w \in \mathbb{R}^l$ and $(u_0,z_0) \in Z_w$. We show that $\log(h_w)$ is bounded at $(u_0,z_0)$.

If $(u_0,z_0) \notin \overline{\mathcal{K}_w}$ we are clearly done by taking an open neighborhood of $(u_0,z_0)$ in $Z_w$ disjoint from $\mathcal{K}_w$ (by Definition~\ref{3.11}). So assume $(u_0,z_0) \in \overline{\mathcal{K}_w}$.

\textbf{Case 1:} Let $\exp(h) \in \mathcal{E}_C|_\mathcal{K}$ (where every function from $\mathcal{E}_C$ is considered as a function on $\mathcal{K}$). Suppose that $h_w$ is not bounded at $(u_0,z_0)$ (i.e. at $u_0$). Then there is a definable curve $\gamma: \textnormal{}]0,1[ \textnormal{} \to \mathcal{K}_w$ such that $\lim_{y \searrow 0} \gamma(y)=(u_0,z_0)$ and $\lim_{y \searrow 0} \vert{h(w,\gamma_u(y))}\vert = \infty$. Note that $\gamma$ is not necessarily compatible with $C_w$ at this point. Let $t(y):=(w,\gamma_u(y))$. Note that $(\gamma_u,\gamma_z)$ runs through $(\hat{\Xi}_C)_w$, that $(\gamma_u,\textnormal{Re}(\gamma_z))$ through $X_w$ (by property $(\omega_1)$) and that $(u_0,\text{Re}(z_0)) \in X_w$. By property $(\omega_3)$ we see that $\lim_{y \searrow 0} \mathcal{L}_C(t(y)) = 0$ (since $\exp(h) \in \mathcal{E}_C$) and hence $\lim_{y \searrow 0} \text{dist}(C_{t(y)},\gamma_z(y))=0$ implying $z_0 \in \mathbb{R}$ and $(u_0,z_0) \in \overline{C_w}$. So by passing to a suitable subcurve of $\gamma$ if necessary we may assume that $\gamma(]0,1[) \subset Z_w$ since $(u_0,z_0) \in Z$ and $Z_w$ is open in $\mathbb{R}^m \times \mathbb{C}$. Since $h$ is $C$-consistent, we see that $\lim_{y \searrow 0} L_C(t(y))=0$ (otherwise $\gamma$ would be compatible with $C_w$ and with Lemma~\ref{lem:C-consistent-enlarged} we obtain $\lim_{y \searrow 0} h(t(y)) \in \mathbb{R}$). But then there is also $1 \leq i<j$ with $C^*:=C_{B,i}<C$ such that $(u_0,z_0) \in \overline{(C^*)_w}$. Let $k:=k_{C^*}$ (i.e. the change index of $\mathcal{Y}_{C^*}$). As in the proof of Case~2 in Claim~\ref{claim-19} we look at two subcases.\\
\textbf{Subcase 1:} $\lim_{y \searrow 0} \vert{\gamma_z(y)-\mu_{C^*,\kappa}(t(y))}\vert>0$ for $\kappa=r$ and all $\kappa \in \{k,r\}$ if $k \geq 0$. Note that $\gamma$ is compatible with $(C^*)_w$ and hence $\lim_{y \searrow 0} \mathcal{L}_{C^*}(t(y))>0$ by Lemma~\ref{lem:C-consistent-enlarged}. Therefore we see with the proof of Subcase~2.1 (and hence also with the proof of Subcase~1.1) in Claim~\ref{claim-19} that 
$$(t(y),\gamma_z(y)) \in \Xi_{C^*} \cap Z \subset \bigcup_{i=1}^{j-1} (Z_{B,i,\mathcal{T}} \text{ }\dot{\cup} \text{ } Z_{B,i,\Gamma})$$
for every $y \in \text{}]0,1[$ after passing to a suitable subcurve of $\gamma$ if necessary which is a contradiction since $\mathcal{K}$ is disjoint to the set to the right.\\
\textbf{Subcase 2:} $\lim_{y \searrow 0} \vert{\gamma_z(y)-\mu_{C^*,\kappa}(t(y))\vert}=0$ for $\kappa=r$ or a $\kappa \in \{k,r\}$ if $k \geq 0$. Then we see with the proof of Subcase~1.1 and Subcase~1.2 in Claim~\ref{claim-19} that there is $\varepsilon>0$ such that after passing to a suitable subcurve of $\gamma$ if necessary we have that $\gamma_z(y) \in B(\mu_{C^*,\kappa}(t(y)),\varepsilon) \setminus \mathbb{R}$ and $B(\mu_{C^*,\kappa}(t(y)),\varepsilon) \setminus \mathbb{R} \subset (\Xi_{C^*})_{(t(y))}$ for every $y \in \text]0,1[$. This implies that $(t(y),\gamma_z(y)) \in \Xi_{C^*} \cap Z$ for every $y \in \text{}]0,1[$ which is the same contradiction as in Subcase~1.\\

\textbf{Case 2:} Assume that $\exp(h) \in (E_{\Phi_C}|_\mathcal{K})^{\text{Re}}$. Then $h=\text{Re}(\Phi_g)|_\mathcal{K}$ for a $g \in \log(E_C)$ where $\Phi_g$ is complex log-exp-analytically $(l,r)$-prepared in $z$ for $r \in \mathbb{N}_0$ and $l \in \{-1, \ldots ,e-1\}$. We show the statement by induction on $l$. If $l=-1$ then $g=0$ and the assertion follows.

$l-1 \to l$: Let 
$$(r,\mathcal{Z},a,\exp(D_0),q,s,V,b,\exp(D),P):=$$
$$(r,\mathcal{Z}_{C,g}|_\mathcal{K},a_{C,g}|_\mathcal{K},\exp(D_{C,0,g})|_\mathcal{K},q_{C,g},s,V_{C,g},b_{C,g}|_\mathcal{K},\exp(D_{C,g})|_\mathcal{K},P_{C,g})$$
where $(b_1, \ldots ,b_s):=(b_{C,1,g}, \ldots , b_{C,s,g})$, $P:=(p_1, \ldots ,p_s)^t$ for $p_i \in \mathbb{Q}^{r+1}$ and $(\exp(D_1), \ldots ,\exp(D_s)):=(\exp(D_{C,1,g}), \ldots , \exp(D_{C,s,g}))$. So we have 
$$\Phi_g(t,z) = a(t) (\sigma\mathcal{Z})^{\otimes q}(t,z) \exp(D_0(t,z))$$
$$V(b_1(t) (\sigma\mathcal{Z})^{\otimes p_1}(t,z) \exp(D_1(t,z)), \ldots ,b_s(t) (\sigma\mathcal{Z})^{\otimes p_s}(t,z) \exp(D_s(t,z)))$$
for $(t,z) \in \mathcal{K}$ where $t=(w,u)$. 
We show that $\vert{\Phi_g}\vert$ is bounded at $(u_0,z_0)$ and we obtain also that $h$ is bounded at $(u_0,z_0)$. Suppose that $\vert{\Phi_g}\vert$ is not bounded at $(u_0,z_0)$. Then there is a definable curve $\gamma: \textnormal{}]0,1[ \textnormal{} \to \mathcal{K}_w$ such that 
$$\lim_{y \searrow 0} \vert{\Phi_g}\vert(w,\gamma(y)) = \infty.$$
Note that $\gamma$ is not necessarily compatible with $C_w$ at this point and that $\gamma$ runs through $(H_{C,\text{small}})_w$ since $\mathcal{K} \subset H_{C,\text{small}}$. Let $t(y):=(w,\gamma_u(y))$.  
With Case~1 we can assume that $\lim_{y \searrow 0} g^*(t(y)) \in \mathbb{R}$ for every $g^* \in \log(\mathcal{E}_C)$ and hence $\lim_{y \searrow 0} a(t(y)) \in \mathbb{R}$ and $\lim_{y \searrow 0} \Theta_i(t(y)) \in \mathbb{R}$ for $i \in \{0, \ldots , r\}$ (since $\exp(a)$ and $\exp(\Theta_j) \in \mathcal{E}_C$ for $j \in \{0, \ldots , r\}$). We can also assume that $\lim_{y \searrow 0} L_C(t(y)) \in \mathbb{R}_{>0}$. Otherwise we find $C^*:=C_{i,B}<C$ with $(u_0,z_0) \in \overline{(C^*)_w}$ such that $\gamma$ is compatible with $(C^*)_w$ which gives the same contradiction as in Case~1.
Hence, we get with the induction hypothesis and the fact that there is $R>1$ such that $V$ is bounded on $\overline{D^s(0,R)}$ and 
$$\vert{b_i(t)(\sigma \mathcal{Z})^{\otimes p_i}(t,z)\exp(D_i(t,z))}\vert \leq R$$
for all $i \in \{1, \ldots ,s\}$ and $(t,z) \in \mathcal{K}$, that
$\lim_{y \searrow 0} \vert{(\sigma\mathcal{Z})^{\otimes q}(t(y))}\vert = \infty$. So by the proofs of 
Lemma~\ref{5.28} and Lemma~\ref{5.31} we see that there is $j \in \{0, \ldots ,r\}$ with 
$$\lim_{y \searrow 0} \vert{\mu_{C,j}(t(y))-\gamma_z(y)}\vert = 0$$
which implies with the definition of $H_{C,\text{small}}$ and Corollaries~\ref{only-r},~\ref{k-and-r} that 
$$\lim_{y \searrow 0} \vert{\gamma_z(y)-\mu_{C,\kappa}(t(y))}\vert = 0$$
for $\kappa=r$ or a $\kappa \in \{k_C,r\}$ if $k_C \geq 0$ since for $j \in \{0, \ldots, r\} \setminus \{k_C,r\}$
$$\lim_{y \searrow 0} \vert{\mu_{C,\kappa}(t(y))-\mu_{C,j}(t(y))\vert} > 0$$ by the proof of Corollary~\ref{5.29}(2). Fix such a $\kappa$. We see $(u_0,z_0) \in X_w$. Let $x_0:=z_0$. By construction of $(H_{C,\textnormal{small}})_w$ we see that
$$\lim_{y \searrow 0} \textnormal{dist}(\mu_{C,\kappa}(t(y)),C_{(t(y))})=0.$$ 
So $(u_0,x_0) \in \overline{C_w}$. Now we see that $\gamma$ is compatible with $C_w$ and therefore with the proof of Subcase~1.1 and Subcase~1.2 in Claim~\ref{claim-19} that there is $\varepsilon>0$ such that after passing to a suitable subcurve if necessary we have that $\gamma_z(y) \in B(\mu_{C,\kappa}(t(y)),\varepsilon) \setminus \mathbb{R}$ and $B(\mu_{C,\kappa}(t(y)),\varepsilon) \setminus \mathbb{R} \subset (\mathcal{T}_{C,\kappa})_{(t(y))}$ for every $y \in \text]0,1[$. This implies that $(t(y),\gamma_z(y)) \in \mathcal{T}_{C^*,\kappa}$ for $y \in \text{}]0,1[$ which is a contradiction since $\mathcal{T}_{C^*,\kappa}$ and $\mathcal{K}$ are disjoint by construction.
\end{proof}

Similarly as in Claim~\ref{claim-23} one sees that $\log(h)$ is locally bounded in $(u,z)$ with reference set $Z$ if $h \in E_{B,j,\mathcal{T}}$ for $B \in \mathcal{B}$ and $j \in \{1, \ldots ,u_B\}$. This finishes the proof of Theorem~\ref{5.90}.
\end{proof}
	
\subsection{The Multivariate Case}
We prove by an induction on the number of variables that the class of real analytic restricted log-exp-analytic functions exhibits parametric global complexification which is Theorem B. The ''highness'' from the unary case above is necessary to enable the induction. 
	
\begin{theorem}
\label{5.91}
Let $m \in \mathbb{N}$. Let $w$ range over $\mathbb{R}^m$ and $v$ over $\mathbb{C}^m$. Let $X \subset \mathbb{R}^n \times \mathbb{R}^m$ be definable. Assume that $X_t$ is open for all $t \in \mathbb{R}^n$. Let $f:X \to \mathbb{R}, (t,w) \mapsto f(t,w),$ be restricted log-exp-analytic in $w$ such that $f_t$ is real analytic for every $t \in \mathbb{R}^n$. Then $f$ has an $m$-ary parametric global complexification $F:Z \to \mathbb{C}, (t,v) \mapsto F(t,v),$ which is restricted log-exp-analytic in $v$.
\end{theorem}

\begin{proof}
	For $l \in \{1, \ldots ,m\}$ let $u$ range over $\mathbb{R}^{m-l}$ and $x$ over $\mathbb{R}^l$. Let $\pi_n:\mathbb{R}^n \times \mathbb{R}^m \to \mathbb{R}^n, (t,u,x) \to t,$ be the projection on the first $n$ real coordinates. We show by induction on $l \in \{1, \ldots ,m\}$ that $f:X \to \mathbb{R}, (t,u,x) \mapsto f(t,u,x),$ has an $l$-ary high parametric global complexification $F:Z \to \mathbb{C}, (t,u,z) \mapsto F(t,u,z),$ with respect to $(u,x)$ such that $F$ is also restricted log-exp-analytic in $(u,z)$.
	
	$l=1$: This follows from Theorem~\ref{5.90}.
	
	$l-1 \to l$: Given $x=(x_1, \ldots ,x_l) \in \mathbb{R}^l$ we set $x'':=(x_2, \ldots ,x_l)$. The set
	$$X_{(t,u,x_1)}=\{x'' \in \mathbb{R}^{l-1} \mid (t,u,x_1,x'') \in X\} \subset \mathbb{R}^{l-1}$$
	is open and the function $f_{(t,u,x_1)}:X_{(t,u,x_1)} \to \mathbb{R}$ is real analytic since it coincides with $(f_t)_{(u,x_1)}$ for $(t,u,x_1) \in \mathbb{R}^n \times \mathbb{R}^{m-l} \times \mathbb{R}$. Let $\mathbb{R}^{n+m-l+1}$ be the parameter space and write $\hat{X}$ if we consider $X$ as a subset of $\mathbb{R}^{n+m-l+1} \times \mathbb{R}^{l-1}$. Moreover we write $\hat{f}$ if we consider $f$ as a function on $\hat{X}$. By the inductive hypothesis we find a definable set $\hat{Y} \subset \mathbb{R}^{n+m-l+1} \times \mathbb{C}^{l-1}$ such that $\hat{Y}_t \subset \mathbb{R}^{m-l+1} \times \mathbb{C}^{l-1}$ is open for $t \in \mathbb{R}^n$ and an $(l-1)$-ary high parametric global complexification $\hat{G}:\hat{Y} \to \mathbb{C}, (t,u,x_1,z'') \mapsto \hat{G}(t,u,x_1,z''),$ of $f$ with respect to $(u,x)$ which is restricted log-exp-analytic in $(u,x_1,z'')$. We set $G_1:=\textnormal{Re}(\hat{G}):\hat{Y} \to \mathbb{R}$ and $G_2:=\textnormal{Im} (\hat{G}):\hat{Y} \to \mathbb{R}$. Note that $G_1|_X=f$, $G_2|_X=0$, and that $G_1$ and $G_2$ are restricted log-exp-analytic in $(u,x_1,z'')$ (according to Definitions~\ref{3.8} ,~\ref{3.20} and Remark~\ref{rem:rest-log-exp-real-imaginary}). Considering $z''$ as a tuple of $2(l-1)$ real variables we set 
	$$Y:=\{(t,u,x_1,z'') \in \hat{Y} \mid (G_j)_t \textnormal{ is real analytic at $(u,x_1,z'')$ for $j \in \{1,2\}$}\}$$ 
	Then $Y$ is definable by a parametric version of Tamm's theorem for restricted log-exp-analytic functions (see Corollary 3.27 in~\cite{28}). We have $X \subset Y \subset \hat{Y}$. Since real analytiticity is an open property $Y_t$ is open in $\mathbb{R}^{m-l+1} \times \mathbb{C}^{l-1}$ for all $t \in \mathbb{R}^n$. 
	
	Again considering $z''$ as a tuple of $2(l-1)$ real variables we find with Theorem~\ref{5.90} a unary high parametric global complexification 
	$$F_j:Z_j \to \mathbb{C}, (t,u,\hat{z},z'') \mapsto F_j(t,u,\hat{z},z''),$$ 
	for $G_j|_Y$ with respect to $(u,x_1,z'')$ which is also restricted log-exp-analytic in $(u,\hat{z},z'')$ where $j \in \{1,2\}$. (Here $\hat{z}$ is an additional variable which ranges over $\mathbb{C}$. Note also that $(Z_j)_t$ is open in $\mathbb{R}^{m-l} \times \mathbb{C}^l$ for $t \in \mathbb{R}^n$ with $X \subset Z_j$ for $j \in \{1,2\}$.) $(*)$

	Let
	$$\tilde{Z}:=\{(t,u,z) \in Z_1 \cap Z_2 \mid (F_1 +iF_2)_{(t,u)} \textnormal{ is holomorphic at }z\}.$$
	Then $\tilde{Z}$ is definable and $F:=(F_1+iF_2)|_{\tilde{Z}}:\tilde{Z} \to \mathbb{C}$ is restricted log-exp-analytic in $(u,z)$ with reference set $Z_1 \cap Z_2$. 
	
	Let
	$$\zeta_{\tilde{Z}}(t,u,x):=\sup\{r \in \mathbb{R}_{\geq 0} \mid x+iQ^l(0,r) \subset \tilde{Z}_{(t,u)}\}.$$

	\begin{claim}
	\label{claim-24}
	Let $(t,u,x) \in X$. Then
	$$\liminf_{(s,y) \to (u,x)} \zeta_{\tilde{Z}}(t,s,y) > 0.$$ 
	\end{claim}
	
	\begin{proof}
	Since $f_t$ is real analytic for every $t \in \mathbb{R}^n$ we find some $0<r$ and an open neighborhood $V$ of $(u,x)$ in $X_t$ such that $f_t(s,y)=p((s,y)-(u,x))$ for $(s,y) \in V$ where $p$ is a convergent power series on $Q^m(0,r)$. We find some $0<r'<r/2$ such that
	$$U:= \textnormal{} Q^{m-l}(u,r') \textnormal{} \times \textnormal{} ]x_1-r',x_1+r'[ \textnormal{} \times (Q^{l-1}(x'',r') + iQ^{l-1}(0,r')) \subset \hat{Y}_t$$
	and that $\hat{G}_t(s,y_1,z'')=p((s,y_1,z'')-(u,x))$ for all $(s,y_1,z'') \in U$. So for $(s,y) \in Q^m((u,x),r')$ we see that
	$$y'' + iQ^{l-1}(0,r') \subset Y_{(t,s,y_1)}.$$
	So by $(*)$ we find some $0<\varepsilon<r'/2$ and an open neighborhood $W$ of $(u,x)$ in $X_t$ with $W \subset Q^m((u,x),r')$ such that for $(s,y) \in W$ we have $y+iQ^l(0,\varepsilon) \subset (Z_j)_{(t,s)}$ for $j \in \{1,2\}$ and $F_{(t,s)}(z)=p_{s-u}(z-x)$ for all $z \in y+iQ^l(0,\varepsilon)$. This yields $y+iQ^l(0,\varepsilon) \subset \tilde{Z}_{(t,s)}$ for all $(s,y) \in W$. Hence
	\begin{align*}
	\liminf_{(s,y) \to (u,x)} &\zeta_{\tilde{Z}}(t,s,y) \geq \varepsilon > 0. \qedhere
	\end{align*}
	\end{proof}
	Let
	$$Z:=\{(t,u,z) \in \mathbb{R}^n \times \mathbb{R}^{m-l} \times \mathbb{C}^l \mid (t,u,\textnormal{Re}(z)) \in X$$
	$$\textnormal{ and } \textnormal{Im}(z) \in Q^l(0,\rho(t,u,\textnormal{Re}(z)))\}$$
	where 
	$$\rho(t,u,x):=\liminf\limits_{(s,y) \to (u,x)} \zeta_{\tilde{Z}}(t,s,y)$$
	for $(t,u,x) \in \mathbb{R}^n \times \mathbb{R}^{m-l} \times \mathbb{R}^l$. Note that $Z_t$ is open for $t \in \mathbb{R}^n$. (Compare with the proof of Claim~\ref{claim-20} in Theorem~\ref{5.90}: The function $\mathbb{R}^{m-l} \times \mathbb{R}^l, (u,x) \mapsto \rho(t,u,x),$ defines a lower semi-continuous function for every $t \in \mathbb{R}^n$.) By Claim~\ref{claim-24} we have $X \subset Z$. So we see that $F|_Z$ is an $l$-ary high parametric global complexification of $f$ with respect to $(u,x)$ which is restricted log-exp-analytic in $(u,z)$. This finishes the proof of Theorem B.
\end{proof}
	
	An immediate consequence of Theorem~\ref{5.91} is Theorem A.
	
	\begin{theorem}
	\label{5.92}
	Let $n \in \mathbb{N}$ and $U \subset \mathbb{R}^n$ be a definable open set. Let $f:U \to \mathbb{R}$ be a real analytic restricted log-exp-analytic function. Then $f$ has a global complexification which is again restricted log-exp-analytic, i.e. there is a definable open $V \subset \mathbb{C}^n$ with $U \subset V$ and a holomorphic restricted log-exp-analytic $F:V \to \mathbb{C}$ such that $F|_U=f$.
	\end{theorem}
	
	\begin{corollary}
	The following properties hold.
	\begin{itemize}
			\item[(1)] Let $X \subset \mathbb{R}^n \times \mathbb{R}^m$ be definable such that $X_t$ is open for every $t \in \mathbb{R}^n$ and let $f:X \to \mathbb{R}, (t,x) \mapsto f(t,x),$ be log-analytic. Suppose that $f_t$ is real analytic for every $t \in \mathbb{R}^n$. Then there is a definable $V \subset \mathbb{R}^n \times \mathbb{C}^m$ with $U \subset V$ such that $V_t$ is open for every $t \in \mathbb{R}^n$ and a restricted log-exp-analytic $F:V \to \mathbb{C}, (t,z) \mapsto F(t,z),$ in $z$ such that $F_t$ is holomorphic for every $t \in \mathbb{R}^n$ and $F|_U=f$.
			\item[(2)] Let $U \subset \mathbb{R}^n$ be open. Let $g:U \to \mathbb{R}$ be a real analytic log-analytic function. Then there is an open definable $V \subset \mathbb{C}^n$ with $U \subset V$ and a holomorphic restricted log-exp-analytic $G:V \to \mathbb{C}$ such that $G|_U=g$.
	\end{itemize}
	\end{corollary}
	
	\begin{proof}
	(1): This follows directly from the fact that a log-analytic function is restricted log-exp-analytic in $x$.
	(2): This follows directly from (1).
	\end{proof}

\section{Discussions and Conclusions}

We have identified real analytic restricted log-exp-analytic functions as a large subclass of $\mathbb{R}_{\text{an,exp}}$-definable functions which have a global complexification which is again restricted log-exp-analytic. We also provided a parametric version of this result. While prior research has focused on globally subanalytic functions (see Kaiser~\cite{17}), this work demonstrates how analoguous results can be obtained for the structure $\mathbb{R}_{\text{an,exp}}$. It supports the thesis that restricted log-exp-analytic functions share many analytical properties with globally subanalytic ones, despite their distinct behavior from the point of o-minimality. \\
Another novel contribution of this paper is the consideration of parametric integrals in $\mathbb{R}_{\text{an,exp}}$ which have been rarely touched so far (compare for example with Kaiser et al.~\cite{20}, Speissegger~\cite{35} or Van den Dries~\cite{12} for a few papers in this direction). We described a simple method how $\mathbb{R}_{\text{an,exp}}$-definability results of special parametric integrals can be obtained. \\

Finally we want to discuss two open questions for further research outgoing from Theorem A and Theorem B.\\


{\bf Question 1}\\
{\it Does a real analytic log-analytic function have a global complexification which is again log-analytic?\\}

This question could have a positive answer simply by observing that a log-analytic function can be piecewise written as $\mathcal{L}_{\textnormal{an}}(^{-1},(\sqrt[n]{ \ldots })_{n=2,3, \ldots },\log)$-terms and therefore, extends piecewise to a holomorphic log-analytic function. But to enable the induction on the number of variables, one needs a preparation theorem for log-analytic functions with log-analytic data only.
But such questions are pretty unsolved at present. One idea could be to adapt the preparation theorem of Cluckers and Miller for constructible functions from~\cite{5},~\cite{6} respectively~\cite{7} to the log-analytic context in some way.

Since we established global complexification for real analytic restricted log-exp-analytic functions and the class of restricted log-exp-analytic functions is a proper subclass of the class of all definable functions (compare with Example~\ref{3.18}) the following question naturally arises.\\

{\bf Question 2}\\
{\it Does $\mathbb{R}_{\textnormal{an,exp}}$ have global complexification?\\}

With our results above this question can be reformulated as follows: Is every real analytic definable function restricted log-exp-analytic? This may be true from the point of analysis: If the global exponential function comes not locally bounded into the game, features like flatness may occur (compare with Example~\ref{3.18}). But neither real analytic functions nor restricted log-exp-analytic functions are flat. \\

We are hopeful that our results and the proofs may serve as stepping stones towards a
better understanding of the role of restricted log-exp-analytic functions in o-minimality, because a lesson leart from these analyses is that this large class of functions may be the key for generalizing many analytical results, which are valid for globally subanalytic functions, to the structure $\mathbb{R}_{\text{an,exp}}$. 

\printindex 

\end{document}